\documentclass[11pt,
]{book}
\usepackage{amssymb,amstext,amsmath,amsthm,latexsym,mathrsfs,mathbbol}
\usepackage{makeidx}      

\makeindex

\newcommand{\nek}{\newcommand}
\nek{\renek}{\renewcommand}

\DeclareMathAlphabet{\cur}{U}{eur}{m}{n}

\nek{\skr}{\mathscr}

3

\nek{\bfsf} {\sffamily\bfseries\upshape}

\nek{\bfit}{\bfseries\itshape}
\nek{\bfsl}{\bfseries\slshape}
\nek{\sfbs}{\mdseries\sffamily\itshape}        
\DeclareMathAlphabet{\bma}{OML}{cmm}{b}{it}

\nek{\vyk} [1] {}

\nek{\imar}[1]{\marginpar[
\flushright\footnotesize\sl
$\mtho\longrightarrow$\\ \vspace{-1ex}{#1}$\mtho\dashv$\vspace*{1ex}]
{
\flushleft\footnotesize\sl
$\mtho\longleftarrow$\\ \vspace{-1ex}{#1}$\mtho\dashv$\vspace*{1ex}}}

\nek{\jmar}[1]{\marginpar[\vspace{5ex}%
\flushright\footnotesize\sl
$\mtho\longrightarrow$\\ \vspace{0ex}{#1}$\mtho\dashv$\vspace*{0ex}]
{\vspace{5ex}%
\flushleft\footnotesize\sl
$\mtho\longleftarrow$\\ \vspace{0ex}{#1}$\mtho\dashv$\vspace*{0ex}}}

\nek{\imae}[1]{\marginpar[
\flushright\footnotesize\vspace{-4ex}%
$\mtho\longrightarrow$\\%
\vspace{-1ex}{#1}$\mtho$\vspace*{2ex}]%
{
\flushleft\footnotesize\vspace{-4ex}%
$\mtho\longleftarrow$\\%
\vspace{-1ex}{#1}$\mtho$\vspace*{2ex}}
}%

\nek{\lam} [1]
{\label{#1}\hspace*{-3pt}\imar{#1}%
}%
\nek{\las} [1]
{\label{#1}\imae{#1}}%

\renek{\las} [1]{\label{#1}}
\renek{\lam} [1]{\label{#1}}
\renek{\imae}[1]{}
\renek{\imar}[1]{}
\renek{\jmar}[1]{}

\renek{\thechapter}{\arabic{chapter}}
\renek{\thesubsection}
{\setcounter{cly}0\thechapter\hspace{0.1ex}\alph{subsection}}

\nek{\parf}{\setcounter{subsection}0\chapter}
\nek{\punk}[1]{\subsection{{\protect\boldmath#1}}%
\setcounter{clt}0}

\renek{\chaptermark}[1]{\markboth{%
{\small\upshape Chapter\ \thechapter\hspace{1em}#1}}
{}%
}%

\renek{\thefigure}{\arabic{figure}}

\renek{\subsectionmark}[1]{\markright{%
{\small\upshape\thesubsection\hspace{1em}#1}}
{}%
}%

\vyk{
\chardef\eurl ='025
\chardef\eurM ='026
\nek{\fla} {{\cur{\eurl}}} 
\nek{\fmu} {{\cur{\eurM}}} 
\chardef\eurF = '011
\nek{\fFi} {{\cur{\eurF}}}
}

\newcounter{enui}

\newcounter{enuf}
\nek{\enufi}{\addtocounter{enuf}{1}}
\nek{\fenu}{
\def\theenumi{$\mtho(\fnsymbol{enuf})$}
\def\labelenumi{\theenumi}
\enufi\itsep
}
\newcounter{enuc}
\nek{\enuci}{\addtocounter{enuc}{1}}
\nek{\cenu}{
\def\theenumi{$\mtho\arabic{enuc}^\circ$}
\def\labelenumi{$\mtho\arabic{enuc}^\circ.$}
\enuci\itsep
}

\newcounter{enuF}
\nek{\enuFi}{\addtocounter{enuF}{1}}
\nek{\Fenu}{
\def\theenumi{\psur(\arabic{enuF})\psur}
\def\labelenumi{\theenumi}
\nek{\ifla}{\addtocounter{enuF}1\itla}
\enuFi\itsep
}

\nek{\itsep}{\itemsep=0.4ex plus 0.15ex minus 0.15ex}
\nek{\tenu}[1]{
\def\theenumi{#1}
\def\labelenumi{\theenumi}\itsep
}
\nek{\tenui}[1]{
\def\theenumii{#1}
\itsep
}

\theoremstyle{plain}

\newtheorem{theore}      {Theorem}       [chapter]       
\newtheorem{corollar}  [theore]{Corollary}
\newtheorem{propo}     [theore]{Proposition}
\newtheorem{lemm}      [theore]{Lemma}
\newtheorem{cla}       [theore]{Claim} 
\newtheorem{clt}       {Claim} [theore]
\newtheorem{cort}  [clt]{Corollary}  
\newtheorem{prot}  [clt]{Proposition}  
\newtheorem{lemt}  [clt]{Lemma}

\newtheorem{cly}       {{\ubf Claim}} 
\newtheorem{cory}    [cly]  {{\ubf Corollary}}  

\theoremstyle{definition}
\newtheorem{defn}      [theore]{Definition}
\newtheorem{prim}      [theore]{Example}
\newtheorem{primt}     {Example}[theore]
        \renek{\theprimt}{\thetheore:\roman{primt}}
        \nek{\sbrospri}{\addtocounter{theore}1\setcounter{primt}0}
\newtheorem{rem}       [theore]{Remark}
\newtheorem{que}       [theore] {Question}  
\newtheorem{remt}  [clt]{Remark}  
\newtheorem{vop}       [theore]{Question}
\newtheorem{upt}       [theore]{Exercise}
\newtheorem{uprt}       [clt]{Exercise}
\newtheorem{deft}  [clt]{Definition}

\newtheorem*{qrF}{{\ubf Proof}}               
\newtheorem*{prF}{{\bf Proof}}               
\newtheorem*{prC}{{\bf Proof of the claim}}  
\newtheorem*{con}{{\bfit Construction\/}{\bf.}}    
\newtheorem*{aq} {{\bf Acknowledgements\/}}          

\nek{\thsp}{\hspace{0.1ex plus \mathsurround}}
\nek{\bpro}{\begin{propo}}
\nek{\epro}{\end{propo}}
\nek{\bcl} {\begin{cla}}
\nek{\ecl} {\end{cla}}
\nek{\bct} {\begin{clt}}
\nek{\ect} {\end{clt}}
\nek{\bcy} {\begin{cly}}
\nek{\ecy} {\end{cly}}
\nek{\bco} {\begin{con}}
\nek{\eco} {\qeDD{Construction}\end{con}}
\nek{\bcor}{\begin{corollar}}
\nek{\ecor}{\end{corollar}}
\nek{\bcoy}{\begin{cory}}
\nek{\ecoy}{\end{cory}}
\nek{\baq} {\begin{aq}}
\nek{\eaq} {\end{aq}}
\nek{\bex} {\begin{prim}}
\nek{\eex} {\qeD\end{prim}}
\nek{\eeX} {\end{prim}}
\nek{\bpri} {\begin{primt}}
\nek{\epri} {\qeD\end{primt}}
\nek{\bdf} {\begin{defn}} 
\nek{\eDf} {\end{defn}}
\nek{\edf} {\qeD\end{defn}}
\nek{\edF} {\end{defn}}
\nek{\ble} {\begin{lemm}}
\nek{\ele} {\end{lemm}}
\nek{\bte} {\begin{theore}}
\nek{\ete} {\end{theore}}
\nek{\bre} {\begin{rem}} 
\nek{\ere} {\qeD\end{rem}} 
\nek{\bqe} {\begin{que}} 
\nek{\eqe} {\qeD\end{que}}
\nek{\bdt}{\begin{deft}\thsp\rm}
\nek{\edt}{\qeD\end{deft}}
\nek{\blt}{\begin{lemt}}
\nek{\elt}{\end{lemt}}
\nek{\bcot}{\begin{cort}}
\nek{\ecot}{\end{cort}}
\nek{\bprt}{\begin{prot}}
\nek{\eprt}{\end{prot}}
\nek{\bup} {\begin{upt}}
\nek{\eup} {\qeD\end{upt}}
\nek{\bqu} {\begin{vop}}
\nek{\equ} {\qeD\end{vop}}
\nek{\bret}{\begin{remt}}
\nek{\eret}{\qeD\end{remt}}
\nek{\bupt} {\begin{uprt}}
\nek{\eupt} {\qeD\end{uprt}}
\nek{\bqf} {\begin{qrF}} 
\nek{\eqf} {\,\hfill{\mtho$\vdash$}\end{qrF}} 

\nek{\bpf} {\begin{prF}} 
\nek{\epf} {\qed\end{prF}} 
\nek{\ePf} {\end{prF}} 
\nek{\bpc} {\begin{prC}} 
\nek{\epc} {\qeDD{Claim}\end{prC}} 
\nek{\qeD} {\hfill$\mtho\Box$}
\nek{\qfD} {\hfill$\mtho\vdash$}
\nek{\qeDD} [1] 
{\hfill\hbox{$\mtho\Box$~({\sl #1\/}\hspace{0.1ex})}}
\nek{\qedd} [1] 
{\hfill\hbox{$\mtho\vdash$~({\sl #1\/}\hspace{0.1ex})}}
\nek{\epF} [1] {\qeDD{#1}\end{prF}} 
\nek{\eqF} [1] {\qedd{#1}\end{qrF}} 

\nek{\bde}{\begin{description}}
\nek{\ede}{\end{description}}
\nek{\ben}{\begin{enumerate}}
\nek{\een}{\end{enumerate}}
\nek{\bit}{\begin{itemize}}
\nek{\eit}{\end{itemize}}
\nek{\bay}{\begin{array}}
\nek{\eay}{\end{array}}
\nek{\bmp}{\begin{minipage}}
\nek{\emp}{\end{minipage}}
\nek{\bus}{\begin{equation}}   
\nek{\eus}{\end{equation}}

\nek{\btb}{\begin{tabular*}}
\nek{\etb}{\end{tabular*}}
\nek{\fF} {{\bf F}}
\nek{\fG} {{\bf G}}
\nek{\Gd} {\fG_\fda}
\nek{\Gds}{\fG_{\fda\fsg}}
\nek{\Fs} {\fF_\fsg}
\nek{\fS} {{\bf S}}
\nek{\fT} {{\bf T}}
\nek{\fH} {{\bf H}}
\nek{\fa} {{\bf a}}
\nek{\fb} {{\bf b}}
\nek{\ZFC}{\text{\ubf ZFC}}
\nek{\ZC} {\text{\ubf ZC}}
\nek{\zhc}{\text{\ubf ZFC}^-}
\nek{\pli}{{\bf I}}
\nek{\plj}{{\phantom{I}\bf I}}
\nek{\pld}{{\bf II}}
\nek{\gai} [2] {\fH^\cZ_{#1#2}}
\nek{\gad} [2] {\fG^{#1}_{#2}} 
\nek{\iai} [1] {\fH_{#1}}
\nek{\iad} [1] {\fG_{#1}}
\nek{\hi} {\fH}
\nek{\hd} {\fG}
\nek{\etc} {{\sl etc}}
\nek{\iesp}{\hspace{0.3ex}}
\nek{\pw} {\hbox{a.\iesp e.}}
\nek{\ie} {\hbox{\sl i.\iesp e.}}
\nek{\eg} {\hbox{\sl e.\iesp g.}}
\nek{\ea} {\hbox{\sl et.\hspace{0.3ex}al.}}
\nek{\noo}
{\text{\sl w.\iesp l.\iesp o.\iesp g.}}
\nek{\pv} {\text{a.\iesp a.}}
\nek{\vrt} {\text{w.\iesp r.\iesp t.}}
\nek{\vva} {{\sl vice versa}}
\nek{\lsc} {\text{\sc lsc}}
\nek{\er} {{ER}}
\nek{\bm} {{BM}}
\nek{\dd}[1]{$\mtho\hspace{0.2ex}{#1}$-\hspace{0.0ex}}
\nek{\dw}{\dd\om}
\nek{\lis}[1] {\mathop{\tt lim\hspace{0.2ex}sup}_{#1}}
\nek{\len}[1] {\mathop{\tt lh}{#1}}
\nek{\Ord}  {{\tt{Ord}}}
\nek{\Exh}  {{\tt{Exh}}}
\nek{\Nul}  {{\tt{Null}}}
\nek{\Mod}  {\mathop{\tt{Mod}}}
\nek{\Aut}  {\mathop{\tt{Aut}}}
\nek{\card} {\mathop{\tt card}}
\nek{\lh}   {\mathop{\tt lh}}
\nek{\pr}   {\mathop{\tt pr}}
\nek{\sr}   {\mathop{\tt sr}}
\nek{\ran}  {\mathop{\tt ran}}
\nek{\dom}  {\mathop{\tt dom}}
\nek{\fld}  {\mathop{\tt field}}
\nek{\otp}  {\mathop{\tt otp}}
\nek{\Max}  {\mathop{\tt Max}}
\nek{\tsup} {\mathop{\tt sup}}
\nek{\tinf} {\mathop{\tt inf}}
\nek{\tmin} {\mathop{\tt min}}
\nek{\tmax} {\mathop{\tt max}}
\nek{\tlim} {\mathop{\tt lim}}
\nek{\tlis} {\mathop{\tt lim\hspace{0.3ex}sup}}
\nek{\tlii} {\mathop{\tt lim\hspace{0.3ex}inf}}
\nek{\Fin}  {{\tt Fin}} 
\nek{\bFin} {{\bf Fin}} 
\nek{\maxi}[1] {\Max^\xi_{#1}}
\nek{\HC} {{\rm HC}}
\nek{\ccc}{{\sc ccc}}
\nek{\al} {\alpha}
\nek{\ba} {\beta}
\nek{\ga} {\gamma}
\nek{\Da} {\Delta}
\nek{\da} {\delta}
\nek{\ka} {\kappa}
\nek{\la} {\lambda}
\nek{\La}{\Lambda}
\nek{\sg} {\sigma}
\nek{\Sg} {\Sigma}
\nek{\vpi}{\varphi}
\nek{\vpy}{\vpi_\iy}
\nek{\vt} {\vartheta}
\nek{\vT} {\Theta}
\nek{\ovt}{{\overline\vt}}
\nek{\ovi}{{\overline\vpi}}
\nek{\ops}{{\overline\psi}}
\nek{\ve} {\varepsilon}
\nek{\om} {\omega}
\nek{\Om} {\Omega}
\nek{\lom}{^{<\om}}
\nek{\za}{\zeta}
\nek{\tpi}{\tau_\vpi}
\nek{\omi} {\om_1}
\nek{\omm} [1] {\om^{\om^{#1}}}
\nek{\bse} {2\lom}
\nek{\nse} {\dN\lom}
\nek{\alo} {{\aleph_0}}
\nek{\sd}   {\mathbin{\Da}}
\nek{\bigd} {\mathbin{\hbox{\large$\mtho\Da$}}}
\nek{\sqe} {\fmu\hspace{0.05ex}}
\newcommand{\fs}[2]{{\bf\iSg}^{#1}_{#2}}
\newcommand{\fp}[2]{{\bf\iPi}^{#1}_{#2}}
\newcommand{\fd}[2]{{\bf\iDa}^{#1}_{#2}}
\newcommand{\iSg}{{\mathchar"7106}}
\newcommand{\iPi}{{\mathchar"7105}}
\newcommand{\iDa}{{\mathchar"7101}}
\newcommand{\is}[2]{\iSg^{#1}_{#2}}
\newcommand{\ip}[2]{\iPi^{#1}_{#2}}
\newcommand{\id}[2]{\iDa^{#1}_{#2}}
\nek{\BBB}{\hspace{0.05ex}}
\nek{\dA}{{\BBB{\mathbb A}\BBB}}
\nek{\dC}{{\BBB{\mathbb C}\BBB}}
\nek{\dF}{{\BBB{\mathbb F}\BBB}}
\nek{\dN}{{\BBB{\mathbb N}\BBB}}
\nek{\dP}{{\BBB{\mathbb P}\BBB}}
\nek{\dQ}{{\BBB{\mathbb Q}\BBB}}
\nek{\dqp}{\dQ^+}
\nek{\dR}{{\BBB{\mathbb R}\BBB}}
\nek{\dS}{{\BBB{\mathbb S}\BBB}}
\nek{\dT}{{\BBB{\mathbb T}\BBB}}
\nek{\dV}{{\BBB{\mathbb V}\BBB}}
\nek{\dZ}{{\BBB{\mathbb Z}\BBB}}
\nek{\dX}{{\BBB{\mathbb X}\BBB}}
\nek{\dY}{{\BBB{\mathbb Y}\BBB}}
\nek{\dyn} {\dY^\dN}
\nek{\dvp} {\dV^+}
\nek{\dn}{2^\dN}
\nek{\dntn}{2^{\dN\ti\dN}}
\nek{\dnqn}{{(2^\dN)}{}^\dN}
\nek{\pnqn}{{\pn}{}^\dN}
\nek{\bn}{\dN^\dN}
\nek{\rn} {\dR^\dN}
\nek{\tm} [1] {2^{\mxi\xi}}
\nek{\nn}{{\dN\ti\dN}}
\nek{\ccs} {}
\nek{\cA}{{\ccs{\skr A}\ccs}}
\nek{\cD}{{\ccs{\skr D}\ccs}}
\nek{\cE}{{\ccs{\skr E}\ccs}}
\nek{\cF}{{\ccs{\skr F}\ccs}}
\nek{\cS}{{\ccs{\skr S}\ccs}}
\nek{\cP}{{\ccs{\skr P}\ccs}}
\nek{\cW}{{\ccs{\skr W}\ccs}}
\nek{\cX}{{\ccs{\skr X}\ccs}}
\nek{\cI} {{\skr I}} 
\nek{\cJ} {{\skr J}} 
\nek{\cO} {{\skr O}} 
\nek{\cZ} {{\skr Z}}
\nek{\zo} {\cZ_0}
\nek{\zw} {\cZ_{\hbox{\small\rm w}}}
\nek{\xn} {\cX^\dN}
\nek{\an} {\dA^\dN}
\nek{\cd} [1] {\cD_{#1}}
\nek{\pws}  [1] {\cP(#1)}
\nek{\cp}  [2] {\cP_{\hspace*{-0.4ex}\tt cnt}^{#1}(#2)}
\nek{\pwf} [1] {\cP_{\hspace*{-0.4ex}\tt fin}(#1)}
\nek{\pwc} [1] {\cP_{\hspace*{-0.4ex}\tt ctbl}(#1)}
\nek{\pnn}{\cP(\nn)}
\nek{\pn}{\cP(\dN)}
\nek{\ps}{\cP(\dS)}
\nek{\pz}{\cP(\dZ)}
\nek{\mm}{{\BBB{\mathfrak M}\BBB}}
\nek{\mn}{{\BBB{\mathfrak N}\BBB}}
\nek{\gE}{{\BBB{\mathfrak E}\BBB}}
\nek{\shi} {{\mathfrak s}}
\newcommand{\gc}  {{\BBB{\mathfrak c}\BBB}}
\newcommand{\goa} {{\BBB{\mathfrak a}\BBB}}
\newcommand{\gob} {{\BBB{\mathfrak b}\BBB}}
\nek{\ilo}[1] {{[0\hspace{0.1ex},\hspace{0.1ex}n_{#1})}}
\nek{\ii} [1] {{[n_{#1}\hspace{0.1ex},\hspace{0.1ex}\infty)}}
\nek{\iv} [2] {{(#1\hspace{0.1ex},\hspace{0.1ex}#2)}}
\nek{\ix} [2] {{[#1\hspace{0.1ex},\hspace{0.1ex}#2]}}
\nek{\ir} [2] {{[#1\hspace{0.1ex},\hspace{0.1ex}#2)}}
\nek{\iry}[1] {{[#1\hspace{0.1ex},\hspace{0.1ex}\iy)}}
\nek{\iq} [2] {\ir{\nu_{#1}}{\nu_{#2}}}
\nek{\iqy}[1] {\ir{\nu_{#1}}\iy}
\nek{\iqo}[1] {\ir0{\nu_{#1}}}
\nek{\iqn}[1] {\iq{#1}{#1+1}}
\nek{\opl} {\oplus}
\nek{\ap}  {\cdot}
\nek{\cj}  {\mathbin{\hspace{0.1ex}\&\hspace{0.1ex}}}
\nek{\dm}  {$$}
\nek{\sus} {\mathopen{\exists\hspace{0.35ex}}}
\nek{\kaz} {\mathopen{\forall\hspace{0.35ex}}}
\nek{\imp} {\Longrightarrow} 
\nek{\eqv} {\Longleftrightarrow} 
\nek{\ti}  {\times} 
\nek{\mo}  {\models} 
\nek{\sq}  {\subseteq}
\nek{\qs}  {\supseteq}
\nek{\su}  {\subset}
\nek{\sneq}{\subsetneqq}
\nek{\we}  {{\mathbin{\hspace*{0.2ex}^\wedge}}}
\nek{\obr} {^{-1}}
\nek{\dif} {\smallsetminus}
\nek{\res} {\mathbin{\restriction}}
\nek{\lef} {\preccurlyeq}
\nek{\gef} {\succcurlyeq}
\nek{\pu}  {\varnothing}
\nek{\iy}  {\infty}
\nek{\piy} {+\iy}
\nek{\nin} {\not\in}
\nek{\limp}{\,\imp\,}
\nek{\leqv}{\,\eqv\,}
\nek{\onto}{\stackrel{{\rm onto}}{\longrightarrow}}
\nek{\ang} [1] {\langle #1\rangle}
\nek{\stk} [2] {\ang{#1\hspace{0.3ex};\hspace{0.1ex}#2}}
\nek{\sis} [2] {\ans{#1}_{#2}}
\nek{\ans} [1] {\{\hspace{0.01ex}#1\hspace{0.01ex}\}}

\nek{\zz} {\linebreak[0]} 
\nek{\ens} [2] {\ans{{#1\hspace{0.5ex}{:}}\zz\hspace{0.5ex}#2}}

\nek{\suh} [1] {[\hspace{0.3pt}#1\hspace{0.3pt}]_{\sq}}
\nek{\itla} {\item\label}

\nek{\aeq} {\mathbin{\|}}
\nek{\laeq}{\,\aeq\,} 

\nek{\tz} {\mathbin{;}}

\nek{\seq}[2] {(#1)_{#2}}

\nek{\kb}[2]{#1^{(#2)}}

\nek{\ssty} {\textstyle}
\nek{\isum} [3] {{{\ssty\ugl\sum_{#1}\,#3}\mid{#2}\ugr}}

\nek{\ifi}  {{\cur{Fin}\hspace*{0.1ex}}}
\nek{\frt}  {{\cur{Fr}\hspace*{0.1ex}}}
\nek{\ibo}  {{\cur{Bou}\hspace*{0.1ex}}}
\nek{\fifi} {\ifi\ti\ifi}
\nek{\fio} {\ifi\ti0}
\nek{\ofi} {0\ti\ifi}

\nek{\xip} {{\xi+1}}
\nek{\etp} {{\eta+1}}

\nek{\vx} [1] {^{(#1)}}

\nek{\dop}[1] {{#1}^{\complement}}

\nek{\skl} {\hbox{\mtho\large$($}}
\nek{\skp} {\hbox{\mtho\large$)$}}

\nek{\ugl} {\hbox{\mtho\large$\langle$}}
\nek{\ugr} {\hbox{\mtho\large$\rangle$}}

\nek{\df} [1] {\dop {#1}}
\nek{\pl} [1] {{#1}^+}

\nek{\bbW} {\hbox{\mtho\boldmath $W$}}

\nek{\bo}{{\bf O}}

\nek{\Zo}  {{\cZ_0}}
\nek{\Eo}  {\rE_{\text{\sf0}}}
\nek{\Fo}  {\rF_0}
\nek{\ovli} [1] {\widehat{#1}}
\nek{\Ec}  {\mathbin{\ovli{{\rE}}}}
\nek{\Eco} {\mathbin{\ovli{\Eo}}}
\nek{\rzo}  {\rZ_{\text{\sf0}}}
\nek{\neo} {\mathbin{\not{{\hspace{-0.4ex}\rE}}_0}}
\nek{\nfo} {\mathbin{\not{{\hspace{-0.4ex}\rF}}_0}}
\nek{\nrE} [1] {\mathbin{\not{{\hspace{-0.4ex}\rE}}_{#1}}}
\nek{\rC}  {\mathbin{\sf C}}
\nek{\rD}  {\mathbin{\sf D}}
\nek{\rR}  {\mathbin{\sf R}}
\nek{\rT}  {\mathbin{\sf T}}
\nek{\rtd} {\rT_2}
\nek{\rZ}  {\mathbin{\sf Z}}
\nek{\rdi} {\mathbin{{\sf D}_I}}
\nek{\rda} {\mathbin{{\sf D}_A}}
\nek{\rF}  {\mathbin{\sf F}}
\nek{\rE}  {\mathbin{\sf E}}
\nek{\rP}  {\mathbin{\sf P}}
\nek{\rS}  {\mathbin{\sf S}}

\nek{\epo} [1]
{\mathbin{{{#1}}{\vphantom{\pn}}^{\text{\bf+}}}}
\nek{\ei}   {\epo{\rE}}
\nek{\rfi}  {\epo{\rF}}

\nek{\nD}  {\mathbin{{\not\hspace{-0.35ex}\sf D}}}
\nek{\nE}  {\mathbin{{\not\hspace{-0.35ex}\sf E}}}
\nek{\nR}  {\mathbin{{\not\hspace{-0.35ex}\sf R}}}
\nek{\nF}  {\mathbin{{\not\hspace{-0.25ex}\sf F}}}
\nek{\reo} {\rE_{\hspace{-1.0pt}\Zo}}
\nek{\rez} {\rE_{\hspace{-1.0pt}\cZ}}
\nek{\nrz} {\mathbin
{\not{{\hspace{-0.4ex}\rE}}_{\hspace{-1.0pt}\cZ}}}
\nek{\nre} {\mathbin
{\not{{\hspace{-0.4ex}\rE}}_{\hspace{-1.0pt}\Zo}}}
\nek{\reff}{\rE_{\fifi}}
\nek{\reb} {\le_{\rm B}}
\nek{\eqb} {\approx_{\rm B}}
\nek{\dzo} {\dd{\Zo}}
\nek{\dde} {\dd{\rE}}
\nek{\ddec}{\dd{\Ec}}
\nek{\dep} {\dd{\rE'}}
\nek{\ddf} {\dd{\rF}}
\nek{\ddv} {\dd{\vt}}
\nek{\ddz} {\dd\cZ} 
\nek{\der} {\er-}
\nek{\dee} {\dd{\rE,\rE'}}
\nek{\Def} {\dd{\rE,\rF}}
\nek{\ddd} [2] {\dd{#1,#2}}

\nek{\fdo}{\hbox{\raisebox{0.2ex}{\mtho\tiny$\bullet$}}}

\nek{\fdt}{\hbox{\raisebox{-0.25ex}{\LARGE\bf.}}}
\nek{\bdot}[1] {\raisebox{-0.07ex}{\mtho$\stackrel{\fdt}{#1}$}}
\nek{\doa} {{\bdot a}} 
\nek{\dox} {{\bdot x}}
\nek{\dog} {\raisebox{-0.28ex}{\mtho$\stackrel{\fdt}g$}}
\nek{\doxl}{\dox_{\tt left}}
\nek{\doxr}{\dox_{\tt right}}
\nek{\rkb} [1] {|#1|_{\rm CB}}
\nek{\rkt} [2] {{^{\om}\hspace*{-1.3pt}|#1|_{#2}}}
\nek{\rko} [2] {{|#1|_{#2}}}
\nek{\rkT} [1] {{^{\om}\hspace*{-1.3pt}|#1|}}
\nek{\rkO} [1] {{|#1|}}
\nek{\kt}  [1] {{^{#1}\hspace*{-0.8pt}T}}
\nek{\bsf} [1] {I_{#1}}
\nek{\fri} [1] {\mathbin{\displaystyle\rE^{\tt fr}_{#1}}}
\nek{\fre} [1] {\frt_{{#1}}}

\nek{\fps} [3] {{\prod_{#3}{#2}\,/\,{#1}}}
\nek{\fss} [3] {{\sum_{#3}{#2}\,/\,{#1}}}
\nek{\fpt} [2] {{\sum{\sis{#2}{}}\,/\,{#1}}}
\nek{\fpd} [2] {{#1}\otimes{#2}}
\nek{\fsm} [2] {{#1}\oplus{#2}}
\nek{\dsu} {\fsm}

\nek{\sto} {[s_0,t_0]}
\nek{\ostp}{[s',t']}
\nek{\ost} {[s,t]}
\nek{\osu} {[s,u]} 
\nek{\otu} {[t,u]}
\nek{\ouv} {[u,v]}

\nek{\zs} {{\tilde s}}
\nek{\zt} {{\tilde t}}
\nek{\zu} {{\tilde u}}
\nek{\zv} {{\tilde v}}
\nek{\zn} {{\tilde n}}
\nek{\zm} {{\tilde m}}

\nek{\ff}[2] {F_{#1}^{#2}}

\nek{\spa} {\dS}
\nek{\spn} {\spa^\dN}
\nek{\spp} {\dY}

\nek{\poq} {\underline}
\nek{\nad} {\overline}
\nek{\nadd}{\widehat}

\nek{\sm} {\text{\mtho\large\boldmath$\fda$}}
\nek{\smp} [2] {\sm_{#1}^{#2-1}}
\nek{\smy} [1] {\sm_{#1}^{\iy}}

\nek{\sui} [1] {\cS_{\ans{#1}}}
\nek{\sun} {{\sui{1/n}}}
\nek{\srn} {{\sui{r_n}}}
\nek{\ern} {\rS_{\ans{r_n}}}
\nek{\erpn} {\rS_{\ans{r'_n}}}
\nek{\eun} {\rS_{\ans{1/n}}}
\nek{\suo} {{\sui0}}
\nek{\eso} {\rE_{\hspace{-1.0pt}\suo}}
\nek{\nso} {\mathbin
{\not{{\hspace{-0.4ex}\rE}}_{\hspace{-1.0pt}\suo}}}

\nek{\gal} [3] {{\tt Gal}^{#1}_{#2}(#3)} 

\nek{\nbd} [1] {{\skr O}_1(#1)}

\nek{\aH} {H^\ast}
\nek{\aB} {B^\ast}

\nek{\ovl} [1] {\overline{#1}} 
\nek{\ovg} [1] {\ovl{g_{#1}}}
\nek{\ovp} [1] {\ovl{\ga_{#1}}}

\nek{\plo} {+1} 

\nek{\yk} [1] {k_{#1}}

\nek{\mtho}{\mathsurround=0mm}
\nek{\msur}{\hspace{-1\mathsurround}}
\nek{\psur}{\hspace{0.3\mathsurround}}
\nek{\dsur}{\hspace{-0.3\mathsurround}}
\nek{\hsur}{\hspace{-0.5\mathsurround}}
\nek{\noi}{\noindent}
\nek{\vom}{\vspace{1mm}}
\nek{\vtm}{\vspace{2mm}}

\nek{\uv}{{\bf V}}

\nek{\wA} {{\widehat A}}
\nek{\ha} {{\hat a}}
\nek{\hb} {{\hat b}}
\nek{\he} {{\hat\ve}}
\nek{\hT} {{\hat t}}
\nek{\hl} {{\hat l}}
\nek{\hs} {{\hat s}}
\nek{\hx} {{\hat x}}
\nek{\hm} {{\hat m}}
\nek{\hn} {{\widehat n}}
\nek{\hk} {{\widehat k}}

\nek{\fo} {{\mathbf 0}}
\nek{\fr} {{\mathbf 1}}

\nek{\vnu} {{\vec \nu}} 
\nek{\dvn} {\dd\vnu} 
\nek{\wnu} {\cW_{\vnu}} 
\nek{\ewn} {\rE_{\vnu}}
\nek{\ret} {\rE_{T}}

\nek{\ci} [1] {I_{#1}}

\nek{\vex}{{\vec x}}
\nek{\vey}{{\vec y}}

\nek{\dns} [1] {d_{\vnu}^{#1}}
\nek{\den} {d_{\vnu}}
\nek{\din} {d_\vnu}

\nek{\poo}{=_{\tt df}}

\nek{\dpi}{d_\vpi}

\nek{\nol}[1] {{\bf 0}_{#1}}
\nek{\edi}[1] {{\bf 1}_{#1}}

\nek{\nrn} {_{\ans{r_n}}}
\nek{\vpr} {\vpi\nrn}
\nek{\dpr} {d\nrn}
\nek{\drn} {\dd{\ans{r_n}}}

\nek{\okr} [2] {{\skr O}_{#1}(#2)}

\nek{\nid}{\gE}

\nek{\zid} [2] {\cI_{#1/#2}}
\nek{\zidef} {\zid\rE\rF}

\nek{\zfo} [2] {\dP_{#1/#2}}
\nek{\zfoef} {\zfo\rE\rF}

\nek{\peo} {\dP_{\Eo}}
\nek{\pep} {\dP'_{\Eo}}
\nek{\ieo} {\cI_{\Eo}}

\SetMathAlphabet{\cur}{bold}{U}{eur}{b}{n}

\renek{\Gd} {\fG_\fda}
\renek{\Fs} {\fF_\fsg}

\nek{\ubf}{\fontseries{b}\selectfont}

\mathchardef\alphA ="710B
\mathchardef\betA ="710C
\mathchardef\gammA ="710D
\mathchardef\deltA ="710E
\mathchardef\vartA ="7123
\mathchardef\kpA   ="7114
\mathchardef\mU    ="7116
\mathchardef\nU    ="7117
\mathchardef\rhO   ="711A
\mathchardef\sigmA ="711B
\mathchardef\xI ="7118

\nek{\fal} {{\cur{\alphA}}}
\nek{\fba} {{\cur{\betA}}}
\nek{\fsg} {{\cur{\sigmA}}}
\nek{\fda} {{\cur{\deltA}}}

\nek{\dds}{\dd\fsg}

\nek{\nmp} {\Longleftarrow}

\nek{\tsc}[1]{\hbox{\footnotesize\sc{#1}}}

\nek{\ddi} {\dd{\cI}}
\nek{\ddj} {\dd{\cJ}}
\nek{\ddij}{\dd{(\cI,\cJ)}}

\nek{\ren}  {\le_{\tsc{c}}}
\nek{\rens} {<_{\tsc{c}}}
\nek{\eqc}  {\sim_{\tsc{c}}}

\nek{\eui}[1] {{\text{\it{EU}}}_{\ans{#1}}}

\nek{\sqa} {\sq^\ast}
\nek{\sua} {\su^\ast}

\nek{\prf} [1] {\hbox{\S\hspace{0.3ex}\ref{#1}}}
\nek{\prff}[1] {\S\S~\ref{#1}}
\nek{\pff} [1] {\S~{#1}}

\nek{\mrn} {\mu\nrn}

\renek{\dop} [1] {\complement #1}
\renek{\df}  [1] {{#1}^\complement}
\nek{\doP}  [1] {{#1}^\complement}

\nek{\ima} {\hbox{\hspace{0.1ex}''}}

\nek{\aprb}{\approx_{\tsc b}}
\nek{\ismb}{\cong_{\tsc b}}

\nek{\incs} {<_{\tsc i}}
\nek{\inc} {\le_{\tsc i}}
\nek{\eqi} {\sim_{\tsc i}}

\nek{\reas} {<_{\tsc a}}
\nek{\rea} {\le_{\tsc a}}
\nek{\eqa} {\sim_{\tsc a}}

\nek{\reaas} {<_{\tsc{aa}}}
\nek{\reaa} {\le_{\tsc{aa}}}
\nek{\eqaa} {\sim_{\tsc{aa}}}

\nek{\rebs} {<_{\tsc b}}
\renek{\reb} {\le_{\tsc b}}
\renek{\eqb} {\sim_{\tsc b}}

\nek{\reBs}{<_{\tsc{bm}}}
\nek{\reB} {\le_{\tsc{bm}}}
\nek{\eqB} {\sim_{\tsc{bm}}}

\nek{\rds} {<^{\Da}_{\tsc{b}}}
\nek{\rd} {\le^{\Da}_{\tsc{b}}}
\nek{\eqd} {\approx^{\Da}_{\tsc{b}}}

\nek{\rbas} {<_{\tsc{b,ba}}}
\nek{\rba} {\le_{\tsc{b,ba}}}
\nek{\eqba} {\approx_{\tsc{b,ba}}}

\nek{\rbasp} {<_{\tsc{b,ba}}^+}
\nek{\rbap} {\le_{\tsc{b,ba}}^+}
\nek{\eqbab} {\approx_{\tsc{b,ba}}^+}

\nek{\nab} [1] {\nabla(#1)}
\nek{\orb} {\le_{\tsc{rb}}}
\nek{\srb} {<_{\tsc{rb}}}
\nek{\erb} {\sim_{\tsc{rb}}}

\nek{\orbpp} {\le_{\tsc{rb}}^{++}}
\nek{\srbpp} {<_{\tsc{rb}}^{++}}
\nek{\erbpp} {\sim_{\tsc{rb}}^{++}}

\nek{\orbp} {\le_{\tsc{rb}}^{+}}
\nek{\srbp} {<_{\tsc{rb}}^{+}}
\nek{\erbp} {\sim_{\tsc{rb}}^{+}}

\nek{\ork} {\le_{\tsc{rk}}}
\nek{\srk} {<_{\tsc{rk}}}
\nek{\erk} {\sim_{\tsc{rk}}}

\nek{\obe} {\le_{\tsc{be}}}
\nek{\sbe} {<_{\tsc{be}}}
\nek{\ebe} {\sim_{\tsc{be}}}

\nek{\rei} {\rE_{\cI}}
\nek{\rej} {\rE_{\cJ}}

\nek{\eeb} {\sim_{\tsc{b}}}

\nek{\obep} {\le_{\tsc{be}}^+}
\nek{\sbep} {<_{\tsc{be}}^+}
\nek{\ebep} {\sim_{\tsc{be}}^+}

\nek{\seb} {<_{\tsc{b}}}

\nek{\supp} {\mathop{\tt supp}}

\nek{\atp} {\mathop{\tt at}^+}
\nek{\atm} {\mathop{\tt at}^-}

\nek{\hv} [2] {||#1||_{#2}}

\nek{\vmu} {{\vec \mu}}

\nek{\bel} [1] {\mathrel{\text{\boldmath\mtho$\ell$}}^{#1}}
\nek{\beL} [1] {\mathrel{\text{\bfsf L}}^{#1}}
\nek{\bem}     {\mathrel{\text{\bfsf m}}}
\nek{\fCo}{\mathbin{\rC_{\text{\sf0}}}}
\nek{\fco}{\mathrel{\text{\bfsf c}_{\text{\sf0}}}}
\nek{\fc} {\mathrel{\text{\bfsf c}}}
\nek{\fvt}{\mathrel{\text{\boldmath\mtho$\vt$}}}

\nek{\beli} {\bel\iy}

\nek{\oin} {{[0,1]^\dN}}

\nek{\iz} [2] {\cI_{#1}^{#2}}
\nek{\jz} [1] {\iz{}{}(#1)}
\nek{\iw} [2] {\cW_{#1}^{#2}}
\nek{\jw} [1] {\iw{}{}(#1)}
\nek{\ib} [2] {\cB^{#1}_{#2}}
\nek{\jb} [1] {\ib{}{#1}}

\nek{\ovu} {{\overline u}}

\nek{\nr}[2]{\nor{#1}_{#2}} 

\nek{\TS}{\textstyle}
\nek{\DS}{\displaystyle}

\nek{\Ba}{B^\ast}

\nek{\resi} [1] {\mathop{\restriction_{#1}}}

\nek{\gP}{{\BBB{\mathfrak P}\BBB}}
\nek{\gF}{{\BBB{\mathfrak F}\BBB}}
\nek{\gJ}{{\BBB{\mathfrak J}\BBB}}
\nek{\dG}{{\BBB{\mathbb G}\BBB}}
\nek{\dH}{{\BBB{\mathbb H}\BBB}}

\nek{\ac} {\cdot} 

\nek{\curle}{\preccurlyeq}
\nek{\cle}{\curle}
\nek{\cge}{\succcurlyeq}
\nek{\cl} {\prec}
\nek{\resic} [1] {\resi{\cle #1}}
\nek{\rec} {\resic}
\nek{\recb} [1] {\rec{(#1)}}
\nek{\rsic} [1] {\resi{\cl #1}}
\nek{\rc} {\rsic}
\nek{\rcb} [1] {\rc{(#1)}}

\nek{\kai} {\forall^\iy\hspace{0.1ex}}
\nek{\exi} {\exists^\iy\hspace{0.1ex}}

\nek{\ovw}{\nad w}

\nek{\il}[2] {\ir{n_{#1}}{n_{#2}}}
\nek{\ia}[2] {\ir{a_{#1}}{a_{#2}}}
\nek{\ij}[2] {\ir{j_{#1}}{j_{#2}}}
\nek{\cB}{{\BBB{\skr B}\BBB}}
\nek{\cG}{{\BBB{\skr G}\BBB}}
\nek{\cL}{{\BBB{\skr L}\BBB}}
\nek{\cN}{{\BBB{\skr N}\BBB}}
\nek{\cU}{{\BBB{\skr U}\BBB}}

\nek{\pnd}{\pn^\dN}
\nek{\dnd}{{(\dn)}{}^\dN}

\nek{\anp} [2] {\ang{#1}^{#2}}

\nek{\ta}{\tau}

\nek{\lev} {\mathop{\tt{lev}}}
\nek{\glu} {\mathop{\tt{dep}}}
\nek{\dia} {\mathop{\tt{diam\hspace{0.15ex}}}}

\nek{\Ei}{\rE_{\text{\sf 1}}}
\nek{\Ed}{\rE_{\text{\sf 2}}}
\nek{\Et}{\rE_{\text{\sf 3}}}
\nek{\Ey}{\rE_\iy}
\nek{\Eya} {\mathbin{(\rE_\iy)^\alo}}

\nek{\npi}{\nu_\vpi}
\nek{\nsi}{\nu_\psi}

\nek{\Ii}{\cI_1}
\nek{\Id}{\cI_2}
\nek{\It}{\cI_3}

\nek{\emb}  {\sqsubseteq_{\tsc{b}}}
\nek{\emn}  {\sqsubseteq_{\tsc{c}}}
\nek{\embi} {\sqsubseteq_{\tsc{b}}^{\rm i}}
\nek{\emni} {\sqsubseteq_{\tsc{c}}^{\rm i}}

\nek{\sio}{\cS_0}
\nek{\esn}{\rS_{\ans{1/n}}}

\nek{\nsn}{\mathbin{\not\hspace*{-0.3ex}\esn}}

\nek{\req}[2]{{\DS|_{#1}^{#2}}}

\nek{\rlq}[1]{\resi{\ge#1}}
\nek{\rmq}[1]{\resi{<#1}}
\nek{\rnq}[1]{\resi{\le#1}}

\nek{\Dij} {\dd{\cI,\cJ}}
\nek{\deo} {\dd{\Eo}}

\nek{\dc}[2] {{\bf W}^{#1}_{#2}}

\nek{\Za}{{\nad Q}}
\nek{\Ya}{{\widetilde Y}}

\nek{\Ga}{\Gamma}
\nek{\rG}  {\mathbin{\sf G}}

\nek{\inva}{{\tt{inv}}}

\nek{\tP}{{P^\ast}}
\nek{\app} {{\hspace{0.2ex}{\cdot}\hspace{0.2ex}}}

\nek{\lland}{\,\land\,}

\nek{\seqv} {\hspace{0.3ex}\Leftrightarrow\hspace{0.3ex}}
\nek{\simp} {\hspace{0.3ex}\Rightarrow\hspace{0.3ex}}

\nek{\eqn} [1] {\equiv_{#1}}
\nek{\eqo} [2] {\xrightarrow{#1,#2}}

\nek{\pit}{\tilde\pi}
\nek{\tve}{\tilde\ve}

\nek{\ef} {\dd{\rE,\rF}}

\nek{\isi} {\cong}
\nek{\isa} {\mathrel{\hspace{0.2ex}\isi^\ast}}

\nek{\ske} [3] {\equiv_{#1#2}^{#3}}
\nek{\skab}[1] {\ske AB{#1}}
\nek{\spab}[1] {\ske{A'}{B'}{#1}}

\nek{\Indent}{\hspace*{3ex}}
\nek{\fsur}{\hspace{0.5\mathsurround}}

\nek{\rdm}  {{\mathsf c}_{\tt max}}

\nek{\co} {\dd{\text{\bfsf c}_{\mathsf0}}}
\nek{\lv} {{\normalfont\scshape lv}-}

\nek{\susi} {\exists^{\iy}\hspace*{0.2ex}}
\nek{\kazi} {\forall^{\iy}\hspace*{0.2ex}}

\nek{\igi}{$\hbox{\mtho\boldmath$\fal$}$}
\nek{\igd}{$\hbox{\mtho\boldmath$\fba$}$}

\nek{\cg}[1] {{\tt Choq}(#1)}
\nek{\cgs}[1]{{\tt Choq}^{\rm s}(#1)}

\nek{\dnn}{\dN^\dN}
\nek{\nnn}{(\dnn){\vphantom{\dN}}^\dN}

\nek{\isg} {{S_\iy}}

\nek{\uset}{{Universal sets}}

\nek{\ler}[2] {\mathbin{\sim^{#1}_{#2}}}
\nek{\aer}[2] {\mathbin{\rE_{#1}^{#2}}}
\nek{\ergx}{\aer\dG\dX}
\nek{\egx} {\ergx}
\nek{\lo} [3] {\cO(#3,#1,#2)}
\nek{\sym}{\ler} 
\nek{\ong} {1_\dG}
\nek{\rr} [2] {\rR^{#1}_{#2}}

\nek{\rav}[1]
{\mathbin{\mathsf{EQ}}_{#1}}

\nek{\toq}{\circle*{0.5}}
\nek{\tob}{\circle*{1.0}}

\nek{\stob}{\circle{3.0}}
\nek{\ktob}{\kras\circle{3.0}}

\nek{\cob}{\circle{1.5}}
\nek{\mtir}{\line(-1,0){2}}

\nek{\bon} [2] {\cO_{#1}(#2)}

\nek{\dnnn} {\dnnp\dN}
\nek{\dnnp} [1] {(\dnn){\vphantom{\dN}}^{#1}}

\nek{\prift}{\sf}
\nek{\Penu}{{\prift$\id11$ Enumeration}}
\nek{\Refl}{{\prift Reflection}}
\nek{\Uset}{{\prift Universal Sets}}
\nek{\Kres}{{\prift Kreisel Selection}}
\nek{\Cenu}{{\prift Countable-to-1 Enumeration}}
\nek{\Cuni}{{\prift Countable-to-1 Uniformization}}
\nek{\Cpro}{{\prift Countable-to-1 Projection}}
\nek{\Bore}{{\prift Borel Extension}}
\nek{\Sepa}{{\prift Separation}}
\nek{\Redu}{{\prift Reduction}}

\nek{\dpf} {\mathord{{\dP^2\hspace*{-0.3ex}}\res\rF}}
\nek{\dpe} {\mathord{{\dP^2\hspace*{-0.3ex}}\res\rE}}

\nek{\ek}[2] {[#1]_{{#2}}}

\nek{\eke}[1] {\ek{#1}{\rE}}
\nek{\ekeo}[1] {\ek{#1}{\Eo}}
\nek{\ekec}[1] {\ek{#1}{\Eco}}
\nek{\ekco}[1] {\ek{#1}{\Ec}}
\nek{\ekfo}[1] {\ek{#1}{\Fo}}
\nek{\ekf}[1] {\ek{#1}{\rF}}
\nek{\ekg}[1] {\ek{#1}{G}}

\nek{\ur}{_{\tt right}}
\nek{\ul}{_{\tt left}}

\nek{\cont}{\hbox{\mtho\large$\gc$}}

\nek{\mem} {\dd\in}

\nek{\bL}{{\bf L}}
\nek{\cli} {{\sc cli}}

\nek{\di} [1] {{#1}^\#}
\nek{\drE} {\mathbin{\di\rE}}
\nek{\drF} {\mathbin{\di\rF}}

\nek{\kon} {\hbox{\mtho\large${\mathfrak c}$}}
\nek{\bk} [1] {{\cur B}_{#1}}

\nek{\wtau}{{\widehat\tau}}

\nek{\gra}[1]{\dd{#1}``grainy''}
\nek{\grap}{``grainy''}

\nek{\xE}[2] {\mathbin{\rR^{#2}_{\ge #1}}}

\nek{\moq} [1] {\Mod_{#1}}
\nek{\mox} {\moq} 
\nek{\loa} [1] {j_{#1}}
\nek{\ism} [1] {\cong_{#1}}
\nek{\izm} [2] {\cong_{#1}^{#2}}
\nek{\aut} [1] {\Aut_{#1}}

\nek{\hfn} {{\rm HF}(\dN)}
\nek{\tce} [1] {{\rm TC}_\ve(#1)}
\nek{\ihf} {\simeq_{\hfn}}

\nek{\symr}{\equiv}
\nek{\rrt} [3] {\symr_{#2#3}^{#1}}
\nek{\rrq} [5] {{#4}\symr_{#2#3}^{#1}{#5}}
\nek{\nrq} [5] {{#4}\not\symr_{#2#3}^{#1}{#5}}
\nek{\rrQ} [5] {{#4}\symr_{#2\,,\,#3}^{#1}{#5}}
\nek{\rro} [5] {\ang{#2,#4}\symr^{#1}\ang{#3,#5}}
\nek{\rrO} [1] {\symr^{#1}}

\nek{\lww} {\cL_{\omi\om}}

\nek{\forc}
{\mathrel{{|\hspace{-0.2ex}|\hspace{-1.1ex}-\hspace{-1.7ex}-}}}

\newlength{\dxii}
\nek{\fC} {{\bf C}}
\nek{\pg} {\fC_\dG}
\nek{\px} {\fC_\dX}

\nek{\gen} {gen.\ }
\nek{\hg}  {loc.\hspace{0.4ex}gen.\ }
\nek{\hgp} {loc.\hspace{0.4ex}gen.}

\nek{\rfy} {\rF^\iy}

\nek{\incl} [1] {\mathop{\text{\sc Int}}\overline{#1}}  

\nek{\PP} {pinned}
\nek{\PPP}{Pinned}

\nek{\bap}{{\bar p}}
\renek{\wtau} {{\widehat p}}

\nek{\lbr} {\linebreak[0]}

\nek{\zO} {\yo}
\nek{\zi} {\yi}
\nek{\zd} {\yd}
\nek{\zT} {\yt}
\nek{\yo} {,\linebreak[0]}
\nek{\yi} {,\linebreak[0]\hspace{0.1ex}}
\nek{\yd} {,\linebreak[0]\:}
\nek{\yt} {,\linebreak[0]\;}
\nek{\zq} {,\linebreak[0]\;\,}

\nek{\prit} [1] {[{{\rm #1}}]}
\nek{\nor} [1] {\|#1\|}

\nek{\fap} {f.\hspace{0.1ex}a.\hspace{0.1ex}p.\hspace{0.1ex}m.}

\nek{\bpr} [1] {\bpf[{{\sl#1}\/}]}

\nek{\eqr} {equivalence relation}

\nek{\fx} {{\bf x}}

\nek{\rzd} {^{\text{\tt red}}}

\nek{\srez} [2] {{\bf S}_{#1}(#2)}

\nek{\qc}  [1] {\resi{> #1}}
\nek{\qec} [1] {\resi{\ge #1}}
\nek{\rme} [1] {\resi{\le#1}}

\nek{\alex} {<_{\text{\tt alex}}}
\nek{\lex} {<_{\text{\tt lex}}}
\nek{\act} {<_{\text{\tt act}}}

\nek{\bbo} {\mathbb 0}

\nek{\fz} {{\mathbf z}}

\nek{\dln} {2\lom}

\nek{\fK} {{\bf K}}

\nek{\Ks} {\fK_\fsg}

\nek{\dva}{{\ans{0,1}}}

\nek{\rtn}{\dR^\dN}
\nek{\ntn}{\dN^\dN}
\nek{\ztn}{\dZ^\dN}

\nek{\snos} [1] {\,\footnote{\ #1}}
\nek{\snom}   {\,\footnotemark}
\nek{\snot} [1] {\footnotetext{\ #1}}

\nek{\rH}  {\mathbin{\sf H}}

\nek{\fras}[2] {\text{\footnotesize$\DS\frac{#1}{#2}$}}
\nek{\fral}[2] {\text{\large$\frac{#1}{#2}$}}

\nek{\renu}{\tenu{{\rm(\roman{enumi})}}}

\nek{\ergg}{\aer\dG\dG}

\nek{\rit} [1] {{\it#1\/}}

\nek{\lap} [1] {``#1''}

\nek{\mto} {\mapsto}

\nek{\pgcr} {\addtocounter{page}1}

\nek{\gs} [2] {g_{#1#2}}

\nek{\bie} [2]
{\mathord{{\text{\large\mtho$\prec$}}#1,#2{\text{\large\mtho$\succ$}}}}
\nek{\wa} {\tilde a}
\nek{\wb} {\tilde b}
\nek{\ws} {\bar s}
\nek{\wt} {\bar t}

\nek{\itlm}[1] {\itla{#1}\hspace*{0.0ex}\imar{#1}}

\nek{\xH} {\mathbf F}
\nek{\xh} [1] {\xH_{>#1}}
\nek{\yh} [1] {\xH_{\le#1}}
\nek{\zh} [1] {\xH_{\ge#1}}

\nek{\hh} [2] {{\DS H^{#1}_{#2}}}
\nek{\ah} [2] {{\DS \cA^{#1}_{#2}}}
\nek{\ar} [2] {{\DS A^{#1}_{#2}}}

\nek{\gru} {\hspace*{0.2ex}}

\nek{\lmt} [1] {\le_{\text{\sc nt}}^{#1}}
\nek{\emt} [1] {\rE^{\vphantom{\xi}#1}_{\text{\sc nt}}}
\nek{\mti} [1] {\cJ_{\text{\sc nt}}^{#1}}
\nek{\wmti}[1] {{\wcI}_{\text{\sc nt}}^{#1}}

\nek{\scw} {<_{\text{\tt cw}}}
\nek{\gcw} {\ge_{\text{\tt cw}}}
\nek{\lcw} {\le_{\text{\tt cw}}}
\nek{\lnt} {\le_{\text{\sc nt}}}
\nek{\ent} {\rE_{\text{\sc nt}}}
\nek{\nti} {\cJ_{\text{\sc nt}}}
\nek{\wnti}{\wcJ_{\text{\sc nt}}}

\nek{\ndi}  {\cI_{\text{\sc nt}}}
\nek{\ndie} {\cE_{\ndi}}

\nek{\bez} {\dif}

\nek{\Eqr} {Equivalence relation}

\nek{\rah} [1] {\text{\tt rnk}({#1})}
\nek{\raf} [1] {\text{\tt rnk}_{#1}}
\nek{\rag} [2] {\text{\tt rnk}_{#1}(#2)}

\nek{\pcw} {\mathbin{+_{\text{\hspace*{-0.5ex}\tt cw}}}}

\nek{\tlw} {2\lom}
\nek{\nw} {\dN^\om}
\nek{\nlw} {\dN\lom}
\nek{\tw} {2^\om}

\nek{\nt} {\text{\ubf{NT}}}

\nek{\wS} {\widehat S}
\nek{\wQ} {\widehat Q}

\nek{\ceq} [3] {{#1}_{#2#3}}
\nek{\cer} [4] {{#1}_{#2#3}^{#4}}

\nek{\rU}  {\mathbin{\sf U}}

\nek{\bai} {\cN}

\nek{\lip} [2] {\text{\sc Emb}(#1,#2)}
\nek{\wid} [1] {\text{\sc Wid}(#1)}
\nek{\vid} [1] {\text{\sc Wid}(#1)}

\nek{\lmto}{\longmapsto}

\nek{\atlh} {\addtolength{\itemsep}{-\dxh}}
\nek{\atli} {\addtolength{\itemsep}{-\dxi}}

\nek{\cT} {{\skr T}}

\newlength{\dxh}
\nek{\cip} {{\sc cip}} 

\nek{\uap} [2] {A^{#1}(#2)}
\nek{\uaq} [3] {A^{#1}_{#3}(#2)}
\nek{\use} [1] {A(#1)}

\nek{\imb} {\text{\hspace*{0.3ex}''\hspace*{-0.2ex}}}

\theoremstyle{definition}
\newtheorem{coj}       [theore]{Conjecture} 
\nek{\bcj} {\begin{coj}}
\nek{\ecj} {\end{coj}}

\nek{\api}{\addtocounter{page}{1}}

\nek{\cntr} [1] {{#1}^{\downarrow}}

\nek{\gam} [3] {\text{\ubf G}_{#1}^{#2}(#3)}

\nek{\ts} {\widehat s}
\nek{\tsg}{\widehat \sg}

\nek{\urav}  [1] {\bus\itsep#1\eus} 

\nek{\efp}{\dd{\rE,{\rF'}}}

\nek{\dfp} {\dd{\rF'}}

\nek{\inw} [2] {\dd{(#1\to#2)}invariant}

\nek{\come} {co-meager\ }

\nek{\bfei} {{\ubf either}}
\nek{\bfor} {{\ubf or}}

\nek{\gla} {Chapter}
\nek{\paf} {Section}

\nek{\grf} [1] {\gla~\ref{#1}}
\nek{\nrf} [1] {\paf~\hbox{\ref{#1}}}

\nek{\dPP} {\dP\ti\dP}
\nek{\dSP} {\dS\ti\dP}
\nek{\xib} {{\bma\xI\hspace{0.1ex}}}
\nek{\sgb} {\bma\sigmA}
\nek{\ddp} {\dd\dP}
\nek{\coll} {\text{\sc Coll}}
\nek{\tal}{\xib\ul}
\nek{\tar}{\xib\ur}
\nek{\dCP} {\dC\ti\dP}
\nek{\dCC} {\dC\ti\dC}
\nek{\sgbl}{\sgb\ul}
\nek{\sgbr}{\sgb\ur}

\nek{\uoa} {\mathbin{{\text{\bfsf u}}^\ast_0}}
\nek{\uo}  {\mathbin{{\text{\bfsf u}}_0}}

\nek{\ld} [1] {\mathbin{{\rD}(#1)}}

\nek{\xd} [2] {{\mathrel{{\td}(#1\hspace*{0.1ex};\hspace*{0.1ex}#2)}}}
\nek{\qd} [3] {{\mathrel{{\td}
(\ang{#1\hspace*{0.1ex};\hspace*{0.1ex}#2}_{#3\in\dN})}}}
\nek{\qqd} [3] {{\mathrel{{\td}
(\ang{#1\hspace*{0.1ex};\hspace*{0.1ex}#2}_{#3})}}}

\nek{\td}  {\mathrel{\rD}}
\nek{\ntd} {\mathrel{{\not\hspace{-0.0ex}\td}}}
\nek{\mast} {\text{\mtho$\ast$}}

\nek{\ja} {\text{\mtho\boldmath$a$}}

\renek{\api}{}

\makeatletter
\if@twoside%
\else\fi%
\makeatother
         
\begin{document}

\frontmatter

\title{{\sc Varia}\\  Ideals and Equivalence Relations,
beta-version}

\author{Vladimir Kanovei}

\date{\today}
\maketitle
                                                                  
\chapter*{Introduction}
\addcontentsline{toc}{chapter}{Introduction}

We present a selection of basic results on Borel reducibility
of ideals and \eqr s, especially those with comparably short
proofs.
The focal point are reducibility/irreducibility results
related to some special equivalences like
$\Eo\yd\Ei\yd\Ed\yd\Et\yd\Ey\yd\Zo,$ and Banach-induced
equivalences $\bel p.$
The bulk of results included in the book were obtained
in the 1990s, but some rather recent theorems are presented
as well, like Rosendal's proof that Borel ideals are cofinal
within Borel equivalences of general form. 
\snos
{{\tt kanovei@mccme.ru} and {\tt vkanovei@math.uni-wuppertal.de}
are my contact addresses.}

\newpage

General set-theoretic notation used.
\bit
\item\msur
$\dN=\ans{0,1,2,\dots}:$ natural numbers. 
$\dN^2=\dN\ti\dN.$ 

\item\msur
$\dnn$ is {\it the Baire space\/}. 
If $s\in\dN\lom$ (a finite sequence of natural numbers) 
then $\bon s\dnn=\ens{x\in\dnn}{s\su x},$ a basic clopen 
nbhd in $\dnn$.
\index{zzN@$\dN$}%
\index{zzons@$\bon s\dnn$}%

\item\msur
$X\sqa Y$ means that the difference $X\dif Y$ is finite.
\index{zzzsqa@$\sqa$}

\item
If a basic set $A$ is fixed then $\dop X=\doP X=A\dif X$ 
for any $X\sq A$.
\index{zzcX@$\dop X$}

\item
If $X\sq A\ti B$ and $a\in A$ then 
$\seq Xa=\ens{b}{\ang{a,b}\in X},$ 
a {\it cross-section\/}.
\index{cross-section}%
\index{zzXa@$\seq Xa$}

\item\msur
$\#X=\#(X)$ is the number of elements of a finite set $X$.
\index{zzznx@$\#(X)$}%
\index{number of elements}

\item\msur
$f\ima X=\ens{f(x)}{x\in X\cap\dom f},$ 
the \dd f{\it image\/} of $X$.  
\index{image}%
\index{zzfiX@$f\ima X$}

\item\msur
$\sd$ is the symmetric difference. 
\index{symmetric difference}%
\index{zzda@$\sd$}

\item\msur
$\susi x\;\dots\;$ means: 
``there exist infinitely many $x$ such that \dots'', \\[1mm]
$\kazi x\;\dots\;$ means: 
``for all but finitely many $x,$ \dots holds''.
\index{zzzex@$\susi x$}%
\index{zzzka@$\kazi x$}
\index{quantifiers!zzzex@$\susi x$}%
\index{quantifiers!zzzka@$\kazi x$}%

\item\msur
$\sis{x_a}{a\in A}$ is the map $f$ defined on $A$ by
$f(a)=x_a\zt\kaz a$.

\item\msur
$\pws X=\ens{x}{x\sq X}$.
\index{zzPX@$\pws X,$ the power set}

\item\msur
$\pwf X=\ens{x}{x\sq X\,\text{ is finite}}$.
\index{zzPfinX@$\pwf X,$ all finite subsets}
\eit

\newpage

{\pagestyle{empty}\small
\def\contentsname{\large\bf Contents}
\tableofcontents
}

\mainmatter

\api

\vyk{\appendix

\renek{\thesection}{\Alph{section}}
\renek{\thesubsection}{\thesection\arabic{subsection}}

\parf{Appendix. Descriptive set theory}
}

\parf{Descriptive set theoretic background}
\label i

We assume that the reader of this book has a basic 
knowledge of descriptive set theory,
both classical and effective,
in Polish spaces 
(recursively presented, in the effective case),
including Borel and projective hierarchy, 
Borel sets and functions, analytic and coanalytic sets,  
and the like.  

A map $f$ 
(between Borel sets in Polish spaces) is 
{\it Borel\/} iff its graph is a Borel set iff all 
\dd fpreimages of open sets are Borel. 
A map $f$ is {\it Baire measurable\/} 
({\it BM\/}, for brevity) iff all 
\dd fpreimages of open sets are Baire measurable sets. 
\index{map!Borel}%
\index{map!Baire measurable, BM}


Apart of the very common knowledge, the whole 
instrumentarium of ``effective'' descriptive set theory 
employed in the study of reducibility of ideals and 
\er s, can be summarized in a rather short list  
of key ``principles''.
In those below, by a {\it recursively presented\/} Polish
\index{space!recursively presented}%
\index{recursively presented (space)}%
space one can understand any product space of the form
$\dN^m\ti(\dnn)^n$ without any harm for applications
below, yet in fact this notion is much wider.

\bre
\label{relrel}
For the sake of brevity, the results below are formulated 
only for the ``lightface'' parameter-free classes
$\is11\zd \ip11\zd \id11,$ 
but they remain true for $\is11(p)\zd \ip11(p)\zd \id11(p)$
for any fixed real parameter $p,$ as well as for the
\lap{boldface} classes $\fs11\zd \fp11\zd \fd11$ of resp.\
analytic, coanalytic, Borel sets.
\ere

\bte
[\Redu, \Sepa]
\lam{rese}
\index{separation@\Sepa}%
\index{principle!separation@\Sepa}%
\index{reduction@\Redu}%
\index{principle!reduction@\Redu}%
If\/ $X\yi Y$ are\/ $\ip11$ sets
(of a recursively presented Polish space)
then there exist\/ {\ubf disjoint} $\ip11$ 
sets\/ $X'\sq X$ and\/ $Y'\sq Y$ with\/ $X'\cup Y'=X\cup Y.$
{\rm The sets\/ $X'\zi Y'$ are said to {\it reduce\/} the pair 
\index{reduce}%
$X\zi Y.$}

If\/ $X\yi Y$ are disjoint\/ $\is11$ sets
then there is a\/ $\id11$ set\/ $Z$ 
with\/ $X\sq Z$ and\/ $Y\cap Z=\pu$.
\rm
Such a set\/ $Z$ is said to {\it separate\/} $X$ from 
\index{separate}%
$Y$.\qed
\ete

\bte
[\Cpro]
\lam{cpro}
\index{countablep@\Cpro}%
\index{principle!countablep@\Cpro}%
If\/ $P$ is a\/ $\id11$ subset of the product\/ $\dX\ti\dY$ 
of two recursively presented Polish spaces, and for any\/ 
$x\in\dX$ the cross-section\/ $P_x=\ens{y}{P(x,y)}$ is at most 
countable then\/ $\dom P$ is a\/ $\id11$ set.
\qed
\ete
 
It follows that images of $\id11$ sets via 
countable-to-1, in particular, 1-to-1 $\id11$ maps are 
$\id11$ sets, while images via arbitrary $\id11$ maps 
are, generally, $\is11$.

\bte
[\Cenu]
\lam{cenu}
\index{countablee@\Cenu}%
\index{principle!countablee@\Cenu}%
If\/ $P\yi\dX\yi\dY$ are as in Theorem~\ref{cpro} 
then there is a\/ $\id11$ map\/ $f:\dom P\ti\dN\to\dY$
such that\/ $P_x=\ens{f(x,n)}{n\in\dN}$ for all\/
$x\in\dom P$.
\qed
\ete
 
\bte
[\Bore]
\lam{bore}
\index{borele@\Bore}%
\index{principle!borele@\Bore}%
If\/ $P$ is a\/ $\is11$ subset of the product\/ $\dX\ti\dY$ 
of two recursively presented Polish spaces, and for any\/ 
$x\in\dX$ the cross-section\/ $P_x=\ens{y}{P(x,y)}$ is at most 
countable then there is a\/ $\id11$ superset\/ $Q\supseteq P$ with  
all cross-sections\/ $Q_x$ at most countable.
Similarly, if\/ $P\sq\dX\ti\dY$ is a uniform\/ $\is11$ set then
there is a uniform\/ $\id11$ superset\/ $Q\supseteq P.$
\qed
\ete

Recall that a set $P\sq\dX\ti\dY$ is \rit{uniform} iff the
cross-section $P_x$ contains at most one point for any $x\in\dX.$
This is the same as a \rit{partial function} $\dX\to \dY$.

\bte
[\Cuni]
\lam{cuni}
\index{countableu@\Cuni}%
\index{principle!countableu@\Cuni}%
If\/ $P\yi\dX\yi\dY$ are as in Theorem~\ref{cpro} 
then\/ $P$ can be uniformized by a\/ $\id11$ set.
\qed
\ete

\bte
[\Kres]
\lam{kres}
\index{kreisel@\Kres}%
\index{principle!kreisel@\Kres}%
If\/ $\dX$ is a recursively presented Polish space, 
$P\sq\dX\ti\dN$ is a\/ $\ip11$ set, and\/ $X\sq\dom P$ 
is a\/ $\id11$ set then there is a\/ $\id11$ function 
$f:X\to\dN$ such that\/ $\ang{x,f(x)}\in P$ for all\/ 
$x\in X$.
\qed
\ete
\bpf
Let $Q\sq P$ be a $\ip11$ set which uniformizes $P.$ 
For any $x\in X$ let $f(x)$ be the only $n$ with 
$\ang{x,n}\in Q.$ 
Immediately, (the graph of) $f$ is $\ip11,$ however, 
as $\ran f\sq\dN,$ we have 
$f(x)=n\eqv \kaz m\ne n\:(f(x)\ne m)$ whenever 
$x\in X,$ which demonstrates that $f$ is $\is11$ as well.
\epf

The next theorem provides a useful enumeration of $\id11$ sets.

\bte
[\Penu]
\lam{penu}
\index{denumera@\Penu}%
\index{principle!denumera@\Penu}%
If\/ $\dX$ is a recursively presented Polish space then 
there exist\/ $\ip11$ sets\/ $C\sq\dN$ and\/ 
$W\sq\dN\ti\dX$ and a\/ $\is11$ set\/ $W'\sq\dN\ti\dX$ such that 
$W_e=W'_e$ for all\/ $e\in C,$ and a set\/ $X\sq\dX$ is\/
$\id11$ iff there is\/ $e\in C$ such that\/ $X=W_e=W'_e.$ \ 
{\rm(Here $W_e=\ens{x}{W(e,x)}$ and similarly $W'_e.$)}\qed
\ete
%
There is a generalization useful for relativised classes 
$\id11(y)$.

\bte
[\sf Relativized \Penu] \ 
\lam{penur}
\index{denumera@\Penu!relativized}%
If\/ $\dX\yi \dY$ are recursively presented Polish spa\-ces 
then there exist\/ $\ip11$ sets\/ $C\sq\dY\ti\dN$ and\/ 
$W\sq\dY\ti\dN\ti\dX$ and a\/ $\is11$ set\/ 
$W'\sq\dY\ti\dN\ti\dX$ such that\/
$W_{ye}=W'_{ye}$ for all\/ $\ang{y,e}\in C$ and, 
for any\/ $y\in\dY,$ a set 
$X\sq\dX$ is\/ $\id11(y)$ iff there is\/  
$e$ such that\/ $\ang{y,e}\in C$ and\/ $X=W_{ye}=W'_{ye}.$ \ 
{\rm(Here $W_{ye}=\ens{x}{W(y,e,x)}$ and similarly $W'_{ye}$.)}
\qed
\ete
%

Suppose that $\dX$ is a recursively presented Polish space.
A set $U\sq\dN\ti\dX,$ 
is a {\it a universal\/ $\ip11$ set\/} if   
\index{set!universal}%
for any $\ip11$ set $X\sq\dX$ there is an 
index $e\in\dN$ with $X=U_e=\ens{x}{\ang{e,x}\in U},$
and a {\it a ``good'' universal\/ $\ip11$ set\/} if in
\index{set!universal!good}%
addition for any other $\ip11$ set $V\sq\dN\ti\dX$
there is a
recursive function $f:\dN\to\dN$ such that
$V_e=U_{f(e)}$ for all $e$.

The notions of universal and ``good'' universal $\is11$
sets are similar.

\bte
[\Uset]
\lam{uset}
\index{universalsets@\Uset}%
\index{principle!universalsets@\Uset}%
For any recursively presented Polish space\/ $\dX$
there exist a ``good'' universal\/ $\ip11$ set
$U\sq\dN\ti\dX$ and a ``good'' universal\/ $\is11$ set
$V\sq\dN\ti\dX.$ 
{\rm(In fact we can take $V=(\dN\ti\dX)\dif U$.)}\qed
\ete

If a ``good'' universal $\ip11$ set $U$ is fixed then
a collection $\cA$ of $\ip11$ sets $X\sq\dX$ is 
{\it $\ip11$ in the codes\/}
\index{in the codes}%
if $\ens{e}{U_e\in\cA}$ is a $\ip11$ set. 
Similarly, if a ``good'' universal $\is11$ set $V$
is fixed then a collection $\cA$ of $\is11$ sets
$X\sq\dX$ is {\it $\ip11$ in the codes\/} if
$\ens{e}{V_e\in\cA}$ is a $\ip11$ set.
These notions quite obviously do not depend on the
choice of ``good'' universal sets.

To show how ``good'' universal sets work, we prove:

\bpro
\label{effred}
Let\/ $\dX$ be a recursively presented Polish space and\/
$U\sq\dN\ti\dX$ a ``good'' universal\/ $\ip11$ set.
Then for any pair of\/ $\ip11$ sets\/ $V,W\sq\dN\ti\dX$
there are recursive functions\/ $f,g:\dN\to\dN$ such
that for any\/ $m,n\in\dN$ the pair of cross-sections\/
$U_{f(m,n)}\zi U_{g(m,n)}$ reduces the pair\/  
$V_m\zi W_n$.
\epro
\bpf
Consider the following $\ip11$ sets in $(\dN\ti\dN)\ti\dX$:
\dm
P=\ens{\ang{m,n,x}}{\ang{m,x}\in V\land n\in\dN},\;\,
Q=\ens{\ang{m,n,x}}{\ang{n,x}\in W\land m\in\dN}.
\dm
By \Redu, there is a pair of $\ip11$ sets $P'\sq P$
and $Q'\sq Q$ which reduce the given pair $P\zi Q.$
Accordingly, the pair $P'_{mn}\zi Q'_{mn}$ reduces
$P_{mn}\zi Q_{mn}$ for any $m,n.$
Finally, by the ``good'' universality there are
recursive functions $f,g$ such that
$P'_{mn}=U_{f(m,n)}$ and $Q'_{mn}=U_{g(m,n)}$ 
for all $m,n$. 
\epf

The following theorem is less 
elementary than the results cited above, but it is 
very useful because it allows to \lap{compress} some  
sophisticated arguments with multiple applications of 
Separation and Kreisel selection.

\bte
[\Refl]
\lam{refl}
\index{reflection@\Refl}%
\index{principle!reflection@\Refl}%
Let\/ $\dX$ be a recursively presented Polish 
space. 
\bde
\item[$\ip11$ {\ubf form.}]
Suppose that a collection $\cA$ of $\ip11$ sets $X\sq\dX$  
is $\ip11$ in the codes.
(In the sense of a fixed ``good'' 
universal $\ip11$ set $U\sq\dN\ti\dX.$)  
Then for any $X\in\cA$ there is a $\id11$ set 
$Y\in\cA$ with $Y\sq X$.

\item[$\is11$ {\ubf form.}]
Suppose that a collection $\cA$ of $\is11$ sets $X\sq\dX$   
is $\ip11$ in the codes.
Then for any $X\in\cA$ there is a $\id11$ set 
$Y\in\cA$ with $X\sq Y$.\qed
\ede
\ete

One of (generally, irrelevant here) consequences of 
this principle is that the set of all codes of a 
properly $\ip11$ set or properly $\is11$ set 
is never $\ip11$.

\api

\parf{Borel ideals}
\las{n-id}

This \gla\ does not mean any broad introduction into
Borel ideals; we rather consider some issues close to the content
of the book, including P-ideals, polishable ideals,
\lsc\ submeasures, summable, density, and Fr\'echet ideals,
and Rudin -- Blass reducibility of ideals.

\punk{Introduction to Borel ideals}
\las{borid}

Recall that an \rit{ideal} on a set $A$ is any non-empty set 
$\cI\sq\pws A$ closed under $\cup$ and 
satisfying ${x\in\cI}\imp{y\in\cI}$ whenever 
$y\sq x\sq A.$ 
Thus, any ideal contains the empty set $\pu.$ 
Usually they consider only {\it nontrivial\/} ideals, 
\ie, those which contain all 
singletons $\ans{a}\yt a\in A,$ and do not contain $A,$ \ie, 
$\pwf A\sq\cI\sneq\pws A$.
\index{ideal}%
\index{ideal!nontrivial}%
\index{ideal!pu@$0=\ans\pu$}%
But sometimes the ideal $\ans\pu,$ whose only element is
the empty set $\pu$ is considered and often denoted by $0$.

If $A$ is a countable set then identifying $\pws A$ with $2^A$
via characterictic functions we equip $\pws A$ with the Polish
product topology.
In this sense, a \rit{Borel} ideal on $A$ is any ideal
\index{ideal!nontrivial}%
which is a Borel subset of $\pws A$ in this topology. 
Let us give several important examples of Borel ideals.
\bit
\item
$\ifi=\ens{x\sq\dN}{x\,\text{ is finite}},$ the 
ideal of all finite sets;

\item
$\Ii=
\ens{x\sq\dN^2}{\ens{k}{\seq xk\ne\pu}\in\ifi},$
where $\seq xa=\ens{b}{\ang{a,b}\in x}$;

\item
$\Id=
\ens{x\sq \dN}{\sum_{n\in x}\frac1{n+1}}<\piy,$ 
the {\it summable ideal\/};

\item
$\It=
\ens{x\sq\dN^2}{\kaz k\:(\seq xk\in\ifi)}$;

\item
$\zo=
\ens{x\sq\dN}
{\tlim_{n\to\piy}\frac{\#(x\cap\ir0n)}n=0},$
\index{zzzo@$\zo$}
\index{ideal!zzzo@$\zo$}
the {\it density ideal\/}.
\eit

For any ideal $\cI$ on a set $A,$ we define 
$\pl\cI=\cP(A)\dif\cI$ (\ddi{\it positive\/} sets) and 
$\df\cI=\ens{X}{\dop X\in\cI}$ ({\it the dual filter\/}). 
Clearly $\pu\ne\df\cI\sq\pl\cI.$ 

If $B\sq A,$ then we put  
$\cI\res B=\ens{x\cap B}{x\in\cI}$.
\index{ideal!restriction $\cI\res B$}%
\index{zzIB@$\cI\res B$}%

\punk{P-ideals, submeasures, polishable ideals}
\las{pideals}

Many important Borel ideals belong to the class of 
P-ideals. 

\bdf
\label{api}
An ideal $\cI$ on $\dN$ is 
a {\it P-ideal\/} if for any sequence 
of sets $x_n\in\cI$ there is a set $x\in\cI$ such that 
$x_n\sqa x$ (\ie, $x_n\dif x\in\ifi$) for all $n$;
\edf

For instance, the ideals $\ifi,\,\Id,\,\It,\,\Zo$ 
(but not $\Ii$!) are P-ideals. 

This class admits several apparently different 
but equivalent characterizations, one of which is 
connected with submeasures. 

\bdf
\lam{df:subm}
A {\it submeasure\/} on a set $A$ is any map 
\index{submeasure}%
$\vpi:\cP(A)\to[0,\piy],$ satisfying $\vpi(\pu)=0,$ 
$\vpi(\ans a)<\piy$ for all $a,$
and $\vpi(x)\le \vpi(x\cup y)\le \vpi(x)+\vpi(y)$.

A submeasure $\vpi$ on $\dN$ is 
{\it lover semicontinuous\/}, 
\index{submeasure!lover semicontinuous}%
\index{submeasure!lsc@\lsc}%
or \lsc\ for brevity, if we have 
$\vpi(x)=\tsup_n\vpi(x\cap\ir0n)$ for all $x\in\pn$. 
\edf
 
To be a {\it measure\/}, a submeasure $\vpi$ has to 
satisfy, in addition, that 
$\vpi(x\cup y)=\vpi(x)+\vpi(y)$ whenever $x,\,y$ 
are disjoint. 
Note that any \dds additive measure is \lsc, 
but if $\vpi$ is \lsc\ then $\vpy$ is not 
necessarily \lsc\ itself.

Suppose that $\vpi$ is a submeasure on $\dN.$ 
Define the {\it tailsubmeasure\/} 
\index{tailsubmeasure}%
\index{submeasure!tailsubmeasure}%
$\vpy(x)=\hv x\vpi=\tinf_n(\vpi(x\cap\iry n)).$ 
The following ideals are considered:
\dm
\bay{rcllll}
\Fin_\vpi &=& \ens{x\in\pn}{\vpi(x)<\piy} &&&;\\[1ex]
\index{zzfinf@$\Fin_\vpi$}%

\Nul_\vpi &=& \ens{x\in\pn}{\vpi(x)=0} &&&;\\[1ex]
\index{zznulf@$\Nul_\vpi$}%

\index{zzexhf@$\Exh_\vpi$}%
\Exh_\vpi &=& \ens{x\in\pn}{\vpy(x)=0}&=
& \Nul_{\vpy}&.\\[0.5ex]
\eay
\dm

\bex
\label{exh:e}
$\ifi=\Exh_\vpi=\Nul_\vpi,$ where $\vpi(x)=1$ for any 
$x\ne \pu.$ 
We also have $\ofi=\Exh_\psi,$ where 
$\psi(x)=\sum_k\,2^{-k}\,\vpi(\ens{l}{\ang{k,l}\in x})$ 
is \lsc.
\eex

It turns out (Solecki, see Theorem~\ref{sol} below) that 
analytic P-ideals are the same as ideals of the form 
$\Exh_\vpi,$ where $\vpi$ is a 
\lsc\ submeasure on $\dN.$ 
This implies that any analytic P-ideal is $\fp03$.


There is one more useful characterization of Borel P-ideals. 
Let $T$ be the ordinary Polish product topology on $\pn.$ 
Then $\pn$ is a Polish group in the sense of $T$ and the 
symmetric difference as the operation, and any ideal 
$\cI$ on $\dN$ is a subgroup of $\pn$. 

\bdf
\lam{polI}
An ideal $\cI$ on $\dN$ is {\it polishable\/} if 
there is a Polish group topology $\tau$ on $\cI$ which 
produces the same Borel subsets of $\cI$ as $T\res\cI$.
\edf

The same Solecki's theorem (Theorem~\ref{sol}) proves 
that, for analytic ideals, to be a P-ideal is the same 
as to be polishable. 
It follows (see Example~\ref{exh:e}) that, for instance, 
$\ifi$ and $\It$ are polishable, but $\Ii$ 
is not. 
The latter will be shown directly after the next lemma.

\ble
\lam{sol:?}
Suppose that an ideal\/ $\cI\sq\pn$ is polishable. 
Then there is a unique Polish group topology\/ 
$\tau$ on\/ $\cI.$ 
This topology\/ refines\/ $T\res\cI$ 
and is metrizable by a\/ \dd\sd invariant metric. 
If\/ $Z\in\cI$ then\/ $\tau\res\cP(Z)$ coincides 
with\/ $T\res\cP(Z).$ 
In addition, $\cI$ itself is\/ \dd TBorel.
\ele
\bpf
Let $\tau$ witness that $\cI$ is polishable. 
The identity map 
$f(x)=x{:}\,\stk\cI\tau\to\stk\pn T$ 
is a \dd\sd homomorphism and is Borel-measurable 
because all \dd{(T\res\cI)}open sets are \dd\tau Borel, 
hence, by the Pettis theorem (Kechris~\cite[??]{dst}), 
$f$ is continuous. 
It follows that all \dd{(T\res\cI)}open subsets of 
$\cI$ are \dd\tau open, and that $\cI$ is 
\dd TBorel in $\pn$ because $1-1$ continuous 
images of Borel sets are Borel. 

A similar ``identity map'' argument shows that 
$\tau$ is unique if exists. 

It is known (Kechris~\cite[??]{dst}) 
that any Polish group topology admits a 
left-invariant compatible metric, which, in this case, 
is right-invariant as well since $\sd$ is an 
abelian operation.

Let $Z\in\pn.$ 
Then $\cP(Z)$ is \dd Tclosed, hence, 
\dd\tau closed by the above, subgroup of $\cI,$ and 
$\tau\res\cP(Z)$ is a Polish group topology on $\cP(Z).$ 
Yet $T\res\cP(Z)$ is another Polish group topology on 
$\cP(Z),$ with the same Borel sets. 
The same ``identity map'' argument proves  
that $T$ and $\tau$ coincide on $\cP(Z)$.
\epf

\bex
\lam{fio}
The ideal $\Ii$ is not polishable. 
Indeed we have $\Ii=\bigcup_n W_n,$ where  
$W_n=\ens{x}{x\sq\ans{0,1,\dots,n}\ti\dN}.$  
Let, on the contrary, $\tau$ be a Polish group topology 
on $\Ii.$ 
Then $\tau$ and the ordinary topology $T$ coincide on 
each set $W_n$ by the lemma, in particular, 
each $W_n$ remains \dd\tau nowhere dense 
in $W_{n+1},$ hence, in $\Ii,$ a contradiction with 
the Baire category theorem for $\tau$.
\eex

\punk{Summable and density ideals}
\las{sumden}

Any sequence $\sis{r_n}{n\in\dN}$ of positive reals $r_n$ 
with $\sum r_n=\piy$ defines the ideal 
\dm
\sui{r_n}=\ens{X\sq\dN}{\sum_{n\in X}r_n<\piy}=
\ens{X}{\mrn(X)<\piy}\,,
\dm 
where $\mrn(X)=\sum_{n\in X}r_n.$
These ideals are called {\it summable ideals\/};
all of them are $\Fs$ in the product Polish topology on $\pn.$
Any summable ideal is easily a P-ideal:  
indeed, $\sui{r_n}=\Exh_\vpi,$ where 
$\vpi(X)=\sum_{n\in X}r_n$ is a \dds additive measure. 
Summable ideals are perhaps the easiest to study among 
all P-ideals.
More on summable ideals see \cite{mat72,maz91,aq}. 

\vyk{
Further entries: 
1) Farah \cite[\pff1.12]{aq} on summable 
ideals under $\obe,$ 
2) Hjorth: \dd\reb structure 
of ideals \dd\reb reducible to summable ideals, in 
\cite{h-ban}.
}

Farah \cite[\pff 1.10]{aq} defines a non-summable $\Fs$ 
P-ideal as follows. 
Let $I_k=\ir{2^k}{2^{k+1}}$ and 
$\psi_k(s)=k^{-2}\tmin\ans{k,\# s}$ for all $k$ and 
$s\sq I_k,$ and then
\dm
\psi(X)=\sum_{k=0}^\iy\psi_k(X\cap I_k)
\quad\hbox{and}\quad
\cI=\Fin_\psi\;;
\dm
it turns out that $\cI$ is an $\Fs$ P-ideal, but not 
summable. 
To show that $\cI$ distincts from any $\sui{r_n},$ 
Farah notes that there is a set $X$ 
(which depends on $\sis{r_n}{}$) 
such that the differences 
$|\mrn(X\cap I_k)-\psi_k(X\cap I_k)|,\msur$ 
$k=0,1,2,\dots\,,$ 
are unbounded. 

There exist other important types of Borel P-ideals.
Any sequence $\sis{r_n}{n\in\dN}$ of positive reals $r_n$ 
with $\sum r_n=\piy$ defines the ideal 
\dm
\eui{r_n}=\left\{x\sq\dN:\tlim_{n\to\piy}
\frac{\sum_{i\in x\cap\ir0n}r_i}
{\sum_{i\in \ir0n}r_i}=0\right\}.
\dm 
These ideals are called 
{\it Erd\"os -- Ulam\/ {\rm(or: EU)} ideals\/}. 
Examples: $\zo=\eui1$ and $\cZ_{\tt log}=\eui{1/n}$. 

This definition can be generalized. 
Let $\supp\mu=\ens{n}{\mu(\ans n)>0},$ for any measure 
$\mu$ on $\dN.$ 
Measures $\mu,\,\nu$ are {\it orthogonal\/} if we have 
$\supp\mu\cap\supp\nu=\pu.$  
Now suppose that $\vmu=\sis{\mu_n}{n\in\dN}$ is a sequence 
of pairwise orthogonal measures on $\dN,$ with finite 
sets $\supp\mu_i.$ 
Define $\vpi_\vmu(X)=\tsup_n\mu_n(X):$ this is a \lsc\ 
submeasure on $\dN.$ 
Let finally 
$\cD_\vmu=\Exh(\vpi_\vmu)=\ens{X}{\hv X{\vpi_\mu}=0}.$ 
Ideals of this form are called {\it density ideals\/} 
by Farah~\cite[\pff 1.13]{aq}. 
This class includes all EU ideals 
(although this is not immediately transparent), 
and some other ideals: for instance, $\It$ is a 
density but non-EU ideal.  
Generally density ideals are more complicated than summables. 
We obtain an even wider class if the requirement, 
that the sets $\supp\mu_n$ are finite, is dropped: 
this wider family includes all summmable ideals, too.

See \cite{jk84}, \cite[\pff1.13]{aq} on density ideals. 

\vyk{
Further entries: 
1) Farah: structure of density ideals under $\obe,$ 
2) Farah: \dd{c_0}equalities, 
3) Relation to Banach spaces: Hjorth, SuGao. 

Which ideals are both summable and density ? 
}

\punk{Operations on ideals and Fr\'echet ideals}
\label{op:id}

Suppose that $A$ is any non-empty set, and $\cJ_a$ is an
ideal on a set $B_a$ for all $a\in A.$
The following two operations on such a family of ideals
are defined. 

\bde
\item[{\ubf Disjoint sum}]
$\sum_{a\in A}\cJ_a$ is the ideal on the set
\index{ideal!disjoint sum $\sum_{a\in A}\cI_a$}%
\index{zzSaAIa@$\sum_{a\in A}\cI_a$}%
$B=\ens{\ang{a,b}}{a\in A\land b\in B_a}$
that consists of all sets $x\sq B$ such that 
$\seq xa\in\cJ_a$ for all $a\in A,$ where 
$\seq xa=\ens{b}{\ang{a,b}\in x}$ (the cross-section).
If the sets $B_a$ are pairwise disjoint then $\sum_{a\in A}\cI_a$ 
can be equivalently defined as the ideal on $B=\bigcup_{a\in A}B_a$
that consists of all sets of the form
$\bigcup_{a\in A}x_a,$ where $x_a\in\cI_a$ for all $a.$

In the case of two or finitely many summands, 
the disjoint sum $\dsu\cI\cJ$ of ideals $\cI,\cJ$ on disjoint
sets $A\yi B$ is equal to
$\ens{x\cup y}{x\in\cI\land y\in\cJ}.$\vom
\index{disjoint sum (of ideals)}%
\index{ideal!disjoint sum $\dsu\cI\cJ$}%
\index{zzI+J@$\dsu\cI\cJ$}%
 
\item[{\ubf Fubini sum and product}]
Suppose in addition that $\cI$ is an ideal on $A.$
The \rit{Fubini sum} $\fss{\cI}{\cJ_a}{a\in A}$ of the ideals
$\cJ_a$ modulo $\cI$ is the ideal on the set 
$B$ (defined as above) which 
consists of all sets $y\sq B$ such that the set 
$\ens{a}{\seq xa\nin \cJ_a}$ belongs to $\cI.$  
\index{Fubini sum}%
\index{ideal!Fubini sum $\fss{\cI}{\cJ_a}{a\in A}$}%
\index{zzSJaI@$\fss{\cI}{\cJ_a}{a\in A}$}%
This ideal obviously coincides with the plain
disjoint sum $\sum_{a\in A}\cJ_a$ whenever $\cI=\ans\pu$.
 
In particular, the \rit{Fubini product} $\fpd\cI\cJ$ of two
\index{Fubini product}%
\index{ideal!Fubini product $\fpd{\cI}{\cJ}$}%
\index{zzIxJl@$\fpd{\cI}{\cJ}$}%
ideals $\cI,\cJ$ on sets resp.\ $A,B$ is equal to
$\fss\cI{\cJ_a}{a\in A},$ 
where $\cJ_a=\cJ\zd\kaz a.$ 
Thus $\fpd\cI\cJ$ consists of all sets $y\sq A\ti B$ 
such that $\ens{a}{(y)_a\nin \cJ}\in \cI$.
\ede

Coming back to the ideals defined in \nrf{borid},
$\Ii$ and $\It$ coincide with resp.\ $\fio$ and $\ofi,$
where, we recall, $0$ denotes the least ideal
$0=\ans\pu$.

The operations of Fubini sum and product allow us to
define the following interesting family of Borel 
ideals (mainly, non-P-ideals) on countable sets.\vom

{\ubf Fr\'echet ideals.\/} 
This family consists of ideals $\frt_\xi,\msur$ 
$\xi<\omi,$ defined by transfinite induction.
\index{ideal!Fr\'echet ideal $\frt\xi$}%
\index{zzFxi@$\frt\xi$}%
We put $\frt_1=\ifi$ and 
$\frt_{\xi+1}=\fpd\ifi{\frt_\xi}$ for all $\xi.$ 
Limit steps cause a certain problem. 
The most natural idea would be to define 
$\frt_\la=\fss{\ifi_\la}{\frt_\xi}{\xi<\la}$ for any limit 
$\la,$ where $\ifi_\la$ is the ideal of all finite subsets 
of $\la,$ or perhaps 
$\frt_\la=\fss{\ibo_\la}{\frt_\xi}{\xi<\la},$ 
where $\ibo_\la$ is the ideal of all bounded subsets 
of $\la,$ or even $\frt_\la=\fss{\ans\pu}{\frt_\xi}{\xi<\la}.$ 
Yet this appears not to be fully satisfactory in \cite{jkl}, 
where they define 
$\frt_\la=\fss{\ifi}{\frt_{\xi_n}}{n\in\dN},$ where 
$\sis{\xi_n}{}$ is a once and for all fixed cofinal increasing 
sequence of ordinals below $\la,$ with understanding that 
the result is independent of the choice of $\xi_n,$ modulo 
a certain equivalence.

\punk{Some other ideals}
\label{o}

We consider two interesting families of Borel 
ideals (mainly, non-P-ideals), united by their 
relation to countable ordinals.  
Note that the underlying sets of these ideals   
are countable sets different from $\dN$.\vom

{\ubf Indecomposable ideals.\/} 
Let $\otp X$ be the order type of $X\sq\Ord.$ 
For any ordinals $\xi,\,\vt<\omi$ define: 
\dm
\iz\vt\xi \,=\, \ens{A\sq \vt}{\otp A<\om^\xi} 
\quad\hbox{(nontrivial only if $\vt\ge\om^\xi$)}\,.
\dm
To see that the sets $\iz\vt\xi$ are really ideals 
note that ordinals of the form $\om^\xi$ and only 
those ordinals are {\it indecomposable\/}, \ie, 
are not sums of a pair of smaller ordinals, hence, 
the set $\ens{A\sq \vt}{\otp A<\ga}$ is 
an ideal iff $\ga=\om^\xi$ for some $\xi.$\vom

{\ubf Weiss ideals.\/} 
Let $\rkb X$ be the {\it Cantor-Bendixson rank\/} 
of $X\sq\Ord,$ 
\ie, the least ordinal $\al$ such that $\kb X\al=\pu.$ 
Here $\kb X\al$ is defined by induction on $\al:$ 
$\kb X0=X,\msur$ 
$\kb X\la=\bigcap_{\al<\la}\kb X\al$ at limit steps 
$\la,$ and finally $\kb X{\al+1}=(\kb X\al)',$ where  
$A',$ the Cantor-Bendixson 
derivative, is the set of all ordinals $\ga\in x$ which 
are limit points of $X$ in the interval topology.  
For any ordinals $\xi,\,\vt<\omi$ define: 
\dm
\iw\vt\xi \,=\, \ens{A\sq \vt}{\rkb A<\om^\xi}  
\quad\text{(nontrivial only if $\vt\ge\omm\xi$)}\,.
\dm
It is less transparent that all $\iw\vt\xi$ are 
ideals 
(see Farah~\cite[\pff1.14]{aq}) 
while $\ens{A\sq \vt}{\rkb A<\ga}$ is not an ideal 
if $\ga$ is not of the form $\om^\xi$.

\vyk{
This title intends to include those interesting ideals 
which have not yet been subject of comprehensive study. 
A common method to obtain interesting ideals is to 
consider a countable set bearing a nontrivial structure, 
as the underlying set. 
In principle, there is no difference between different 
countable set as which of them is taken as the 
underlying set for the ideals considered. 
Yet if the set bears a nontrivial structure 
(\ie, more than just countability) 
then this 
gives additional insights as which ideals are 
meaningful. 
This is already transparent for the ideals defined 
in \nrf{non-P}.

We give two examples. \vom
}

{\bfit Ideals on finite sequences.\/} 
The set $\nse$ of all finite sequences of natural 
numbers is countable, yet its own order structure is 
quite different from that of $\dN.$ 
We can exploit this in several ways, for instance, 
with ideals of sets $X\sq\nse$ which intersect every 
branch in $\nse$ by a set which belongs to a given 
ideal~on~$\dN$.

Further entry: Farah~\cite{f-tsir,f-bas,f-co} 
on {\it Tsirelson ideals\/}.

Nowhere dense ideal etc

\punk{Reducibility of ideals}
\las{reli}

There are different methods of reduction of an ideal
$\cI$ on a set $A$ to an ideal $\cJ$ on a set $B,$ where
the reducibility means that $\cI$ is in some sense simpler
(in non-strict way) than $\cJ.$

\bde
\item[{\ubf Rudin--Keisler order\/}:]
$\cI\ork\cJ$ iff there exists a function $\ba:B\to A$ 
(a {\it Rudin -- Keisler\/} reduction)
\index{reducibility!Rudin -- Keisler $\ork$}
such that $x\in\cI\eqv b\obr(x)\in\cJ$. 

\item[{\ubf Rudin--Blass order\/}:]
$\cI\orb\cJ$ iff there is a \poq{finite-to-one} 
function $\ba:B\to A$ 
\index{reducibility!Rudin -- Blass $\orb$}%
(a {\it Rudin--Blass\/} reduction) 
with the same property.

A version: $\cI\orbp\cJ$ allows $\ba$ to be defined on a 
\index{reducibility!Rudin -- Blass!modified $\orbp$}%
proper subset of $B,$ in other words, we have pairwise 
disjoint finite non-empty sets $w_a=\ba\obr(\ans{a})\yt a\in A,$ 
such that $x\in\cI\leqv w_x=\bigcup_{a\in x}w_a\in\cJ.$

Another version $\cI\orbpp\cJ,$ applicable in the case when
\index{reducibility!Rudin -- Blass!modified $\orbpp$}%
$A=B=\dN,$ requires that in addition  
the sets $w_a$ satisfy $\tmax w_a<\tmin w_{a+1}$.
\ede

\ble
\lam{sumi}
Suppose that\/ $r_n\ge 0,\msur$ $r_n\to0,$ and\/ 
$\sum_nr_n=\piy.$ 
Then any\/ summable ideal\/ $\cI$ satisfies\/ 
$\cI\orbpp\sui{r_n}.$
\ele
\bpf
Let $I=\sui{p_n},$ where $p_n\ge0$ 
(no other requirements !). 
Under the assumptions of the lemma we can associate a 
finite set $w_n\sq\dN$ to any $n$ so that 
$\tmax w_n<\tmin w_{n+1}$ and 
$|r_n-\sum_{j\in w_n}r_i|<2^{-n}.$ 
\epf

Another type of reducibility is connected with
\dd\sd homomorphisms.

Suppose that $\cI,\cJ$ are ideals on sets resp.\ $A,B.$
The power sets $\pws A\yd\pws B$ can be considered as groups
with $\sd$ as the operation and $\pu$ as the neutral element.
Then a \rit{\dd\sd homomorphism}
\index{Dhomomorphism!\dd\sg homomorphism}%
is any map $\vt:\pws A\to\pws B$ such that
$\vt(x)\sd\vt(y)=\vt(x\sd y)$ for all $x\yi y\sq A.$

The quotient $\pws A/\cI$ consists of \dd\cI classes
$\ek x\cI=\ens{x\sd a}{a\in\cI}$ of sets $x\sq A;$
it is endowed by the group operation
$\ek x\cI\sd \ek y\cI=\ek{x\sd y}\cI.$
Similarly $\pws B/\cJ.$
For a map $\vt:\pws A\to\pws B$ to induce in obvious way
a group homomorphism of $\pws A/\cI$ to $\pws B/\cJ,$ it is
necessary and sufficient that\vtm

(1)
$(\vt(x)\sd\vt(y))\sd\vt(x\sd y)\in \cJ$ for all
$x\yi y\sq A,$ \ and\,\vtm

(2)
$x\in\cI\eqv \vt(x)\in\cJ$ for all
$x\sq A.$\vtm

\noi
Let us call any such a map
\rit{an\/ \ddij approximate \dd\sd homomorphism}. 
\index{Dhomomorphism!\dd\sd homomorphism!approximate}%
\bde
\item[{\ubf Borel \dd\sd reducibility\/}:]
\index{reducibility!Borel \dd\sd reducibility $\rd$}%
$\cI\rd\cJ$ iff there is a Borel \ddij approximate  
\dd\sd homo\-mor\-phism $\vt:\pws A\to\pws B$.
\ede

Note that if a map $\ba:B\to A$ witnesses, say,
$\cI\ork\cJ$ then the map $\vt(x)=\ba\obr(x)$ obviously
witnesses $\cI\rd\cJ.$ 

\bde
\item[{\ubf Isomorphism\/}] $\cI\isi\cJ$ of ideals 
\index{ideal!isomorphism $\isi$}
$\cI,\cJ$ on sets resp.\ $A\yi B$ 
means that there is a 
\index{isomorphism!icj@$\cI\isi\cJ$}%
\index{zzicj@$\cI\isi\cJ$}%
bijection $\ba:A\onto B$ such that   
${x\in \cI}\leqv{\ba\ima x\in\cJ}$ for all $x\sq A$.
\ede

The following notion belongs to a somewhat different
category since it does not allow to really define $\cI$
in terms of $\cJ.$

\bde
\item[{\ubf Reducibility via inclusion\/}
{\rm(see \cite{jkl})}:] 
$\cI\inc\cJ$ iff there is a map $\ba:B\to A$ such  
\index{reducibility!via inclusion $\inc$}%
that $x\in\cI\imp \ba\obr(x)\in\cJ.$ 
(Note $\imp$ instead of $\eqv$!)
\ede
In particular if $\cI\sq\cJ$ (and $B=A$) then $\cI\inc\cJ$ via 
$\ba(a)=a.$ 
It follows that this order is not fully compatible with the
Borel reducibility $\reb$.

\vyk{
Sometimes they use a weaker definition: let
\index{isomorphism!icja@$\cI\isa\cJ$}%
\index{zzicja@$\cI\isa\cJ$}%
$\cI\isa\cJ$ mean that there are sets $A'\in\df\cI$ 
and $B'\in\df\cJ$ such that 
${\cI\res A'}\isi{\cJ\res B'}.$ 
Yet this implies $\cI\isi\cJ$ in most usual cases, 
the only notable exception (among nontrivial ideals), 
is produced by the ideals $\cI=\ifi$ and 
$\cJ=\dsu\ifi{\pn}\isi\ens{x\sq\dN}{x\cap D\in\ifi},$ 
where $D$ is an infinite and coinfinite set~\footnote
{\ Kechris~\cite{rig} called ideals $\cJ$ of this kind 
{\it trivial variations of\/ $\ifi$}.}
: then $\cI\isa\cJ$ but not $\cI\isi\cJ$.

\punk{Remarks}

\imar{check this subsection once again}%
The following shows simple relationships between different 
reducibilities:   %
\dm
{\cI\orb\cJ}\simp{\cI\ork\cJ}\simp{\cI\obe\cJ}
\simp{\cI\obep\cJ}\simp{\cI\odl\cJ}\simp{\cI\reb\cJ}.
\dm
For instance if 
$b:\dN\to\dN$ witnesses $\cI\ork\cJ$ then 
$\vt_b(X)=b\obr(X)$ witnesses $\cI\obe\cJ.$ 
Note that any $\vt_b$ is an exact Boolean algebra 
homomorphism $\pn\to\pn;$ moreover, it is known that 
any BM Boolean algebra homomorphism $\pn\to\pn$ is 
$\vt_b$ for an appropriate $b:\dN\to\dN.$ 
{\it Approximate\/} homomorphisms are liftings of 
homomorphisms into quotients of $\pn,$ thus, any 
\ddj approximate $\vt:\pn\to\pn$ induces the map 
$\Theta(X)=\ens{\vt(X)\sd Y}{Y\in\cJ},$ which is a 
homomorphism $\pn\to\pn/\cJ.$ 
Farah~\cite{aq}, and Kanovei and Reeken \cite{kr} 
demonstrated that in some important cases 
(of ``nonpatological'' P-ideals and, generally, for all 
Fatou, or Fubini, ideals) 
we have ${\cI\ork\cJ}\eqv{\cI\obe\cJ}.$ 
On the other hand ${\cI\ork\cJ}\nmp{\cI\obe\cJ}$ fails 
for rather artificial P-ideals. 

The right-hand end is the most intrigueing: 
is there a pair of Borel ideals $\cI,\,\cJ$ such that  
$\cI\reb\cJ$ but not $\cI\odl\cJ\,?$ 
If we actually have the equivalence then the 
whole theory of Borel reducibility for Borel ideals 
can be greatly simplified 
because reduction maps which are \dd\sd homomorphisms 
are much easier to deal with. 

}

\parf{Introduction to \eqr s}
\las{somi}

Recall that an \rit{\eqr}
(\er, for brevity)
on a set $A$ is any
\index{equivalence relation, ER}%
reflexive, transitive, and symmetric binary relation on $A.$

\bit
\item
If $\rE$ is an \er\ on a set $X$ then

$\eke y=\ens{x\in X}{y\rE x}$ for any $y\in X$ 
\index{zzye@$\ek y\rE$}%
\index{equivalence class}%
(the \dde{\it class\/} of $x$) and

$\eke Y=\bigcup_{y\in Y}\eke y$ 
(the \dde{\it saturation\/} of $Y$) for $Y\sq X.$ 
\index{saturation}
\index{zzYe@$\eke Y$}%

A set $Y\sq X$ is \dde{\it invariant\/} if 
$\eke Y=Y$. 
\index{set!einvariant@\dde invariant}%

\item 
If $\rE$ is an \er\ on a set $X$ then a set $Y\sq X$ is 
{\it pairwise \dde equivalent, {\rm resp.}, 
pairwise \dde inequivalent\/}, 
\index{set!pairwiseeeqv@pairwise \dde equivalent}%
\index{set!pairwiseeineqv@pairwise \dde inequivalent}%
if $x\rE y,$ resp., $x\nE y$ holds
for  all $x\ne y$ in $Y$.  

\item
If $X,\,Y$ are sets and $\rE$ any binary relation then 
$X\rE Y$ means that we have both 
\index{zzXEY@$X\rE Y$}%
$\kaz x\in X\:\sus y\in Y\:(x\rE y)$ and 
$\kaz y\in Y\:\sus x\in X\:(x\rE y)$.
\eit


\punk{Some examples of Borel \eqr s}
\las{someq}

Let $\rav X$ denote the equality on a set $X,$ 
\index{equivalence relation, ER!DN@$\rav\dN=\alo$}%
\index{equivalence relation, ER!D2N@$\rav\dn=\cont$}%
\index{zzequN@$\rav\dN=\alo$}%
\index{zzequ2N@$\rav\dn=\cont$}%
considered as an \eqr\ on $X.$
\index{equivalence relation, ER!equality@equality $\rav X$}%
\index{zzequX@$\rav X$}%
This is the most elementary type of \er s. 
A much more diverse family consists of \eqr s $\rei$ generated
by Borel ideals.
\bit   
\item
If $\cI$ is an ideal on a set $A$ then $\rei$ denotes  
an \eqr\ on $\pws A,$ 
defined so that $x\rei y$ iff ${x\sd y}\in\cI$. 
\index{equivalence relation, ER!ei@$\rei$}%
\index{zzei@$\rei$}%
\eit
Equivalently, $\rei$ can be considered as an \eqr\ on
$2^A$ defined so that $f\rei g$ iff ${f\sd g}\in\cI,$
where $f\sd g=\ens{a\in A}{f(a)\ne g(a)}.$
Note that $\rei$ is Borel provided so is $\cI.$
We obtain the following important \eqr s:\snos
{The notational system we follow is not the only one used
in modern texts.
For instance $\Ei\yi\Ed\yi\Et$ are sometimes denoted
differently, see \eg\ \cite{su9as}.}
\newlength{\eoe}
\settowidth{\eoe}{$\Eo=\rE_{\ifi}$}
\nek{\eoebox}[1]{\hspace*{0ex}\makebox[\eoe][l]{#1}}
\bde
\item 
[\eoebox{$\Eo=\rE_{\ifi}$}]
is an \er\ on $\pn,$ 
\index{equivalence relation, ER!E0@$\Eo$}%
\index{zzE0@$\Eo$}%
and $x\Eo y$ iff $x\sd y\in\ifi$.

\item 
[\eoebox{$\Ei=\rE_{\Ii}$}]
\parbox[t]
{0.85\textwidth}
{is an \er\ on $\pws{\dN\ti\dN},$ 
\index{equivalence relation, ER!E1@$\Ei$}%
\index{zzE1@$\Ei$}%
and $x\Ei y$ iff $\seq xk=\seq yk$ for all but finite 
$k,$ where, we recall, $\seq xk=\ens{n}{\ang{k,n}\in x}$ 
for $x\sq \dN\ti\dN.$} 

\item 
[\eoebox{$\Ed=\rE_{\Id}$}]
is an \er\ on $\pn,$ 
\index{equivalence relation, ER!E2@$\Ed$}%
\index{zzE2@$\Ed$}%
and $x\Ed y$ iff $\sum_{k\in x\sd y}k\obr<\iy.$ 

\item
[\eoebox{$\Et=\rE_{\It}$}]
is an \er\ on $\pws{\dN\ti\dN},$ 
\index{equivalence relation, ER!E3@$\Et$}%
\index{zzE3@$\Et$}%
and $x\Et y$ iff $\seq xk\Eo\seq yk\zd\kaz k.$

\item 
[\eoebox{$\rzo=\rE_{\Zo}$}]
is an \er\ on $\pn,$ 
\index{equivalence relation, ER!Z0@$\rzo$}%
\index{zzZ0@$\rzo$}%
and $x\rzo y$ iff 
\index{zzz0@$\rzo$}%
$\tlim_{n\to\iy}\frac{\#((x\sd y)\cap \ir0n)}n=0.$ 
\ede

Alternatively, $\Eo$ can be viewed as an \eqr\ on 
$\dn$ defined as $a\Ei b$ iff $a(k)=b(k)$ 
for all but finite $k.$ 
Similarly, $\Ei$ can be viewed as a \er\ on 
$\pnqn,$ or even on $\dnqn,$
defined as $x\Ei y$ iff $x(k)=y(k)$ 
for all but finite $k,$ for all $x,y\in \pn^\dN,$
while $\Et$ can be viewed as a \er\ on 
$\pnqn,$ or on $\dnqn,$ defined as
$x\Et y$ iff $x(k)\Eo y(k)$ for all $k$.

The next group includes \eqr s induced by actions of
(the additive groups of) some Banach spaces ---
see below on group actions.
The following Banach spaces are well known from textbooks:
\dm
\bay{rclclcl} 
\index{Banach space!lp@$\bel p$}%
\index{zzlp@$\bel p$}%
\index{Banach space!linfty@$\bel\iy$}%
\index{zzlinfty@$\bel\iy$}%
\index{Banach space!c@$\fc$}%
\index{zzc@$\fc$}%
\index{Banach space!c0@$\fco$}%
\index{zzc0@$\fco$}%
\bel p &=& \ens{x\in\rn}{\sum_n|x(n)|^p < \iy}\;\;(p\ge1);&& 
\nr x p &=& (\sum_n |x(n)|^p)^{\frac1p}\,;\\[0.8\dxii]
\index{norm!lpnorm@\dd{\bel p}norm $\nr\cdot p$}%
\index{zzpx@$\nr\cdot p$}%

{\bel\iy} &=& \ens{x\in\dR^\dN}{\tsup_n|x(n)| < \iy};&&
\nr x\iy &=& \tsup_n |x(n)|\,;\\[0.8\dxii]
\index{norm!linfnorm@\dd{\bel\iy}norm $\nr\cdot\iy$}%
\index{zzinf@$\nr\cdot\iy$}%

\fc &=& \ens{x\in\dR^\dN}{\tlim_n x(n)<\iy\,\text{ exists}};&&
\nor x &=& \tsup_n |x(n)|\,;\\[0.8\dxii]

\fco &=& \ens{x\in\dR^\dN}{\tlim_n x(n)=0};&&
\nor x &=& \tsup_n |x(n)|\,.\\[0.8\dxii]
\vyk{
\fC[a,b] &=& \ens{f\in\dR^{[a,b]}}{f\,\text{ is continuous}};&&
\nor f &=& \tmax_{a\le t\le b} |f(t)|\,;\\[0.8\dxii]

\fC(K) &=& \ens{f\in\dR^K}{f\,\text{ is continuous}};&&
\nor f &=& \tmax_{t\in K} |f(t)|\,;\\[0.8\dxii]

\beL p[a,b] &=&  
\ens{f\in\dR^{[a,b]}}{\int_a^b|f(t)|^p dt<\iy};&&
\nor f &=& (\int_a^b|f(t)|^p dt)^{\frac1p}.
}%
\eay
\dm
Note that ${\bel p}\yd{\fc}\yd{\fco}$ are separable spaces 
while ${\bel\iy}$ is non-separable.
The domain of each of these spaces 
consists of infinite sequences 
$x=\sis{x(n)}{n\in\dN}$ of reals, and is a subgroup of 
the group $\rn$ (with the componentwise addition). 
The latter can be naturally equipped with the Polish 
product topology, so that 
${\bel p}\yd{\bel\iy}\yd{\fc}\yd{\fco}$ are Borel subgroups 
of $\rn.$ 
(But not topological subgroups since the distances are different. 
The metric definitions as in ${\bel p}$ or ${\bel\iy}$ 
do not work for $\rn$.)

Each of the four mentioned Banach spaces induces
an \rit{orbit equivalence} --- 
a Borel \eqr\ on $\rn$ also denoted by, resp.,  
${\bel p}\yd{\bel\iy}\yd{\fc}\yd{\fco}.$  
\index{equivalence relation, ER!lp@$\bel p$}%
\index{zzlp@$\bel p$}%
\index{equivalence relation, ER!linfty@$\bel\iy$}%
\index{zzlinfty@$\bel\iy$}%
\index{equivalence relation, ER!c@$\fc$}%
\index{zzc@$\fc$}%
\index{equivalence relation, ER!c0@$\fco$}%
\index{zzc0@$\fco$}%
For instance, $x\bel p y$ if and only if 
$\sum_k|x(k)-y(k)|^p<\piy$ (for all $x\yi y\in\rn$).

\vyk{
It is known (see \grf{BAN})
that ${\bel1}\eqb\Ed$ and
${\bel p}\rebs{\bel q}$ whenever $1\le p<q,$ in 
particular, ${\bel1}\eqb\Ed\rebs{\bel q}$ for any $q>1.$
On the other hand, $\fco\eqb\rzo,$ where
\index{equivalence relation, ER!density 0}%
\index{equivalence relation, ER!z0@$\rzo$}%
$\rzo$ is the \rit{density 0} \eqr:
}

Another important \eqr\ is 
%
\bde
\item 
[$\rtd,$] often called 
\index{equivalence relation, ER!t2@$\rtd$}%
\index{zzt2@$\rtd$}%
\rit{the equality of countable sets of reals}, is an \er\
defined on $\nnn$ so that  
$g\rtd h$ iff $\ran g=\ran h$ (for $g\yi h\in\nnn$). 
\ede
%
There is no reasonable way to turn $\pwc\dnn,$ the set 
of all at most countable subsets of $\dnn,$ into a 
Polish space, in order to directly define the equality
of countable sets of reals in terms of $\rav{\cdot}.$ 
However, nonempty members of $\pwc\dnn$ can be
identified with equivalence classes in $\nnn/\rtd.$ 
(See \grf{groas} on the whole series of
\eqr s $\rT_\al\yt \al<\omi$.)

We finish with another important \eqr,

\bde
\item  
[$\Ey,$]
\rit{the universal countable Borel \er\/}. 
\index{equivalence relation, ER!Einfty@$\Ey$}%
\index{zzEinfty@$\Ey$}%
The countability here means that all \dde equivalence
classes $\eke x$ are at most countable sets.
The notion of universality will be explained below.
\ede
See Example~\ref{ex:ac4} on an exact definition of $\Ey$.

\punk{Operations on \eqr s}
\las{opeer}

The following operations over \eqr s are in 
part parallel to the operations on ideals in \nrf{op:id}.
Suppose that $A$ is any non-empty and 
{\it at most countable\/} set, and $\rF_a$ is an
\eqr\ on a set $X_a$ for all $a\in A.$
The following operations on such a family of \er s
are defined. 

\ben
\tenu{(o\arabic{enumi})} 
\itla{oe:beg}
\lam{cun}
{\ubf Union} $\bigcup_{a\in A}\rF_a$ 
(if it results in an \eqr)
and {\ubf intersection} $\bigcap_{a\in A}\rF_a$
(it always results in an \eqr)
--- in the case when all $\rF_a$ are \er s
on one and the same set $X=X_a\zd\kaz a$. 

\itla{cdun}
{\ubf Countable disjoint union}
$\bigvee_{a\in A}\rF_a$
is an \er\ $\rE$ on the set $X=\bigcup_a(\ans a\ti X_a)$ 
defined as follows: $\ang{a,x}\rE \ang{b,y}$ iff 
$a=b$ and $x\rE_a y.$ 

If the sets $X_a$ are pairwise disjoint
then we can equivalently define an \eqr\ 
$\rE=\bigvee_a\rF_a$ on the set  $Y=\bigcup_aX_a$ so that 
$x\rE y$ iff $x\yi y$ belong to the same $X_a$ 
and $x\rF_a y$.  

\itla{prod} 
{\ubf Product} $\prod_{a\in A}\rF_a$ 
\index{product of ers@product of \er s}%
is an \er\ $\rE$ on the cartesian product  
$\prod_{a\in A}X_a$ defined so that 
$x\rE y$ iff $x(a)\rF_a y(a)$ for all $a\in A.$

In particular the product $\rF_1\ti\rF_2$ of \er s $\rE\yi\rF$
on sets resp.\ $X_1\yi X_2$ is an \er\ $\rE$ on $X_1\ti X_2$
defined so that $\ang{x_1,x_2}\rE\ang{y_1,y_2}$ iff
$x_1\rF_1 y_1$ and $x_2\rF_2 y_2$.

If $X_a=X$ and $\rF_a=\rF$ for all $a$ then the power
notation $\rF^A$ can be used instead of $\prod_{a\in A}\rF_a$.

\itla{fp}
The {\ubf Fubini product} (ultraproduct) 
\index{Fubini product!of ers@of \er s}%
\index{ultraproduct of ers@ultraproduct of \er s}%
$\fps{\cI}{\rF_{a}}{a\in A}$ 
modulo an ideal $\cI$ on $A$ 
is the \er\ on the product space 
$\prod_{a}X_a$ defined as follows: 
$x\rE y$ iff the set $\ens{a}{x(a)\nF_a y(a)}$
belongs to $\cI.$  

If $X_a=X$ and $\rF_a=\rF$ for all $a$ then the ultrapower
notation $\rF^A/\cI$ can be used instead of
$\fps{\cI}{\rF_{a}}{a\in A}$.

\itla{oe:end}
\lam{cp}
{\ubf Countable power\/} of an \eqr\ $\rF$ on a set $X$ is  
an \er\ $\rfi$ defined on the set $X^\dN$ as follows: 
$x\rfi y$ iff 
$\ens{[x(k)]_{\rE}}{k\in\dN}=\ens{[y(k)]_{\rE}}{k\in\dN},$ 
so that for any $k$ there 
is $l$ with $x(k)\rF y(l)$ and for any $l$ there 
is $k$ with $x(k)\rF y(l)$.
\een

\sbrospri

\bex
\lam{rtx1}
In these terms, the \eqr s $\Ei$ and $\Et$ coincide with
resp.\ $(\rav\dn)^\dN/\ifi$ and ${\Eo}^\dN$ modulo obvious
bijections between the spaces considered.
Generally, the operations on ideals introduced
in \nrf{op:id} transform
in some regular way into operations on the corresponding \eqr s.
For instance $\rE_{\fss{\cI}{\cJ_a}{a\in A}}$ is equal to
$\fps{\cI}{\rE_{\cJ_a}}{a\in A},$ while $\rE_{\fpd\cI\cJ}$ is
equal to ${(\rE_{\cJ})}^A/\cI,$ where $A$ is the domain of $\cI.$

Accordingly, $\rE_{\sum_a\cJ_a}$ is equal to $\prod_a\rE_{\cJ_a}.$
In particular if $\cI,\cJ$ are ideals on disjoint sets $A\yi B$
then $\rE_{\dsu\cI\cJ}$ is equal to $\rei\ti\rej$.
\eex

\bex
\lam{rtx2}
The \eqr\ $\rtd$ defined in \nrf{someq}
coincides with $\epo{\rav\dnn}$.
\eex

Iterating these operations, we obtain a lot of interesting 
\eqr s starting just with very primitive ones. 

\bex
\lam{rtx}
Iterating the operation of countable power,
H.~Friedman~\cite{frid} defines the sequence of 
\er s $\rT_\xi\yt 1\le\xi<\omi,$   
as follows~\footnote
{\ Hjorth~\cite{h} uses $\rF_\xi$ instead of $\rT_\xi$.}. 
Let $\rT_1=\rav\dnn,$ the equality relation on $\dnn.$ 
Put $\rT_{\xi+1}=\epo{\rT_\xi}$ for all $\xi<\omi.$ 
\index{equivalence relation, ER!Txi@$\rT_\xi$}%
\index{zzTxi@$\rT_\xi$}%
If $\la<\omi$ is a limit ordinal, then put 
$\rT_\la=\bigvee_{\xi<\la}\rT_\xi.$
The definition for the second term $\rtd$ is equivalent with
the separate definition of $\rtd$ in \nrf{someq}
by \ref{rtx2}.
\eex

\vyk{
Note that this transfinite sequence contains 
$\rT_2=\epo{{\epo{\rav\dN}}},$ an \eqr\
also separately defined in 
\nrf{someq}.
The two definitions do not literally coincide with each other.
However the two versions of $\rtd$ are Borel reducible to
each other (see Example~\ref{rtd=rtd}) and also to
$\epo{\rav\dn},$
and hence mutually interchangable in major issues. 
}

\vyk{
In addition to the families of \eqr s introduced by
Definition~\ref{typesER}, some more complicated families
will be considered below, including \er s induced by
Polish group actions, turbulent \er s,
\er s classifiable by countable structures, pinned \er s,
and some more.
}

\punk{Orbit \eqr s of group actions}
\las{groa1}

An {\it action\/} of a group $\dG$ on a space $\dX$ is 
\index{action}%
any map $\goa:{\dG\ti\dX}\to\dX,$ usually written as 
$\goa(g,x)=g\app x,$ such that 
\index{action!g.x@$g\app x$}%
\index{zzg.x@$g\app x$}%
$1)\msur$ $e\app x=x,$ and 
%
$2)\msur$ $g\app(h\app x)=(gh)\app x.$  
Then, for any $g\in\dG,$ the map $x\mapsto g\app x$ is 
a bijection $\dX$ onto $\dX$ with $x\mapsto g\obr\app x$
being the inverse map.
A {\it\dd\dG space\/} is a pair $\stk\dX \goa,$ where $\goa$ is an 
\index{Gspace@{}\dd\dG space}%
action of $\dG$ on $\dX\,;$ 
in this case $\dX$ itself is also called a \dd\dG space, 
and the {\it orbit \er\/},  
\index{equivalence relation, ER!orbit er@orbit \er\ $\ergx$}%
\index{equivalence relation, ER!induced@induced by an action}%
or {\it\er\ induced by the action\/}, 
\index{zzEGX@$\ergx$}%
$\aer \goa\dX=\ergx$ is defined on $\dX$ so that $x\ergx y$ iff 
there is $g\in\dG$ with $y=g\app x.$
{}\dd\ergx classes are the same as 
\dd\dG orbits, \ie,
\dm
[x]_\dG=[x]_{\ergx}=\ens{y}{\sus g\in\dG\;(g\app x=y)}\,.
\dm

\vyk{
A {\it homomorphism\/} (or \dd\dG homomorphism) 
of a \dd\dG space $\dX=\stk\dX \goa$ into a \dd\dG space
$\dY=\stk\dY\gob$ is any map $F:\dX\to\dY$ compatible 
with the actions in the sense that $F(g\app x)=g\app F(x)$
(more exactly, $F(\goa(g,x))=\gob(g,F(x))$)
for any $x\in\dX$ and $g\in\dG.$ 
A $1-1$ homomorphism is an {\it embedding\/}. 
\index{homomorphism!of actions}%
\index{embedding!of actions}%
\index{isomorphism!of actions}%
An embedding $\onto$ is an {\it isomorphism\/}. 
Note that a homomorphism $\stk\dX\goa\to\stk\dY\gob$ is a 
reduction of the induced \eqr{}
$\aer \goa\dX$ to $\aer\gob\dY,$ but not conversely.
}

Recall that a {\it Polish group\/} is a group whose underlying
set is a Polish space and the operations are continuous. 
\index{group!Polish}%
\index{group!Borel}%
A {\it Borel group\/} is a group whose underlying set is a 
Borel set (in a Polish space) and the operations are Borel maps. 
%
\index{group!Polishable}%
A Borel group is {\it Polishable\/} if there is a Polish 
topology on the underlying set which 1) produces the same 
Borel sets as the 
original topology and 2) makes the group Polish.

\bit
\itsep
\item
If both $\dX$ and $\dG$ are Polish and the action continuous, 
then $\stk\dX\goa$ (and also $\dX$) is called a {\it Polish\/} 
\dd\dG space. 
\index{Gspace@{}\dd\dG space!Polish}%
\index{Gspace@{}\dd\dG space!Borel}%
If both $\dX$ and $\dG$ are Borel and the action is a 
Borel map, then $\stk\dX\goa$ (and also $\dX$) is called a 
{\it Borel\/} \dd\dG space.
\eit

\sbrospri

\bex
\lam{ex:ac1}
Any ideal $\cI\sq\pn$ is a group with $\sd$ as the operation.
\index{group!ideal as \dd\sd group}%
We cannot expect this group to be Polish in the product
topology inherited from $\pn$
(indeed, $\cI$ would have to be $\Gd$).
However if $\cI$ is a P-ideal then it is Polishable   
(see \nrf{pideals}),
in other words, $\stk\cI\sd$
is a Polish group in an appropriate Polish topology
compatible with the Borel structure of $\cI.$
Given such a topology, the \dd\sd action of
\imar{correct?}%
(a P-ideal) $\cI$ on $\pn$ is Polish, too.
\eex

\bex
\lam{ex:ac2}
$\dG=\pwf\dN$ is a countable subgroup of
$\stk\pn\sd.$
Define an action of $\dG$ on $\dn$ as follows:
$(w\app x)(n)=x(n)$ whenever $n\nin w$ and
$(w\app x)(n)=1-x(n)$ otherwise.
The orbit \eqr\ $\ergx$ of this action is obviously $\Eo$.
This action is Polish
(given $\dG=\pwf\dN$ the discrete topology)
and {\it free\/}:
\index{action!free}%
$x=w\app x$ implies $w=\pu$ (the neutral element of $\dG$)
for any $x\in\dn.$
\eex

Consider any Borel pairwise \dd\Eo inequivalent set
$T\sq\dn.$ 
Then $w\app T\cap T=\pu$ for any $\pu\ne w\in\pwf\dN$
by the above. 
It easily follows that $T$ is meager in $\dn.$
(Otherwise $T$ is co-meager on a basic clopen set
$\bon s\dn=\ens{x\in\dn}{s\su x},$ where
$s\in\dln.$
Put $w=\ans{n},$ where $n=\lh s.$
Then $w\in\dG$ maps
$T\cap\bon{s\we0}\dn$ onto $T\cap\bon{s\we1}\dn.$
Thus $w\app T\cap T\ne\pu$ -- contradiction.)
We conclude that $\dG\app T=\bigcup_{w\in\dG}w\app T$
is still a meager subset of $\dn$ in this case, and hence
$T$ cannot be a full (Borel) transversal for $\Eo.$

\bex
\lam{ex:ac3}
The {\it canonical\/ {\rm(or {\it shift\/})} action\/}
\index{action!canonicalofGonXG@canonical of $\dG$ on $X^\dG$}%
\index{action!shift}%
of a group $\dG$ on a set of the form $X^\dG$ ($X$ any set)
is defined as follows:
$g\app\sis{x_f}{f\in\dG}=\sis{x_{g\obr f}}{f\in\dG}$ for any
element $\sis{x_f}{f\in\dG}\in X^\dG$ and any $g\in\dG.$
This is easily a Polish action provided $\dG$ is countable,
$X$ a Polish space, and $X^\dG$ given the product topology.
The \eqr\ on $X^\dG$ induced by this action is denoted by
$\rE(\dG,X)$.
\index{equivalence relation, ER!EGX@$\rE(\dG,X)$}%
\index{zzEGX@$\rE(\dG,X)$}%
\index{equivalence relation, ER!induced by the shift action}%
\eex

\bex
\lam{ex:ac4}
The free group of two generators $F_2$ consists of finite
\index{group!free with two generators, $F_2$}%
\index{zzF2@$F_2$}%
irreducible words composed of the symbols
$a\yi b\yi a\obr\yi b\obr,$ including the empty 
word (the neutral element).
The group operation is the concatenation of words
(followed by reduction, if necessary, \eg\ $ab\cdot b\obr a=aa$).

The shift action of $F_2$ on the compact space
$2^{F_2}$ is defined in accordance with the general scheme of
Example~\ref{ex:ac3}, so that if $x\in 2^{F_2}$ and $w\in F_2$
then $(w\app x)(u)=x(w\obr u)$ for all $u\in F_2.$ 
Put, for $x\yi y\in2^{F_2},$ $x\Ey y$ iff 
$x=w\app y$ for some $w\in F_2.$
Thus $\Ey$ is $\rE(F_2,2)$ in the sense of \ref{ex:ac3}.
\index{equivalence relation, ER!einf@$\Ey$}%
\index{zzeinf@$\Ey$}%
\eex

\bex
\lam{ex:ac5}
Come back to Banach spaces ${\bel\iy}\yi{\bel p}\yi{\fc}\yi{\fco}$
discussed in \nrf{someq}.
Each of them can be considered as a Polish group in the sense of
componentwise addition in $\rn.$
Each of them canonically acts on $\rn$ also by componentwise
addition.
For the sake of brevity,
the orbit equivalence relations of these actions, \ie\
$\aer{{\bel\iy}}{\rn}\zd
\aer{{\bel p}}{\rn}\zd
\aer{{\fc}}{\rn}\zd
\aer{{\fco}}{\rn},$
are denoted by the same symbols resp.\
${\bel\iy}\yi{\bel p}\yi{\fc}\yi{\fco}$.
\eex

\bex
\lam{ex:ac6}
The group $\isg$ of all permutations of $\dN$
\index{group!sinfty@$\isg$}
(that is, all bijections $f:\dN\onto\dN,$ with the superposition
as the group operation)
is a Polish group
in the Polish product topology of $\dnn.$
It acts on any set of the form $X^\dN$ as follows:
for any $g\in\isg$ and $x\in X^\dN,$
$(g\app x)(k)=x(g\obr(k))$ for all $k,$ or equivalently
$(g\app x)(g(k))=x(k)$ for all $k.$
Formally, $g\app x=xg\obr$ the the sense of the superposition
in the right-hand side.

Take $X=\dnn.$
Note that $\nnn$ with the product topolody is a Polish space
and the above action is Polish.
Its orbit equivalence $\aer\isg\nnn$ is quite similar to $\rtd,$
but in fact not equal.
Indeed if $x\yi y\in \nnn$ satisfy $x(0)=x(1)=y(0)=u$ and
$x(k)=u(l)=v$ for all $k\ge2\yt l\ge 1,$ where $u\ne v\in \dnn,$
then $x\rtd y$ holds while $x\aer\isg\nnn$ fails.
Still Lemma~\ref{l:rtd} will prove that $\rtd$ and $\aer\isg\nnn$
are Borel equivalent.
\eex

\punk{Borel and Polish actions}
\las{groa2}

The next theorem (too difficult to be proved here) 
shows that the type of the group is the essential 
component in the difference between Polish and Borel actions: 
roughly, any Borel action of a Polish group $\dG$ is a Polish 
action of $\dG$. 

\bte
[{{\rm\cite[5.2.1]{beke}}}]
\lam{b2c}
Suppose that\/ $\dG$ is a Polish group and\/ $\stk\dX a$ is 
a Borel\/ \dd\dG space. 
Then\/ $\dX$ admits a Polish topology which 1) produces the 
same Borel sets as the original topology, and 2) makes 
the action to be Polish.\qeD
\ete

If $\stk\dX\goa$ is a Borel \dd\dG space 
(and $\dG$ is a Borel group) 
then $\ergx$ is easily a $\fs11$ \eqr\ on $\dX.$  
Sometimes $\ergx$ is even Borel: for instance, 
when $\dG$ is a countable group and the action is Borel, 
or if $\dG=\cI\sq\pn$ is a Borel ideal, considered as 
a group with $\sd$ as the operation, which acts on 
$\dX=\pn$ by $\sd$ --- thus ${\aer\dG\pn}=\rE_{\cI}$ is a
Borel relation because $x \aer\dG\pn y$ is
equivalent to $x\sd y\in\cI.$ 
Several much less trivial cases when $\ergx$ is Borel are 
described in \cite[Chapter 7]{beke}, 
for instance, if all \dd{\ergx}classes are Borel sets of 
bounded rank then $\ergx$ is Borel \cite[7.1.1]{beke}.
Yet rather surprisingly equivalence classes generated by
Borel actions are always Borel.

\bte
[{{\rm see \cite[15.14]{dst}}}]
\lam{borcl} 
If\/ $\dG$ is a Polish group and\/ $\stk\dX\goa$ is 
a Borel\/ \dd\dG space  
then every equivalence class of\/ $\ergx$ is Borel.
\ete

The first notable case of this theorem was established by
Scott~\cite{sco} in the course of the proof that 
for any countable order type $t$
(not necessarily well-ordered)
the set of all sets $x\sq\dQ$ of order type $t$ is Borel in
$\pws \dQ.$

\bpf
It can be assumed, by Theorem~\ref{b2c}, that the action 
is continuous.
Then for any $x\in\dX$ the {\it stabilizer\/} 
$\dG_x=\ens{g}{g\app x=x}$
is a closed subgroup of $\dG.$~\footnote 
{\ Kechris~\cite[9.17]{dst} gives an independent proof. 
Both $\dG_x$ and its topological closure, say, $G'$ are 
subgroups, moreover, $G'$ is a closed subgroup, hence, we can 
assume that $G'=\dG,$ in other words, that $\dG_x$ is dense 
in $\dG,$ and the goal is to prove that $\dG_x=\dG.$ 
By a simple argument, $\dG_x$ is 
either comeager or meager in $\dG.$ 
But a comeager subgroup easily coincides with the whole group, 
hence, assume that $\dG_x$ is meager (and dense) in 
$\dG$ and draw a contradiction. 

Let $\sis{V_n}{n\in\dN}$ be a basis of the topology of $\dX,$ 
and $A_n=\ens{g\in\dG}{g\app x\in V_n}.$ 
Easily $A_nh=A_n$ for any $h\in\dG_x.$ 
It follows, because $\dG_x$ is dense, that every $A_n$ is 
either meager or comeager. 
Now, if $g\in\dG$ then $\ans g=\bigcap_{n\in N(g)}A_n,$ 
where $N(g)=\ens{n}{g\app x\in V_n},$ thus, at least one of 
sets $A_n$ containing $g$ is meager. 
It follows that $\dG$ is meager, contradiction.} 
\imar{Quotient spaces~?}%
We can consider $\dG_x$ as continuously acting on $\dG$
by $g\app h=gh$ for all $g,h\in\dG.$
Let $\rF$ denote the associated orbit \er. 
Then every \ddf class $\ek g\rF=g\,\dG_x$
is a shift of $\dG_x,$ 
hence, $\ek g\rF$ is closed.
On the other hand, the saturation $\ek \cO\rF$
of any open set $\cO\sq\dG$ is obviously open. 
Therefore, by Lemma~\ref{transv}\ref{transv4} below,
$\rF$ admits a Borel transversal $S\sq\dG.$  
Yet $g\longmapsto g\app x$ is a Borel $1-1$ map 
of a Borel set $S$ onto $\ek x\rE,$ hence,
$\ek x\rE$ is Borel by \Cpro.
\epf

It follows that \poq{not} all $\fs11$ \er s are orbit \er s 
of Borel actions of Polish groups: indeed, take a non-Borel 
$\fs11$ set $X\sq\dnn,$ define $x\rE y$ if either $x=y$ or 
$x\yi y\in X,$ this is a $\fs11$ \er\ with a non-Borel class 
$X$.

\imar{$\fs11$ or Borel \er\ not induced by Borel grp~?}%

\imar{Borel \er\ not induced by Polish grp~?}%

\vyk{
\punk{Forcings associated with pairs of  \eqr s}
\las{f2}

The range of applications of this comparably new topic
is not yet clear, but at least it offers interesting
technicalities.

\bdf[{{\rm Zapletal~\cite{zap}}}]
\lam{df:f2}
Suppose that $\rE$ is a Borel \eqr\ on a Polish space
$\dX,$ and $\rF\rebs\rE$ is another Borel \eqr.

$\zidef$ is the collection of all Borel
\index{ideal!ief@$\zidef$}%
\index{zzief@$\zidef$}%
sets $X\sq\dX$ such that ${\rE\res X}\reb\rF.$
Clearly $\zidef$ is an ideal in the algebra of all
Borel subsets of $\dX.$
The associated forcing $\zfoef$  
\index{forcing!pef@$\zfoef$}%
\index{zzpef@$\zfoef$}%
consists of all Borel sets $X\sq\dX\zt X\nin\zidef$.
\edf

For instance, the ideal $\zid{\rav\dn}{\rav\dN}$  consists
of all countable Borel sets $X\sq\dn,$ therefore 
$\zfo{\rav\dn}{\rav\dN}$ contains all
\poq{un}countable Borel sets $X\sq\dn$ and is equal to the
Sacks forcing.
The ideal $\zid\Eo{\rav\dn}$  consists
of all Borel sets $X\sq\dn$
such that $\Eo\res X$ is non-smooth
(since smoothness is equivalent to being $\reb\rav\dn$).
See \nrf{zeo} on the associated forcing $\zfo\Eo{\rav\dn}$.
}
 
\vyk{
 
$C^n[a,b]$ funk
iz [a,b] v R, nepreryvn proizvodn do poryadka n vkl
s normoj  \sum_{k=0}^n \max_x |f^n(x)|,  a\le x\le b

$C^n(I_m)$ funk
iz m-mernogo kuba v R, nepreryvn proizvodn do poryadka n vkl
s ravnomernoj normoj po vsem proizvodn

$L_p[a,b]$ funk
izmerimy na [a,b], sunmmiruemy so stepen'yu p
norma kak l_p

This group includes, for any $p\ge 1,$ 
the \er\ $\bel p$ on $\rtn$ defined 
\imar{Are tsirelsons Banach ?\\[1ex]
Is any ctble Borel $\rE\reb{\fco}$ necessarily 
$\reb\Eo$?\\[1ex] 
Are $\bel p$ and $\fco$ comparable with tsirelsons~?}%
as follows: $x\bel p y$ iff $\sum_k|x_k-y_k|^p<\piy$ 
for all sequences 
$x=\sis{x_k}{k\in\dN}$ and $y=\sis{y_k}{k}$ in $\rtn.$ 
Clearly $\bel p$ is the orbit relation of the 
(additive group of the) Banach 
space $\ell^p=\ens{x\in\rtn}{\sum_k|x_k|^p<\piy}$ 
acting on $\rtn$ by componentwise addition. 
}

\api

\parf{Borel reducibility of \eqr s}
\las{n:ers}

There are several reasonable ways to compare \eqr s  
in terms of existence of a 
{\it reduction\/}, that is, a map of certain kind which 
allows to derive one of the \er s from the other one.  
The Borel reducibility $\reb$ is the key one.
The plan of this \gla\ is to define $\reb$ and present
a diargam which displays mutual \dd\reb reducibility of
the \eqr s introduced in \nrf{someq} (the key \eqr s).
The proof of related reducibility/irreducibility claims
will be the mail content of the remainder
of the book.

\punk{Borel reducibility}
\las{bored}

Suppose that $\rE$ and $\rF$ are \eqr s on Borel sets $X\yi Y$
in some Polish spaces.
We define
\bde
\item
[$\rE\reb\rF$]
({\it Borel reducibility\/} of $\rE$ to $\rF$)\,
\index{reducibility!Borel}%
iff there is a Borel map 
\index{zzElebF@$\rE\reb\rF$}%
$\vt: X\to Y$ (called {\it reduction\/}) 
\index{reduction}%
such that ${x\rE y}\eqv {\vt(x)\rF\vt(y)}$ for all 
$x\yi y\in X$;

\item
[$\rE\eqb\rF$]
iff $\rE\reb\rF$ and $\rF\reb\rE$ 
({\it Borel bi-reducibility\/}, or
\rit{Borel equivalence});  
\index{reducibility!Borel bi-reducibility $\eqb$}%
\index{reducibility!Borel equivalence $\eqb$}%
\index{zzEeqF@$\rE\eqb\rF$}%

\item 
[$\rE\rebs\rF$]
iff $\rE\reb\rF$ but not $\rF\reb\rE$ 
\index{zzElesF@$\rE\rebs\rF$}%
\index{reducibility!Borel strict $\rebs$}%
({\it strict Borel reducibility\/}).
\ede
If $\rE\reb\rF$ (resp.\ $\rE\rebs\rF\yt \rE\eqb\rF$)
then $\rE$ is said to be
\rit{Borel reducible}
(resp.\ \rit{Borel strictly reducible},
\rit{Borel equivalent} or \rit{bi-reducible})
to $\rF.$ 
\index{equivalence relation, ER!Borel reducible}%
\index{equivalence relation, ER!Borel reducible!strictly}%
\index{equivalence relation, ER!Borel bi-reducible}%
\index{equivalence relation, ER!Borel equivalent}%

\bde
\item  
[$\rE\emb\rF$]
iff there is a Borel 
\index{zzEemF@$\rE\emb\rF$}%
\index{embedding@Borel embedding $\emb$}%
{\it embedding\/}, \ie, a $1-1$ reduction;

\item 
[$\rE\aprb\rF$]
iff $\rE\emb\rF$ and $\rF\emb\rE$ 
(a rare form, \cite[\pff0]{sinf});
\index{zzEemeF@$\rE\aprb\rF$}%

\item 
[$\rE\embi\rF$]
iff there is a Borel 
{\it invariant\/} embedding, \ie, an embedding $\vt$ 
\index{embedding!invariant}%
\index{zzEemiF@$\rE\embi\rF$}%
such that $\ran\vt=\ens{\vt(x)}{x\in\dX}$ is an 
\index{set!invariant}%
\ddf{\it invariant\/} set 
(meaning that the \ddf{\it saturation\/} 
$\ekf{\ran\vt}=\ens{y'}{\sus x\:(y\rF\vt(x))}$ 
equals $\ran\vt$);
\ede 
Sometimes they write ${\dX/\rE}\reb{\dY/\rF}$ instead 
of $\rE\reb\rF$.

\bde
\item[{\ubf Borel reducibility of ideals\/}:]
$\cI\reb\cJ$ iff $\rei\reb\rej.$ 
Thus it is required that there is a Borel map 
$\vt:\pws A\to\pws B$ such that ${x\sd y}\in \cI$ 
iff ${\vt(x)\sd \vt(y)}\in \cJ.$ 
(Here $\cI,\cJ$ are ideals on countable sets $A,B.$) 
%
\ede

In the domain of ideals, $\reb$ is weaker than all
reducibilities of more special nature discussed in
\nrf{reli} --- in the sense that, for instance, each of 
$\cI\orb\cJ$ and $\cI\rd\cJ$ implies $\cI\reb\cJ.$
The only exception is the reducibility via inclusion
$\inc$ --- it does not imply $\reb.$
Indeed we have $\sui{1/n}\sq\Zo$ while the summable 
ideal $\sui{1/n}$ and the density-0 ideal 
$\Zo$ are known to be \dd\reb incomparable, see below.~\snos
{Some examples of this kind were recently found in the class
of Borel countable \eqr s, see \cite{ad,sim}.}

It would be interesting to figure out exact relationship
between $\reb$ and the \dd\sd reducibility $\rd.$
If the next questions answers in the negative
then the whole theory of
Borel reducibility for Borel ideals can be greatly simplified 
because reduction maps which are \dd\sd homomorphisms 
are much easier to deal with. 

\bqu
\lam{dred}
Is there a pair of Borel ideals $\cI,\cJ$ such that  
$\cI\reb\cJ$ but not $\cI\rd\cJ\,?$ 
\equ

\vyk{

\bde
\item[{\bfsl Additive reducibility\/}:]
$\rE\rea\rF$ if there is an additive reduction $\rE$ 
to $\rF$.
$\rE\reaa\rF$ if there is an asymptotically 
additive reduction $\rE$ to $\rF$.
\ede

\ble[{{\rm Farah~\cite{f-co}}}] 
\lam{arb} 
Suppose that\/ $\cI$ and\/ $\cJ$ are Borel ideals on\/ 
$\dN.$ 
Then\/ $\cI\orbpp\cJ$ iff\/ $\rei\rea\rej$.
\ele
(By definition $\rei$ and $\rej$ are \er s on $\pn,$ 
yet we can consider them as \er s on 
$\dn=\prod_{k\in\dN}\ans{0,1},$ 
as usual, which yields the intended 
meaning for $\rei\rea\rej.$)
\bpf
If $\cI\orbpp\cJ$ via a sequence of finite sets $w_i$ with 
$\tmax w_i<\tmin w_{i+1}$ then we put $n_0=0$ and 
$n_i=\tmin w_i$ for $k\ge1,$ 
so that $w_i\sq \il{i}{i+1},$ and, for any $i,$ 
put $H_i(0)=\il{i}{i+1}\ti\ans0$ and let $H_i(1)$ be 
the characteristic function of $w_i$ within $\il{i}{i+1}.$ 
Conversely, if $\rei\rea\rej$ via a sequence 
$0=n_0<n_1<n_2<\dots$ and a family of maps 
$H_i:\ans{0,1}\to2^{\il{i}{i+1}}$ then 
$\cI\orbpp\cJ$ via the sequence of sets  
$w_i=\ens{k\in\il{i}{i+1}}{H_i(0)(k)\ne H_i(1)(k)}$.
\epf
}

The remainder of the book will be concentrated on the
Borel reducibility/irreducibility theorems.
The following rather elementary lemma gives a couple
of examples.

\ble
\lam{l:rtd}
{\rm(i)} $\rav\dnn\eqb\epo{\rav\dN}.$\;
{\rm(ii)} $\rtd\eqb\aer\isg\nnn$
{\rm(see Example~\ref{ex:ac6})}.
\ele
\bpf
(i)
By definition $\epo{\rav\dN}$ is an \er\ on $\dnn$ and
$x\epo{\rav\dN}y$ holds iff $\ran x=\ran y.$ 
Thus the map $\vt(x)=\chi_{\ran x}$
(the characterictic function)
witnesses that $\epo{\rav\dN}\reb{\rav\dnn}.$
To prove the converse put, for $x\in\dnn,$ 
\dm
r(x)=\ans{x(0)\,,\;x(0)+x(1)+1\,,\;x(0)+x(1)+x(2)+2\,,\;\dots}\,;
\dm
then $\vt(x)=\chi_{r(x)}$ witnesses
${\rav\dnn}\reb\epo{\rav\dN}.$

(ii)
Suppose that $x\yi y\in\nnn.$
Then $x\rtd y$ means that
\dm
\kaz k\:\sus l\:(x(k)=y(l))
\quad\text{and}\quad
\kaz l\:\sus k\:(x(k)=y(l)),
\dm
while $x\aer\isg\nnn y$ means that there is a bijection
$f:\dN\onto\dN$ such that $x(k)=y(f(k))$ for all $k.$
The latter condition is, generally speaking, stronger, but
the two are equivalent provided for any $k$ there exist
infinitely many indices $l$ such that $x(k)=x(l)$ and the
same for $y.$ 
It follows that the map $\vt:\dnn\to\dnn,$ defined so that
$\vt(x)=x'$ iff $x'(2^n(2k+1)-1)=x(k)$ for all $n\yi k,$
is a Borel reduction of ${\rtd}$ to ${\aer\isg\nnn}.$

A Borel reduction $\vt$ of ${\aer\isg\nnn}$ to ${\rtd}$
can be defined as follows:
$\vt(x)=x',$ where $x'(k)=n_x(k)\we x(k)$ for all $k,$
$n_x(k)$ is the number of all $l$ satisfying $x(l)=x(k)$
(including $l=k$)
or $0$ if there exist infinitely many of such $l,$
and $n\we a$ for $a\in\dnn$ is defined as the only element
of $\dnn$ such that $(n\we a)(0)=n$ and
$(n\we a)(j+1)=a(j)$ for all $j.$
\epf

\punk{Borel, continuous, Baire measurable, additive reductions}
\las{b-c}

The Borel reducibility and related notions in \nrf{bored}
admit weaker Baire measurable (\bm, for brevity) versions,
which claims that the reduction postulated to exist is only
a \bm, not necessarily Borel, map.
(Recall that a map is \rit{Baire measurable} if the preimages
of open sets are sets with the Baire property.)
Those versions will be denoted with a subscript BM
instead of B.
Also there are stronger continuous versions, that will be
denoted with a subscript C. 
Thus

\bde
\item  
[$\rE\reB\rF\yt\rE\eqB\rF\yt\rE\reBs\rF$] 
\index{reducibility!Baire measurable@Baire measurable $\reB$}%
\index{zzErenF@$\rE\reB\rF$}%
\index{zzEemnF@$\rE\eqB\rF$}%
\index{zzEemnF@$\rE\reBs\rF$}%
mean the reducibility, resp., bi-reducibility,
strict reducibility 
by {\ubf Baire measurable} maps.

\item  
[$\rE\ren\rF\yt\rE\eqc\rF\yt\rE\rens\rF$] 
\index{reducibility!continuous@continuous $\ren$}%
\index{zzErenF@$\rE\ren\rF$}%
\index{zzErensF@$\rE\rens\rF$}%
\index{zzEemnF@$\rE\eqc\rF$}%
mean the reducibility, resp., bi-reducibility,  
strict reducibility 
by {\ubf continuous} maps.
\ede

It is known that a Baire measurable map defined on a
Polish space is continuous on a comeager set.
Thus \bm\ reducibility is the same as a Borel, or even
continuous reducibility on a comeager set.
On the other hand, according to the following result
of Just~\cite{justm} and Louveau~\cite{L94},
continuous reducibility on the full
domain can sometimes be derived from Borel reducibility.

\ble
\lam{l:bc}
{} \ 
If\/ $\cI$ is a Borel ideal on a countable\/ $A,$
$\rE$ an \eqr\ on a Polish space\/ $\dX,$ and\/ 
$\rei\reB\rE,$ then\/ 
$\rei\ren{\rE\ti\rE}$ 
{\rm(via a continuous reduction)}. 
In addition there is a set\/ $X\sq A\yt X\nin\cI$ such
that\/ ${\rE_{\cI\res X}}\ren\rE,$
where\/ $\cI\res X=\cI\cap\pws X$. 
\ele

Here $\rE\ti\rE$ is an \eqr\ on $\dX\ti\dX$ defined 
so that $\ang{x,y}$ and $\ang{x',y'}$ are equivalent 
iff both $x\rE x'$ and $y\rE y'.$
Note that $\rE\ti\rE\ren \rE$
holds for various \eqr s $\rE,$
and in such a case the condition 
$\rei\ren{\rE\ti\rE}$ in the theorem can be replaced by
$\rei\ren\rE.$


\bpf
We have to define continuous maps  
$\vt_0\zd\vt_1:\pws A\to\dX$ such that, for any\/ 
$x\yi y\in\pn,$ 
$x\sd y\in\cI$ iff both\/ $\vt_0(x)\rE\vt_0(y)$ 
and\/ $\vt_1(x)\rE\vt_1(y)$.
Suppose \noo\ that $A=\dN.$ 
Let $\vt:\pn\to\dX$ witness that $\rei\reB\rE.$  
Then $\vt$ is continuous on a dense $\Gd$ set 
$D=\bigcap_i D_i\sq\pn,$ all $D_i$ being dense open, and 
$D_{i+1}\sq D_i.$  
A sequence $0=n_0<n_1<n_2<\dots$ and, for any $i,$ a set 
$u_i\sq \il{i}{i+1}$ can be easily defined, by
induction on $i,$ so that 
$x\cap \il{i}{i+1}=u_i\imp x\in D_i.$~\snos
{Sets like $u_i$ are called {\it stabilizers\/}, 
\index{stabilizer}%
they are of much help in study of Borel ideals.}  
Let  
\dm
\textstyle
N_1=\bigcup_i \il{2i}{2i+1}\,\yt \; 
N_2=\bigcup_i \il{2i+1}{2i+2}\,\yt \; 
U_1=\bigcup_i u_{2i}\,\yt \; U_2=\bigcup_i u_{2i+1}\,.
\dm
Now set 
$\vt_1(x)=\vt((x\cap N_1)\cup U_2)$ and 
$\vt_2(x)=\vt((x\cap N_2)\cup U_1)$ for $x\sq\dN$.

To prove the second claim let $X$ be that one of the sets
$N_1\yd N_2$ which does not belong to $\cI.$
(Or any of them if neither belongs to $\cI.$)
Let say $X=N_1\nin\cI.$
Then the map $\vt_1$ proves ${\rE_{\cI\res X}}\ren\rE$. 
\epf

The following question should perhaps be answered in the 
negative in general and be open for some particular cases.

\bqu
Suppose that $\rE\reb \rF$ are Borel \er s. 
Does there always exist a {\it continuous\/} reduction ? 
\equ

There is a special useful type of continuous reducibility,
actually a \lap{clone} of the Rudin--Blass order of ideals
considered in \nrf{reli}. 

Suppose that $X=\prod_{k\in\dN}X_k$ and 
$Y=\prod_{k\in\dN}Y_k,$ the sets $X_i\yi Y_i$ are finite,
$0=n_0<n_1<n_2<\dots,$ and 
$H_i:X_i\to\prod_{n_i\le k<n_{i+1}}Y_k$ for any $i.$ 
Define 
\dm
\Psi(x)=H_0(x(0))\cup H_1(x(1))\cup H_2(x(2))\cup\dots\in Y
\dm
for each $x\in X.$ 
Maps $\Psi$ of this kind are called {\it additive\/} 
(Farah~\cite{f-co}). 
More generally, if, in addition, $0=m_0<m_1<m_2<\dots,$ and 
$H_i:\prod_{m_i\le j<m_{i+1}}X_j\to
\prod_{n_i\le k<n_{i+1}}Y_k$ for any $i,$ 
then define 
\dm
\Psi(x)=H_0(x\res\ir{m_0}{m_1})
\cup H_1(x\res\ir{m_1}{m_2})
\cup H_2(x\res\ir{m_2}{m_3})\cup\dots\in Y
\dm
for each $x\in X.$ 
Farah~\cite{f-co} calls maps $\Psi$ of this kind 
{\it asymptotically additive\/}. 
All of them are continuous functions $X\to Y$ 
in the sense of the product Polish topology.
(Recall that $X_i\yi Y_i$ are finite.)

Suppose now that $\rE$ and $\rF$ are \er s on resp.\ 
$X=\prod_kX_k$ and $Y=\prod_kY_k$.

\bde
\item[{\ubf Additive reducibility\/}:]
$\rE\rea\rF$ if there is an additive reduction of $\rE$ 
\index{reducibility!additive@additive $\rea$}%
\index{zzEleaF@$\rE\rea\rF$}%
to $\rF.$
As usual $\rE\eqa\rF$ means that simultaneously
\index{zzEeqaF@$\rE\eqa\rF$}%
$\rE\rea\rF$ and $\rF\rea\rE,$ while
$\rE\reas\rF$ means that $\rE\rea\rF$ but not $\rF\rea\rE$.
\index{zzEleasF@$\rE\reas\rF$}%

A version: $\rE\reaa\rF$ if there exists an asymptotically 
\index{reducibility!additive!asymptotically $\reaa$}%
additive reduction of $\rE$ to $\rF$.
\ede

The additive reducibility
coincides with $\orbpp$ on the domain of Borel ideals:

\ble[{{\rm Farah~\cite{f-co}}}] 
\lam{arb} 
Assume that\/ $\cI$ and\/ $\cJ$ are Borel ideals on\/ 
$\dN.$ 
Then\/ $\cI\orbpp\cJ$ iff\/ $\rei\rea\rej$.
\ele

By definition $\rei$ and $\rej$ are \eqr s on $\pn,$ 
however we can consider them as \er s on 
$\dn=\prod_{k\in\dN}\ans{0,1},$ 
as usual, which yields the intended 
meaning for the relation $\rei\rea\rej.$

\bpf
If $\cI\orbpp\cJ$ via a sequence of finite sets $w_i$ with 
$\tmax w_i<\tmin w_{i+1}$ then we put $n_0=0$ and 
$n_i=\tmin w_i$ for $k\ge1,$ 
so that $w_i\sq \il{i}{i+1},$ and, for any $i,$ 
put $H_i(0)=\il{i}{i+1}\ti\ans0$ and let $H_i(1)$ be 
the characteristic function of $w_i$ within $\il{i}{i+1}.$ 
Conversely, if $\rei\rea\rej$ via a sequence 
$0=n_0<n_1<n_2<\dots$ and a family of maps 
$H_i:\ans{0,1}\to2^{\il{i}{i+1}}$ then 
$\cI\orbpp\cJ$ via the sequence of sets  
$w_i=\ens{k\in\il{i}{i+1}}{H_i(0)(k)\ne H_i(1)(k)}$.
\epf

\punk{Diargam of Borel reducibility of key \eqr s}
\las{brb}

The diagram on page~\pageref{p-p} begins, at the low end, 
with cardinals $1\le n\in\dN\yt \alo\yt \cont,$ 
naturally identified with the \eqr\  
of equality on resp.\
finite (of a certain number $n$ of elements), countable,
uncountable Polish spaces.
As all uncountable Polish spaces are Borel isomorphic, 
the \eqr s $\rav\dX,$ $\dX$ a Polish space, are 
characterized, modulo $\reb,$
or even modulo Borel isomorphism between the domains,
by the cardinality of the domain,
which can be any finite $1\le n<\om,$ or 
$\alo,$ or $\cont=2^\alo.$


\begin{figure}[h]
\vspace*{1mm}

\setlength{\unitlength}{1mm} 
\begin{picture}(130,110)(0,0) 

\label{p-p}

\put(80,0){\tob}
\put(82,-1){$1$}
\put(80,0){\line(0,1){8}}
\put(80,8){\tob}
\put(82,7){$2=\rav{\ans{1,2}}$}
\put(80,9) {\toq}
\put(80,10){\toq}
\put(80,11){\toq}
\put(80,12){\toq}
\put(80,13){\toq}
\put(80,14){\toq}
\put(80,15){\toq}
\put(80,16){\tob}
\put(82,15){$n=\rav{\ans{1,2,...,n}}$}
\put(52,15){($1\le n<\alo$)}
\put(80,17){\toq}
\put(80,18){\toq}
\put(80,19){\toq}
\put(80,20){\toq}
\put(80,21){\toq}
\put(80,22){\toq}
\put(80,23){\toq}
\put(80,24){\tob}
\put(82,23){$\alo=\rav\dN$}
\put(80,24){\line(0,1){8}}
\put(80,32){\tob}
\put(82,31){$\cont=\rav{\dn}$}
\put(80,32){\line(0,1){8}}
\put(80,40){\tob}
\put(82,38){$\Eo$}

\put(80,40){\line(-3,1){30}}
\put(50,50){\line(-2,1){40}}
\put(10,70){\tob}
\put(10,70){\line(3,4){25}}
\put(4.5,69){$\Ei$}

\put(80,40){\line(-3,2){20.9}}
\put(55,56.7){\line(-3,2){20.3}}
\put(35,70){\tob}


\put(37,69){$\Ed\eqb{\bel1}$}
\put(55,55){\framebox{\sf ?}}

\put(80,40){\line(1,2){6.9}}
\put(95,70){\line(-1,-2){5.7}}
\put(95,70){\line(-1,1){19.5}}
\put(95,70){\tob}
\put(82.5,55){\framebox{\sf ctble}}
\put(96,69){$\Ey$}

\put(80,40){\line(3,2){45}}
\put(125,70){\tob}
\put(123,66){$\Et$}

\put(125,70){\line(0,1){14}}
\put(125,100){\line(0,-1){11}}
\put(125,100){\tob}
\put(117,102){$\rzo\eqb{\fco}$}
\put(118,86){\framebox{\sf \co eqs}}

\put(95,70){\line(0,1){30}}
\put(93,102){$\rtd$}
\put(95,100){\tob}
\put(95,100){\line(1,-1){30}}

\newlength{\borlen}
\settowidth{\borlen}{\small the {\sf non-P}\ }

\put(5,92){{\bmp{\borlen}{\small border of\\ the {\sf non-P}\\ domain}\emp}}

\put(35,70){\line(0,1){10}}
\put(35,85){\line(0,1){18}}
\put(35,103){\line(3,-1){40.5}}
\put(35,103){\tob}
\put(32.0,104){$\bel\iy$}
\put(31.5,81.2){\framebox{$\bel p$}}

\vyk{
\put(95,100){\mtir}
\put(91,100){\mtir}
\put(87,100){\mtir}
\put(80.5,99){$\ni$}
\put(79,100){\mtir}
\put(75,100){\mtir}
}


\linethickness{0.2mm}

\qbezier(9,68)(28,75)(39,109)

\put(17,90){\vector(1,-1){8.8}}

\end{picture}

\caption{\small Reducibility between the key \eqr s. 
Connecting lines here indicate
Borel reducibility of lower ERs to upper ones.
}%
\end{figure}
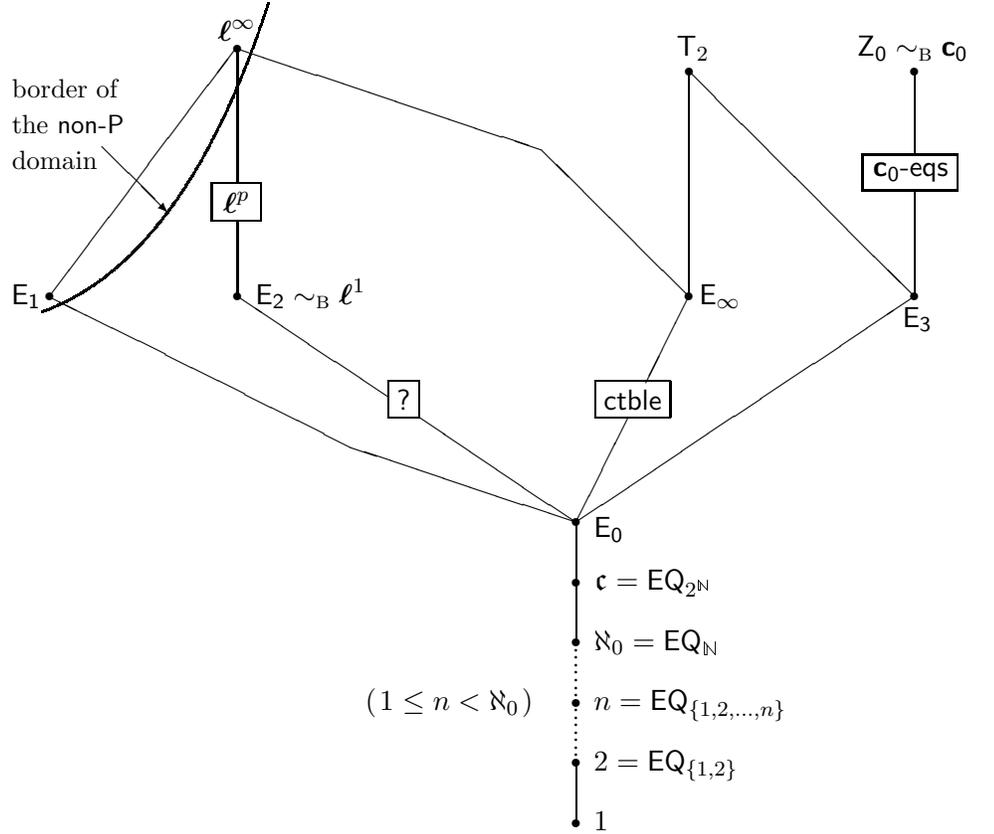
%

  
The linearity breaks above $\Eo:$ 
each one of the four \eqr s $\Ei\zd\Ed\zd\Et\zd\Ey$ of the 
next level is strictly \dd{\rebs}bigger than $\Eo,$ 
and they are pairwise \dd\reb incomparable 
with each other.

The framebox \fbox{\sf?} points on an
interesting open problem (Question~\ref{4dq} below).
The framebox \fbox{\sf\co eqs} denotes 
\co{\it equalities\/}, a family of Borel \er s introduced 
by Farah~\cite{f-co}, all of them are \dd\reb between 
$\Et$ and ${\fco}\eqb\rzo,$ and there is 
continuum-many \dd\reb incomparable among them.

The \lap{{\sf non-P} domain} denotes the family of all Borel
\er s that cannot be induced by a Polish action.
$\Ei$ belongs to this family, and it is conjectured that
$\Ei$ is a \dd\reb least \er\ in this family. 
Solecki~\cite{sol,sol'} proved this conjecture for \er s  
generated by Borel ideals: for instance for a Borel ideal
$\cI$ to be 
\imar{$\Ei$ and Polish grps action problem}%
not a P-ideal it is necessary and sufficient that 
$\Ei\reb\rei.$
See \grf{idI1} for more details.

Finally, the framebox \fbox{\sf ctble} denotes the family of
all Borel countable \er s
(meaning that equivalence classes are at most countable);
all of them are Borel reducible to $\Ey.$
The following theorem of Adams -- Kechris \cite{ak} shows that
this is quite a rich family. 

\bte
[{\rm not to be proved here}]
\lam{akth}
There is continuum many 
pairwise \dd\reb incomparable countable Borel \eqr s.
\ete

A somewhat weaker result that implies the existence of
continuum many pairwise \dd\reb incomparable
(not necessarily countable) Borel \eqr s will be established
by Theorem~\ref{coAB}.

\punk{Reducibility and irreducibility on the diagram}
\las{rir}

Here we discuss, without going into technicalities,
the structure of the diagram on page~\pageref{p-p}  and
related theorems.

Recall that any straight line on the diagram indicates
Borel reducibility of the \er\ at the lower end to the
\er at the upper end.
Some of these reducibility claims are witnessed by a
simple and obvious reductions.
Slightly less obvious are reductions of $\Ey$ and $\Et$
to $\rtd$ and $\Et$ to $\fco,$ see lemmas \ref{e3'}, \ref{e3"}.
Finally, to prove that $\Ei\yi\Ey,$ and all of $\bel p$
(including ${\bel1}\eqb\Ed$),
are Borel reducible to $\bel\iy,$ we apply 
Rosendal's theorem in \cite{roz} saying that $\bel\iy$ is a
\dd\reb largest $\Fs$.

That $\Ed\eqb{\bel1}$ and ${\fco}\eqb\rzo$ see lemmas \ref{co=d}
and \ref{l1=s}.

See Theorem~\ref{hd-t} on the equivalence
${{\bel p}\reb{\bel q}}\eqv{p\le q}$.

It is the most interesting question whether 
the diagram on page~\pageref{p-p} is complete in the
sense that there is no Borel reducibility interrelations
between the \er s mentioned in the diagram except for those
explicitly indicated by straight lines.
Some of these irreducibility claims are trivial by a simple
cardinality argument:
clearly an \er\ $\rE$ having strictly more equivalence classes
than $\rF$ is not Borel reducible to $\rF.$

However this argument is not applicable in more complicated
cases, beginning with the irreducibility claim
$\Eo\not\reb \rav\dn:$
each of the two relations has exactly continuum-many
classes.
Here we have to employ the borelness.
Suppose towards the contrary that $\vt:\dn\to\dn$ is a Borel
reduction of $\Eo$ to $\rav\dn.$
Then the pre-image $\ens{x}{\vt(x)=y}$ of any $y\in\dn$ is
countable (or empty).
We conclude,
using some classical theorems of descriptive set theory
(theorems \ref{cuni} and \ref{cpro})
that there is a Borel set
$T\sq\dn$ having exactly one element with each \dd\Eo class.
%
But this contradicts the borelness of $T$ --- see a short
argument after Example~\ref{ex:ac2}.~\snos
{Alternatively, one can derive $\rav\dn\reb\Eo$   
from an old result of Sierpi\'nski~\cite{sie}: 
any linear ordering of 
all \dd\Eo classes yields a Lebesgue non-measurable set of 
the same descriptive complexity as the given ordering.}

As for the rest of the diagram, to establish its completeness 
one has to prove the following irreducibility claims:\vtm

$
\bay[t]{clclllll}
(1) & \Ei &\not\reb\;:& \Ed, & \rtd, & \fco;\\[0.7\dxii]

(2) & \bel\iy &\not\reb\;:& \Ei, & \Ed, & \rtd, & \fco;\\[0.7\dxii]

(3) & \Ed &\not\reb\;:& \Ei, & \rtd, & \fco;\\[0.7\dxii]

(4) & \Ey &\not\reb\;:& \Ei, & {\Ed}
, &
{\fco} &
 \text{--- \ubf\ this group contains open problems};\\[0.7\dxii]

(5) & \Et &\not\reb\;:& \bel\iy;\\[0.7\dxii]

(6) & \rtd &\not\reb\;:& \bel\iy, & \fco;\\[0.7\dxii]

(7) & \fco &\not\reb\;:& \bel\iy, & \rtd.\\[0.7\dxii]
\eay
$\vtm

Beginning with (1), we note that $\Ei$ is not Borel
reducible to any \eqr\ induced by a Polish action
by Theorem~\ref{e1pga}
(Kechris -- Louveau).
On the other hand, $\Ed\zi\rtd\zi{\fco}$ obviously
belong to this category of \er s.

(2) follows from (1) and (3) since $\Ei\reb{\bel\iy}$ and
$\Ed\reb{\bel\iy}$.

The result $\Ed\not\reb{\fco}$ in (3) is Hjorth's 
Theorem~\ref{t:>s}\ref{t>s2}. 
The result $\Ed\not\reb\Ei$ (Corollary~\ref{ne<e1})
will be established by
reference to Kechris' Theorem~\ref{rig1} on the structure
of ideals Borel reducible to $\Ei$.

The results $\Ed\not\reb\rtd$ and ${\fco}\not\reb\rtd$ in
(3) and (7) are proved below in \grf{groat}
(Corollary~\ref{-t2}); this will involve turbulence
theory by Hjorth and Kechris.

The result of (5) is Lemma~\ref{e3}.
It also implies ${\fco}\not\reb{\bel\iy}$ in (7).

(6) was obtained by Hjorth, see \grf{rtd}.

This leaves us with (4).
We don't know how to prove $\Ey\not\reb\Ei$
easily and directly.
There are two indirect ways.
The first one is to apply some results in the theory of
countable and hyperfinite \eqr s --- see Corollary~\ref{eiey}.
The second one is based on theorems \ref{kelu}
(3rd dichotomy) and \ref{e1pga} --- see Corollary~\ref{eiey*}.

\bqe 
\lam{liy2co}
Is $\Ey$ Borel reducible to $\fco$? to $\bel1$ or any
other $\bel p$?
\eqe

A related question whether $\Ey$ is Borel reducible to $\Et$
answers in the negative on the base of 6th dichotomy theorem
by Corollary~\ref{etey}.

The irreducibility results in (1) -- (7)
can be partitioned into two rather distinct categories.
The first group consists of those having proofs that
involve only common methods of descriptive set theory,
as the proof of $\Eo\not\reb\rav\dn$ outlined above. 
This includes such results as
$\Ed\not\reb{\fco},$ 
${\bel\iy}\not\reb{\fco},$ 
$\Et\not\reb{\bel\iy},$ 
${\fco}\not\reb{\bel\iy},$ and also
$\Ed\not\reb\Ei$ as a transitional claim between the first
and second groop: it refers to Theorem~\ref{rig1}, a special
result on the \dd\reb structure of ideals below $\Ii,$
rather complicated but still based on classics of
descriptive set theory.

Note that some results in this group belong to the earliest
of this type.
For instance Just proved that $\Ed$ is
mutually \dd\reb irreducible with $\rzo$ \cite{just}
and with $\rE_{\fpd\ifi\ifi}$ \cite{justm}.
According to \cite[1.4]{hypsm} the irreducibility claim
$\Ei\not\reb\Ey$ goes back to even earlier paper \cite{fr}.

The other group consists of irreducibility results that
involve (as far as we know) methods that definitely go 
beyond common tools of descriptive set theory.
This includes such resuts as
$\Ei\not\reb{\Ed},$ 
$\Ei\not\reb{\rtd},$ 
$\Ei\not\reb{\fco},$
based on the fact that
$\Ei$ is not reducible to a Polish action (Theorem~\ref{e1pga}),
$\Ed\not\reb{\rtd}$ and ${\fco}\not\reb{\rtd}$ based on the
turbulence theory,
$\Ey\not\reb\Ei$ and $\Ey\not\reb\Et$ based on resp.\ 3th and 6th
dichotomy theorems (see the next \paf), and finally 
$\rtd\not\reb{\bel\iy}$ and $\rtd\not\reb{\fco}$ based on the
theory of \rit{pinned} equivalence relations (\grf{rtd}).

\punk{Dichotomy theorems}
\las{dihot}

Another general problem related to the diagram is the
\dd\reb structure of certain domains, for instance,
\dd\reb intervals between adjacent \eqr s.
Some results in this direction are known as dichotomy
theorems because of their distinguished dichotomical form.

\bde
\itsep
\item[{\ubf 1st dichotomy\/}] 
(
Theorem~\ref{1dih} below).  
{\it 
Any Borel, even any\/ $\fp11$ \eqr\/ $\rE$ \bfei\/ 
has at most countably many equivalence classes, formally, 
$\rE\reb\rav\dN,$ \bfor\/ satisfies\/
$\rav\dn\reb\rE$.\/}
\ede

Thus not only the strict \dd\rebs interval between the
\er s $\alo=\rav\dN$ and $\cont=\rav\dn$ is empty, but
the union of the lower \dd\reb cone of the former and
the upper \dd\reb cone of the latter cover the whole
family of Borel \eqr s!
The same is true for the next \dd\rebs interval:

\bde
\item[{\ubf 2nd dichotomy\/}] 
(
Thm~\ref{2dih}).  
{\it 
Any Borel \er\/ $\rE$ satisfies\/ \bfei\/  
$\rE\reb\cont$ \bfor\/ $\Eo\reb\rE.$\/} 
\ede

What is going on in the \dd\rebs intervals
between $\Eo$ and the \eqr s
${\Ei}\yd{\Ed}\yt{\Et}$?
The following dichotomy theorems provide some answers. 

\bde
\itsep
\item[{\ubf 3rd dichotomy\/}] 
(%
Theorem~\ref{kelu}). 
{\it Any \eqr\/ $\rE\reb\Ei$ satisfies\/ \bfei\ 
$\rE\reb\Eo$ \bfor\/ $\rE\eqb\Ei.$} 

\item[{\ubf 4th dichotomy\/}] 
(
Theorem~\ref{h}).  
{\it Any \eqr\/ $\rE\reb\Ed$ \bfei\ is 
essentially countable \bfor\ satisfies \/ $\rE\eqb\Ed$.}
\ede

An \eqr\ $\rE$ is \rit{essentially countable} iff it is Borel
reducible to a Borel countable
(\ie, with at most countable equivalence classes) \er.
The \bfei\ case in 4th dichotomy remains mysterious.
This is marked by the framebox 
\framebox{\sf ?} on the diagram. 

\bqe 
\lam{4dq}
In 4th dichotomy, can the \bfei\ case be
improved to $\reb\Eo$?
\eqe

The fifth dichotomy theorem is a bit more special, it will
be addressed below.

\bde
\itsep
\item[{\ubf 6th dichotomy\/}] 
(
Theorem~\ref{6dih}).  
{\it Any \eqr\/ $\rE\reb\Et$ satisfies\/ \bfei\/
$\rE\reb\Eo$ \bfor\/ $\rE\eqb\Et$.}
\ede

On the other hand, the interval
between $\Eo$ and $\Ey$ contains all countable Borel \er s
and among them plenty of pairvise \dd\eqb inequivalent \er s
by Theorem~\ref{akth}.

It was once considered \cite{hk:nd} as a plausible
hypothesis that any Borel \er\ which is not $\reb\Ey,$ \ie, 
not essentially countable, satisfies
$\rE_i\reb\rE$ for at least one $i=1,2,3.$ 
This turns out to be not the case: 
Farah~\cite{f-bas,f-tsir} and Velickovic~\cite{vel:tsir}
found an independent family of uncountable Borel \er s, 
based on {\it Tsirelson ideals\/}, \dd\reb incomparable 
with $\Ei\zi\Ed\zi\Et$.

\bqe
\lam{min}
It there any reasonable ``basis'' of Borel \er s
above $\Eo$?
\eqe

\punk{Borel ideals in the structure of Borel reducibility}
\las{bibe}

Some of \eqr s on the diagram are obviously generated by
Borel ideals, for some other ones this is not clear.
This leads to the question what is the place of Borel ideals
in the whole structure od Borel \eqr s.
The answer obtained in the studies of last years can be
formulated as follows:
Borel ideals are \dd\reb cofinal, but rather rare, in the
\dd\reb structure of Borel \er s.
We prove the following theorem, the cofinality claim of which
is due to Rosental~\cite{roz}
(Theorem~\ref{M} in \grf{rosen})
while on the other claim see Corollary~\ref{t2i}.

\bte
\label{troz}
For any Borel \eqr\/ $\rE$ there exists a Borel ideal\/
$\cI\sq\pn$ such that\/ $\rE\reb\rei.$
On the other hand there is no Borel ideal\/
$\cI$ such that\/ $\rtd\eqb\rei.$ 
\ete

\api

\parf{``Elementary'' stuff}
\las{BAN}

This \gla\ gathers proofs of some
reducibility/irreducibility results related to the diagram
on page~\pageref{p-p}, elementary in the sense that they do not
involve any special concepts.
Some of them are really simple, as \eg\ some lemmas on $\Et$ and
$\rtd$ in \nrf{e3t2} or the equivalences ${\fco}\eqb\rzo$ and
$\Ed\eqb{\bel1}$ in \nrf{e2i}, while
some other quite tricky.
The latter category includes Hjorth's theorem on the irreducibility
of nontrivial summable ideals to $\fco$ in \nrf{summ}, interrelations
in the family of \eqr s $\bel b$ in \nrf{hd1}, and the
\dd\reb universality of $\bel\iy$ in the class of all $\Fs$ \eqr s
in \nrf{liy}.
\imar{add v russkij text}

\punk{$\Et$ and $\rtd$}
\las{e3t2}

These \eqr s, together with ${\fco}\eqb\rzo,$ are the only
non-$\fs02$ equivalences explicitly mentioned on the
diagram.

\ble
\lam{e3}
$\Et$ is Borel irreducible to\/ ${\bel\iy}.$ 
\ele
\bpf
Suppose towards the contrary that $\vt:\dntn\to\rtn$
is a Borel reduction of $\Et$ to\/ ${\bel\iy}.$~\snos
{Recall that, for $x,y\in\dntn,$ $x\Et y$ means  
$\seq xk\Eo\seq yk\zd\kaz k,$ where
$\seq xk\in\dn$ is defined by $\seq xk(n)=x(k,n)$
for all $n$ while $a\Eo b$ means that
$a\sd b=\ens m{a(m)\ne b(m)}$ is finite.}
Since obviously ${\bel\iy}\eqb{\bel\iy}\ti{\bel\iy},$ 
Lemma~\ref{l:bc}
reduces the general case to the case of continuous $\vt.$
Define $\fo,\fr\in\dn$ by 
$\fo(n)=0\zd\fr(n)=1\zd\kaz n.$ 
Define $\bbo\in\dntn$  by
$\fo(k,n)=0$ for all $k,n,$ thus $\seq\bbo k=\fo\zd\kaz k.$
Finally, for any $k$ define $\fz_k\in\dn$ by
$\fz_k(n)=1$ for $n<k$ and $\fz_k(n)=0$ for $n\ge k$.

We claim that there are increasing sequences
of natural numbers $\sis{k_m}{}$ and $\sis{j_m}{}$ 
such that $|\vt(x)(j_m)-\vt(\bbo)(j_m)|>m$
for any $m$ and any $x\in\dntn$ satisfying
\dm
\seq x{k}=
\left\{
\bay{rl}
\fz_{k_i} & \text{whenever }\,i<m\,\text{ and }\,
k=k_i\\[0.7\dxii]

\fo & \text{for all }\,k<k_m\,
\text{not of the form }\,k_i.
\eay
\right.
\dm
To see that this implies contradiction
define $x\in\dntn$ so that
$\seq x{k_i}=\fz_{k_i}\zd\kaz i$ and
$\seq x{k}=\fo$ whenever $k$ does
not have the form $k_i.$
Then obviously $x\Et\bbo,$ but
$|\vt(x)(j_m)-\vt(\bbo)(j_m)|>m$ for all $m,$ hence
$\vt(x)\bel\iy\vt(\bbo)$ fails, as required.

We put $k_0=0.$ 
To define $j_0$ and $k_1,$ consider $x_0\in\dntn$ defined
by $\seq{x_0}0=\fr$ but $\seq{x_0}k=\fo$ for all $k\ge1.$
Then $x_0\Et \bbo$ fails, and hence
$\vt(x_0)\bel\iy\vt(\bbo)$ fails either.
Take any $j_0$ with $|\vt(x_0)(j_0)-\vt(\bbo)(j_0)|>0.$
As $\vt$ is continuous, there is a number $k_1>0$ such that
$|\vt(x)(j_0)-\vt(\bbo)(j_0)|>0$
holds for any $x\in\dntn$ with $\seq x0=\fz_{k_1}$
and $\seq xk=\fo$ for all $0<k<k_1$. 

To define $j_1$ and $k_2,$ consider $x_1\in\dntn$
defined so that $\seq{x_1}0=\fz_{k_1},$   
$\seq{x_1}k=\fo$ whenever $0<k<k_1,$ and
$\seq{x_1}{k_1}=\fr.$
Once again there is a number $j_1$ with
$|\vt(x_1)(j_1)-\vt(\bbo)(j_1)|>1,$
and a number $k_2>k_1$ such that
$|\vt(x)(j_1)-\vt(\bbo)(j_1)|>1$
for any $x\in\dntn$ with $\seq x0=\fz_{k_1},$
$\seq x{k_1}=\fz_{k_1},$
and $\seq xk=\fo$ for all $0<k<k_1$ and $k_1<k<k_2$.

Et cetera.
\epf

\ble
\lam{e3'}
$\Et$ is Borel reducible to both\/ $\rtd$ and\/ $\fco$.
\ele
\bpf
(1)
If $a\in\dn$ and $s\in\dln$ then define $sx\in\dn$ by
$(sx)(k)=x(k)+_2s(k)$ for $k<\lh s$ and
$(sx)(k)=x(k)$ for $k\ge\lh s.$
If $m\in\dN$ then $m\we x\in\dn$ denotes the concatenation.
In these terms, if $x,y\in\dntn$ then obviously
\dm
x\Et y\leqv
\ens{m\we{(s\seq xm)}}{s\in\dln\yt m\in\dN}=
\ens{m\we{(s\seq ym)}}{s\in\dln\yt m\in\dN}.
\dm
Now any bijection $\dln\ti\dN\onto\dN$ yields a Borel
reduction of $\Et$ to $\rtd$.

(2)
To reduce $\Et$ to $\fco$ consider a Borel map
$\vt:\dntn\to\rtn$ such that
$\vt(x)(2^n(2k+1)-1)=n\obr \seq xn(k)$.
\epf

\ble
\lam{e3"}
Any countable Borel ER is Borel reducible to\/ $\rtd$.
\ele
\bpf
Let $\rE$ be a countable Borel ER on $\dn.$
It follows from \Cenu\ that there is a Borel map
$f:\dn\ti\dN\to\dn$ such that 
$[a]_{\rE}=\ens{f(a,n)}{n\in\dN}$ for all $a\in\dn.$ 
The map $\vt$ sending any $a\in\dn$ to $x=\vt(a)\in\dntn$
such that $\seq xn=f(a,n)\zd\kaz n,$ is a reduction required.
\epf

See further study on $\rtd$ in \grf{rtd}, where
it will be shown that $\rtd$ is not Borel reducible to
a big family of \eqr s that includes
${\fco}\zi{\bel p}\zi{\bel\iy}\zi{\Ei}\zi{\Ed}\zi{\Et}\zi{\Ey}.$
On the other hand, the \eqr s in this list,
with the exception of ${\Et}\zi{\Ey},$ 
are not Borel reducible to $\rtd$ --- this follows from the
turbulence theory presented in \grf{groat}.

\punk{Discretization and generation by ideals}
\las{e2i}

Some \eqr s on the diagram on page~\pageref{p-p} are
explicitly generated by ideals, like
$\rE_i\yt i=0,1,2,3.$
Some other \er s are defined differently.
It will be shown below (\grf{rosen})
that {\ubf any} Borel \er\ $\rE$ is Borel reducible to a
\er\ of the form $\rei,$ $\cI$ a Borel ideal.
On the other hand, $\fco\zi\bel1\zi\bel\iy$ turn out to
be Borel equivalent to some meaningful Borel ideals. 
Moreover, these \eqr s admit ``discretization'' by means
of restriction to certain subsets of $\rtn.$

\bdf
\lam{dx}
We define 
$\dX=\prod_{n\in\dN}X_n=\ens{x\in\rtn}{\kaz n\:(x(n)\in X_n)},$
where
$X_n=\ans{\frac0{2^n},\frac1{2^n},\dots,\frac{2^n}{2^n}}$.
\index{zzXXn@$\dX=\prod_nX_n$}%
\edf

\ble
\lam{redx}
${\fco}\reb {{\fco}\res{\dX}}$ and\/
${\bel p}\reb {{\bel p}\res{\dX}}$ for any\/ $1\le p<\iy.$

On the other hand, ${\bel\iy}\reb {{\bel\iy}\res{\ztn}}$.
\ele
\bpf
We first show that ${\fco} \reb {{\fco}\res{\oin}}.$
Let $\pi$ be any bijection of $\dN\ti\dZ$ onto $\dN.$ 
For $x\in\rtn,$ define $\vt(x)\in\oin$ as follows. 
Suppose that $k=\pi(n,\eta)$ ($\eta\in\dZ$). 
If $\eta\le x(n)<\eta+1$ then let $\vt(x)(k)=x(n).$ 
If $x(n)\ge \eta+1$ then put $\vt(x)(k)=1.$ 
If $x(n)<\eta$ then put $\vt(x)(k)=0.$ 
Then $\vt$ is a Borel reduction of
${\fco}$ to ${\fco}\res{\oin}.$
Now we prove that ${{\fco}\res{\oin}}\reb {{\fco}\res{\dX}}.$
For $x\in\oin$ define $\psi(x)\in\dX$ so that $\psi(x)(n)$
the largest number of the form $\frac i{2^n}\zd 0\le i\le{2^n}$
smaller than $x(n).$
Then obviously $x\fco {\psi(x)}$ holds for any $x\in\oin,$
and hence $\psi$ is a Borel reduction of ${\fco}\res{\oin}$ to
${\fco}\res{\dX}$.

Thus ${\fco}\reb{{\fco}\res{\dX}},$ 
and hence in fact ${\fco} \eqb {{\fco}\res\dX}.$

The argument for $\bel1$ is pretty similar.
The result for $\bel\iy$ is obvious:
given $x\in\rtn,$ replace any $x(n)$ by the largest integer
value $\le x(n)$.

The version for $\bel p\yi1<p<\iy,$ needs some comments in
the first part (reduction to $\oin$).
Note that if $\eta\in\dZ$ and
$\eta-1\le x(n)<\eta<\za\le y(n)<\za+1$
then the value $(y(n)-x(n))^p$ in the distance 
$\nor{y-x}_p=(\sum_n|y(n)-x(n)|^p)^{\frac1p}$
is replaced by $(\za-\eta)+(\eta-x(n))^p+(y(n)-\za)^p$
in $\nor{\vt(y)-\vt(x)}_p.$
Thus if this happens infinitely many times then both
distances are infinite, while otherwise this case can be
neglected.
Further, if $\eta-1\le x(n)<\eta\le y(n)<\eta+1$
then $(y(n)-x(n))^p$ in $\nor{y-x}_p$
is replaced by $(\eta-x(n))^p+(y(n)-\eta)^p$
in $\nor{\vt(y)-\vt(x)}_p.$
However
$(\eta-x(n))^p+(y(n)-\eta)^p\le (y(n)-x(n))^p\le
2^{p-1} ((\eta-x(n))^p+(y(n)-\eta)^p),$
and hence these parts of the sums in $\nor{y-x}_p$
and $\nor{\vt(y)-\vt(x)}_p$ differ from each other by
a factor between $1$ and $2^{p-1}.$
Finally, if $\eta\le x(n)\zi y(n)<\eta+1$ for one and the
same $\eta\in\dZ$ then the term $(y(n)-x(n))^p$ in
$\nor{y-x}_p$ appears unchanged in $\nor{\vt(y)-\vt(x)}_p.$
Thus totally $\nor{y-x}_p$ is finite iff so is
$\nor{\vt(y)-\vt(x)}_p$.
\epf

\ble
[{{\rm Oliver \cite{odis}}}]
\lam{co=d}
${\fco}\eqb{\rzo}.$ \
{\rm(Recall that $\rzo=\rE_{\zo}$.)}
\ele
\bpf
Prove that ${\fco}\reb \rzo.$
It suffices,
by Lemma~\ref{redx}, to define a Borel reduction 
${{\fco}\res\dX}\to {\rzo},$
\ie, a Borel map $\vt:\dX\to\pn$ such 
that ${x\fco y}\leqv{\vt(x)\sd\vt(y)\in\zo}$ 
for all $x\yi y\in\dX.$ 
Let $x\in\dX.$ 
Then, for any $n,$ we have $x(n)=\fras{k(n)}{2^n}$ 
for some natural $k(n)\le 2^n.$ 
The value of $k(n)$ determines the intersection 
$\vt(x)\cap\ir{2^n}{2^{n+1}}:$
for each $j<2^n,$ we define $2^n+j\in\vt(x)$ iff $j<k(n).$
Then   
$x(n)=\frac{\#(\vt(x)\cap\ir{2^n}{2^{n+1}})}{2^n}$ 
for any $n,$ and moreover  
$|y(n)-x(n)|=
\fras{\#([\vt(x)\sd\vt(y)]\cap\ir{2^n}{2^{n+1}})}{2^n}$ 
for all $x,y\in\dX$ and $n.$
This easily implies that $\vt$ is as required.

To prove $\rzo\reb {\fco},$ 
we have to define a Borel map $\vt:\pn\to\rtn$ such 
that ${x\sd x\in\zo}\leqv{\vt(x)\fco\vt(x)}.$ 
Most elementary ideas like $\vt(x)(n)=\frac{\#(x\cap\ir0n)}n$ 
do not work, the right way is based on the following 
observation: for any sets $s\yi t\sq\ir0n$ to satisfy 
$\#(s\sd t)\le k$ it is necessary and sufficient that 
$|\#(s\sd z)-\#(t\sd z)|\le k$ 
for any $z\sq\ir0n.$ 
To make use of this fact, let us fix an enumeration 
(with repetitions) $\sis{z_j}{j\in\dN}$ 
of all finite subsets of $\dN$ such that 
\dm
\ens{z_j}{2^n\le j<2^{n+1}}
\quad=\quad\hbox{all subsets of $\ir0n$}
\dm
for every $n.$ 
Put, for any $x\in\pn$ and $2^n\le j<2^{n+1},$ 
$\vt(x)(j)=\frac{\#(x\cap z_j)}n.$ 
Then $\vt:\pn\to[0,1]^\dN$ is a required reduction. 
\epf

Recall that for any sequence of reals $r_n\ge0,$
$\ern$ is an \eqr\ on $\pn$ generated by
the ideal $\srn=\ens{x\sq\dN}{\sum_{n\in x}r_n<\piy}$.

\ble
[{{\rm Attributed to Kechris in \cite[2.4]{h-ban}}}]
\lam{l1=s}
If\/ $r_n\ge 0\yt r_n\to0\yt\sum_nr_n=\piy$ then\/
$\ern\eqb {\bel1}.$
In particular,\/ $\Ed=\eun$ satisfies\/
$\Ed\eqb{\bel1}.$
\ele
\bpf
To prove $\ern\reb {\bel1},$ 
define $\vt(x)\in\rtn$ for any $x\in\pn$ as 
follows: $\vt(x)(n)=r_n$ for any $n\in x,$ and
$\vt(x)(n)=0$ for any other $n.$ 
Then ${x\sd y\in\sui{r_n}}\leqv{\vt(x)\bel1\vt(y)},$ 
as required.

To prove the other direction, it suffices
to define a Borel reduction of 
${\bel1\res\dX}$ to $\ern.$
We can associate a (generally, infinite) set $s_{nk}\sq\dN$
with any pair of $n$ and $k<2^n,$ so that 
the sets $s_{nk}$ are pairwise disjoint and 
$\sum_{j\in s_{nk}}r_j=2^{-n}.$ 
The map   
$\vt(x)=\bigcup_n\bigcup_{k<2^nx(n)}s_{nk}\yt x\in\dX,$
is the reduction required.
\epf

\punk{Summables irreducible to density-0}
\las{summ}

The \dd\reb independence of $\bel1$ and $\fco,$ two
best known ``Banach'' \eqr s, is quite important.
In one direction it is provided by \ref{t>s2} of
the next theorem.
As for the other direction, Lemma~\ref{e3}
contains an even stronger irreducibility claim.


Is there any example of Borel ideals $\cI\reb\cJ$ 
which do not satisfy $\cI\rd\cJ$? \  
Typically the reductions found 
to witness $\cI\reb\cJ$ are \dd\sd homomorphisms, and 
even better maps. 
The following lemma proves that Borel reduction yields 
\dd\orbpp reduction in quite a representative case.
Suppose that $\cI,\cJ$ are ideals over $\dN.$
Let us say that $\cI\orbpp\cJ$ {\it holds exponentially\/}
\index{reducibility!Rudin -- Blass!$\orbpp$ exponentially}%
if there exist a sequence of natural numbers $k_i$ with 
and $k_{i+1}\ge 2k_i$ and a sequence of sets
$w_i\sq\ir{k_i}{k_{i+1}}$ that withesses $\cI\orbpp\cJ$ 
--- in other words, the equivalence
$A\in\cI\leqv w_A=\bigcup_{k\in A}w_k\in\cJ$
holds for any $A\sq\dN$. 

\bte
\lam{t:>s}
Suppose that\/ $r_n\ge 0\yt r_n\to0\yt\sum_nr_n=\piy.$
Then
\ben
\renu
\itsep
\itla{t>s1}
{\rm(Farah~\cite[2.1]{f-tsir})} \ 
If\/ $\cJ$ is a Borel P-ideal and\/ $\sui{r_n}\reb\cJ$
then we have\/ $\sui{r_n}\orbpp\cJ$ exponentially$;$

\itla{t>s2}
{\rm(Hjorth~\cite{h-ban})} \
$\sui{r_n}$ is not 
Borel-reducible to\/ $\zo$.
\een
\ete
\bpf
\ref{t>s1}
Let a Borel map $\vt:\pn\to\pn$ witness $\sui{r_n}\reb\cJ.$ 
Let, according to Theorem~\ref{sol}, $\nu$ be a \lsc\ 
submeasure on $\dN$ with $\cJ=\Exh_\nu.$ 
The construction makes use of stabilizers. 
%
Suppose that $n\in\dN.$ 
If $u\yi v\sq\ir0n$ then
$(u\cup x)\sd(v\cup x)=u\sd v\in\sui{r_n}$ 
for any $x\sq\ir n\piy,$ therefore  
$\vt(u\cup x)\sd\vt(v\cup x)\in\cJ.$ 
It follows, by the choice of the submeasure $\nu,$
that for any $\ve>0$ 
there are numbers $n'>k>n$ and a set $s\sq\ir n{n'}$ 
such that    
\dm
\nu\skl(\vt(u\cup s\cup x)\sd\vt(v\cup s\cup x))
\cap\ir k\iy\skp\,<\,\ve
\dm
holds for all $u\yi  v\sq\ir0n$ 
and all generic $x\sq\ir{n'}\iy$.

\bre
\lam{ff:zhc} 
In the course of the proof, \lap{generic} means Cohen-generic
over a fixed countable transitive model $\mm$ of $\zhc,$ 
\index{theory!zfhc@$\zhc$}%
\index{zfhc@$\zhc$}%
\index{zzzfhc@$\zhc$}%
the theory containing all of $\ZFC$ minus the Power Set axiom 
but plus the axiom: 
\lap{for every set $X,$ the countable power set
$\pwc X=\ens{y\sq X}{\card y\le\alo}$ exists}.\snos
{In fact generic points are precisely those which avoid
certain meager $\Fs$ sets, but the genericity assumption
allows us not to specify those sets explicitly,
giving instead a reference to $\mm$ where all
relevant meager $\Fs$ sets have to be coded.}

Note that Cohen-generic extensions of such a model are still
models of $\zhc.$ 

We require that in addition $\mm$ contains all relevant
real-type objects,
together with codes of all relevant Borel sets.
In particular, in the case considered, $\mm$ contains the
sequence $\sis{r_n}{n\in\dN}$ and also contains Borel codes
of the ideal $\cI$ and of the map $\vt$.
\ere

This allows us to define an increasing sequence of 
natural numbers 
$0=k_0=a_0<b_0<k_1<a_1<b_1<k_2<\dots$ and, for any 
$i,$ a set $s_i\sq\ir{b_i}{a_{i+1}}$ such that, for all 
generic $x\yi x\sq\ir {a_{i+1}}\iy$ and all $u\yi v\sq\ir0{b_i},$
we have 
\ben
\tenu{(\arabic{enumi})}
\itla{stab1}\msur
$\nu\skl(\vt(u\cup s_i\cup x)\sd\vt(v\cup s_i\cup x))
\cap\ir{k_{i+1}}\iy\skp<2^{-i}$; 

\itla{stab2}\msur
$\skl\vt(u\cup s_i\cup x)\sd\vt(u\cup s_i\cup y)\skp
\cap\ir0{k_{i+1}}=\pu$;

\itla{stabg}
any $z\sq\dN,$ satisfying $z\cap\ir{b_i}{a_{i+1}}=s_i$ 
for infinitely many $i,$ is generic;

\itla{stabk}\msur
$k_{i+1}\ge 2k_i$ for all $i$;
\een
and in addition, under the assumptions on $\sis{r_n}{}$, 
\ben
\tenu{(\arabic{enumi})}
\addtocounter{enumi}4
\itla{stab3}
there is a set $g_i\sq\ir{a_i}{b_i}$ such that 
$|r_i-\sum_{n\in g_i}r_n|<2^{-i}$.
\een
It follows from \ref{stab3} that 
$a\mapsto g_a=\bigcup_{i\in a}g_i$ 
is a reduction of $\sui{r_n}$ to $\sui{r_n}\res N,$ where 
$N=\bigcup_i\ir {a_i}{b_i}.$
Let $S=\bigcup_i s_i;$ note that $S\cap N=\pu.$

Put $\xi(z)=\vt(z\cup S)\sd \vt(S)$ for any $z\sq N.$ 
\imar{why $\sd\vt(S)$ added?}%
Then, for any sets $x\yi y\sq N,$ 
\dm
{x\sd y \in\sui{r_n}}\eqv
{\vt(x\cup S)\sd\vt(y\cup S)\in\cJ}\eqv
{\xi(x)\sd\xi(y)\in\cJ,}
\dm
thus $\xi$ reduces $\sui{r_n}\res N$ to $\cJ.$ 
Now put $w_i=\xi(g_i)\cap\ir{k_{i}}{k_{i+1}}$ 
and $w_a=\bigcup_{i\in a} w_i$ for $a\in\pn.$  
We assert that the map $i\mapsto w_i$ proves 
$\sui{r_n}\orbpp\cJ.$ 
In view of the above, 
it remains to show that $\xi(g_a)\sd w_a\in\cJ$ for 
any $a\in\pn$.

As $\cJ=\Exh_\nu,$ it suffices to demonstrate that 
$\nu\skl w_i\sd(\xi(g_a)\cap\ir{k_{i}}{k_{i+1}})\skp
<2^{-i}$ 
for all $i\in a$ while 
$\nu(\xi(g_a)\cap\ir{k_{i}}{k_{i+1}})<2^{-i}$ 
for $i\nin a.$ 
After dropping the common term $\vt(S),$    
it suffices to check that 
\ben
\tenu{(\alph{enumi})}
\itla{111}
$\nu\skl (\vt(g_i\cup S)\sd\vt(g_a\cup S))
\cap\ir{k_{i}}{k_{i+1}}\skp<2^{-i}$ 
for all $i\in a$ while 

\itla{222}
$\nu\skl(\vt(S)\sd\vt(g_a\cup S))\cap
\ir{k_{i}}{k_{i+1}}\skp<2^{-i}$ 
for $i\nin a.$ 
\een
Note that any set of the form $x\cup S,$ where $x\sq N,$ 
is generic by \ref{stabg}. 
It follows, by \ref{stab2}, that we can assume, in \ref{111} 
and \ref{222}, that $a\sq[0,i],$ \ie, resp.\ $\tmax a=i$ 
and $\tmax a<i.$ 
We can finally apply \ref{stab1}, with 
$u=a\cup\bigcup_{j<i}s_j,\msur$ $x=\bigcup_{j>i}s_j,$ and 
$v=u_i \cup\bigcup_{j<i}s_j$ if $i\in a$ while 
$v=\bigcup_{j<i}s_j$ if $i\nin a$.

\ref{t>s2}
Otherwise $\sui{r_n}\orbpp\zo$ exponentially by \ref{t>s1}. 
Let this be witnessed by $i\mapsto w_i$ and a sequence of 
numbers $k_i,$ so that $k_{i+1}\ge 2k_i$ and 
$w_i\sq\ir{k_i}{k_{i+1}}.$ 
If $d_i=\frac{\#(w_i)}{k_{i+1}}\to0$ 
then easily $\bigcup_iw_i\in\zo$ by the choice of $\sis{k_i}{}.$ 
Otherwise there is a set $a\in\sui{r_n}$ such that 
$d_i>\ve$ for all $i\in a$ and one and the same $\ve>0,$ 
so that $w_a=\bigcup_{i\in a}w_i\nin\zo.$ 
In both cases we have a contradiction with the assumption 
that the map $i\mapsto w_i$ witnesses $\sui{r_n}\orbpp\zo$.
\epf

\bqu
\lam{far-s}
Farah~\cite{f-tsir} points out that
Theorem~\ref{t:>s}\ref{t>s1} also holds for $\ofi$
(instead of $\sui{r_n}$) and asks for which other ideals 
it is true.
\equ

\punk
{The family $\bel p$}
\las{hd1}

It follows from the next theorem that Borel reducibility
between \eqr s $\bel p\yt 1\le p<\iy,$ is fully determined
by the value of $p$.

\bte[{{\rm Dougherty -- Hjorth \cite{dh}}}]
\lam{hd-t}
If\/ $1\le p<q<\iy$ then\/ ${\bel p}\rebs{\bel q}$.
\ete
\bpf
{\ubf Part 1\/}:  
show that ${\bel q}\not\reb{\bel p}.$  

By Lemma~\ref{redx}, it suffices to prove that 
${\bel q\res\dX}\not\reb{\bel p\res \dX}.$ 
Suppose, on the contrary, that $\vt:\dX\to\dX$
is a Borel reduction of
${\bel q\res\dX}$ to ${\bel p\res \dX}.$
Arguing as in the proof of Theorem~\ref{t:>s},
we can reduce the general case to the case when
there exist increasing sequences of numbers
$0=j(0)<j(1)<j(2)<\dots$ and $0=a_0<a_1<a_2<\dots$ 
and a map $\tau:\dY\to\dX,$ where
$\dY=\prod_{n=0}^\iy X_{j(n)},$
which reduces ${\bel q\res\dY}$ to $\bel p\res \dX$ 
and has the form
$\tau(x)=\bigcup_{n\in\dN}t^{x(n)}_n,$
where $t^r_n\in \prod_{k=a_n}^{a_{n+1}-1}X_k$ for any
$r\in X_{j_n}.$
(See Definition~\ref{dx}.)

{\it Case 1\/}: 
there are\/ $c>0$ and a number\/ $N$ such that\/ 
$\nr{\tau^1_n-\tau^0_n}p\ge c$ for all\/ $n\ge N.$ 
Since $p<q,$ there is a non-decreasing sequence of natural 
numbers $i_n\le {j_n}\yt n=0,1,2,\dots,$ 
such that $\sum_n2^{p(i_n-j_n)}$ 
diverges but $\sum_n2^{q(i_n-j_n)}$ converges. 
({\it Hint\/}: $i_n\approx j_n-p\obr\log_2n$.) 

Now consider any $n\ge N.$ 
As $\nr{\tau^1_n-\tau^0_n}p\ge c$ and because $\nr{\dots}p$ 
is a norm, there exists a pair of rationals $u(n)<v(n)$ in 
$X_{j_n}$ with $v(n)-u(n)=2^{i_n-j_n}$ and 
$\nr{\tau^{v(n)}_n-\tau^{u(n)}_n}p\ge c\,2^{i_n-j_n}.$ 
In addition, put $u(n)=v(n)=0$ for $n<N.$ 
Then the \dd{\bel q}distance between the infinite
sequences $u=\sis{u(n)}{n\in\dN}$ 
and $v=\sis{v(n)}{n\in\dN}$ is equal to 
$\sum_{n=N}^\iy 2^{q(i_n-j_n)}<+\iy,$
while the \dd{\bel p}distance between $\tau(u)$ and
$\tau(v)$ is non-smaller than
$\sum_{n=N}^\iy c^p\,2^{p(i_n-j_n)}=\iy.$ 
But this contradicts the assumption that $\tau$
is a reduction.

{\it Case 2\/}: 
otherwise.
Then there is a strictly  increasing sequence 
$n_0<n_1<n_2<\dots$ with 
$\nr{\tau^1_{n_k}-\tau^0_{n_k}}p\le 2^{-k}$ for all $k.$ 
Let now $x\in \dY$ be the constant $0$ while $y\in\dY$ 
be defined by $y(n_k)=1\zd\kaz k$ and $y(n)=0$ 
for all other $n.$ 
Then $x\bel q y$ fails
($|y(n)-x(n)|\not\to 0$)
but $\tau(x)\bel p \tau(y)$ 
holds, contradiction.\vom

{\ubf Part 2\/}: 
show that ${\bel p}\reb{\bel q}.$ 

It suffices to prove that ${\bel p\res\oin}\reb{\bel q}$
(Lemma~\ref{redx}).
We \noo\ assume that $q<2p$:
any bigger $q$ can be approached in several steps. 
For $\vec x=\ang{x,y}\in\dR^2,$ let 
$\nr{\vec x}h=(x^h+y^h)^{1/h}.$ 

\ble
\lam{hd+1}
For any\/ $\frac12<\al<1$ there is a continuous 
map\/ $K_\al:[0,1]\to[0,1]^2$ and positive real numbers\/ 
$m<M$ such that for all\/ $x<y$ in $[0,1]$ we have\/
$m(y-x)^\al\le\nr{K_\al(y)-K_\al(x)}2\le M(y-x)^\al$. 
\ele 
\bpf[{{\rm Lemma}}]
The construction of such a map $K$ can be easier described
in terms of fractal geometry rather than by an analytic
expression. 
Let $r=4^{-\al},$ so that $\frac14<r<\frac12$ and 
$\al=-\log_4r.$ 
Starting with the segment $[(0,0)\,\yi (1,0)]$
of the horisontal axis of the cartesian plane, we replace it
by four smaller segments of length $r$ each
(thin lines on Fig.~2, left).
Each of them we replace by four segments of
length $r^2$ (thin lines on Fig.~2, right). 
And so on, infinitely many steps.
The resulting curve $K$ is parametrized by giving the
vertices of the polygons values equal to multiples of $4^{-n},$
$n$ being the number of the polygon.
For instance, the vertices of the left polygon on Fig.~2 are
given values $0,\frac14,\frac12,\frac34,1.$
%

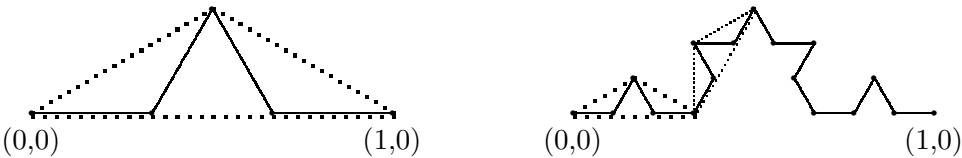
\begin{figure}[h]
\label{fig2}


\setlength{\unitlength}{0.8mm} 
\begin{picture}(130,35)(2,0) 

\label{f-f}

\put(5,4) {(0,0)}
\put(65,4) {(1,0)}
\put(95,4) {(0,0)}
\put(155,4) {(1,0)}

\put(10,10){\tob}
\put(30,10){\tob}
\put(50,10){\tob}
\put(70,10){\tob}
\put(40,27.32){\tob}

\put(10,10){\line(1,0){20}}
\put(50,10){\line(1,0){20}}

\put(100,10){\tob}
\put(120,10){\tob}
\put(140,10){\tob}
\put(160,10){\tob}

\put(106.67,10){\tob}
\put(113.33,10){\tob}
\put(146.67,10){\tob}
\put(153.33,10){\tob}

\put(130,27.32){\tob}
\put(123.33,15.77){\tob}
\put(126.67,21.55){\tob}
\put(133.33,21.55){\tob}
\put(136.67,15.77){\tob}

\put(110,15.77){\tob}
\put(150,15.77){\tob}

\put(120,21.55){\tob}
\put(140,21.55){\tob}

\put(100,10){\line(1,0){6.67}}
\put(140,10){\line(1,0){6.67}}
\put(120,10){\line(-1,0){6.67}}
\put(160,10){\line(-1,0){6.67}}


\qbezier(30,10)(35,18.66)(40,27.32)
\qbezier(50,10)(45,18.66)(40,27.32)

\qbezier(106.67,10)(108.33,12.88)(110,15.77)
\qbezier(113.33,10)(111.67,12.88)(110,15.77)
\qbezier(146.67,10)(148.33,12.88)(150,15.77)
\qbezier(153.33,10)(151.67,12.88)(150,15.77)

\qbezier(120,10)(121.67,12.88)(123.33,15.77)
\qbezier(140,10)(138.33,12.88)(136.67,15.77)

\qbezier(120,21.55)(121.67,18.66)(123.33,15.77)
\qbezier(140,21.55)(138.33,18.66)(136.67,15.77)

\qbezier(120,21.55)(123.33,21.55)(126.67,21.55)
\qbezier(140,21.55)(136.67,21.55)(133.33,21.55)

\qbezier(133.33,21.55)(131.67,24.44)(130,27.32)
\qbezier(126.67,21.55)(128.33,24.44)(130,27.32)

\linethickness{0.5mm}



\linethickness{0.4mm}

\qbezier[20](10,10)(25,18.66)(40,27.32)
\qbezier[20](70,10)(55,18.66)(40,27.32)
\qbezier[30](10,9.3)(40,9.3)(70,9.3)

\qbezier[10](100,9.3)(110,9.3)(120,9.3)
\qbezier[6](100,10)(105,12.88)(110,15.77)
\qbezier[6](120,10)(115,12.88)(110,15.77)

\linethickness{0.2mm}

\qbezier[15](120,10)(120,15)(120,21.55)
\qbezier[15](120,21.55)(125,24.44)(130,27.32)
\qbezier[20](120.1,9.3)(125.1,17.96)(130.1,26.62)

\end{picture}

\caption{$r=\frac13$, \ left: step 1, \ right: step 2}

\end{figure}

%
Note that the curve $K:[0,1]\to[0,1]^2,$
approximated by the polygons, is
bounded by certain triangles built on the sides of the
polygons.
For instance, the whole curve lies inside the triangle
bounded by dotted lines in Fig.~2, left.
(The dotted line that follows the basic side
$[(0,0)\,\yi (1,0)]$ of the triangle is drawn  
slightly below its true position.)
Further, the parts $0\le t\le\frac14$ and
$\frac14\le t\le\frac12$ of the curve 
lie inside the triangles bounded by (slightly different)
dotted lines in Fig.~2, right.
And so on.
Let us call those triangles {\it bounding triangles\/}.
\vyk{
It is quite obvious geometrically that for any $t$
in $[\frac{i-1}{4^n},\frac{i}{4^n}]$ ($1\le i\le4^n$)
the euclidean distance $\nr{K_r(t)-C^i_n}2$ between
$K(t)$ and the center $C^i_n$ of the straight interval
$[K(\frac{i-1}{4^n})\,\yi K(\frac{i}{4^n})]$
(inside the unit square)
satisfies $\nr{K_r(t)-C^i_n}2\le \frac{r^n}2.$
}

To prove the inequality of the lemma, consider any pair of
reals $x<y\in[0,1].$
Let $n$ be the least number such that $x,y$ belong to
non-adjacent intervals, resp.,
$[\fras{i-1}{4^n},\fras{i}{4^n}]$ and
$[\fras{j-1}{4^n},\fras{j}{4^n}],$ with $j>i+1.$
Then $4^{-n}\le |y-x|\le 8\cdot 4^{-n}.$

The points $K(x)$ and $K(y)$ then belong to one and the
same side or adjacent sides of the \dd{n-1}th polygon.
Let $C$ be a common vertice of these sides.
It is quite clear geometrically that the euclidean distances
from $K(x)$ and $K(y)$ to $C$ do not exceed $r^{n-1}$
(the length of the side), thus 
$\nr{K(x)-K(y)}2\le 2\,r^{n-1}.$

Estimation from below needs more work.
The points $K(x)\zd K(y)$ belong to the bounding triangles
built on the segments, resp.,
$[K(\frac{i-1}{4^n})\,\yi K(\frac{i}{4^n})]$ and
$[K(\frac{j-1}{4^n})\,\yi K(\frac{j}{4^n})],$ and
obviously $i+1<j\le i+8,$ so that there exist at most six
bounding triangles between these two.
Note that adjacent bounding triangles meet each other at
only two possible angles (that depend on $r$ but not on $n$),
and taking it as geometrically evident that non-adjacent
bounding triangles are disjoint, we conclude that there is
a constant $c>0$ (that depends on $r$ but not on $n$)
such that the distance between two non-adjacent
bounding triangles of rank $n,$ having at most $6$
bounding triangles of rank $n$ between them, does not exceed
$c\app r^n.$
In particular, $\nr{K(x)-K(y)}2\ge c\app r^{n}.$
Combining this with the inequalities above, we conclude that
$m(y-x)^\al\le\nr{K(y)-K(x)}2\le M(y-x)^\al,$
where $m=\frac c{8^\al}$ and $M=\frac2r$
(and $\al=-\log_4r$). 
\epF{Lemma}

Coming back to the theorem, let $\al=p/q,$ 
and let $K_\al$ be as in the lemma. 
Let $x=\ang{x_0,x_1,x_2,\dots}\in\oin.$ 
Then $K_\al(x_i)=\ang{x'_i,x''_i}\in[0,1]^2.$ 
We put $\vt(x)=\ang{x'_0,x''_0,x'_1,x''_1,x'_2,x''_2,\dots}.$ 
Prove that $\vt$ reduces $\bel p\res\oin$ to $\bel q$.

Let $x=\sis{x_i}{i\in\dN}$ and $y=\sis{y_i}{i\in\dN}$ 
belong to $\oin;$ we have to prove that $x\bel p y$ iff 
$\vt(x)\bel q\vt(y).$ 
To simplify the picture note the following: 
\dm
2^{-1/2}\nr w 2\le \tmax\ans{w',w''}\le \nr w q\le \nr w 1\le
2\nr w2
\dm 
for any $w=\ang{w',w''}\in\dR^2.$ 
The task takes the following form:
\dm
\sum_i(x_i-y_i)^p<\iy
\leqv 
\sum_i {\nr{K_\al(x_i)-K_\al(y_i)}2}^q<\iy\,.
\dm
Furthermore, by the choice of $K_\al,$ this converts to 
\dm
\sum_i(x_i-y_i)^p<\iy
\leqv 
\sum_i (x_i-y_i)^{\al q}<\iy\,,
\dm
which holds because $\al q=p.$ 
\epF{Theorem~\ref{hd-t}}

\punk{$\bel\iy$: maximal $\Ks$}
\las{liy}

Recall that $\Ks$ denotes the class of all \dd\fsg compact
sets in Polish spaces.
Easy computations show that this class contains, among
others, the \eqr s ${\Ei}\zd{\Ey}\zd{\bel p}\zd1\le p\le\iy,$
considered as sets of pairs in corresponding Polish spaces.
Note that if $\rE$ a $\Ks$ equivalence on a Polish space $\dX$
then $\dX$ is $\Ks$ as well 
since projections of compact sets are compact.
Thus $\Ks$ \er s on Polish spaces is one and the same as 
$\fs02$ \er s on $\Ks$ Polish spaces.

\bte
\lam{eiy}
Any\/ $\Ks$ \eqr\ on a Polish space, in particular,
${\Ei}\zd{\Ey}\zd{\bel p},$
is Borel reducible to\/ ${\bel\iy}.$~\snos
{The result for $\bel p$ is due to Su Gao \cite{sudis}.
He defines
$d_p(x,s)=(\sum_{k=0}^{\lh s-1} |x(k)-s(k)|^p)^{\frac1p}$
for any $x\in\rtn$ and $s\in\dQ\lom$
(a finite sequence of rationals). 
Easily the \dd{\bel p}distance
$(\sum_{k=0}^{\iy} |x(k)-y(k)|^p)^{\frac1p}$
between any pair of $x,y\in\rtn$ is finite iff 
there is a constant $C$ such that
$|d_p(x,s)-d_p(y,s)|<C$ for all $s\in\dQ\lom.$
This yields a reduction required.}
\ete
\bpf[{{\rm from Rosendal~\cite{roz}}}]
Let $\dA$ be the set of all \dd\sq increasing sequences
$a=\sis{a_n}{n\in\dN}$ of subsets $a_n\sq\dN$ --- a closed
subset of the Polish space $\pn^\dN.$ 
Define an \er\ $\rH$ on $\dA$ by
\dm
\sis{a_n}{}\rH\sis{b_n}{}
\quad\text{ iff }\quad
\sus N\:\kaz m\:
({a_m\sq b_{N+m}}\land {b_m\sq a_{N+m}}).
\dm

{\it Claim 1\/}:
$\rH\reb{\bel\iy}.$
This is easy.
Given a sequence $a=\sis{a_n}{n\in\dN}\in\pn^\dN,$ define
$\vt(a)\in \dN^{\dN\ti\dN}$
by $\vt(a)(n,k)$ to be the least $j\le k$ such that $n\in a_j,$
or $\vt(a)(n,k)=k$ whenever $n\nin a_k.$
Then $\sis{a_n}{}\rH\sis{b_n}{}$ iff there is $N$ such that 
$|\vt(a)(n,k)-\vt(b)(n,k)|\le N$ for all $n\zi k$.

{\it Claim 2\/}:
any $\Ks$ equivalence $\rE$ on a Polish space $\dX$ is Borel
reducible to $\rH.$ 
As a $\Ks$ set, $\rE$ has the form $\rE=\bigcup_nE_n,$ where
each $E_n$ is a compact subset of $\dX\ti\dX$
(not necessarily an \er)
and $E_n\sq E_{n+1}.$
We can \noo\ assume that each $E_n$ is reflexive and
symmetric on its domain $D_n=\dom{E_n}=\ran{E_n}$
(a compact set), in particular,
$x\in D_n\imp{\ang{x,x}\in E_n}.$
\vyk{
Furthermore, both $\dX$ and
$D=\rav\dX=\ens{\ang{x,x}}{x\in\dX}$ are $\Ks;$
let $D=\bigcup_nD_n,$ where $D_n$ are compact 
and $D_n\sq D_{n+1}.$
}%
Define $P_0=E_0$ and 
\dm
P_{n+1}=P_n\cup E_{n+1}\cup {P_n^{(2)}},
\;\text{ where }\,
P_n^{(2)}=\ens{\ang{x,y}}{\sus z\:
({\ang{x,z}\in P_n}\land{\ang{z,y}\in P_n})},
\dm
by induction.
Thus all $P_n$ are still compact subsets of $\dX\ti\dX,$
moreover, of $\rE$ since $\rE$ is an \eqr, and
$E_n\sq P_n\sq P_{n+1},$ therefore $\rE=\bigcup_nP_n.$

Let $\ens{U_k}{k\in\dN}$ be a basis for the topology of
$\dX.$
Put, for any $x\in\dX,$
$\vt_n(x)=\ens{k}{U_k\cap R_n(x)\ne\pu},$
where $R_n(x)=\ens{y}{\ang{x,y}\in R_n}.$
Then obviously $\vt_n(x)\sq\vt_{n+1}(x),$
and hence $\vt(x)=\sis{\vt_n(x)}{n\in\dN}\in\dA.$
Then $\vt$ reduces $\rE$ to $\rH.$

Indeed if $x\rE y$
then $\ang{y,x}\in P_n$ for some $n,$ and
for all $m$ and $z\in\dX$ we have 
${\ang{x,z}\in R_m}\imp{\ang{y,z}\in R_{1+\tmax\ans{m,n}}}.$
In other words, $R_m(x)\sq R_{1+\tmax\ans{m,n}}(y)$
and hence
$\vt_m(x)\sq \vt_{1+\tmax\ans{m,n}}(y)$
hold for all $m.$
Similarly, for some $n'$ we have
$\vt_m(y)\sq \vt_{1+\tmax\ans{m,n'}}(y)\zd \kaz m.$
Thus ${\vt(x)}\rH{\vt(y)}.$

Conversely, suppose that ${\vt(x)}\rH{\vt(y)},$ thus,
for some $N,$ we have 
$R_m(x)\sq R_{N+m}(y)$ and $R_m(y)\sq R_{N+m}(x)$
for all $m$ and $y.$
Taking $m$ big enough for $P_m$ to contain
$\ang{x,x},$ we obtain $x\in R_{N+m}(y),$ so that
immediately $x\rE y$.
\epf

\api

\parf{Smooth, hyperfinite, countable}
\las{parf:1st}

This \gla\ is related to the domain $\reb\Ey$ 
in the diagram on page \pageref{p-p}.
The following types of \eqr s are relevant to this
domain:

\bdf 
\lam{typesER}
A Borel \eqr\ $\rE$ on a (Borel) set $\dX$ is:\vtm 

\noi --
{\it countable\/}, \ 
if every \dde class
$\eke x=\ens{y\in X}{x\rE y}\yt x\in X,$ 
\index{equivalence relation, ER!countable}%
is at most countable;\vtm
\imar{Is any ctble $\fs11$ \er\ actually Borel ?}%

\noi --
{\it essentially countable\/}, \ 
if $\rE\reb\rF,$ where $\rF$  
\index{equivalence relation, ER!essentially countable}%
is a countable Borel \er;\vtm

\noi --
{\it finite\/}, \
if every \dde class $\eke x=\ens{y\in\dX}{x\rE y}\yt x\in X,$ 
\index{equivalence relation, ER!finite}%
is finite;\vtm

\noi --
{\it hyperfinite\/}, \
if $\rE=\bigcup_n\rF_n$ for an increasing sequence of Borel 
\index{equivalence relation, ER!hyperfinite}%
finite \er s $\rF_n$;\vtm

\noi --
{\it smooth\/}, \
if $\rE\reb\rav\dn$;\vtm 
\index{equivalence relation, ER!smooth}%

\noi --
{\it hypersmooth\/}, \
if $\rE=\bigcup_n\rF_n$ for an increasing sequence of  
\index{equivalence relation, ER!hyperfinite}%
smooth \er s $\rE_n$.\qed
\eDf

After a few rather simple results on smooth \eqr s, we
proceed to countable equivalences.
We prove in \nrf{cer} that every countable Borel \er\ is
Borel reducible to $\Ey,$ and hence the whole domain
$\reb\Ey$ is equal to the class of essentially countable
Borel \er s.

Then we turn to hyperfinite \eqr s, a very interesting subclass
of countable Borel equivalences.
A typical hyperfinite equivalence is $\Eo$ --- in fact the
\dd\reb largest, or \rit{universal} hyperfinite \er. 
Hyperfinite \er s admit several different characterizations 
--- some of them are presented by Theorem~\ref{thf}.

The \eqr\ $\Ey$ turns out to be (countable but) non-hyperfinite
by Theorem~\ref{femo}.
It follows that $\Eo\rebs\Ey$ strictly.

We finish with two separate theorems.
One of them, Theorem~\ref{cudf}, asserts that, given a
countable \eqr\ satisfying $\rF+\rF\reb\rF,$ the property
\lap{being Borel reducible to $\rF$} is \dd\fsg additive. 
Theorem~\ref{jnmt} shows that $\ifi$ is the \dd\reb least
ideal.

\vyk{
Theorem~\ref{1dih} is the 1st dichotomy theorem;
it asserts that
$\rav\dn$ is the \dd\reb least among all Borel \eqr s with
uncountably many equivalence classes.

Finally, we present a characterization, in terms of the 
existence of transversals, of those Borel sets $X$ for which
\imar{where is this?}
$\Eo\res X$ is smooth.
}

\punk{Smooth and below}
\las{<smooth}

By definition an \eqr\ $\rE$ is smooth iff there is a Borel map
$\vt:X\to\dn$ such that the
equivalence ${x\rE y}\eqv{\vt(x)=\vt(y)}$ holds for all
$x\yi y\in X=\dom\rE.$ 
In other words, it is required that the equivalence classes
can be counted by reals (here: elements of $\dn$) in Borel way.
An important subspecies of smooth \eqr s consists of 
those having a Borel {\it transversal\/}: 
a set with 
\index{transversal}%
exactly one element in every equivalence class. 

\ble
\label{transv}
\ben
\renu
\itsep
\itla{transv1}
Any Borel \er\ that has a Borel transversal is smooth$;$
\imar{transv}

\itla{transv2}
any Borel finite (with finite classes) \er\ admits a
Borel transversal$;$

\itla{transv3}
any Borel countable smooth \er\ admits a Borel transversal;

\itla{transv4}
any Borel \er\ $\rE$ on a Polish space\/ $\dX,$ such that 
every\/ \dde class is closed and the saturation\/ 
$\ek\cO\rE$ of every open set\/ $\cO\sq\dX$ is Borel, 
admits a Borel transversal, hence, is smooth$.$~\snos
{Srivastava~\cite{sri} proved the result for 
\er s with $\Gd$ classes,  
which is the best possible as $\Eo$ is a Borel \er, whose 
classes are $\Fs$ and saturations of open sets are even open, 
but without any Borel transversal. 
See also \cite[18.20~iv)]{dst}.}

\itla{transv5}
$\Eo$ is not smooth.

\itla{transv6}
there exists a smooth \er\ $\rE$ that does not have a Borel
transversal.
\een
\ele
\bpf
\ref{transv1}
Let $T$ be a Borel transversal for $\rE.$   
The map $\vt(x)=$
``the only element of $T$ \dde equivalent to $x$''
reduces $\rE$ to $\rav T.$
\vyk{
\snos
{To see that a smooth \er\ does not necessarily have a 
Borel transversal  
take a closed set $P\sq\dnn\ti\dnn$ with $\dom P=\dnn,$ 
not uniformizable by a Borel set, and let  
$\ang{x,y}\rE\ang{x',y'}$ iff both $\ang{x,y}$ and 
$\ang{x',y'}$ belong to $P$ and $x=x'$.}
}

\ref{transv2}
Consider the set of the \dd<least elements of \dde classes, 
where $<$ is a fixed Borel linear order on the 
domain of $\rE$.

\ref{transv3}
Use \Cuni\ (Theorem~\ref{cuni}).

\ref{transv4}
Since any uncountable Polish space is a continuous image of 
$\dnn,$ 
we can assume that $\rE$ is an \eqr\ on $\dnn.$ 
Then, for any $x\in\dnn,$
the equivalence class $\ek x\rE$ is a closed subset of 
$\dnn,$ naturally identified with a tree, say, $T_x\sq\dN\lom.$ 
Let $\vt(x)$ denote the leftmost branch of $T_x.$ 
Then $x\rE\vt(x)$ and ${x\rE y} \imp{\vt(x)=\vt(y)},$ so that 
it remains to show that $Z=\ens{\vt(x)}{x\in\dnn}$ is Borel. 
Note that 
\dm
{z\in Z}\,\leqv\,
\kaz m\;\kaz s,\,t\in \dN^m\;
\skl
{s<_{\text{\tt lex}} t}\land {z\in\cO_t}\limp 
{\ek z\rE\cap\cO_t=\pu}
\skp,
\dm
where $<_{\text{\tt lex}}$ is the lexicographical order on 
$\dN^m$ 
and $\cO_s=\ens{x\in\dnn}{s\su x}.$ 
However $\ek x\rE\cap\cO_t=\pu$ iff $x\nin \ek{\cO_t}\rE$ 
and $\ek{\cO_t}\rE$ is Borel for any $t$.

\ref{transv5}
Otherwise $\Eo$ has a Borel transversal $T$ by \ref{transv3},
which is a contradiction, see an argument after
Example~\ref{ex:ac2}.

\ref{transv6}
Take a closed set $P\sq\dnn\ti\dnn$ with $\dom P=\dnn$ 
that is not uniformizable by a Borel set,  
\imar{reference?}%
and define $\ang{x,y}\rE\ang{x',y'}$ iff both $\ang{x,y}$ and  
$\ang{x',y'}$ belong to $P$ and $x=x'$. 
\epf

\vyk{

\punk{Assembling smooth \eqr s}
\las{cud}

If $\rE$ and $\rF$ are smooth \er s on disjoint sets, resp.,  
$X$ and $Y,$ then easily $\rE\cup\rF$ is a smooth \er\ on 
$X\cup Y.$ 
The question becomes less clear when we have a Borel \er\ 
$\rE$ on a Polish space $X\cup Y$ such that both $\rE\res X$
and $\rE\res Y$ are smooth but the sets $X,Y$ not necessarily
\dd\rE invariant in $X\cup Y$ if even disjoint; is $\rE$ smooth? 
We answer this in the positive, even in the case of countable 
unions.

\bte
\lam{cud1}
Let\/ $\rE$ be a Borel \eqr\ on a Borel set\/ 
$X=\bigcup_kX_k,$ with all\/ $X_k$ also Borel. 
Suppose that each\/ $\rE\res X_k$ is smooth. 
Then\/ $\rE$ is smooth.
\ete
\bpf\footnote
{\ The shortest proof is to note that otherwise 
$\Eo\reb\rE$ by the 2-nd dichotomy, easily leading to 
contradiction by a Baire category argument.
Yet we prefer to give a direct proof.
Note that even in the case when the sets $X_k$ are
pairwise disjoint, most obvious ideas like
``to define $\vt(x)$ take the least $k$ such that
$X_k$ intersects $\eke x$ and apply $\vt_k$'' do not
work.}
It suffices to prove the following.
If $\rE$ is a Borel \eqr\ on the union $X\cup Y$ of two
disjoint Borel sets $X,Y,$ and both ${\rE}\res X$
and ${\rE}\res Y$
are smooth via Borel maps, resp., $f:X\to P$ and $g:Y\to Q$
($P,Q$ being disjoint Borel) 
then $\rE$ itself 
is also smooth, and this can be witnessed by a Borel map
$h$ defined on $X\cup Y$ so that $h\res X=f.$ 

Let us prove this claim.
It follows from the smoothness that the set 
\dm
F=\ens{\ang{q,p}}{\sus x\in X\;\sus y\in Y\;
\skl{f(x)=p}\land{g(y)=q}\land{x\rE y}\skp}
\sq Q\ti P
\dm
is a $\fs11$ bijection.
By Theorem~\ref{bore} there is a Borel bijection 
$\Phi\sq Q\ti P$ such that $F\sq\Phi.$ 
Now the $\fp11$ set 
\dm
W=\ens{\ang{q,p}\in \Phi}{\kaz x\in X\;\kaz y\in Y\;
\skl f(x)=p\land g(y)=q\imp x\rE y\skp}\,,
\dm
satisfies $F\sq W\sq \Phi,$ hence there is 
a Borel bijection $\Psi$ with $F\sq\Psi\sq W.$  
The sets $B=\dom\Psi$ and $A=\ran\Psi$ are Borel subsets 
of resp.\ $Q,P,$  
and it follows from the construction that  
$\Psi\cap(\dom F\ti \ran F)=F.$ 
Finally, put $h(x)=f(x)$ for $x\in X,$
$h(y)=\Psi(g(y))$ whenever $y\in Y$ and $g(y)\in B,$ and
$h(y)=y$ whenever $y\in Y$ and $g(y)\nin B$.
\vyk{
\dm
D=\Psi\cup\ens{\ang{p,p}}{p\in P\dif A}
\cup\ens{\ang{q,q}}{q\in Q\dif B}\,.
\dm
Then, for any $x\in X$ there is a unique 
$h(x)=\ang{p,q}\in D$ with $p=f(x),$ 
correspondingly, for any $y\in Y$ there is unique 
$h(y)=\ang{p,q}\in D$ with $q=g(y),$ and if $y\rE z$ 
then $h(x)=h(y)=\ang{f(x),g(y)},$ hence $h$ 
witnesses that $\rE$ is smooth.
}%
\epf

}

\punk{Countable \eqr s}
\las{cer}


This class of \eqr s
is a subject of ongoing intence study. 
We present here the following important theorem
(\cite[Thm 1]{femo1}, \cite[1.8]{djk}) 
and a few more results below, leaving
\cite{jkl,gab,kemi} as sources of further information
regarding countable \eqr s.

\bte
\lam{femo}
Any Borel countable \er\/ $\rE$ on a Polish space\/ $\dX\,{:}$
\ben
\tenu{{\rm(\roman{enumi})}}
\itsep
\itla{femo1}
is induced by a Polish action of a countable group\/
$\dG$ on\/ $\dX\,;$   

\itla{femo2}
satisfies\/ $\rE\reb \Ey=\rE(F_2,2),$ where\/
$F_2$ is the free group with two generators and\/
$\rE(F_2,2)$ is the \er\ induced by the 
shift action of\/ $F_2$ on\/ $2^{F_2}.$
\een
\ete
\bpf
\ref{femo1}
We \noo\ assume that $\dX=\dn.$
According to \Cenu\
(Theorem~\ref{cenu}, in a relativized version, if necessary, see
Remark~\ref{relrel}),
there is a sequence of Borel maps $f_n:\dn\to\dn$
such that $\eke a=\ens{f_n(a)}{n\in\dN}$ for each
$a\in\dn.$
Put $\Ga'_n=\ens{\ang{a,f_n(a)}}{a\in\dN}$
(the graph of $f_n$)
and $\Ga_n=\Ga'_n\dif\bigcup_{k<n}\Ga'_k.$
The sets $P_{nk}=\Ga_n\cap{\Ga_k}\obr$ form a partition 
of (the graph of) $\rE$ onto countably many Borel
injective sets.
Further define $\Da=\ens{\ang{a,a}}{a\in\dn}$ and let
$\sis{D_m}{m\in\dN}$ be an enumeration of all
non-empty sets of the form $P_{nk}\dif\Da.$
Intersecting the sets $D_m$ with the rectangles of the
form
\dm
R_s=\ens{\ang{a,b}\in\dn\ti\dn}{s\we0\su a\land s\we1\su b}
\quad\text{and}\quad
{R_s}\obr,
\dm
we reduce the general case to the case when
$\dom{D_m}\cap\ran{D_m}=\pu\zd\kaz m.$

Now, for any $m$ define $h_m(a)=b$ whenever either
$\ang{a,b}\in D_m$ or $\ang{a,b}\in {D_m}\obr,$ or
$a=b\nin \dom{D_m}\cup\ran{D_m}.$
Clearly $h_m$ is a Borel bijection $\dn\onto\dn.$
Thus $\sis{h_m}{m\in\dN}$ is a family of Borel
automorphisms of $\dn$ such that
$\eke a=\ens{h_m(a)}{m\in\dN}.$
It does not take much effort to expand this system to a Borel
action of $F_\om,$ the free group
with countably many generators $a_1,a_2,a_3,\dots,$ 
\index{group!free with $\alo$ generators, $F_\om$}%
\index{zzFw@$F_\om$}%
on $\dn,$ whose induced \eqr\ is $\rE$.

\ref{femo2}
First of all, by \ref{femo1}, $\rE\reb\rR,$ where $\rR$
is induced by a Borel action $\app$ of $F_\om$ on $\dn.$
The map $\vt(a)=\sis{g\obr\app a}{g\in F_\om}\yt a\in\dn$
is a Borel reduction of $\rR$ to $\rE(F_\om,\dn).$
If now $F_\om$ is a subgroup of a countable group $H$ then
$\rE(F_\om,\dn)\reb\rE(H,\dn)$ by means of the map sending
any $\sis{a_g}{g\in F_\om}$ to $\sis{b_h}{h\in H},$ where
$b_g=a_g$ for $g\in F_\om$ and $b_h$ equal to any fixed
$b'\in\dn$ for $h\in H\dif F_\om.$
As $F_\om$ admits an injective homomorphism into $F_2$~\snos
{Indeed, let $F$ be the subgroup of $F_2$ generated by all  
elements of the form $\al_n=a^nb^n$ and $\al_n\obr=b^{-n}a^{-n}.$ 
The map sending any $a_n$ to $\al_n$ and accordingly
$a_n\obr$ to $\al_n\obr$ is an isomorphism of $F_\om$ onto $F$.}
we conclude that $\rE\reb\rE(F_2,\dn).$

It remains to define a Borel reduction of $\rE(F_2,\dn)$ to
$\rE(F_2,2).$
The inequality $\rE(F_2,\dn)\reb\rE(F_2,2^{\dZ\dif\ans0})$
is clear.
Further $\rE(F_2,2^{\dZ\dif\ans0})\reb\rE(F_2\ti\dZ,3),$
by means of the map sending any
$\sis{a_g}{g\in F_2}\msur$ $(a_g\in2^{\dZ\dif\ans0})$ to
$\sis{b_{gj}}{g\in F_2,\;j\in\dZ},$ where
$b_{gj}= a_g(j)$ for $j\ne0$ and $b_{g0}= 2.$
Further, for any group $G$ it holds
$\rE(G,3)\reb\rE(G\ti \dZ_2,2)$
by means of the map sending every element 
$\sis{a_g}{g\in G}\;\,(a_g=0,1,2)$ to
$\sis{b_{gi}}{g\in G,\;i\in\dZ_2},$ where
\dm
b_{gi}=
\left\{
\bay{rclcl}
0,&\text{if}& a_g=0 &\text{\ubf or}&
a_g=1\,\text{ and }\, i=0,\\[0.7\dxii]

1,&\text{if}& a_g=2 &\text{\ubf or}&
a_g=1\,\text{ and }\, i=1.
\eay
\right.
\dm
Thus $\rE(F_2,\dn)\reb\rE(F_2\ti\dZ\ti\dZ_2,2).$
However, $F_2\ti\dZ\ti\dZ_2$ admits a homomorphism
into the group $F_\om,$ and then into $F_2$ by the above,
so that $\rE(F_2,\dn)\reb\rE(F_2,2),$ as required.
\epf

We add here a technical lemma,
attributed to Kechris in \cite{h-ban},
that will be used in \grf{summI}.
Recall that equivalences Borel reducible to Borel countable
ones are called \rit{essentially countable}.
The lemma shows that maps much weaker than reductions
lead to the same class.

\ble
\lam{KT}
Suppose that\/ $A\yi X$ are Borel sets, 
$\rE$ is a\/ Borel \er\ on\/ $A,$ and\/ 
$\rho:A\to X$ is a Borel map satisfying the 
following$:$ first, the\/ \dd\rho image of any\/ 
\dde class is at most countable, 
secong,\/ \dd\rho images of different\/ 
\dde classes are pairwise disjoint. 
Then\/ $\rE$ is an essentially countable \eqr.
\ele
\bpf
The relation: 
$x\rR y$ iff $x\yi y\in Y$ belong to the \dd\rho image 
of one and the same \dde class in $A,$ 
is a $\fs11$ \eqr\ on the set $Y=\ran \vt.$ 
Moreover, 
\dm
\rR\sq P=\ans{\ang{x,y}:\neg\;\sus a\yi b\in A\;
(a\nE b\land x=\rho(a)\land y=\rho(b))}\,,
\dm
where $P$ is $\fp11.$
Thus there is a Borel set $U$ with $\rR\sq U\sq P.$
In particular, $U\cap(Y\ti Y)=\rR.$ 
As all \dd\rR equivalence classes are at most countable, 
we can assume that all cross-sections of $U$ are 
at most countable, too.   

To prove the lemma it suffices to find a Borel \eqr\ $\rF$ with 
$\rR\sq\rF\sq U.$ 
Say that a set $Z\sq X$ is \rit{stable} if $U\cap(Z\ti Z)$ 
is an \eqr.
For example, $Y$ is stable. 
We observe that the set 
$D_0=\ans{y:Y\cup\ans{y}\,\text{ is stable}}$
is $\fp11$ and satisfies $Y\sq D_0,$ hence, there is a 
Borel set $Z_1$ with $Y\sq Z_1\sq D_0.$ 
Similarly, 
\dm
D_1=\ans{y'\in Z_1:Y\cup\ans{y,y'}\,
\hbox{ is stable for any }\,y\in Z_1}
\dm
is $\fp11$ and satisfies $Y\sq D_1$ by the definition 
of $Z_1,$ so that there is a Borel set $Z_2$ with 
$Y\sq Z_2\sq D_1.$ 
Generally, we define 
\dm
D_n=\ans{y'\in Z_n:Y\cup\ans{y_1,\dots,y_n,y'}\,
\hbox{ is stable for all }\,y_1,\dots,y_n\in Z_n}\,
\dm
find that $Y\sq D_n,$ and choose a Borel set $Z_n$ 
with $Y\sq Z_n\sq D_n.$ 
Then, by the construction, $Y\sq Z=\bigcap_nZ_n,$ 
and, for any finite 
$Z'\sq Z,$ the set $Y\cup Z'$ is stable, so that 
$Z$ itself is stable, and we can take 
$\rF=U\cap(Z\ti Z)$.
\epf

\punk{Hyperfinite \eqr s}
\las{hfer}


The class of Borel hyperfinite \eqr s has been a topic of
intense study since 1970s.
Papers \cite{djk,jkl,kemi} give a comprehensive account of  
the results obtained regarding hyperfinite relations,
with further references.

\bte
[{{\rm Theorems 5.1 and, partially, 7.1 in \cite{djk} 
and 12.1(ii) in \cite{jkl}}}]
\lam{thf}
The following are equivalent for a Borel \eqr\/ $\rE$ on\/ 
a Polish space\/ $\dX:$
\ben
\tenu{{\rm(\roman{enumi})}}
\itla{thf1}\msur
$\rE\reb\Eo$ and\/ $\rE$ is countable$;$

\itla{thf2}\msur
$\rE$ is hyperfinite$;$

\itla{thf3}\msur
$\rE$ is  hypersmooth and countable$;$

\itla{thf4}
there is a Borel set\/ $X\sq\dnd$ such that\/ $\Ei\res X$ 
is a countable \er\ and\/ $\rE$ is isomorphic, 
via a Borel bijection of\/ $\dX$ onto\/ $X,$  
to\/ $\Ei\res X\;;$~\snos
{This transitional condition refers to $\Ei,$ here an \eqr\ on
$\dnd$ defined so that $x\Ei y$ iff $x(n)=y(n)$ for all but
finite $n$.}
         
\itla{thf5}\msur
$\rE$ is induced by a Borel action of\/ $\dZ,$ the additive 
group of the integers.

\itla{thf6} 
there exists a pair of Borel \er s\/ $\rF\yi\rR$ of type\/ 
$2$ such that\/ $\rE=\rF\lor\rR.$~\snos 
{An \eqr\ $\rF$ is {\it of type $n$\/} if any \ddf class 
contains at most $n$ elements. 
$\rF\lor\rR$ denotes the least \er\ which includes 
$\rF\cup\rR$.}
\een
\ete

Note that all Borel finite \eqr s are smooth by
Lemma~\ref{transv}. 
Accordingly, all hyperfinite \er s are hypersmooth. 
On the other hand, all finite and hyperfinite \er s are 
countable, of course. 
It follows from the theorem that, conversely, 
every hypersmooth countable \eqr\ is hyperfinite.

The theorem also shows that $\Eo$ is a universal 
hyperfinite equivalence.   
(To see that $\Eo$ is hyperfinite, define $x\rF_n y$ 
iff $x\sd y\sq[0,n)$ for $x,\,y\sq\dN.$)

Some other characterizations of hyperfinite \eqr s are
known.
For instance, for a Borel \er\ $\rE$ to be hyperfinite it
is necessary and sufficient that there is a Borel partial
order $\prec $ on the domain of $\rE$ that orders each \dde class
$\eke x$ linearly and with the order type being a suborder
of $\dZ,$ the integers positive and negative.
Thus the \dd\prec order type of $\eke x$ has to be either finite,
or $\om,$ or $\om^\ast$ (the inverse $\om$), or $\om^\ast+\om,$
the order type of $\dZ$ itself.
On this see the references given above.

\bpf
It does not seem possible to prove the theorem by a simple
cyclic argument. The structure of the proof will be the
following:
\dm
\bay{clclclclcllll}
\ref{thf1}&\imp&\ref{thf3}&\imp&\ref{thf4}&\imp& 
\ref{thf5}&\imp&\ref{thf1}&;\\[1\dxii]

\ref{thf5}&\imp&\ref{thf6}&\imp&\ref{thf5}&;\\[1\dxii]

\ref{thf5}&\imp&\ref{thf2}&\imp&\ref{thf3}&.\\[1\dxii]
\eay
\dm

The implications 
$\ref{thf2}\imp\ref{thf3}$ and $\ref{thf1}\imp\ref{thf3}$ 
are quite elementary. 

$\ref{thf3}\imp\ref{thf4}.$ 
Let $\rE=\bigcup_n\rF_n$ be a countable and hypersmooth \er\
on a space $\dX,$ all $\rF_n$ being smooth (and countable), 
and $\rF_{n}\sq\rF_{n+1},\:\kaz n.$ 
We may assume that $\dX=\dn$ and $\rF_0=\rav{\dn}.$ 
Let $T_n\sq\dX$ be a Borel transversal for $\rF_n$ 
(recall Lemma~\ref{transv}\ref{transv3}).
Now let $\vt_n(x)$ be the only element of $T_n$ with 
$v\rF_n{\vt_n(x)}.$ 
Then $x\mapsto\sis{\vt_n(x)}{n\in\dN}$ is a $1-1$ Borel map 
$\dX\to \dnd$ and ${x\rE y}\eqv{\vt(x)\Ei\vt(y)}.$ 
Take $X$ to be the image of $\dX$.

$\ref{thf4}\imp\ref{thf5}.$
Let $X$ be as indicated. 
For any \dd\dN sequence $x$ and $n\in\dN,$ let 
$x\qc n=x\res{(n,\iy)}.$ 
It follows from (the relativized versions of)
\Cpro\ and \Cenu\ (theorems \ref{cpro} and \ref{cenu})
that for any $n$ the set
$X\qc n=\ens{x\qc n}{x\in X}$ is Borel 
and there is a countable family of Borel functions 
$g^n_i:X\qc n\to X\zT i\in\dN,$ such that the set  
$X_\xi=\ens{x\in X}{x\qc n=\xi}$ is equal to 
$\ens{g^n_i(\xi)}{i\in\dN}$ for any $\xi\in X\qc n.$ 
Then it holds   
$\ens{g^n_i(\xi)(n)}{i\in\dN}=\ens{x(n)}{x\in X_\xi}.$ 
        
For any $x\in\dnd$ let $\vpi(x)=\sis{\vpi_n(x)}{n\in\dN},$ 
where $\vpi_n(x)$ is the least number $i$ such that 
$x(n)=f^n_i(x)(n);$ thus, $\vpi(x)\in\dnn.$
Let $\mu(x)$ be the sequence
\dm
\vpi_0(x),\vpi'_0(x),
\vpi_1(x)+1,\vpi'_1(x)+1,\dots,
\vpi_n(x)+n,\vpi'_n(x)+n,\dots,
\dm
where $\vpi'_n(x)=\tmax_{k\le n}\vpi_k(x).$ 
Easily if $x\ne y\in X$ satisfy $x \Ei y,$ \ie, 
$x\qc n=y\qc n$ for some $n,$ then 
$\vpi(x)\qc n= \vpi(y)\qc n$ but 
$\vpi(x)\ne \vpi(y)\zT \mu(x)\ne\mu(y),$ 
and $\mu(x)\qc m= \mu(y)\qc m$ for some $m\ge n.$ 
  
Let $\alex$ be the anti-lexicographical partial order on 
$\dnn,$ \ie, $a\alex b$ iff there is $n$ such that 
$a\qc n=b\qc n$ and $a(n)<b(n).$  
For $x,\,y\in X$ define $x<_0 y$ iff $\mu(x)\alex \mu(y).$ 
It follows from the above that $<_0$ linearly orders every 
\dd\Ei class $\ek x\Ei\cap X$ of $x\in X.$ 
Moreover, it follows from the definition of $\mu(x)$ that 
any \dd\alex interval between some $\mu(x)\alex \mu(y)$ 
contains only finitely many elements of the form $\mu(z).$ 
(For $\vpi$ this would not be true.)
We conclude that any class $\ek x\Ei\cap X\zT x\in X,$ is 
linearly ordered by $<_0$ similarly to a subset of $\dZ,$ 
the integers. 
That $<_0$ can be converted to a required Borel action of 
$\dZ$ on $X$ is rather easy 
(however those \dd\Ei classes in $X$ ordered similarly 
to $\dN,$ the inverse of $\dN,$ or finite, should be 
treated separately). 

$\ref{thf5}\imp\ref{thf2}.$
Assume \noo\ that $\dX=\dn.$
An increasing sequence of \er s $\rF_n$ whose union 
is $\rE$ is defined separately on each \dd\rE class 
$C;$ they ``integrate'' 
into Borel \er s $\rF_n$ defined on the whole of $\dn$ 
because the action allows 
to replace quantifiers over a \dde class $C$ by quantifiers 
over $\dZ.$

Let $C$ be any \dde class of $x\in X.$ 
Note that if an element $x_C\in C$ can be chosen in some 
Borel-definable way then we can define $x\rF_n y$ iff there 
exist integers $j,\,k\in\dZ$ with $|j|\le n\zT |k|\le n,$ 
and $x=j\app x_C\zT y=k\app x_C.$  
This applies, for instance, when $C$ is finite, thus, we 
can assume that $C$ is infinite. 
Let $\lex$ be the lexicographical ordering of $\dn,$ 
and $\act$ be the partial order induced by the action, 
\ie, $x\act y$ iff $y=j\app x\zT j>0.$ 
By the same reason we can assume that neither of  
$a=\tinf_{\lex}C$ and $b=\tsup_{\lex}C$ belongs to $C.$ 
Let $C_n$ be the set of all $x\in C$ with 
$x\res n\ne a\res n$ and $x\res n\ne b\res n.$ 
Define $x\rF_n y$ iff $x,\,y$ belong to one and the same 
\dd\lex interval in $C$ lying entirely within $C_n,$ or 
just $x=y.$ 
In our assumptions, any $\rF_n$ has finite classes, and for 
any two $x,\,y\in C$ there is $n$ with $x\rF_n y$. 

$\ref{thf5}\imp\ref{thf1}.$
This is more complicated. 
A preliminary step is to show that $\rE\reb\rE(\dZ,\dn),$ 
where $\rE(\dZ,\dn)$ is the orbit equivalence induced by the 
shift action of $\dZ$ on $(\dn)^\dZ$: 
$k\app\sis{x_j}{j\in\dZ}=\sis{x_{j-k}}{j\in\dZ}$ for 
$k\in\dZ.$ 
Assuming \noo\ that $\rE$ is a \er\ on $\dn,$ we obtain a 
Borel reduction of $\rE$ to $\rE(\dZ,\dn)$ by 
$\vt(x)=\sis{j\app x}{j\in\dZ},$ where $\app$ is a Borel  
action of $\dZ$ on $\dn$ which induces $\rE.$ 
Then Theorem~7.1 in \cite{djk} proves that 
$\rE(\dZ,\dn)\reb\Eo$.

$\ref{thf6}\imp\ref{thf5}.$
Suppose that $\rE=\rF\lor\rR,$ where $\rF\zd\rR$ are 
type-2 \eqr s on $\dn.$
Let a \ddf pair be any pair $\ans{a,b}$
in $\dn$ such that $a\rF b.$
Let a \ddf singleton be any $x\in\dn$
\ddf equivalent only to itself.
Then any $x\in\dn$ is either a \ddf singleton or a member
of a unique \ddf pair.

Fix an arbitrary $x\in\dn.$
We now define an oriented chain $\to$ on the equivalence
class $\eke x.$
For any \ddf pair $\ang{a,b}$ in $\dX$
define $a\to b$ whenever $a\lex b,$ where $\lex$ is
the lexicographical order on $\dn.$  
If $\ans{a\lex b}$ and $\ans{a'\lex b'}$ are different
\ddf pairs then define $b\to a'$ whenever either 
$b\rR a'$ or $b\rR b'.$
(These two conditions are obviously incompatible.) 
If $c$ is a \ddf singleton then define 
$b\to c$ whenever $b\rR c,$
and $c\to a$ whenever $c\rR a.$
If finally $c\ne d$ are \ddf singletons then 
define $c\to d$ whenever $c\rR d$ and $c\lex d.$

If $\eke x$ has no endpoints in the sense of $\to$ then
either
\dm
\eke x=\ans{\dots\to a_{-2}\to a_{-1}\to a_0\to a_1\to a_2\to\dots}
\dm
is a bi-infinite chain or
$\eke x=\ans{a_{1}\to a_{2}\to a_3\to \dots\to a_n\to a_1}$
is a finite cyclic chain.
In the first subcase we straightforwardly define
an action of $\dZ$ on $\eke x$ by
$1\ap a_n=a_{n+1}\zd \kaz n\in\dZ.$
In the second subcase put $1\ap a_k=a_{k+1}$ for $k<n,$
and $1\ap a_n=a_1.$
If $\eke x=\ans{a_{1}\to a_{2}\to a_3\to \dots\to a_n}$ is a
chain with two endpoints then the action is defined the same
way.
If finally $\eke x$ is a chain with just one endpoint, say
$\eke x=\ans{a_0\to a_1\to a_2\to\dots},$ then put
$1\ap a_{2n}=a_{2n+2},$ $1\ap a_{2n+3}=a_{2n+1},$ and
$1\ap a_{1}=a_{0}$. 

$\ref{thf5}\imp\ref{thf6}.$
The authors of \cite{jkl} present a short proof which 
refers to several difficult theorems on hyperfinite \er s. 
Here we give an elementary proof.

Let $\rE$ be induced by a Borel action of $\dZ.$ 
We are going to define $\rF$ and $\rR$ on any \dde class 
$C=\eke x.$
If we can choose an element $x_C\in C$ in some uniform 
Borel-definable way then a rather easy construction is 
possible, which we leave to the reader.
This applies, for instance, when $C$ is finite, hence, 
let us assume that $C$ is infinite. 
Then the linear order $\act$ on $C$ induced by the action 
of $\dZ$ is obviously similar to $\dZ.$ 
Let $\lex$ be the lexicographical ordering of $\dn=\dom\rE.$ 

Our goal is to define $\rF$ on $C$ so that every \ddf class 
contains exactly two (distinct) elements. 
The ensuing definition of $\rR$ is then rather simple. 
(First, order pairs $\ans{x,y}$ of elements of $C$ in 
accordance with the \dd\act lexicographical ordering of 
pairs $\ang{\tmax_{\act}\ans{x,y},\tmin_{\act}\ans{x,y}},$ 
this is still similar to $\dZ.$  
Now, if $\ans{x,y}$ and $\ans{x',y'}$ are two \ddf classes, 
the latter being the next to the former in the sense just 
defined, and $x\act y\zT x'\act y',$ then define $y\rR x'.$)

Suppose that $W\sq C.$ 
An element $z\in W$ iz {\it lmin\/} (locally minimal)
in $W$ if it is \dd\lex smaller than both of its 
\dd\act neighbours in $W.$
Put $W_{\text{lmin}}=\ens{z\in W}{z\,\text{ is lmin in }\,W}.$
If $C_{\text{lmin}}$ is \poq{not} unbounded in $C$ in both 
directions 
then an appropriate choice of $x_C\in C$ is possible. 
(Take the \dd\act least or \dd\act largest point in 
$C_{\text{lmin}},$ 
or if $C_{\text{lmin}}=\pu,$ so that, for instance, 
$\act$ and $\lex$ coincide on $C,$ we can choose something 
like a \dd\lex middest element of $C.$)
Thus, we can assume that $C_{\text{lmin}}$ is unbounded in 
$C$ in both directions. 

Let a {\it lmin-interval\/} be any \dd\act semi-interval 
$[x,x')$ between two consecutive elements $x\act x'$ 
of $C_{\text{lmin}}.$ 
Let $[x,x')=\ans{x_0,x_1,\dots,x_{m-1}}$ be the enumeration 
in the \dd\act increasing order ($x_0=x$). 
Define $x_{2k}\rF x_{2k+1}$ whenever $2k+1<m.$ 
If $m$ is odd then $x_{m-1}$ remains unmatched. 
Let $C^1$ be the set of all unmatched elements. 
Now, the nontrivial case is when $C^1$ is unbounded in $C$ in 
both directions. 
We define $C^1_{\text{lmin}},$ as above, and repeat the same 
construction, extending $\rF$ to a part of $C^1,$ with, perhaps, 
a remainder $C^2\sq C^1$ where $\rF$ remains indefined. 
{\sl Et cetera\/}.

Thus, we define a decreasing sequence 
$C=C^0\supseteq C^1 \supseteq C^2 \supseteq \dots$ of subsets 
of $C,$ and the equivalence relation $\rF$ on each difference 
$C^n\dif C^{n+1}$ whose classes contain exactly two points each, 
and the nontrivial case is when every $C^n$ is \dd\act unbounded 
in $C$ in both directions. 
(Otherwise there is an appropriate choice of $x_C\in C.$) 
If $C^\iy=\bigcap_nC^n=\pu$ then $\rF$ is defined on $C$ and we 
are done. 
If $C^\iy=\ans{x}$ is a singleton then $x_C=x$ chooses an 
element in $C.$ 
Finally, $C^\iy$ cannot contain two different elements as 
otherwise one of $C^n$ would contain two \dd\act neighbours 
$x\act y$ which survive in $C^{n+1},$ which is easily impossible.
\epf

\punk{Non-hyperfinite countable equivalence relations}
\las{no-hypf}

It follows from Theorem~\ref{thf}\ref{thf1},\ref{thf2}
that hyperfinite \eqr s form an initial segment,
in the sense of $\reb,$ within the collection of all
Borel countable \eqr s.
Let us show that this is a proper initial segment,
that is, not all Borel countable \eqr s
are hyperfinite.

\bte
\lam{Eynon}
The \eqr\ $\Ey$ is not hyperfinite, in particular\/
$\Eo\rebs\Ey$.
\ete
\bpf
We present the original proof of this result given in \cite{slst}.
There is another, more complicated proof, based on the fact
that a certain property called {\it amenability\/} holds for
all hyperfinite \eqr s and associated groups like $\stk\dZ+,$
but fails for $\Ey$ and the group $F_2$ --- see
\cite{k-amen,jkl} and references there for details.

Given a pair of bijections $f,g:\dn\onto\dn,$ we define an
action $\goa_{fg}$ of the free group $F_2$ with two generators
$a,b$ on $\dn$ as follows:
if $w=a_1a_2\dots a_n\in F_2$  then
$\goa_{fg}(w,x)=w\app x=
h_{a_1} ( h_{a_2} (\dots(h_{a_n}(x)) \dots)),$
where
$h_a=f,$ $h_{a\obr}=f\obr,$ $h_b=g,$ $h_{b\obr}=g\obr.$
Separately  $\La\app x=x,$ where $\La,$ the empty word, is the
neutral element of $F_2.$ 
The maps $f,g$ are \rit{independent},
\index{maps!independent}%
iff the action is free, that is, for any $x,$
$w\app x=x$ implies $w=\La$.    

To prove the theorem we define a free action 
of $F_2$ on $\dn$ by 
\rit{Lipschitz homeomorphisms}, \ie\  
\index{Lipschitz homeomorphism}%
those homeomorphisms $f:\dn\onto\dn$ satisfying 
${{x\res n}={y\res n}}\eqv {f(x)\res n=f(y)\res n}$
for all $n$ and $y\in\dn.$
Such an action can be extended to any
set $2^n=\ens{s\in\dln}{\lh s=n}$ so that
$w\app{(x\res n)}={(w\app x)}\res n$ for all $x\in\dn.$

\vyk{
Any pair of bijections $f,g:\dn\onto\dn$ induces  
an action of $F_2,$ the free group with two 
generators $a\zi b,$ on $\dn$ by
$a\app x=f(x),$ $b\app x=g(x),$ and generally for any 
element $w=a_1a_2\dots a_n\in F_2,$
$w\app x=h_w(x),$ where
$h_w(x)=
h_{a_1}(h_{a_2}(\dots h_{a_n}(x)\dots))$
and $h_a=f,$ $h_{a\obr}=f\obr,$ $h_b=g,$ $h_{b\obr}=g\obr.$
Separately, $h_\La$ is the identity, where
$\La,$ the empty word, is the neutral element of $F_2.$
The action is free iff $f,g$ are \rit{independent}, meaning 
that $h_w(x)\ne x$ whenever $x\in\dn$ and $w\in F_2\yt w\ne\La.$
}%

\ble
\lam{eyl}
There exists an independent pair of Lipschitz homeomorphisms\/
$f,g:\dn\onto\dn.$ 
\ele
\bpf
Define $f\res 2^n$ and $g\res 2^n$ by induction on $n.$
We'll take care that
\urav{%
\label{kond}
\lh{f(s)}=\lh{g(s)}=\lh s,\quad
f(s)\su f(s\we i),\quad\text{and}\quad g(s)\su g(s\we i)
}%
for all $s\in\dln$ and $i=0,1.$
Fix a linear ordering of length $\om,$
of the set of all pairs $\ang{w,s}\in F_2\ti\dln$ such that
$w\ne\La.$

Put $f(\La)=g(\La)=\La$ ($n=0$) and
$f(\ang i)=g(\ang i)=\ang{1-i}\yt i=0,1.$

To carry  out the step $n\to n+1,$
suppose that the values $f(s)\zi g(s),$ and subsequently  
$w\app s$ for all $w\in F_2,$ have been defined for all
$w\in F_2$ and $s\in\dln$ with $\lh s\le n.$
%
Let $\ang{w_n,s_n}$ be the least pair
(in the sense of the ordering mentioned above)
such that $k=\lh {s_n}\le n,$ there is $t\in 2^n$ with
$s_n\sq t$ and $w_n\app t=t,$ and $u\app s_n\ne v\app s_n$
for all initial subwords\snos
{$\La$ and $w$ itself are considered as initial subwords 
of any word $w\in F_2$.}
$u\ne v$ of $w_n$ --- except for the case when $u=\La$ and
$v=w_n$ or vice versa.
(Pairs $\ang{w,s}$ of this kind do exist:
as $2^n$ is finite, for any $s\in 2^n$ there
is $w\in F_2\bez\ans\La$ such that $w\app s=s.$)

We put $T_n=\ens{t\in 2^n}{s_n\sq t\land w_n\app t=t}.$
The sets
\dm
C_t=\ens{u\app t}
{u\,\text{ is an initial subword of }\,w_n}\,,
\quad t\in T_n\,,
\dm
are pairwise disjoint.
Indeed if $u\app t_1=v\app t_2=t',$
where $u,v$ are initial subwords of $w_n,$ then $u\ne v$
as otherwise $t_1=u\obr\app t'=v\obr\app t')=t_2.$
But then $u\app s_n=v\app s_n$
(as $t_1,t_2$ extend $s_n$),
which contradicts the choice of $s_n$.

Consider any $t\in T_n.$
The word $w_n$ has the form $a_0\,a_1\,\dots\,a_{m-1}$
for some $m\ge1,$
where all $a_\ell$ belong to $\ans{a,b,a\obr,b\obr}.$
Then $C_t=\ans{t_0,t_1,\dots,t_m},$ where $t_0=t$  
and $t_{\ell+1}=a_\ell\app t_\ell\yd\kaz\ell.$
Easily $t_m=w_n\app t=t=t_0,$ but $t_\ell\ne t_{\ell'}$ 
whenever $\ell<\ell'<m.$  
We define  
$a_0\app {(t_0\we i)}=t_1\we(1-i)$ for $i=0,1,$ but
$a_\ell\app{(t_\ell\we i)}=t_{\ell+1}\we i$ whenever $1\le\ell<m.$
Then easily $w_n\app{(t\we i)}=t\we(1-i)\ne t.$

Note that this definition of {\ubf some} of the values
of $a\app r,b\app r,a\obr\app r,b\obr\app r\yd r\in 2^{n+1},$
is self-consistent.\snos
{The inconsistency would have appeared in the case
$a_{m-1}\obr=a_0.$
Then $a_0\app{(t_0\we i)}=t_1\we(1-i)$
while $a_{m-1}\obr\app{(t_m\we i)}=t_{m-1}\we i,$ and $t_0=t_m.$
However $a_{m-1}\obr\ne a_0,$ since otherwise
$a_0\obr s_n={(a_0\,\dots\,a_{m-2})}\app s_n,$
contrary to the choice of $s_n$.}
Thus it remains consistent on the union of all \lap{cycles}
$C_t\yt t\in T_n.$
It follows that the action of $f$ and $g$ can be defined on
$2^{n+1}$ so that \eqref{kond} holds, while the values of 
$a_\ell\app{(t_\ell\we i)}$ coincide with the abovedefined ones
within each cycle $C_t\yt t\in T_n.$
Then $w_n\app{(t\we i)}\ne t\we i$ for all $t\in T_n\yt i=0,1.$ 
It follows that there can be no pair
$\ang{w_{n'},s_{n'}}\yt n'>n,$  equal to $\ang{w_n,s_n}$.

This definition results in a pair of Lipschitz homeomorphisms
$f,g$ of $\dn.$
To check the independence, suppose towards the contrary that 
$x\in\dn\yt w\in F_2\yt w\ne\La,$ and $w\app x=x,$
and there is no shorter word $w$ of this sort.
Then there exists $k\in\dN$ such that $s=x\res k$
satisfies $u\app s\ne v\app s$ for all initial subwords
$u\ne v$ of $w$ except for the case $u=\La$ and $v=w$
(or vice versa).
The pair $\ang{w,s}$ is equal to $\ang{w_n,s_n}$ for some $n\ge k.$
Then the set $T_n$ contains the element $t=x\res n.$
Put $i=x(n).$
Then by definition
$w\app{(t\we i)}={(w\app t)}\we (1-i)=t\we (1-i)\ne t\we i,$
contary to the assumption $w\app x=x$.
\epF{Lemma}

Fix a pair of independent Lipschitz homeomorphisms
$f,g:\dn\onto\dn.$
Define the action $\goa(w,x)=w\app x$ as above.
This Polish (even \lap{Lipschitz}) action of $F_2$ on $\dn$
induces a Borel countable
\eqr{} $x\rE y$ iff $\sus w\in F_2\:(y=w\app x).$
Let us show that $\rE$ is not hyperfinite.

Suppose towards the contrary that $\rE=\bigcup_n\rF_n$ where
$\sis{\rF_n}{n\in\dN}$ is a \dd\sq increasing sequence of finite 
Borel \eqr s.
For any $x$ let $n_x$ be the least $n$ such that  
$\ans{f(x)\yo g(x)\yo f\obr(x)\yo g\obr(x)}$
is a subset of $\ek x{\rF_n}.$
Then there exist a number $n$ and a closed
$X\sq\dn$ such that $n_x\le n$ for all $x\in X,$
and $\mu(X)\ge 3/4,$
where $\mu$ is the uniform probability measure on $\dn.$

Define the subtree $T=\ens{x\res m}{x\in X\land m\in\dN}$ of $\dln.$
We claim that
the set $U$ of all pairs $\ang{w,s}\in F_2\ti\dln$ such that
$\lh w=\lh s$ and
$u\app s\in T$ for any initial subword $u$ of $w$
(including $\La$ and $w$)
is infinite. 

To prove this fact fix $\ell\in\dN$ and find 
$\ang{w,s}\in U$ such that $\lh s=\lh w\ge\ell.$ 
By the independence of $f,g,$ we have $w\app x\ne x$ for all
$w\in W=\ans{a,b,a\obr,b\obr}$ and $x\in\dn,$ in addition
$w\app x\ne {w'}\app x$ for any $w\ne w'$  in $W.$
Then by K\"onig that there is a number $m\ge \ell$
such that $w\app s\ne s$ and $w\app s\ne {w'}\app s$ for all
$w\ne w'$ in $W$ and all $s\in 2^m.$
Note that the graph
\dm
\Ga=\ens{\ans{s,t}}
{s,t\in 2^m\land\sus w\in W\:(w\app s=t)}
\dm 
on $2^m$ has exactly $2\cdot 2^m$ edges:
indeed, by the choice of $m$ for every $s\in 2^m$
there exist exactly $4$ different nodes $t\in 2^m$ such that
$\ans{s,t}\in \Ga.$ 

Consider the subgraph
$G=\ens{\ans{s,t}\in\Ga}{s,t\in T}.$
The intersection $T\cap2^m$ contains at least $\frac34\cdot 2^m$
elements (as $X$ is a set of measure $\ge 3/4$),
accordingly the difference $2^m\dif T$ contains at most 
$\frac 14\cdot 2^m$ elements.
Thus comparably to $\Ga$ the subgraph $G$ loses at most
$4\cdot\frac 14\cdot 2^m=2^m$ edges.
In other words, $G,$ a graph with $\le2^m$ nodes,
has at least $2\cdot 2^m-2^m=2^m$ edges.

Now we apply the following combinatorial fact.

\ble
\lam{eyu}
Any graph\/ $G$ on a finite set\/ $Y,$ containing
not more nodes than edges, has a cycle with at least three
nodes.
\vyk{\ie,  
there exist pairwise different nodes\/
$y_0,y_1,\dots,y_{k-1}\in Y\;\,(k\ge 3)$ 
such that\/ $\ans{y_i,y_{i+1}}\in G$
for all\/ $i<k,$ and in addition\/ $\ans{y_{k-1},y_{0}}\in G$.
}%
\ele
\bpf[\rm Sketch]
Otherwise $Y$ contains an endpoint, that is, an 
element $y\in Y$ such that $\ans{y,y'}\in G$ holds 
for at most one $y'\in Y\bez\ans y.$
This allows to use induction on the number of nodes.
\epf

Thus $G$ contains a cycle 
$s_0,s_1,\dots,s_{k_1},s_k=s_0.$  
Here $k\ge 3,$ all $s_k$ belong to $T\cap 2^m,$
$s_i\yd i<k,$ are pairwise different,  
and for any $i<k$ there exists $a_i\in W=\ans{a,b,a\obr,b\obr}$
such that $a_i\app s_i=s_{i+1}.$ 
The word $u=a_0a_1\dots a_{k-1}$ is irreducible as 
otherwise $s_{i-1}=s_{i+1}$ for some $0<i<k.$
Moreover the word $uu$
(the concatenation of two copies of $u$)
is irreducible, too, as otherwise $s_{1}=s_{k-1}.$
Therefore $u^m$ (the concatenation of $m$ copies of $u$)
is irreducible as well, and so is its initial subword
$w=u^m\res m.$
It follows that $\ang{w,s_0}\in U,$ as required.

As $U$ is infinite, by K\"onig it contains
an infinite branch, \ie\ there is an (irreducible) word
$w\in{\ans{a,b,a\obr,b\obr}}^\dN$ and $x\in\dn$
such that $\ang{w\res m,x\res m}\in U$ for all $m.$
Then clearly ${(w\res m)}\app x\in X$ for all $m,$ and hence
$x\rF_n {({(w\res m)}\app x)}$ by induction on $m.$
Finally ${(w\res m)}\app x\ne {(w\res m')}\app x$ holds whenever
$m\ne m'$ by the independence of $f,g.$
Thus the equivalence class $\ek x{\rF_n}$ is infinite,
contradiction.

Thus $\rE$ is a countable non-hyperfinite \eqr.
Recall that $\rE\reb\Ey$ by Theorem~\ref{femo}.
Thus $\Ey$ itself is non-hyperfinite as well by 
the equivalence $\ref{thf1}\eqv\ref{thf2}$
of Theorem~\ref{thf}.
\epF{Theorem}

\punk{Assembling countable equivalence relations}
\las{assem:hypf}

The following theorem shows that in certain cases the notion
of being Borel reducible to a given countable Borel \eqr\ is
\dd\fsg additive.
The sum $\rF+\rF$ means the union of two Borel isomorphic
copies of $\rF$
defined on a pair of disjoint (and \ddf disconnected) Borel
sets (in one and the same Polish space).

\bte
\lam{cudf}
Let\/ $\rF$ be a countable Borel \er\ satisfying\/
$\rF+\rF\reb\rF,$ and\/
$\rE$ be a Borel \er\ on a Borel set\/ 
$X=\bigcup_kX_k,$ with all\/ $X_k$ also Borel. 
Suppose that\/ $\rE\res X_k\reb\rF$ for each\/ $k.$ 
Then\/ $\rE\reb\rF$.
\ete
\bpf
It obviously suffices to prove that if $\rE$ is a Borel
\eqr\ defined on the union $X\cup Y$ of disjoint Borel
sets $X,Y,$
$\rF$ is a countable Borel \eqr\ defined on the union
$P\cup Q$ of disjoint Borel sets $P,Q,$ \ddf disconnected in
the sense that $p\nF q$ for all $p\in P\yt q\in Q,$
and $f,g$ are Borel reductions
of resp.\ ${{\rE}\res X}\yd{{\rE}\res Y}$ 
to resp.\ ${{\rF}\res P}\yd{{\rE}\res Q}$ 
then there is a Borel reduction
$h$ of $\rE$ to $\rF.$
%
As $X,Y$ are {\ubf not} assumed to be \dde disconnected,
the key problem is to define $h(y)$ in the case when
$y\in Y$ satisfies $g(y)\in \ran U,$ where
\dm
U=\ens{\ang{p,q}\in P\ti Q}{\sus x\in X\:\sus y\in Y\,
({x\rE y}\land f(x)=p\land f(y)=q)} 
\dm
is a $\fs11$ set.
As $f,g$ are reductions to $\rF,$
$U$ is a subset of the $\fp11$ set
\dm
W=\ens{\ang{p,q}\in P\ti Q}
{\kaz \ang{p',q'}\in U\:({p\rF p'}\eqv{q\rF q'})}\,.
\dm
Therefore by \Sepa\ there is an intermediate Borel set
$V\yd U\sq V\sq W.$

The set $U$ is $1-1$ modulo $\rF$ in the sense that
the equivalence ${p\rF p'}\eqv{q\rF q'}$ holds for any two
pairs $\ang{p,q}$ and $\ang{p',q'}$ in $U.$
The set $V$ does not necessarily have this property.
To obtain a Borel subset of $V$ and still superset of $U,$
$1-1$ modulo $\rF,$ note that $U$ is a subset of the $\fp11$
set
\dm
R=\ens{\ang{p',q'}\in V}
{\kaz \ang{p,q}\in V\,({p\rF p'}\eqv{q\rF q'})}\,.
\dm
It follows that there exists a Borel set $S$ with
$U\sq S\sq R.$
Clearly $S$ is $1-1$ modulo $\rF$ together with $R.$
Since $\rF$ is a countable \eqr, it follows by \Cpro\
and \Cenu\ (Theorems \ref{cpro} and \ref{cenu})
that the set $Z=\ran S$ is Borel and there is a Borel
map $\vt:Z\to P$ such that $\ang{\vt(q),q}\in S$ 
for every $q\in Z.$

In particular, we have $\ran U\sq Z$ and $p\rF\vt(q)$
for all pairs $\ang{p,q}\in U.$ 
In addition, it can be \noo\ assumed that $Z$ is
\ddf invariant, \ie\ $q\in Z\land q'\rF q\imp q'\in Z.$
(Indeed consider the set
$Z'=\ekf Z=\ens{q'}{\sus q\in Z\:(q\rF q')}.$
Note that $\rF$ is the orbit equivalence of a Polish action
of a countable group by Theorem~\ref{femo}.
It follows that there exists a countable system
$\sis{\ba_n}{n\in\dN}$ of Borel isomorphisms of the set
$P\cup Q=\dom\rF$ such that $Z'=\bigcup_{n}\ens{\ba_n(q)}{q\in Z}.$
It follows that $Z'$ is Borel by \Cpro, and by \Cenu\ there
is a Borel map $\za:Z'\to Z$ such that $\za(q')\rF q'$ for
all $q'\in Z'.$
Replace $Z,\vt$ by $Z'$ and the map $\vt'(q')=\za(\vt(q')).$) 

This allows us to define a Borel reduction of $\rE$ to $\rF$
as follows.
Naturally, put $h(x)=f(x)$ for $x\in X.$
If $y\in Y$ and $g(y)\nin Z$ then put $h(y)=g(y),$ while
in the case $g(y)\in Z$ we define $h(y)=\vt(g(y)).$ 
\epf

The condition $\rF+\rF\reb\rF$ holds for many naturally
arising \eqr s $\rF.$
(In fact it is not clear how to cook up a Borel equivalence
not satisfying this reduction.)
In particular it holds for $\rF=\Eo$ and the equalities
$\rF=\rav X.$

\bcor
\lam{cudc}
Suppose that\/ $\rE$ be a Borel \er\ on a Borel set\/ 
$X=\bigcup_kX_k,$ with all\/ $X_k$ also Borel.
If\/ ${\rE}\res X_k$ is smooth (resp.\ hyperfinite)
for all\/ $k$ then\/ $\rE$  itself is smooth
(resp.\ hyperfinite.
If\/ ${{\rE}\res X_k}\reb\Eo$ for all\/ $k$ 
then\/ $\rE\reb\Eo.$
\qed
\ecor

\punk{$\ifi$ is the \dd\reb least ideal!}
\las{ifil}

The proof of the following useful result is based on
a short argument involved in many other results, including 
several proofs in \grf{BAN}.

\bte
\label{jnmt}
\ben
\renu
\itsep
\itla{jnmt1}
{\rm\cite{jn,mat75,tal}}
If\/ $\cI$ is a\/ {\rm(nontrivial)} ideal on\/ $\dN,$ 
\imar{jnmt}
with the Baire property in the topology of $\pn,$ 
then\/ $\ifi\orbpp$ and\/ $\orb\cI\,;$ 

\itla{jnmt2}
however\/ $\rav\dn\rebs\Eo$ strictly, thus\/ $\rav\dn$
is not\/ \dd\eqb equivalent to an \eqr\ of the form\/
$\rei\,;$

\itla{jnmt3}
if\/ $\cI\orbp\cJ$ are Borel ideals, and there is an 
infinite set\/ $Z\sq\dom\cI$ such that\/ 
$\cI\res Z=\pwf Z,$ 
then\/ $\cI\orb\cJ$.
\een
\ete
\bpf
\ref
{jnmt1}
First of all $\cI$ must be meager in $\pn.$ 
(Otherwise $\cI$ would be comeager somewhere, 
easily leading to contradiction.) 
Thus, all $X\sq\dN$ ``generic''~\snos
{That is, Cohen generic in the sense a certain countable
family of dense open subsets of $\pn$.} 
do not belong to $\cI.$ 
Now it suffices to define non-empty finite sets $w_i\sq\dN$ 
with $\tmax w_i<\tmin w_{i+1}$ such that any union of 
infinitely many of them is ``generic''. 
Clearly the following observation yields the result: 
if $D$ is an open dense subset of $\pn$ and $n\in\dN$ then 
there is $m>n$ and a set $u\sq[n,m]$ with $m,\,n\in u$ 
such that any $x\in\pn$ satisfying $x\cap[n,m]=u$ belongs 
to $D$.

Thus we have $\ifi\orbpp\cI.$ 
To derive $\ifi\orb\cI$ cover each $w_k$ by a 
finite set $u_k$ such that $\bigcup_{k\in \dN}u_k=\dN$ 
and still $u_k\cap u_l=\pu$ for $k\ne l$.

\ref{jnmt2}
That $\rav\dn\reb\Eo$ is witnessed by any perfect set 
$X\sq\dn$ which is a
{\it partial\/}
\index{transversal!partial}%
transversal for $\Eo$ 
(\ie, any $x\ne y$ in $X$ are \dd{\Eo}inequivalent). 
On the other hand, $\rav\dn$ is smooth but $\Eo$
is non-smooth by Lemma~\ref{transv}\ref{transv5}.

\ref{jnmt3}
Assume \noo\ that $\cI,\cJ$ are ideals over $\dN.$
Let pairwise disjoint finite sets $w_k\sq\dN$ 
witness $\cI\orbp\cJ.$ 
Put $Z'=\dN\dif Z\yt X=\bigcup_{k\in Z}w_k,$ 
and $Y=\bigcup_{k\in Z'}w_k.$ 
The reduction via $\sis{w_k}{}$ reduces $\pwf Z$ to 
$\cJ\res X$ and $\cI\res Z'$ to $\cJ\res Y.$ 
Keeping the latter, replace the former by a 
\dd\orb like reduction of $\pwf z$ to $\cJ\res Y',$ where 
$Y'=\dN\dif Y,$ which exists by Theorem~\ref{jnmt}.
\epf



\vyk{
It follows from the 2nd dichotomy that if $X\sq \dn$ is a 
Borel set then either $\Eo\res X$ is smooth or 
${\Eo\res X}\eqb \Eo.$ 
Let is show that the existence of a full Borel transversal 
is necessary and sufficient for the smoothness. 
Recall that a \dde{\it transversal\/} is a set 
\index{transversal}%
\index{transversal!partial}%
having exactly one element in common with every \dde class. 
Let a {\it partial transversal\/} be a set 
having $\le1$ element in common with every \dde class, in 
other words, it is required that the set is 
pairwise \dde inequivalent.

\bte
\lam{baby}
Let\/ $X\sq\dn$ be\/ $\is11.$ 
Then exactly one of the following two conditions is 
satisfied$:$ 
\ben
\tenu{{\rm(\Roman{enumi})}}
\itla{bab1}
there is a\/ $\is11$ transversal\/ $S$ for\/ $\Eo\res X\;;$ 

\itla{bab2}
there is a closed set\/ $F\sq X$ such that\/ 
$\Eo\eqb{\Eo\res F}$.
\een
\ete 
\bpf
To see that \ref{bab1} and \ref{bab2} are incompatible, prove

\blt
\lam{transd}
Any\/ $\is11$ partial transversal\/ $S$ for\/ $\Eo$ can be 
covered by a\/ $\id11$ partial transversal\/ $D$.
\elt
\bpf
$P=\ens{x\in\dn}{x\nin\ekeo{S\dif\ans{x}}}$ 
is a $\ip11$ set and $S\sq P:$ 
$P$ is the union of $S$ and all \dd\Eo classes 
disjoint from $S.$  
If $U$ is a $\id11$ set with $S\sq U\sq P$ then 
$D=\ens{x\in U}{\card(\ekeo{x}\cap U)=1}$ 
is as required,  
that $D$ remains $\id11$ follows from the countability  
of \dd\Eo classes.
\epF{Lemma}

Now, if $F$ is as in \ref{bab2} then 
$\Eo\res F$ has a $\id11$ transversal by the lemma, 
hence, $\Eo$ itself has a Borel transversal, then 
$\Eo$ is smooth, contradiction.

To see that one of $\ref{bab1}$ and $\ref{bab2}$ holds, 
define ${x\Eco y}$ iff 
$x\yi y$ belong to the same \dd\Eo invariant $\id11$ sets. 
The rest depends on the relationships between the given 
$X$ and the $\is11$ set  
$H_0=\ens{x}{\ekeo x\sneq\ekco x}$.\vom

{\sl Case 1\/}: 
$X\cap H_0=\pu.$ 
Following the proof of Lemma~\ref{d2eff}, we obtain a 
$\id11$ set $Y$ with $X\sq Y\sq \doP{H_0}$ and a $\id11$ 
reduction $\vt:Y\to\dn$ of $\Eo\res Y$ to $\rav\dn.$ 
The set $P=\ens{\ang{a,y}}{y\in Y\land \vt(y)=a}$ is then a 
$\id11$ set whose cross-sections $P_a=\vt\obr(a)$ are just 
classes of $\Eo\res Y.$ 
Then by \Cenu\ of \nrf{eff} there is $\id11$ set 
$T\sq\dN\ti Y$ such that each $T_n=\ens{y}{T(n,y)}$ is a 
transversal for $\Eo\res Y$ and $Y=\bigcup_n T_n.$ 
Now it takes some efforts, which we leave for the reader 
(in particular \Kres\ is applied) 
to assemble a $\is11$ transversal for $\Eo\res X$ from 
parts of the transversals $T_n$.\vom

{\sl Case 2\/}: 
$Z=X\cap H_0\ne\pu.$ 
Then the argument of Case 2 of the proof of Theorem~\ref{2dih} 
(begin with $Z$ rather than with the whole $H$!) shows that 
there is a perfect set $F\sq Z$ such that 
$\Eo\reb{\Eo\res F}$.\vom

\epF{Theorem~\ref{baby}}
}

\api

\parf{The 1st and 2nd dichotomy theorems}
\las{parf2d}

The following two results are known as the first, or
Silver, and 2nd, or 
\lap{Glimm--Effros}, dichotomy theorems.

\bte
[{{\rm Silver~\cite{sil}}}]
\lam{1dih}
Any\/ $\fp11$
(therefore any Borel)
\eqr\/ $\rE$ on\/ $\dnn$ \bfei\ has at 
most countably many equivalence classes\/ \bfor\ admits a 
perfect set of pairwise\/ \dde inequivalent reals.

In other words, \bfei\/ $\rE\reb\rav\dN$ \bfor\/  
$\rav\dn\emn\rE$.
\ete

\bte
[{{\rm Harrington, Kechris, Louveau \cite{hkl}}}]
\lam{2dih}
If\/ $\rE$ is a Borel \eqr\ then \bfei\/ 
$\rE$ is smooth \bfor\/ $\Eo\emn\rE$.
\ete

Recall that $\emn$ in the \bfor\ part means
the reducibility via a continuous injective map.
Obviously $\emn$ implies $\reb,$ and hence it follows
from the first theorem that the union of the lower
\dd\reb cone of $\rav\dN$ and the upper \dd\reb cone
of $\rav\dn$ fully covers the whole class of Borel \eqr s.
As smoothness means simply $\rE\reb\rav\dn,$
it follows
from the second theorem that the union of the lower
\dd\reb cone of
$\rav\dn$ and the upper \dd\reb cone of $\Eo$ fully covers
the whole class of Borel \eqr s.

The proofs of these theorems  follow below in this \gla.
They make heavy use of methods of effective descriptive set
theory, in particular, the Gandy -- Harrington topology.
We begin with a brief introduction into this technical tool.

This \gla\ ends with an introduction into an interesting
forcing notion that consists of all uncountable Borel sets
$X\sq\dn$ such that $\Eo\res X$ is not smooth.

\punk{The Gandy -- Harrington topology}
\las{bairef}

The following notion is similar to the Choquet property but 
somewhat more convenient to provide the 
nonemptiness of countable intersections of pointsets. 

\bdf
\lam{genb}
A family $\cF$ of sets in a topological space
is {\it Polish--like\/} if there exists a 
\index{Polish--like}%
countable collection $\ens{\cD_n}{n\in\dN}$ of dense subsets 
$\cD_n\sq\cF$ such that we have $\bigcap_nF_n\ne\pu$ 
whenever $F_0\qs F_1\qs F_2\qs\dots$ is a decreasing sequence 
of sets $F_n\in\cF$ which intersects every $\cD_n.$ \ 

Here, a set $\cD\sq\cF$ is {\it dense\/} if 
\index{set!dense}%
$\kaz F\in\cF\:\sus D\in\cD\:(D\sq F).$
\edf

For instance if $\cX$ is a Polish space then the collection 
of all its non-empty closed sets is Polish--like, for take 
$\cD_n$ to be all closed sets of diameter $\le n\obr.$
We'll make use of the following technical fact:

\bte
[{{\rm see \eg\ Kanovei~\cite{umnE}, Hjorth~\cite{h-ban}}}]
\lam{s11pol}
The collection\/ $\cF$ of all non-empty\/ $\is11$ subsets 
of\/ $\dnn$ is Polish--like.\qeD
\ete

\bpf
For any $P\sq\dnn\ti\dnn$ define 
$\pr P=\ens{x}{\sus y\:P(x,y)}$ 
(the projection). 
If $P\sq\dnn\ti\dnn$ and $s\yi t\in\dN\lom$ then let 
$P_{st}=\ens{\ang{x,y}\in P}{s\su x\land t\su y}.$  
Let $\cD(P,s,t)$ be the collection of all $\is11$ sets 
$\pu\ne X\sq\dnn$ such that either $X\cap \pr P_{st}=\pu$ 
or $X\sq\pr P_{s\we i\,\yi t\we j}$ for some $i\yi j.$ 
(Note that in the ``or'' case $i$ is unique but $j$ may be 
not unique.) 
Let $\ens{\cD_n}{n\in\dN}$ be an arbitrary enumeration of all 
sets of the form $\cD(P,s,t),$ where $P\sq\dnn\ti\dnn$ is 
$\ip01.$ 
Note that in this case all sets of the form $\pr P_{st}$ are 
$\is11$ subsets of $\dnn,$ therefore, $\cD(P,s,t)$ is easily a 
dense subset of $\cF,$ so that all $\cD_n\sq \cF$ are dense. 

Now consider a decreasing sequence $X_0\qs X_1\qs \dots$ of 
non-empty $\is11$ sets $X_k\sq\dnn,$ which intersects every 
$\cD_n\,;$ prove that $\bigcap_nX_n\ne \pu.$ 
Call a set $X\sq\dnn$ {\it positive\/} 
if there is $n$ such that $X_n\sq X.$ 
For any $n,$ fix a $\ip01$ set $P^n\sq\dnn\ti\dnn$ such that 
$X_n=\pr P^n.$ 
For any $s\yi t\in\dN\lom,$ if $\pr P^n_{st}$ is positive 
then, by the choice of the sequence of $X_n,$ 
there is a unique $i$ and some
$j$ such that  $\pr P^n_{s\we i\,\yi t\we j}$ is also positive.  
It follows that there is a unique $x=x_n\in\dnn$ and some  
$y=y_n\in\dnn$ (perhaps not unique) such that  
$\pr P^n_{x\res k\,\yi y\res k}$ is positive for any $k.$ 
As $P^n$ is closed, we have $P^n(x,y),$ hence, $x_n=x\in X_n$.

It remains to show that $x_m=x_n$ for $m\ne n.$ 
To see this note that if both $P_{st}$ and 
$Q_{s't'}$ are positive then either $s\sq s'$ 
or $s'\sq s$.
\epf

The collection of all non-empty $\is11$ subsets of $\dnn$ 
is a base of the {\it Gandy -- Harrington topology\/}, which has  
many  remarkable applications in descriptive set theory. 
This topology is not Polish, even not metrizable at all, 
yet it shares the following important property of Polish 
topologies:

\bcor
\lam{s11ba}
The Gandy -- Harrington topology is Baire, that is every 
comeager set is dense.
\ecor
\bpf
This can be proved using Choquet property of the topology, 
see \cite{hkl}, however, the Polish--likeness 
(Theorem~\ref{s11pol}) also immediately yields the result.
\epf

\punk{The first dichotomy theorem.}
\las{>smooth}
  
Beginning the proof\snos
{We present a forcing-style proof of Miller~\cite{am}, with 
some simplifications. 
See \cite{mw} for another proof, based on the Gandy -- 
Harrington topology.
In fact both proofs involve very similar 
combinatorial arguments.}
of Theorem~\ref{1dih}, let us fix a $\fp11$ \eqr\ $\rE$ on $\dnn.$
Then $\rE$ belongs to $\ip11(p)$ for some parameter $p\in\dnn.$ 
As usual, we can suppose that $\rE$ is in fact a lightface
$\ip11$ relation; the case of an arbitrary $p$ does not differ
in any essential detail.\vom

{\sl Case 1\/}: every $x\in\dnn$ belongs to a $\id11$ 
pairwise \dde equivalent set $X.$ 
(A set $X$ is pairwise \dde equivalent iff all elements of $X$ are
\index{set!pairwiseEequiv@pairwise \dde equivalent}%
\dde equivalent to each other, in other words, 
the saturation $\eke X$ is an equivalence class.)  
Then $\rE$  has at most countably many equivalence classes.\vom

{\sl Case 2\/}: otherwise. 
Then the set $H$ (\rit{the domain of nontriviality})
of all $x\in\dnn$ which do {\ubf not} belong to a 
$\id11$ pairwise \dde equivalent set is non-empty.


\bcl
\lam{dih1:1}
$H$ is\/ $\is11.$ 
Any\/ $\is11$ set\/ $\pu\ne X\sq H$ 
is not pairwise\/ \dde equivalent.
\ecl
\bpf
We make use of an enumeration of $\id11$ sets provided by
Theorem~\ref{penu}.
Suppose that $x\in\dnn.$ 
Then obviously $x\in H$ iff for any $e\in\dN:$ 
{\ubf if} $e$ codes a $\id11$ set, say, $W_e\sq\dnn$ and 
$x\in W_e$ 
{\ubf then} $W_e$ is not \dde equivalent. 
The {\ubf if} part of this characterization is $\ip11$ while 
the {\ubf then} part is $\is11$.

If $X\ne\pu$ is a pairwise \dde equivalent $\is11$ set then 
$B=\bigcap_{x\in X}\eke x$ is a $\ip11$ \dde equivalence class  
and $X\sq B.$  
By \Sepa\ (Theorem~\ref{rese}),
there is a $\id11$ set $C$ with $X\sq C\sq B.$ 
Then, if $X\sq H$ then $C\sq H$ is a $\id11$ pairwise 
\dde equivalent set, a 
contradiction to the definition of $H$.
\epF{Claim}

Let us fix a countable transitive model 
\index{modelm@model $\mm$}%
$\mm$ of $\zhc$ (see Remark~\ref{ff:zhc}).
We suppose that $\mm$ is an elementary 
submodel of the universe \vrt\ all analytic formulas\snos
{\label{suf}  
Being an elementary submodel is useful to guarantee that 
relations like the inclusion orders of $\is11$ sets are 
absolute for $\mm$.}.
We consider the set 
$\dP=\ens{X\sq\dnn}{X\,\text{ is non-empty and }\,\is11}$ 
\index{zzP@$\dP$}%
as a forcing to extend $\mm$ 
(smaller sets are stronger conditions) --- 
the {\it Gandy -- Harrington forcing\/}. 
\index{Gandy -- Harrington!forcing}%
\index{forcing!Gandy -- Harrington}%
Obviously $\dP\nin$ and $\not\sq\mm,$ of course, but clearly 
$\dP$ can be adequately coded in $\mm,$ say, via a 
universal $\is11$ set.

\bcor
[{{\rm from Theorem~\ref{s11pol}}}]
\lam{d12}
If\/ $G\sq\dP$ is a\/ \dd\dP generic, over\/ $\mm,$ set, 
then\/ $\bigcap G$ contains a single real, denoted\/ 
$x_G$.\qeD
\ecor

Reals of the form $x_G,$ $G$ as in the corollary, 
are called \dd\dP{\it generic\/} (over $\mm$).
Let $\dox$ be the name for $x_G$ in the machinery of
forcing $\dP.$ 
Then any condition $A\in\dP$ forces that $\dox\in A$.

The forcing product $\dP^2$ consist of all rectangles
$X\ti Y$ with 
\index{zzpu2@$\dP^2$}%
$X\yi Y\in\dP.$  
It follows from the above by the product forcing lemmas 
that any set $G\sq\dP^2$ \dd{\dP^2}generic over $\mm$ 
produces a pair of reals 
(a \dd{\dP^2}{\it generic pair\/}),  
say, $x^G\ul$ and $x^G\ur,$ so that  
$\ang{x^G\ul,x^G\ur}\in W$ for any $W\in G.$  
Let $\doxl$ and $\doxr$ be their names.
The following is the key fact:

\ble
\lam{d14}
$H\ti H$ \dd{\dP^2}forces\/ $\doxl \nE\doxr$.
\ele
\bpf
Otherwise there is a condition $X\ti Y\in\dP^2$ 
with $X\cup Y\sq H$ that \dd{\dP^2}forces $\doxl\rE\doxr,$ 
and hence any \dd{\dP^2}generic pair 
$\ang{x,y}\in X\ti Y$ satisfies $x\rE y.$ 
By the product forcing lemmas for any pair of \dd\dP generic 
$x'\yi x''\in X$ there is $y\in Y$ such that both 
$\ang{x,y}$ and $\ang{x',y}$ are \dd{\dP^2}generic pairs, 
therefore 

\bit
\item[$(\ast)$] 
$x'\rE x''$ holds for any points  $x'\yi x''\in X$
separately \dd\dP generic over\/ $\mm$.
\eit

\vyk{
The $\is11$ set $P=X^2\dif\rE$ is non-empty.
(Otherwise $X\ne\pu$ is a \dde equivalent $\is11$ 
subset of $H.$ 
Then $B=\bigcap_{x\in X}\eke x$ is a $\ip11$ set and 
an \dde equivalence class, and $X\sq B\sq H.$ 
By \Sepa, there is a $\id11$ set $C$ with $X\sq C\sq B.$ 
Then $C\sq H$ is a $\id11$ \dde equivalent set, a 
contradiction to the definition of $H.$)
}

Note that the set $\dP_2$ of all non-empty $\is11$ 
\index{zzpd2@$\dP_2$}%
subsets of $\dnn\ti\dnn$ is just a copy of $\dP$ 
(not of $\dP^2$!) as a forcing.
In particular, if a set $G\sq\dP_2$ is \dd{\dP_2}generic
over $\mm$ then there is a unique pair of reals 
(\dd{\dP_2}{\it generic pair\/}) 
$\ang{x^G\ul,x^G\ur}$ which belongs to every $W$ in $G,$ 
and in this case, 
{\it both\/ $x^G\ul$ and $x^G\ur$ are\/ \dd\dP generic\/},  
because if $G\sq\dP_2$ is \dd{\dP_2}generic then 
the sets $G'$ and $G''$ of all projections of sets $W\in G$ 
to resp.\ 1st and 2nd co-ordinate, are easily 
\dd\dP generic.
Now let $G\sq\dP_2$ be a \dd{\dP_2}generic set, over $\mm,$  
containing the $\is11$ set $P=X^2\dif\rE.$  
(Note that $P\ne\pu$ by Lemma~\ref{dih1:1}.)
Then $\ang{x^G\ul,x^G\ur}\in P,$ hence $x^G\ul\nE x^G\ur.$ 
However, as we observed, both $x^G\ul$ and $x^G\ur$ are 
\dd\dP generic elements of $X$ 
(because $P\sq X\ti X$), thus $x^G\ul\rE x^G\ur$ by $(\ast),$
contradiction. 
\epF{Lemma~\ref{d14}}

Now to accomplish the proof of the theorem
let us fix enumerations 
$\sis{\cD(n)}{n\in\dN}$ and $\sis{\cD^2(n)}{n\in\dN}$ 
of all dense subsets of resp.\ $\dP$ and $\dP^2$ which 
are coded in $\mm.$ 
Then there is a system 
$\sis{X_u}{u\in 2\lom}$ of sets $X_u,$ satisfying 
\ben
\tenu{(\roman{enumi})}
\itla{d1i.}
$X_u\in\dP,$ moreover, $X_\La\sq H$; 

\itla{d1i,}
$X_u\in\cD(n)$ whenever $u\in 2^n;$
\imar{d1i,}

\itla{d1d}
$X_{u\we i}\sq X_u$ for all $u\in2\lom$ and $i=0,1$;

\itla{d1t}
if $u\ne v\in 2^n$ then $X_u\ti X_v\in\cD^2(n)$.
\een
It follows from \ref{d1i,} that, for any $a\in \dn,$ the 
set $\ens{X_{a\res m}}{m\in\dN}$ is \dd\dP generic over $\mm,$ 
hence, $\bigcap_mX_{a\res m}$ is a singleton by 
Corollary~\ref{d12}.
Let $x_a$ be its only element.
The map $a\mapsto x_a$ is continuous because the diameters 
of sets $X_u$ converge to $0$ uniformly with $\lh u\to 0$ by 
\ref{d1i.}. 
In addition, by \ref{d1t} and Lemma~\ref{d14},  
$x_a\nE x_b$ holds for any pair $a\ne b,$ in particular, 
$x_a\ne x_b,$ hence,  
we have a perfect \dde inequivalent set 
$Y=\ens{x_a}{a\in\dn}$.\vom

\qeDD{Theorem~\ref{1dih}}

\punk{The second dichotomy theorem}
\las{ghc}

Beginning the proof of Theorem~\ref{2dih}
(it will be completed in \nrf{d2splitk}),
we suppose, as usual, that $\rE$ is a lightface $\id11$ 
\eqr\ on $\dnn.$
Similarly to Theorem~\ref{1dih},
the proof employs the Gandy -- Harrington topology, but
is considerably more complicated.


Consider an auxiliary \eqr\ $x\Ec y$ iff $x\yi y\in\dnn$ 
belong to the same \dde invariant $\id11$ sets. 
(A set $X$ is \dde{\it invariant\/} iff $X=\eke X.$) 
\index{set!einvariant@\dde invariant}%
Easily $\rE\sq\Ec.$  
In fact it follows from the next lemma that
$\Ec$ is equal to the closure of $\rE$ in the Gandy -- 
Harrington topology.


\ble
\lam{d20}
If\/ $\rF$ is a\/ $\is11$ \er\ on\/ $\dnn,$ and 
$X,Y\sq\dnn$ are disjoint\/ \ddf in\-var\-iant\/ $\is11$ 
sets, then there is an\/ \ddf invariant\/ $\id11$ 
set\/ $X'$ separating\/ $X$ from\/~$Y.$ 
\ele
\bpf
By \Sepa, for any $\is11$ set $A$ with $A\cap Y=\pu$ there is 
a $\id11$ set $A'$ with $A\sq A'$ and $A'\cap Y=\pu$ --- note 
that then $\ekf{A'}\cap Y=\pu$ because $Y$ is \ddf invariant. 
It follows that that there is a sequence 
$X=A_0\sq A'_0\sq A_1\sq A'_1\sq \dots,$  
where 
$A'_i$ are $\id11$ sets, accordingly, $A_{i+1}=\ekf{A'_i}$ 
are $\is11$ sets, and $A_i\cap Y=\pu.$ 
Then $X'=\bigcup_nA_n=\bigcup_nA'_n$ and 
is an \ddf invariant Borel set which separates $X$ from $Y.$  
To ensure that $X'$ is $\id11$ we have to maintain the choice
of sets $A_n$ in effective manner. 

Let $U\sq \dN\ti\dnn$ be a ``good'' universal $\is11$ set.
(We make use of Theorem~\ref{uset}.) 
Then there is a recursive $h:\dN\to\dN$ such that
$\ekf{U_n}=U_{h(n)}$ for each $n.$
Moreover, applying Lemma~\ref{effred}
(to the complement of $U$ as a ``good'' universal $\ip11$ set,
and with a code for $Y$ fixed),
we obtain a pair of recursive functions $f\zi g:\dN\to\dN$
such that for any $n,$ if $U_n\cap Y=\pu$ then
$U_{f(n)}\zi U_{g(n)}$ are complementary $\is11$ sets
(hence, either of them is $\id11$)
containing, resp., $U_n$ and $Y.$ 
A suitable iteration of $h$ and $f,g$ allows us
to define a sequence $X=A_0\sq A'_0\sq A_1\sq A'_1\sq \dots$
as above effectively enough for the union of those sets
to be $\id11$. 
\epF{Lemma}

\ble
\lam{d21}
$\Ec$ is a \/ $\is11$ relation.
\ele
\bpf
Let $C\sq\dN$ and $W\yi W'\sq\dN\ti\dnn$ be as in
Theorem~\ref{penu}. 
The formula $\inva(e)$ saying that $e\in C$ and 
the set $W_e=W'_e$ is \dde invariant, that is,
\dm
e\in C\lland 
\kaz a\yi b\:
\skl {a\in W_e\land b\nin W'_e}\imp {a\nE b}\skp
\dm
--- is obviously a $\ip11$ formula.
On the other hand, $x\Ec y$ iff
\dm
\kaz e\; 
\skl
\inva(e)\limp 
(x\in W_e\imp y\in W'_e) \land (y\in W_e\imp x\in W'_e)
\skp.
\eqno\square\:(\hbox{\sl Lemma})
\dm
\ePf

Let us return to the proof of Theorem~\ref{2dih}. 
We have two cases.\vom

{\ubf Case 1\/}:
$\rE=\Ec,$ that is $\rE$ is Gandy -- Harrington closed. 
The next lemma shows that in this assumption we obtain the
\bfei\ case in Theorem~\ref{2dih}.

\ble
\lam{d2eff}
If\/ $\rE=\Ec$ then there is a\/ $\id11$ reduction of\/ 
$\rE$ to\/ $\rav\dn$.
\ele
\bpf
Let $C\sq\dN$ and $W\yi W'\sq\dN\ti\dnn$ be as in
Theorem~\ref{penu}. 
By \Kres\ (Theorem~\ref{kres})
there is a $\id11$ function $\vpi:X^2\to C$ such 
that $W_{\vpi(x,y)}=W'_{\vpi(x,y)}$ is a \dde invariant 
$\id11$ set containing $x$ but not $y$ whenever $x\yi y\in X$ 
are \dde inequivalent. 
Then $R=\ran \vpi$ is a $\is11$ subset of $C,$ hence, by \Sepa, 
there is a $\id11$ set $N$ with $R\sq N\sq C.$ 
The map $\vt(x)=\ens{n\in N}{x\in D_n}$ is a $\id11$ reduction 
of $\rE$ to $\rav{\dn}$. 
\epF{Lemma and Case 1}

{\ubf Case 2\/}: $\rE\sneq\Ec.$ 
Then the $\is11$ set 
$H=\ens{x}{\eke x\sneq\ekco x}$ 
(the union of all \ddec classes containing more than one 
\dde class)
is non-empty.
We are going to prove that this leads to the \bfor\ case
in Theorem~\ref{2dih}.
This will take some space.
We begin with a couple of technical lemmas.
The first of them says that the property $\rE\sneq\Ec$
holds hereditarily within the key domain $H$.

\ble
\lam{e=oe}
If\/ $X\sq H$ is a\/ $\is11$ set then\/ $\rE\sneq\Ec$ on\/ $X$.
\ele
\bpf
Suppose that ${\rE\res X}={\Ec\res X}.$ 
Then $\rE=\Ec$ on $Y=\eke X$ as well. 
(If $y\yi y'\in Y$ then there are $x\yi x'\in X$ such 
that $x\rE y$ and $x'\rE y',$ so that if $y\Ec y'$ then 
$x\Ec x'$ by transitivity, hence, $x\rE x',$ and $y\rE y'$ 
again by transitivity.) 
It follows that $\rE=\Ec$ on an even bigger set, 
$Z=\ekco X.$ 
(Otherwise the $\is11$ set 
$Y'=Z\dif Y=\ens{z}{\sus x\in X\:(x\Ec y\land{x\nE y})}$ 
is non-empty and \dde invariant, together with $Y,$ hence
by Lemma~\ref{d20} 
there is a \dde invariant $\id11$ set $B$ with $Y\sq B$ 
and $Y'\cap B=\pu$, which implies that 
no point in $Y$ is \dd{\Ec}equivalent to a point in $Y',$ 
contradiction.) 
Then by definition $Z\cap H=\pu$.
\epF{Lemma}

\ble
\lam{d27}
If\/ $A,B\sq H$ are non-empty\/ $\is11$ sets with\/ 
$A\rE B$ then there exist non-empty\/ {\ubf disjoint}\/ $\is11$ 
sets\/ $A'\sq A$ and\/ $B'\sq B$ still satisfying\/  
$A'\rE B'$.
\ele

Recall that $A\rE B$ means that $\eke A=\eke B$.

\bpf
We assert that there are points $a\in A$ and $b\in B$ 
with $a\ne b$ and $a\rE b.$ 

(Otherwise $\rE$ is the equality on $X=A\cup B.$ 
Prove that then $\rE=\Ec$ on $X,$   
a contradiction to Lemma~\ref{e=oe}. 
Take any $x\ne y$ in $X.$ 
Let $U$ be a clopen set containing $x$ but not $y.$ 
Then $A=\eke {U\cap X}$ and $C=\eke {X\dif U}$ are two 
disjoint \dde invariant $\is11$ sets containing resp.\ $x\yi y.$ 
Then $x\Ec y$ fails by Lemma~\ref{d20}.) 

Thus let $a\yi b$ be as indicated. 
Let $U$ be a clopen set containing $a$ but not $b.$ 
Put $A'=A\cap U\cap\eke{\doP U}$ and 
$B'=B\cap{\doP U}\cap\eke U$.
\epF{Lemma}

\punk{Restricted product forcing}
\las{d2forc}

In continuation of the proof of Theorem~\ref{2dih} (Case 2),
we come back to the forcing notions $\dP$ and $\dP_2$ 
introduced in \nrf{>smooth}.
Let us fix a countable model $\mm$ of $\zhc$ chosen as in  
\nrf{>smooth}.

Let $\dpe$ be the collection of all sets of the form 
\index{zzp2e@$\dpe$}%
$X\ti Y,$ where $X\yi Y\sq\dnn$ are non-empty $\is11$ sets and 
$X\rE Y$ (which means here that $\eke X=\eke Y$).  
\index{zzXEY@$X\rE Y$}%
Easily $\dP_2\sq\dpe\sq\dP^2$. 
The forcing 
$\dpe$ is not really a product, yet if $X\ti Z\in\dpe$ and 
$\pu\ne X'\sq X$ is $\is11$ then $Z'=Z\cap\eke{X'}$ is $\is11$ 
and $X'\ti Z'\in\dpe.$ 
It follows that any set $G\sq\dpe,$ \dd\dpe generic over $\mm,$
still produces a pair of \dd\dP generic sets 
$G\ul=\ens{\dom P}{P\in G}$ and 
$G\ur=\ens{\ran P}{P\in G},$ therefore produces a pair 
of \dd\dP generic reals $x^G\ul$ and $x^G\ur,$ whose names 
will be $\doxl$ and $\doxr$ as above.


\ble
\lam{d23}
In the sense of the forcing\/ $\dpe,$ 
any condition\/ $P=X\ti Z$ in\/ $\dpe$
forces\/ $\ang{\doxl,\doxr}\in P$ 
and forces\/ $\doxl\Ec\doxr,$ but\/  
$H\ti H$ forces\/ $\doxl\nE\doxr$.
\ele
\bpf
To see that $\doxl\Ec\doxr$ is forced suppose otherwise. 
Then, by the definition of $\Ec,$ 
there is a condition $P=X\ti Z\in\dpe$ and an  
\dde invariant $\id11$ set $B$
such that $P$ forces  $\doxl\in B$ but $\doxr\nin B.$ 
Then easily $X\sq B$ but $Z\cap B=\pu,$ a contradiction with 
$\eke X=\eke Z$.

To see that $H\ti H$ forces $\doxl\nE\doxr$ suppose towards 
the contrary that some $P=X\ti Z\in\dpe$ with $X\cup Z\sq H$ 
forces $\doxl\rE\doxr,$ thus,

\ben
\tenu{(\arabic{enumi})}
\itla{d23-}
$x\rE z$ holds 
for every \dd\dpe generic pair $\ang{x,z}\in P$. 
\een

\bcl
\lam{d24}
If\/ $x,y\in X$ are\/ \dd\dP generic over\/ $\mm,$ 
and\/ $x\Ec y,$ then\/ $x\rE y$.
\ecl
\bpf
We assert that 
\ben
\tenu{(\arabic{enumi})}
\addtocounter{enumi}1
\itla{d23a}
${x\in A}\eqv {y\in A}$ holds
for each \dde invariant $\is11$ set $A$. 
\een
Indeed, if, say, $x\in A$ but $y\nin A$ then by the 
genericity of $y$ there is a $\is11$ set $C$ with 
$y\in C$ and $A\cap C=\pu.$ 
As $A$ is \dde invariant, Lemma~\ref{d20} yields an 
\dde invariant $\id11$ set $B$ such that $C\sq B$ but 
$A\cap B=\pu.$ 
Then $x\nin B$ but $y\in B,$ a contradiction to 
$x\Ec y$.

Let 
$\sis{\cD_n}{n\in\dN}$ be an enumeration 
of all dense subsets of $\dpe$ which are coded in $\mm.$ 
We define 
two sequences $P_0\qs P_1\qs\dots$ and 
$Q_0\qs Q_1\qs\dots$ of conditions $P_n=X_n\ti Z_n$ and 
$Q_n=Y_n\ti Z_n$ in $\dpe,$ so that $P_0=Q_0=P,$ 
$x\in X_n$ and $y\in Y_n$ for any $n,$ and finally 
$P_n\yi Q_n\in\cD_{n-1}$ for $n\ge1.$ 
If this is done then we have a real $z$ 
(the only element of $\bigcap_nZ_n$) 
such that both $\ang{x,z}$ and $\ang{y,z}$ are \dd\dpe generic, 
hence, $x\rE z$ and $y\rE z$ by \ref{d23-}, hence, 
$x\rE y$. 

Suppose that $P_n$ and $Q_n$ have been defined.  
As $x$ is generic, there is (we leave details for the reader) 
a condition $P'=A\ti C\in \cD_{n}$ and $\sq P_n$ such that  
$x\in A.$ 
Let $B=Y_n\cap\eke A:$ then $y\in B$ by \ref{d23a}, and 
easily $\eke B=\eke C=\eke A$ 
(as $\eke{X_n}=\eke{Z_n}=\eke{Y_n}$), thus, $B\ti C\in\dpe,$ 
so there is a condition $Q'=V\ti W\in\cD_n$ and 
$\sq B\ti C\sq Q_n$ such that $y\in V.$ 
Put $Y_{n+1}=V,\msur$ $Z_{n+1}=W,$ and $X_{n+1}=A\cap\eke W$.
\epF{Claim}

It follows that $\rE=\Ec$ on\/ $X.$  
(Otherwise $S=\ens{\ang{x,y}\in X^2}{x\Ec y\land x\nE y}$ 
is a non-empty $\is11$ set,  
and any \dd{\dP_2}generic pair $\ang{x,y}\in S$ implies a 
contradiction to Claim~\ref{d24}.
Recall that $\dP_2=$ all non-empty $\is11$ subsets of 
$\dnnp2.$) 
But this implies $X\cap H=\pu$ by Lemma~\ref{e=oe}, 
contradiction.
\vyk{
To accomplish the proof of Lemma~\ref{d23}, note that 
\dm
X\sq C'=\ens{c\in\dnn}{\kaz x\in X\:(x\Ec c\imp x\rE c)}\,,
\quad\text{and}\quad C'\:\hbox{ is }\,\ip11\,,
\dm 
hence, there is a $\id11$ set $B'$ with $X\sq B'\sq C'.$ 
Further, 
\dm
X\sq C=\ens{c\in B'}{\kaz b\in B'\:(b\Ec c\imp b\rE c)}\,,
\dm 
hence, there is a $\id11$ set $B$ with $X\sq B\sq C,$ 
and by definition ${\rE\res B}={\Ec\res B}.$ 
It easily follows that $\rE\yi \Ec$ coincide also on $W=\eke B,$ 
hence, $W$ cannot intersect $H.$ 
However $\pu\ne X\sq B\cap H,$ contradiction.
}
\epF{Lemma~\ref{d23}}

\vyk{
\blt
\lam{d27}
$H\ti H$ forces, in the sense of\/ $\dP_2,$ that\/ 
$\doxl\ne\doxr$.
\elt
\bpf
Otherwise there is a condition in $\dP_2,$ of the form 
$X\ti X,$ where $\pu\ne X\sq H$ is a non-empty set, such 
that ${\rE}\res X$ is the equality. 
Let us prove that then $\rE=\rF$ on $X:$ then we will have 
a contradiction as in the end of the proof of Lemma~\ref{d23}. 
Take any $x\ne y$ in $X.$ and
To show $x\nF y$ let $U$ be a clopen set containing $x$ but 
not $y.$ 
Then $A=\eke {U\cap X}$ and $C=\eke {X\dif U}$ are two 
disjoint \dde invariant $\is11$ sets containing resp.\ $x\yi y.$ 
Then by Lemma~\ref{d20} there is a \dde invariant $\is11$ set 
$B$ containing $x$ but not $y,$ as required.
\epF{Lemma}
}

\punk{Splitting system}
\las{d2split}

The \bfor\ case of Theorem~\ref{2dih}, that is $\Eo\emn\rE,$
means that $\Eo$ has a continuous \lap{copy} of the form
${\rE}\res X,$ $X$ being a closed set in $\dnn.$
To obtain such a set, we define a splitting system of sets
in $\dP$ satisfying certain requirements.

Let us fix 
enumerations 
$\sis{\cD(n)}{n\in\dN},\msur$ $\sis{\cD_2(n)}{n\in\dN},\msur$  
$\sis{\cD^2(n)}{n\in\dN}$ 
of all dense subsets of resp.\ $\dP\yd \dP_2\yd\dpe,$ which 
belong to the model $\mm$ fixed above.
We assume that 
$\cD(n+1)\sq\cD(n)\yt \cD_2(n+1)\sq\cD_2(n),$ and 
$\cD^2(n+1)\sq\cD^2(n).$ 
If $u\yi v\in 2^m$ (binary sequences of length $m$) have the 
form $u=0^k\we0\we w$ and $v=0^k\we1\we w$ for some $k<m$ and 
$w\in 2^{m-k-1}$ then we call $\ang{u,v}$ a {\it crucial pair\/}. 
It can be proved by induction on $m$ that $2^m$ is 
a connected tree 
(\ie, a connected graph without cycles) of crucial pairs, with 
sequences beginning with $1$ as the endpoints of the graph. 

We define 
a system of sets $X_u$ ($u\in2\lom$) and 
$\rR_{uv}\,\yd \ang{u,v}$ being a crucial pair, so that 
the following requirements are satisfied:

\ben
\tenu{(\roman{enumi})}
\itla{d2i1.}
$X_u\in\dP,$ moreover, $X_\La\sq H$; 

\itla{d2i1,}
$X_u\in \cD(n)$ for any $u\in2^n$; 

\itla{d2i2}
$X_{u\we i}\sq X_u$ for all $u$ and $i$;

\itla{d2i3}
$\rR_{uv}\in\dP_2,$ moreover, $\rR_{uv}\in \cD_2(n)$
for any crucial pair  $\ang{u,v}$ in $2^n$; 

\itla{d2i4}
$\rR_{uv}\sq \rE$ and $X_u \rR_{uv} X_v$ 
for any crucial pair $\ang{u,v}$ in $2^n$; 

\itla{d2i5} 
$\rR_{u\we i\,\yi v\we i}\sq\rR_{uv}$;

\itla{d2i6}
if $u,v\in2^n$ and $u(n-1)\ne v(n-1)$ then 
$X_u\ti X_v\in \cD^2(n)$ and also $X_u\cap X_v=\pu$.
\een
Note that \ref{d2i4} implies that $X_u\rE X_v$ for any crucial 
pair $\ang{u,v},$ hence, also for any pair in 
$2^n$ because any $u,v\in2^n$ are connected by a 
unique chain of crucial pairs. 
It follows that $X_u\ti X_v\in\dpe$ for any pair of 
$u,v\in 2^n,$ for~any~$n$.

Assume that such a system has been defined.  
Then for any $a\in\dn$ the sequence 
$\sis{X_{a\res n}}{n\in\dN}$ is \dd\dP generic over $\mm$
by \ref{d2i1,},
therefore $\bigcap_n X_{a\res n}=\ans{x_a},$
where $x_a$ is \dd\dP generic,  
and the map $a\mapsto x_a$ is continuous since diameters 
of $X_u$ converge to $0$ uniformly with $\lh u\to 0$ by 
\ref{d2i1.}, and is $1-1$ by the last condition~of~\ref{d2i6}. 

Let $a,b\in\dn.$ 
If ${a\neo b}$ 
then, by \ref{d2i6}, $\ang{x_a,x_b}$ is a \dd\dpe generic 
pair, hence, $x_a\nE x_b$ by Lemma~\ref{d23}.
Now suppose that $a\Eo b,$ prove that then $x_a\rE x_b.$
We can suppose that $a=w\we 0\we c$ and $b=w\we0\we c,$ 
where $w\in2\lom$ and $c\in\dn$ 
(indeed if $a\Eo b$ then $a\yi b$ can be connected by a finite 
chain of such special pairs). 
Then $\ang{x_a,x_b}$ is \dd{\dP_2}generic, actually, the only 
member of the intersection 
$\bigcap_n \rR_{w\we0\we(c\res n)\,\yi w\we1\we(c\res n)}$ 
by \ref{d2i3} and \ref{d2i4}, in particular, $x_a\rE x_b$ 
because we have $R_{uv}\sq{\rE}$ for all $u\yi v$.

Thus we have 
a continuous $1-1$ reduction of $\Eo$ to $\rE.$ 
\vom

\qeDD{Case 2 in Theorem~\ref{2dih} modulo the construction}

\punk{Construction of a splitting system}
\las{d2splitk}

Thus it remains to define a splitting system of sets
satisfying \ref{d2i1.} -- \ref{d2i6}.

Let $X_\La$ be any set in $\cD(0)$ such that $X_\La\sq H.$ 

Now suppose that $X_s$ and $\rR_{st}$ have been defined for all 
$s\in2^n$ and all crucial pairs in $2^n,$ and extend the 
construction on $2^{n+1}.$ 
Temporarily, define $X_{s\we i}=X_s$ and 
$\rR_{s\we i\,\yi t\we i}=\rR_{st}:$ this leaves 
$\rR_{0^n\we0\,\yi 0^n\we 1}$ still undefined, so we put 
$\rR_{0^n\we0\,\yi 0^n\we 1}={\rE}\cap {(X_{0^n}\ti X_{0^n})}.$ 
Note that the system of sets $X_u$ and 
relations $\rR_{uv}$ defined this way at level $n+1$ 
satisfies all requirements of \ref{d2i1.} -- \ref{d2i6} 
except for the requirements of membership in the dense sets
in \ref{d2i1,}, \ref{d2i3}, \ref{d2i6}
--- say in this case that the system is ``coherent''. 
It remains to produce a still ``coherent'' system of smaller 
sets and relations which also satisfies the 
membership in the dense sets. 
This will be achieved in several steps.\vom

{\sl Step 1\/}: achieve that $X_u\in\cD(n+1)$ for any 
$u\in 2^{n+1}.$ 
Take any particular $u_0\in 2^{n+1}.$ 
There is, by the density, 
$X'_{u_0}\in\cD(n+1)$ and $\sq X_{u_0}.$ 
Suppose that $\ang{u_0,v}$ is a crucial pair. 
Put
$\rR'_{u_0,v}=\ens{\ang{x,y}\in\rR_{u_0,v}}{x\in X'_{u_0}}$ 
and $X'_v=\ran \rR'_{u_0,v}.$ 
This shows how the change spreads along the whole set 
$2^{n+1}$ viewed as the tree of crucial pairs. 
Finally we obtain a coherent system with the additional 
requirement that $X'_{u_0}\in\cD(n+1).$ 
Do this consecutively for all $u_0\in 2^{n+1}.$ 
The total result -- we re-denote it as still $X_u$ and $\rR_{uv}$ 
-- is a ``coherent'' system with $X_u\in\cD(n+1)$ for all $u.$ 
Note that still $X_{0^n\we 0}=X_{0^n\we 1}$ and 
\dm
\rR_{0^n\we0\,\yi 0^n\we 1}=
{\rE}\cap (X_{0^n\we 0}\ti X_{0^n\we 1})\,.
\eqno(\ast)
\dm

{\sl Step 2\/}: achieve that 
$X_{s\we0}\ti X_{t\we1}\in\cD^2(n+1)$ for all $s\yi t\in 2^{n+1}.$ 
Consider a pair of $u_0=s_0\we0$ and $v_0=t_0\we 1$ in 
$2^{n+1}.$ 
By the density there is a set 
$X'_{u_0}\ti X'_{v_0}\in\cD^2(n+1)$ and $\sq X_{u_0}\ti X_{v_0}.$ 
By definition we have $X'_{u_0}\rE X'_{v_0},$ but, due to 
Lemma~\ref{d27} we can maintain that $X'_{u_0}\cap X'_{v_0}=\pu.$
The two ``shockwaves'', from the changes at nodes $u_0$ and $v_0,$
as 
in Step 1, meet only at the pair $0^m\we0\yi 0^m\we1,$ where the 
new sets satisfy $X'_{0^m\we0}\rE X'_{0^m\we1}$ just because 
\dde equivalence is everywhere preserved though the 
changes. 
Now, in view of $(\ast),$ we can define 
$\rR'_{0^n\we0\,\yi 0^n\we 1}=
{\rE}\cap (X'_{0^n\we 0}\ti X'_{0^n\we 1}),$ 
preserving $(\ast)$ as well. 
When all pairs are considered, we will be left with a coherent
system 
of sets and relations, re-denoted as $X_u$ and $\rR_{uv},$ 
which satisfies the \dd{\cD(n+1)}requirements in \ref{d2i1,} 
and \ref{d2i6}.\vom

{\sl Step 3\/}: achieve that $\rR_{uv}\in\cD_2(n+1)$ for any 
crucial pair at level $n+1,$ and also that 
$X'_{0^n\we 0}\cap X'_{0^n\we 1}=\pu.$ 
Consider any crucial pair $\ang{u_0,v_0}.$ 
If this is not $\ang{0^n\we 0,0^n\we 1p}$ then let 
$\rR'_{u_0v_0}\sq\rR_{u_0v_0}$ be any set in $\cD_2(n+1).$ 
If this is $u_0=0^n\we 0$ and $v_0=0^n\we 1$ then first we 
choose (Lemma~\ref{d27}) disjoint non-empty $\is11$ sets 
$U\sq X_{0^n\we 0}$ and $V\sq X_{0^n\we 1}$ still with 
$U\rE V,$ and only then a set 
$\rR'_{u_0v_0}\sq\rE\cap(U\ti V)$ which belongs to 
$\in\cD_2(n+1).$ 
In both cases, put $X'_{u_0}=\dom \rR'_{u_0v_0}$ and 
$X'_{v_0}=\ran \rR'_{u_0v_0}.$ 
It remains to spread the changes, along the chain of crucial 
pairs, to the left of $u_0$ and  to the right of $v_0,$ exactly 
as in Case 1. 
Executing such a reduction for all crucial pairs $\ang{u_0,v_0}$
at level $n+1$ one by one, we end up with a system of sets 
fully satisfying \ref{d2i1.} -- \ref{d2i6}.\vtm

\qeDD{Theorem~\ref{2dih}}

\punk{A forcing notion associated with $\Eo$}
\las{zeo}

We here consider a forcing notion $\peo$ 
\index{forcing!peo@$\peo$}%
\index{zzpeo@$\peo$}%
that consists of all Borel sets
$X\sq\dn$ such that $\Eo\res X$ is non-smooth.
A related ideal $\ieo$ (this time an ideal on $\dn$)
\index{ideal!ieo@$\ieo$}%
\index{zzieo@$\ieo$}%
consists of all Borel sets $X\sq\dn$ such that $\Eo\res X$
is smooth.
Alternatively for a Borel $X\sq\dn$ to be in $\ieo$ it is
necessary and sufficient that ${\Eo}\res X$ has a Borel
transversal --- this is by Lemma~\ref{transv}.

Forcings like $\peo,$ that is those defined in the form of
a collection of all Borel sets $X$ such that a given Borel 
\eqr\ $\rE$ satisfies $\rE\reb{\rE\res X},$
are still work in progress, their applications not
yet established.

\ble
\label{bay}
\ben
\tenu{{\rm(\roman{enumi})}}
\itla{bay1}
$\ieo$ is a\/ \dd\fsg additive ideal.
\imar{bay}
Let\/  $X\sq\dn$ be a Borel set. 

\itla{bay2}
$X$ belongs to\/ $\peo$ iff\/
${\Eo}\emn{\Eo\res X}$ (by a continuous injection).

\itla{bay3}
$X$ belongs to\/ $\ieo$ iff\/
${{\Eo}\res X}$ admits a Borel transversal. 
\een
\ele
\bpf
\ref{bay1} immediately follows from Corollary~\ref{cudc}.
In \ref{bay2},
if $X\in\peo$ then ${\Eo}\emn{\Eo\res X}$ by 
Theorem~\ref{2dih}, while if ${\Eo}\emn{\Eo\res X}$ 
then $\Eo\res X$ is not smooth since $\Eo$ itself is 
not smooth by Lemma~\ref{transv}\ref{transv5}.
In \ref{bay3},
if ${{\Eo}\res X}$ admits a Borel transversal then it
is smooth by Lemma~\ref{transv}\ref{transv1} and hence
$X$ belongs to $\ieo.$
To prove the converse apply Lemma~\ref{transv}\ref{transv3}.
\epf

Note that any $X\in\peo$ contains a closed subset $Y\sq X$
also in $\peo$ by Theorem~\ref{2dih}.
(Apply the theorem for $\rE=\Eo\res X.$
As $\Eo\res X$ is not smooth, we have ${\Eo}\emn{\Eo\res X},$
by  a continuous reduction $\vt.$
Take as $Y$ the full image of $\vt.$
$Y$ is compact, hence closed.)
Such sets $Y$ can be chosen in a special family.

\bdf
\lam{df:zap}
Suppose that two binary sequences $u^0_n\ne u^1_n\in \dln$
of equal length $\lh{u^0_n}=\lh{u^1_n}\ge1$
are chosen for each $n,$ together with one more
sequence $u_0\in\dln.$
Define
$\vt(a)=
u_0\we u^{a(0)}_0\we u^{a(1)}_1\we\dots$
for any $a\in\dn.$
Easily $\vt$ is a continuous injection $\dn\to\dn,$
$Y=\ran \vt$ is a closed set in $\dn,$ $\vt$ witnesses
$\Eo\emn{\Eo\res Y},$ and hence $Y\in\peo.$

Let $\pep$ denote the collection of all sets $Y$
definable in such a form.
\edf


The next theorem gives a necessary and sufficient condition for
a Borel set $X\sq\dn$ in the class $\id11$ to belong to $\peo,$
in terms related to the Gandy -- Harrington forcing.
Relativization to $\id11(p)$ for an arbitrary parameter $p\in\dn$
is obvious.
The theorem also proves the density of the subset $\pep$ of
much more transparent \lap{conditions} in $\peo$. 

\bte
\lam{gho:t}
Suppose that\/ $X\sq\dn$ is a\/ $\id11$ set.
Then\/ $X\in\peo$ iff\/ $X$ is not covered by the union of all
pairwise\/ \dd\Eo inequivalent\/ $\id11$ sets.
In addition,
\ben
\renu
\itla{ghot1}
{\rm({{\rm Zapletal~\cite{zap}}})}
$\pep$ is a dense subset of\/ $\peo:$
\imar{ghot1}
for any\/ $X\in\peo$ there exists\/ $Y\in\pep$
such that\/ $Y\sq X\,;$

\itla{ghot2}
{\rm({{\rm Zapletal~\cite{zap}}})}
$\peo$ forces that the \lap{old} continuum\/
\imar{ghot2}
$\gc$ remains uncountable.
\een
\ete
\bpf
{\it The \lap{if} claim.}
This is easy.
It is quite clear that ${\Eo}\res Y$ is smooth whenever
$Y$ is a Borel pairwise \dd\Eo inequivalent $\id11$ set.
However countable unions preserve smoothness by
Corollary~\ref{cudc}.

{\it The \lap{only if} claim.}
Suppose that
\vyk{
${\Eo}\res X$ is smooth.
It follows from Lemma~\ref{transv}\ref{transv3} that then
${\Eo}\res X$ has a Borel transversal, and hence $X$ is a
union of countably many pairwise \dd\Eo inequivalent
Borel sets.
Yet we need $\id11$ sets here.
Suppose towards the contrary that
}%
$X$ is not covered by
the union $U$ of all
pairwise \dd\Eo inequivalent $\id11$ sets.
As in the proofs of the 1st and 2nd dichotomy theorems above,
$U$ is a $\ip11$ set, and hence $A=X\dif U$ is a non-empty
$\is11$ set.

The key property of $A$ is that it does not intersect any
pairwise \dd\Eo inequivalent $\is11$ set.
(To prove this one has to establish that any
pairwise \dd\Eo inequivalent $\is11$ set can be covered 
by a pairwise \dd\Eo inequivalent $\id11$ set.)
It follows that
\ben
\tenu{{\mtho(\fnsymbol{enumi})}}
\itla{ift}
any non-empty $\is11$ set $Y\sq A$ is not pairwise
\dd\Eo inequivalent, \ie\ it contains
a pair of points $x\ne y$ with $x\Eo y$.
\een

For any sequences $r,w\in\dln$ with $\lh r\le\lh w,$
define $r\ap w\in\dln$ (the \dd rshift of $w$)
so that $\lh{(r\ap w)}=\lh w$ and 
$(r\ap w)(k)=1-w(k)$ whenever $k<\lh r$ and $r(k)=1,$ and
$(r\ap w)(k)=w(k)$ otherwise.
Clearly $r\ap{(r\ap w)}=w.$
Similarly define $r\ap a\in\dn$ for $a\in\dn,$ and
$r\ap X=\ens{r\ap a}{a\in X}$ for any set $X\sq\dn.$

We are going to define sequences $u\in\dln$ and
$u^0_n\ne u^1_n\in\dln\;(n\in\dN)$ such that
$\lh{u^0_n}=\lh{u^1_n},$ as in
Definition~\ref{df:zap}, and also a system of $\is11$ sets
$X_s\in\peo\;\,(s\in\dln)$ satisfying the following:
\ben
\tenu{(\arabic{enumi})}
\itla{eo1}\msur
$X_\La\sq A,$ $X_{s\we i}\sq X_s,$ and $\dia{X_s}\le 2^{-\lh s}$;

\itla{eo1+} a condition in terms of the
Gandy -- Harrington forcing,
similar to \ref{d1i,} in \nrf{>smooth} or
\ref{d2i1,} in \nrf{d2split}, such that, as a 
consequence, $\bigcap_nX_{a\res n}\ne\pu$ 
for any $a\in\dn$;

\itla{eo2}\msur
$X_s\sq\cO_{w_s},$ where
$w_s=u_0\we u^{s(0)}_0\we u^{s(1)}_1\we\dots
\we u^{s(k-1)}_{k-1}\in\dln,$
$k=\lh s,$ and $\cO_w=\ens{a\in\dn}{w\su a}$ for $w\in\dln;$

\itla{eo3} 
if $s,t\in 2^n$ for some $n$ then $X_t=w_t\ap w_s\ap X_s$.
\een
Then define the map $\vt$ as in Definition~\ref{df:zap}.
The set $Y=\ran\vt=\bigcap_n\bigcup_{s\in2^n}X_s\sq X$
belongs to $\pep,$ hence to $\peo,$ proving that $X\in\peo$
as well.

This argument also proves claim \ref{ghot1} of the theorem.
Indeed suppose that $X\sq\dn$ is a $\id11$ set.
(As usual the relativization to any $\id11(p)$ is routine.)
It follows from the \lap{if} claim of the theorem that 
$X\not\sq U,$ and hence we are in the domain of the
\lap{only if} claim, thus there is a subset $Y\sq X\yt Y\in\pep.$ 

It remains to carry out the construction of sets $X_s$.\vom

{\it Step 0\/}.
We put $X_\La=A$ and let $u_0\in\dln$
be the largest sequence such that $X_\La\sq \cO_{u_0}.$
Let $\ell_0=\lh{u_0}$.\vom

{\it Step 1\/}.
Here we define $u^i_0$ and $X_{\ang i}$ for $i=0,1.$
It follows from \ref{ift} above that there exist points
$x'\ne y'\in X_\La$ such that $x'\Eo y'.$
This means that there exist two different sequences
$u^0_0\ne u^1_0$ of equal length $\lh u^0_0=\lh u^1_0$
such that
$u_0\we u^0_0\su x'\yt u_0\we u^1_0\su y',$
and $x'(k)=y'(k)$ for all $k\ge\ell_1=\ell_0+\lh u^i_0.$
Put $w_{\ang0}=u_0\we u^0_0$ and $w_{\ang1}=u_0\we u^0_1.$
Then the sets
\dm
\bay{rcl}
X_{\ang0} &=&
\ens{x\in X_\La}{w_{\ang0}\su x\land
\sus y\in X_\La\:(w_{\ang1}\su y\land x\Eo y)},
\quad\text{and}\\[0.8\dxii]

X_{\ang1} &=&
\ens{y\in X_\La}{w_{\ang1}\su y\land
\sus x\in X_\La\:(w_{\ang0}\su x\land x\Eo y)}
\eay
\dm
are still nonempty $\is11$ sets
(containing resp.\ $x\yi y$), and they satisfy \ref{eo2}
and \ref{eo3}.

Finally replace $X_{\ang0}$ by a suitable smaller $\is11$
set $X'_{\ang0}$ in order to fulfill \ref{eo1+}, and put
$X'_{\ang1}=w_{\ang0}\ap w_{\ang1}\ap X'_{\ang0}.$
Now choose suitable smaller $\is11$
set $X''_{\ang1}\sq X'_{\ang1}$ in order to fulfill \ref{eo1+},
and put $X''_{\ang0}=w_{\ang1}\ap w_{\ang0}\ap X''_{\ang1}.$
Re-denote the sets $X''_{\ang0}\yd X''_{\ang1}$ again by
$X_{\ang0}\yd X_{\ang1}$.
\vom

{\it Step 2\/}.
Here we define $u^i_1$ for $i=0,1$ and $X_s$ for $s\in\dln$
with $\lh s=2.$
Once again there exist points
$x'\ne y'\in X_{\ang0}$ such that $x'\Eo y'.$
This means that there exist two different sequences
$u^0_1\ne u^1_1$ of equal length $\lh u^0_1=\lh u^1_1$
such that
$u_0\we u^0_0\we u^0_1\su x'\yt u_0\we u^1_0\we u^1_1\su y',$
and $x'(k)=y'(k)$ for all $k\ge\ell_2=\ell_1+\lh u^i_1.$
Put $w_{\ang{i,j}}=u_0\we u^i_0\we u^j_1$ for 
$i,j\in\ans{0,1}.$
Then the sets
\dm
\bay{l}
X_{\ang{0,0}} =
\ens{x\in X_\La}{w_{\ang{0,0}}\su x\land
\sus y\in X_\La\:(w_{\ang{0,1}}\su y\land x\Eo y)},
\quad\text{and}\\[0.8\dxii]

X_{\ang{0,1}} =
\ens{y\in X_\La}{w_{\ang{0,1}}\su y\land
\sus x\in X_\La\:(w_{\ang{0,0}}\su x\land x\Eo y)}
\eay
\dm
are still nonempty $\is11$ sets
satisfying \ref{eo2} and \ref{eo3}.
There is no need in an additional split of $X_{\ang1}$
in order to define the sets $X_{\ang{1,0}}\zd X_{\ang{1,1}}:$
just put $X_{\ang{1,0}}=w_{\ang0}\ap w_{\ang1}\ap X_{\ang{0,0}}$ 
and $X_{\ang{1,1}}=w_{\ang0}\ap w_{\ang1}\ap X_{\ang{0,1}}.$

It remains to shrink the sets $X_{\ang{i,j}}$ in several
(that is, four) rounds in order to fulfill \ref{eo1+},
applying the actions of $w_{\ang{i,j}}$ as required by
\ref{eo3} to define intermediate sets.\vom

{\it Steps $\ge3$\/}.
Suppose that all sets $X_s\yt s\in 2^n,$ have been
suitably defined.
Let $\sq\in2^n$ be the sequence of $n$ zeros.
We define sets $X_{\sg\we 0}$ and $X_{\sg\we1}$ by
splitting $X_\sg$ as above, and then split every other
$X_s$ applying $w_\sg\ap w_s$.\vom

The construction results in a system of sets and sequences
satisfying requirements \ref{eo1}, \ref{eo2}, \ref{eo3},
as required.\vom

\ref{ghot2}
It suffices to prove the same result for the subforcing $\pep.$ 
Given a sequence of dense sets $D_n\sq\pep,$ we carry out a
splitting construction similar to the one given above,
with the following amendments.
First, each set $X_s$ belongs to $D_{\lh s},$ hence to $\pep,$
therefore is a closed set in $\dn.$ 
Second, condition \ref{eo1+} is abolished, of course.
That any set $X\in\pep$ satisfies \ref{ift}
(that is, it contains a pair of points $x\ne y$ with $x\Eo y$)
is obvious.
\epf

We observe that $\peo$ as a forcing is somewhat closer to
Silver rather than Sacks forcing.
The property of minimality of the generic real, common to
both Sacks and Silver forcings, holds for $\peo$ as well,
the proof
resembles known arguments, but in addition the following
is applied: if $X\in\peo$ and $f:X\to\dn$ is a Borel
\dd\Eo invariant map (that is, ${x\Eo y}\imp{f(x)=f(y)}$)
then $f$ is constant on a set $Y\in\peo\yt Y\sq X.$\snos
{Suppose, for the sake of brevity, that $X=\dn.$
For any $n,$ the set $Y^0_n=\ens{a}{f(a)(n)=0}$ is Borel and
\dd\Eo invariant.
It follows that $Y^0_n$ is either meager or comeager.
Put $b(n)=0$ iff $Y^0_n$ is comeager.
Then $D=\ens{a}{f(a)=b}$ is comeager.
A splitting construction as in the proof of Theorem~\ref{gho:t}
yields a set $Y\in\peo\yt Y\sq D$.}

\api

\parf{Ideal $\Ii$ and P-ideals}
\las{idI1}

By definition the ideal $\fio=\Ii$ consists of all sets 
$x\sq\pnn$ such that all, except for finitely many, 
cross-sections 
$\seq xn=\ens{k}{\ang{n,k}\in x}$ are empty. 
This \gla\ contains proofs of some key results
related to this ideal.
First of all we show following Kechris that there exist
essentially only three types of ideals Borel reducible
to $\Ii,$ two of them being $\ifi$ and $\Ii$ itself.
Then a proof of Solecki's theorem, that characterizes
P-ideals in terms of \lsc\ submeasures and polishability
and shows that $\Ii$
is the least Borel non-polishable ideal, will be given.

\punk{Ideals below $\Ii$}
\las{<i1}

\vyk{
Suppose that $\cI,\cJ$ are ideals on countable sets resp.\
$A\yi B.$ 
Then $\cI$ is called a \rit{trivial variation\/} of $\cJ$
\index{trivial variation}%
iff there is an infinite set $D\sq B$ such that 
$\cI\res D\isi\cJ$~\snos
{Recall that $\cI\isi\cJ$ means the isomorphism of ideals
$\cI,\cJ$ via a bijection between the underlying sets.
$\cI\res D=\ens{x\cap D}{x\in \cI},$ and this is equal
to $\cI\cap\pws D.$}
while $\cI\res{(A\bez D)}=\pws{A\bez D}.$
In other words it is required that $\cI$ is isomorphic
to the disjoint sum $\fsm{\cJ}{\pws{A\bez D}}.$
If the difference $A\bez D$ is infinite then let us call
$\cI$ a \rit{strict trivial variation\/} of $\cJ;$
in this case simply $\cI\cong {\dsu{\cJ}{\pws{\dN}}}.$
If the difference $A\bez D$ is finite then
$\dsu{\cJ}{\pws{A\bez D}}$ is isomorphic to $\cJ$ itself for
typical ideals $\cJ,$ in particular, for the ideals 
$\cJ=\ifi\yi \Ii\yi\It$.

Note that a strict trivial variation of $\Ii$ or $\It$ is
isomorphic to resp.\ $\Ii\yi\It,$ while a strict trivial
variation of $\ifi$ is clearly not isomorphic to $\ifi$.
}

Recall that $\cI\isi\cJ$ means the isomorphism of ideals
$\cI,\cJ$ via a bijection between the underlying sets.
The ideal $\fsm\ifi\pn$
\index{ideal!fin+pn@$\fsm\ifi\pn$}%
\index{zzfin+pn@$\fsm\ifi\pn$}%
(the disjoint sum in the sense of \nrf{op:id})
in the next theorem is isomorphic to the ideal
$\ifi_{\textsc{odd}}=\ens{x\sq\dN}{x\cap{2\dN}\in\ifi},$
where $2\dN=$ all odd numbers.\snos
{\label{trvar}%
Ideals isomorphic to any of $\cI\yi\fsm\cI\pn$ were called
\index{trivial variation}%
\index{ideal!trivial variation}%
\rit{trivial variations} of $\cI$ in \cite{rig}.}

\bte
[{{\rm Kechris~\cite{rig}}}]
\lam{rig1}
If\/ $\cI$ is a Borel (nontrivial) ideal on\/ $\dN$ and\/
$\cI\reb\Ii$ then\/ $\cI$ is isomorphic to one of the
following three ideals$:$
$\Ii\yt\ifi\yt\fsm\ifi\pn$.
\ete

Thus there exist only three different ideals 
Borel reducible to $\Ii,$ they are $\ifi,$ the disjoint sum 
$\ifi\oplus\pn,$ and $\Ii$ itself.

\vyk{
\bupt
Prove that any trivial variation of $\cI_1$ is isomorphic 
to $\cI_1$ 
while any trivial variation of $\ifi$ is isomorphic 
either to $\ifi$ or to the disjoint sum 
$\ifi\oplus\pn,$ \eg, realized in the form of 
$\ens{x\sq\dN}{x\cap{\textsc{odd}}\in\ifi}$. 
\eupt
}

\bpf
We begin with another version of the method 
used in the proof of Theorem~\ref{jnmt}. 
Suppose that $\sis{\cB_k}{k\in\dN}$ is a fixed 
system of Borel subsets of $\pn.$ 
(It will be specified later.) 
Then there exists an increasing sequence of integers 
$0=n_0<n_1<n_2<\dots$ and sets $s_k\sq\il k{k+1}$ such 
that

\ben
\tenu{(\arabic{enumi})}
\itla{ri1}
any $x\sq\dN$ with $\kai k\:(x\cap\il k{k+1}=s_k)$ 
is ``generic''\snos
{\label{genmean}  
We mean, Cohen generic over a certain countable
family of dense open subsets of $\pn$ that depends on
the choice of the family of sets $\cB_k$.}~;

\itla{ri2}
if $k'\ge k$ and $u\sq\ilo{k'}$ then $u\cup s_{k'}$ 
decides $\cB_k$ in the sense that either any 
``generic'' $x\in\pn$ with $x\cap\ilo{k'+1}=u\cup s_{k'}$ 
belongs to $\cB_k$ or any ``generic'' $x$ with 
$x\cap\ilo{k'+1}=u\cup s_{k'}$ does not belong to $\cB_k$.
\een

Now put $\cD_0=\ens{x\cup S_1}{x\sq Z_0}$ and  
$\cD_1=\ens{x\cup S_0}{x\sq Z_1},$ where
\dm
\textstyle
S_0=\bigcup_ks_{2k}\sq Z_0=\bigcup_k\il{2k}{2k+1}\,,
\quad 
S_1=\bigcup_ks_{2k+1}\sq Z_1=\bigcup_k\il{2k+1}{2k+2}.
\dm
Clearly any $x\in\cD_0\cup\cD_1$ is ``generic'' by 
\ref{ri1}, hence it follows from \ref{ri2} that 
\ben
\tenu{(\arabic{enumi})}
\addtocounter{enumi}2
\itla{ri3}
each $\cB_k$ is clopen on both $\cD_0$ and $\cD_1$.
\een

As $\cI\reb\Ii,$ 
it follows from Lemma~\ref{l:bc} 
(and the trivial fact that $\cI_1\oplus\cI_1\cong\cI_1$) 
that there exists a {\it continuous\/} reduction 
$\vt:\pn\to\pws{\dN\ti\dN}$ of $\cI$ to $\Ii.$
Thus $\rE_\cI$ is 
the union of an increasing sequence of (topologically) 
closed \eqr s $\rR_m\sq\rE_\cI$ just because $\cI_1$ 
admits such a form. 
We now require that $\sis{\cB_k}{}$ includes  
all sets 
$B^m_l=\ens{x\in\pn}{\kaz s\sq\ir0l\;x\rR_m(x\sd s)}.$ 
Then by \ref{ri3} and the compactness of $\cD_i$ 
for any $l$ there is $m(l)\ge l$ satisfying 
\ben
\tenu{(\arabic{enumi})}
\addtocounter{enumi}3
\itla{ri4}
$\kaz x\in\cD_0\cup\cD_1\:\kaz s\sq\ir0l\;
\skl x\rR_{m(l)}(x\sd s)\skp$.
\een

To prove the theorem it suffices to obtain a sequence 
$x_0\sq x_1\sq x_2\sq\dots$ of sets $x_k\in\cI$ with 
$\cI=\bigcup_n\pws{x_n}:$ 
that in this case $\cI$ is as required is an easy exercise. 
As any topologically closed ideal is easily $\pws x$ for 
some $x\sq\dN,$ it suffices to show that $\cI$ is a union 
of a countable sequence of closed subideals. 
It suffices to demonstrate this fact separately for 
$\cI\res Z_0$ and $\cI\res Z_1.$ 
Prove that 
$\cI\res Z_0$ is a countable union of closed subideals, 
ending the proof of the theorem.

If $m\in \dN$ and $s\sq u\sq Z_0$ are finite then let 
\dm
I^m_{us}=\ens{A\sq Z_0}{\kaz x\in\cD_0\:
(x\cap u=s\imp {\skl x\cup(A\dif u)\skp\rR_m x})}\,.
\dm

\ble
\lam{rig1l}
Sets\/ $I^m_{us}$ are closed topologically and 
under\/ $\cup,$ and\/ $I^m_{us}\sq\cI$.
\ele
\bpf
$I^m_{us}$ are topologically closed because so are $\rR_m$. 

Suppose that $A\yi B\in I^m_{us}.$  
To prove that $A\cup B\in I^m_{us},$ let 
$x\in\cD_0$ satisfy $x\cap u=s.$ 
Then $x'=x\cup(A\dif u)\in\cD_0$ satisfies $x'\cap u=s,$ 
too, hence, as $B\in I^m_{us},$ we have 
$\skl{x'\cup(B\dif u)}\skp\rR_m x',$ thus, 
$\skl{x\cup((A\cup B)\dif u)}\skp\rR_m x'.$ 
However $x'\rR_m x$ just because $A\in I^m_{us}.$ 
It remains to recall that $\rR_m$ is a \er.

To prove that any $A\in I^m_{us}$ belongs to 
$\cI$ take $x=s\cup S_1.$ 
Then we have ${x\cup(A\dif u)}\rR_m x,$ thus, 
$A\in\cI$ as $s$ is finite and $\rR_m\sq\rE_\cI$.
\epF{Lemma}

\ble
\lam{rig1L}
$\cI\res Z_0=\bigcup_{m\yi u\yi s}I^m_{us}$.
\ele
\bpf
Let $A\in\cI,\msur$ $A\sq Z_0.$ 
The sets 
$Q_m=\ens{x\in\cD_0}{\skl x\cup A\skp\rR_m x}$ 
are closed and satisfy $\cD_0=\bigcup_m Q_m.$ 
It follows that one of them has a non-empty interior 
in $\cD_0,$ thus, there exist finite sets 
$s\sq u\sq Z_0$ and some $m_0$ with 
\dm
\kaz x\in\cD_0\:
(x\cap u=s\imp {\skl x\cup A\skp\rR_{m_0} x})\,.
\dm
This is not exactly what we need, however, 
by \ref{ri4}, there exists a number 
$m=\tmax\ans{m_0,m(\tsup u)}$ 
big enough for 
\dm
\kaz x\in\cD_0:
\skl x\cup A\skp\rR_{m} \skl x\cup(A\dif u)\skp\,.
\dm
It follows that $A\in I^m_{su},$ as required.
\epF{Lemma}

Let $J^m_{su}$ be the hereditary hull of $I^m_{su}$ 
(all subsets of sets in $I^m_{su}$). 
It follows from Lemma~\ref{rig1l} that any $J^m_{su}$ 
is a topologically closed subideal of $\cI\res Z_0,$ 
however, $\cI\res Z_0$ is the union of those 
ideals by Lemma~\ref{rig1L}, as required.
\epf

\bcor
\lam{ne<e1}
The \er s\/ $\Ed$ and\/ $\Et$ are Borel irreducible to\/ $\Ei.$
It follows that they are Borel irreducible to\/ $\Eo,$ and hence\/
$\Eo\rebs\Ed$ and\/ $\Eo\rebs\Et$.
\ecor
\bpf
It is quite clear that
neither $\Id$ nor $\It$ belong to the types of ideals mentioned in 
Theorem~\ref{rig1}.
\epf

That $\Eo\rebs\Ei$ strictly, and even that $\Ei$ is not
essentially countable (formally $\Ei\not\reb\Ey$),
will be established by Lemma~\ref{eoei'} below.

\punk{$\Ii$ and P-ideals}
\las{sol: prf}

The next theorem claims that the ideal $\Ii$ is the
\dd\orb least among all Borel 
ideals which are not P-ideals.  
That it is the \dd\reb least in this family will be shown in
the next \gla.

Recall that \rit{analytic} means $\fs11$ while the notions of
polishable ideals and P-ideals were introduced in \grf{n-id}.

\bte
\lam{sol}
The following conditions are equivalent for any ideal
on\/ $\dN:$
\ben
\tenu{{\rm(\roman{enumi})}}
\itla{s2}\msur
$\cI$ has the form\/ $\Exh_\vpi,$ where\/ $\vpi$ is a 
\lsc\ submeasure on\/ $\dN\;;$

\itla{s3}\msur
$\cI$ is a polishable ideal\/$;$ 

\itla{s1}\msur
$\cI$ is an analytic P-ideal$;$

\itla{s5}\msur
$\cI$ is an analytic ideal such that all countable unions of\/ 
\ddi small sets are \ddi small,
{\rm where a set $X\sq\pn$ is \ddi{\it small\/} if there
\index{set!Ismall@\ddi small}%
is $A\in\cI$ such that 
$X\res A=\ens{x\cap A}{x\in X}\sq\pws A$ 
is meager in $\pws A$\/}$;$

\itla{s4}\msur
$\cI$ is an analytic ideal satisfying\/ $\cI_1\not\orb \cI\;;$ 

\itla{s6}\msur
$\cI$ is an analytic ideal satisfying\/ $\cI_1\not\reb \cI$. 
\een
\ete

By the way it follows that all analytic P-ideals actually
belong to $\fp03,$ simply because any ideal of type \ref{s2} is 
easily $\fp03$.

\bcor
\label{rb=b}
If\/ $\Ii$ is a Borel ideal then\/ $\Ii\orb\cI$ iff\/
$\Ei\reb \rei$.\qed
\ecor

\bcor
\lam{<pideal}
Suppose that\/ $\cJ$ is an analytic P-ideal. 
Then any ideal\/ $\cI$ satisfying\/ $\cI\reb\cJ$
is an analytic P-ideal, too.
\ecor
\bpf
Use equivalence $\ref{s6}\eqv\ref{s1}$ of the theorem. 
\epf

\bpf[Theorem]
We begin with the proof of the equivalence of the first five
conditions, the result of Solecki~\cite{sol,sol'}.
First of all, comparably simple (but tricky in some points)
equivalences $\ref{s2}\eqv\ref{s3}$ and 
$\ref{s1}\eqv\ref{s4}\eqv\ref{s5}$ and implication
$\ref{s2}\imp\ref{s1}$ will be established. 
The hard part will be the implication $\ref{s5}\imp\ref{s2}$
that follows in \nrf{sol-hard}.
The last condition \ref{s6} (Kechris and Louveau \cite{hypsm})
will be addes to the equivalence
by Lemma~\ref{soll} based on several complicated theorems in
the next \gla.\vom 

$\ref{s2}\imp\ref{s3}$ 
If $\vpi(\ans n)>0$ for all $n$ then the required metric 
on $\cI=\Exh_\vpi$ can be defined by 
$\dpi(x,y)=\vpi(x\sd y).$ 
Then any set $U\sq\cI$ open in the sense of the ordinary 
topology (the one inherited from $\pn$) is \dd\dpi open, 
while any \dd\dpi open set is Borel in the ordinary sense.
In the general case we assemble the required metric of 
$\dpi$ on the domain $\ens{n}{\vpi(\ans n)>0}$ and the 
ordinary Polish metric on $\pn$ on the complementary 
domain.\vom

$\ref{s3}\imp\ref{s2}$ 
Let $\tau$ be a Polish group topology on $\cI,$ 
generated by a \dd\sd invariant compatible metric $d.$ 
It can be shown (Solecki~\cite[p.\ 60]{sol'}) that 
$\vpi(x)=\tsup_{y\in\cI,\:y\sq x}\,d(\pu,x)$ is 
a \lsc\ submeasure with $\cI=\Exh_\vpi.$ 
The key observation is that for any $x\in\cI$ the 
sequence $\sis{x\cap\ir0n}{n\in\dN}$ \dd dconverges 
to $x$ by the last statement of Lemma~\ref{sol:?}, which 
implies both that $\vpi$ is \lsc\ 
(because the supremum above can be restricted to 
finite sets $y$) 
and that $\cI=\Exh_\vpi$ 
(where the inclusion $\supseteq$ needs another 
``identity map'' argument).\vom

$\ref{s2}\imp\ref{s1}$ 
That any $\cI=\Exh_\vpi,$ $\vpi$ being \lsc, is 
a P-ideal, is an easy exercise:  
if $x_1\yi x_2\yi x_3\yi \dots\in\cI$ then define 
an increasing sequence of numbers $n_i\in x_i$ with 
$\vpi(x_i\cap\iry{n_i})\le2^{-n}$ and put 
$x=\bigcup_i(x_i\cap\iry{n_i})$.\vom

$\ref{s1}\limp \ref{s4}$ \ 
This is because $\cI_1$ easily does not satisfy \ref{s1}. \vom 

$\ref{s4}\imp\ref{s5}$ \ 
Suppose that sets $X_n\sq\pn$ are \ddi small, so that 
$X_n\res A_n$ is meager in $\pws{A_n}$ for some $A_n\in\cI,$  
but $X=\bigcup_n X_n$ is not \ddi small, and prove 
$\cI_1\orb \cI.$ 
Arguing as in the proof of Theorem~\ref{jnmt}, we use 
the meagerness to find, for any $n,$ a sequence of 
pairwise disjoint non-empty 
finite sets $w^n_k\sq A_n,\msur$ $k\in\dN,$ and 
subsets $u^n_k\sq w^n_k,$ such that 
\ben
\tenu{(\alph{enumi})}
\itla{sa1}
if $x\sq \dN$ and $\exi k\:(x\cap w^n_k=u^n_k)$ 
then $x\nin X_n$. 
\een
Dropping some sets $w^n_k$ away and reenumerating the 
rest, we can strengthen the disjointness to the following: 
$w^n_k\cap w^m_l=\pu$ unless both $n=m$ and $k=l$. 

Now put $w^n_{ij}=w^n_{2^i(2j+1)-1}.$  
The sets $\ovw_{ij}=\bigcup_{n\le i} w^n_{ij}$ 
are still pairwise disjoint, and satisfy the following 
two properties: 

\ben
\tenu{(\alph{enumi})}
\addtocounter{enumi}1
\itla{sa2}\msur
$\bigcup_j\ovw_{ij}\sq A_0\cup\dots\cup A_i,$ hence, $\in \cI,$ 
for any $i$;

\itla{sa3}
if a set $Z\sq\dN\ti\dN$ does not belong to $\cI_1,$ 
\ie, $\exi i\:\sus j\:(\ang{i,j}\in Z),$ then 
$\kaz n\:\exi k\:(w^n_k\sq\ovw_Z),$ 
where $\ovw_Z=\bigcup_{\ang{i,j}\in Z}\ovw_{ij})$.
\een
We assert that the map $\ang{i,j}\mapsto\ovw_{ij}$ 
witnesses $\cI_1\orbp\cI.$ 
(Then a simple argument, as in the proof of 
Theorem~\ref{jnmt}, gives $\cI_1\orb\cI.$) 

Indeed if $Z\sq\dN\ti\dN$ belongs to $\cI_1$ then 
$\ovw_Z\in\cI$ by \ref{sa2}. 
Suppose that $Z\nin\cI_1.$ 
It suffices to show that $X_n\res \ovw_Z$ is meager in 
$\pws{\ovw_Z}$ for any $n.$ 
Note that by \ref{sa3} the set 
$K=\ens{k}{w^n_k\sq\ovw_Z}$ is infinite and in fact 
$\ovw_Z\cap A_n=\bigcup_{k\in K}w^n_k.$ 
Therefore, any $x\sq \ovw_Z$ satisfying 
$x\cap w^n_k=u^n_k$ for infinitely many $k\in K,$ 
does not belong to $X_n$ by \ref{sa1}. 
Now the meagerness of $X_n\res \ovw_Z$ is clear.\vom

$\ref{s5}\imp\ref{s1}$ \ 
This also is quite easy: if a sequence of sets $Z_n\in\cI$ 
witnesses that $\cI$ is not a P-ideal, then the union of 
\ddi small sets $\cP(Z_n)$ is not \ddi small.

\punk{The hard part}
\las{sol-hard}

We finally prove $\ref{s5}\imp\ref{s2},$ the hard part of 
Theorem~\ref{sol}. 
A couple of definitions precede the key lemma.  
\bit
\item
Let $C(\cI)$ be the collection of all hereditary 
(\ie, ${y\sq x\in K}\imp {y\in K}$) compact 
\ddi large sets $K\sq\pn.$ \
(By definition a set $K\sq\pn$ is \ddi\rit{large}
\index{set!Ilarge@\ddi large}%
iff it is not \ddi small in the sense of \ref{s5} of
Theorem~\ref{sol}.)
\eit
Note that if $K\sq\pn$ is hereditary and compact then for
$K\in C(\cI)$ it is necessary and sufficient that for
any $A\in\cI$ there is $n$ such that $A\cap{\ir n\iy}\in K$.

\bit
\item
Given sets $X\yi Y\sq\pn,$ let 
$X+Y=\ens{x\cup y}{x\in X\land y\in Y}$.
\eit

\ble
\lam{sl2}
Assume that\/ $\cI$ is of type \ref{s5} of
Theorem~\ref{sol}.
Then there is a countable sequence of sets\/
$K_m\in C(\cI)$ such that for any set\/ $K\in C(\cI)$
there exist numbers\/ $m\yi n$ with\/ $K_m+K_n\sq K$. 
\ele
\bpf
As $\cI$ is a $\fs11$ subset of $\pn,$ there exists 
a continuous map $f:\bn\onto\cI.$ 
For any $s\in\nse,$ we define 
\dm
N_s=\ens{a\in\bn}{s\su a}\quad\hbox{and}\quad
B_s=f\ima N_s \quad\text{(the \dd fimage of $N_s$)}\,. 
\dm
Consider the set 
$T=\ens{s}{B_s\,\hbox{ is \ddi large}}.$ 
As $\cI$ itself is clearly \ddi large, $\La\in T.$ 
On the other hand, the assumption \ref{s5} easily implies 
that $T$ has no endpoints and no isolated branches, hence, 
$P=\ens{a\in\bn}{\kaz n\:(a\res n\in T)}$ is a perfect 
set. 
Moreover, $F_s=f\ima(P\cap N_s)$ is \ddi large for 
any $s\in T$ because $B_s\dif F_s$ is a countable 
union of \ddi small sets.

Now consider any set $K\in C(\cI).$ 
By definition, if $x\yi y\in\cI$ then $z=x\cup y\in\cI,$ 
thus, $K\res z$ is not meager in $\pws z,$ 
hence, by the compactness, 
$K\res z$ includes a basic nbhd of $\pws z,$ hence, 
by the hereditarity, there is a number $n$ such that 
$Z\cap\iry n\in K.$ 
We conclude that $P^2=\bigcup_n Q_n,$ where each 
$Q_n=\ens{\ang{a,b}\in P^2}
{(f(a)\cup f(b))\cap\iry n\in K}$ 
is closed in $P$ because so is $K$ and $f$ is continuous. 
Thus, there are $s\yi t\in T$ such that 
$P^2\cap (N_s\ti N_t)\sq Q_n,$ in other words, 
$(F_s+F_t)\res\iry n\sq K,$ hence, 
$(\nadd{F_s}+\nadd{F_t})\res\iry n\sq K,$ where 
$\nadd{\vphantom|\dots}$ denotes 
the topological closure of the hereditary hull.
Thus we can take, as $\sis{K_m}{},$ all sets of the 
form $K_{sn}=\nadd{F_s}\res n$. 
\epf

As $C(\cI)$ is obviously a filter, we can transform  
(still in the assumption that $\cI$ is of type \ref{s5})
the sequence of sets given by the lemma into 
a \dd\sq decreasing sequence of sets $K_n\in C(\cI)$ 
such that
\ben
\tenu{(\arabic{enumi})}
\itla{ck1} 
for any $K\in C(\cI)$ there is $n$ with $K_n\sq K$,
\een
and $K_{n+1}+K_{n+1}\sq K_n$ for any $n.$ 
Taking any other term of the sequence, we can sharpen 
the latter requirement to
\ben
\tenu{(\arabic{enumi})}
\addtocounter{enumi}1
\itla{ck2}
for any $n:$ $K_{n+1}+K_{n+1}+K_{n+1}\sq K_n$.
\een

This is the starting point for the construction of a 
\lsc\ submeasure $\vpi$ with $\cI=\Exh_\vpi.$ 
Assuming that, in addition, $K_0=\pn,$ let, for any 
$x\in\pwf\dN$, 
\dm
\bay{lcll}
\vpi_1(x) &=& \tinf\left\{\,2^{-n}:x\in K_n\,\right\}&,
\hbox{ and}\\[1ex]

\vpi_2(x) &=& \tinf\left\{\,\textstyle
\sum_{i=1}^{m}\vpi_1(x_i):m\ge 1\land 
x_i\in\pwf\dN\land x\sq\bigcup_{i=1}^{m}x_i\,
\right\}&.
\eay
\dm
Then set $\vpi(x)=\tsup_n\vpi_2(x\cap\ir0n)$ for 
any $x\sq\dN.$ 
A routine verification shows that $\vpi$ submeasure 
and that $\cI=\Exh_\vpi.$ 
(See Solecki~\cite{sol'}. 
To check that any $x\in\Exh_\vpi$ belongs to $\cI$ 
we use the following observation: $x\in\cI$ iff for 
any $K\in C(\cI)$ there is $n$ such that  
$x\cap\ir n\iy\in K$.)

\epF{Theorem~\ref{sol} without \ref{s6}}

\api

\parf{Equivalence relation $\Ei$}
\las{Ei:er}

The ideal $\Ii$ naturally defines the \er\ $\Ei=\rE_{\Ii}$ 
on $\pnn$ so that  
$x\Ei y$ iff $x\sd y\in\Ii.$ 
\index{equivalence relation, ER!E1@$\Ei$}%
\index{zzE1@$\Ei$}%
We can as well consider $\Ei$ as an \eqr\ on $\dnd,$ or even
on $\dX^\dN$ for any uncountable Polish space $\dX,$ 
defined as $x\Ei y$ iff $x(k)=y(k)$ for all but finite $k.$ 

The following notation will be quite useful in our study 
of subsets of spaces of the form $X^\dN.$
If $x$ is a function defined on $\dN$ then, for any $n,$ let 
\dm
x\rmq n= x\res{[0,n)}\,,\;\; 
x\rme n= x\res{[0,n]}\,,\;\;
x\qc n = x\res{(n,\iy)}\,,\;\; 
x\qec n= x\res{[n,\iy)}\,.
\dm 
\index{zzx<n@$x\rmq n$}%
\index{zzx<=n@$x\rme n$}%
\index{zzx>n@$x\qc n$}%
\index{zzx>=n@$x\qec n$}%
For any set $X$ of \dd\dN sequences, let 
$X\rmq n=\ens{x\rmq n}{x\in X},$ and similarly for 
${\le}\zd{>}\zd{\ge}.$
If $\xi\in X\qc n$ then let 
$\srez X\xi=\ens{x(n)}{x\in X\land x\qc n=\xi}.$
\index{zzSxxi@$\srez X\xi$}%

\punk{$\Ei:$ hypersmoothness and non-countability}
\las{einc}

Recall that a hypersmooth \eqr\ is a countable increasing 
union of Borel smooth \er s. 
This \paf\ contains a several results which describe the
relationships between hypersmooth and countable \eqr s.
The following lemma shows that $\Ei$ is universal in this 
class.

\ble
\lam{hs<e1}
For a Borel \er\/ $\rE$ to be hypersmooth it is necessary 
and sufficient that\/ $\rE\reb\Ei$.
\ele
\bpf
Let $\dX$ be the domain of $\rE.$ 
Assume that $\rE$ is hypersmooth, \ie, $\rE=\bigcup_n\rE_n,$ 
where $x\rE_n y$ iff $\vt_n(x)=\vt_n(y),$ each 
$\vt_n:\dX\to\dn$ is Borel, and $\rE_n\sq\rE_{n+1},\:\kaz n.$ 
Then $\vt(x)=\sis{\vt_n(x)}{n\in\dN}$ witnesses 
$\rE\reb\Ei.$ 
Conversely, if $\vt:\dX\to\dnd$ is a Borel reduction 
of $\rE$ to $\Ei$
then the sequence of \er s $x\rE_n y$ iff 
$\vt(x)\qec n=\vt(y)\qec n$ witnesses that $\rE$ 
is hypersmooth.
\epf


\bcor
\lam{eiey}
$\Ey\not\reb\Ei$.
\ecor
\bpf
Otherwise $\Ey$ is a hypersmooth \eqr\ by Lemma~\ref{hs<e1}.
But $\Ey$ is countable as well. 
It follows that $\Ey\reb\Eo$ by Theorem~\ref{thf}.
This contradicts Theorem~\ref{Eynon}.
\epf

The following result is given in \cite{hypsm} with a 
reference to earlier papers.


\ble
\label{eoei'}
\ben
\tenu{{\rm(\roman{enumi})}}
\itsep
\itla{eoei'1}
$\Ei$ is not essentially countable, that is, there is no Borel 
\imar{eoei'}
countable\/ {\rm(with at most countable classes)} 
\er\ $\rE$ such that\/ $\Ei\reb\rE$.

\itla{eoei'2}
$\Eo\rebs\Ei,$ in other words,\/ $\ifi\rebs\Ii$.
\een
\ele
\bpf
\ref{eoei'1}
(A version of the argument in \cite{hypsm}, 1.4 and 1.5.) 
Let $\dX$ be the domain of $\rE,$ and 
$\vt:\dnd\to\dX$ a Borel map satisfying  
${x\Ei y}\limp{\vt(x)\rF\vt(y)}.$  
Then $\vt$ is continuous on a dense $\Gd$ set $D\sq\dnd.$ 
We begin with a few definitions.
Let us fix a countable transitive model $\mm$ of $\zhc$
(a big enough fragment of $\ZFC,$ see Remark~\ref{ff:zhc}), 
which contains codes for $D\zd {\vt\res D}\zd \dX$. 

We are going to define, for any $k,$ a pair of points
$a_k\ne b_k\in\dn,$ a number $\ell(k)$ and a tuple 
$\tau_k\in(\dn)^{\ell(k)}$ such that 
\ben
\tenu{(\arabic{enumi})}
\itla{^1}
both $x=\ang{a_0}\we\tau_0\we \ang{a_1}\we\tau_1\we\dots$ and 
$y=\ang{b_0}\we\tau_0\we \ang{b_1}\we\tau_1\we\dots$ are  
elements of $\dnd$ Cohen generic over $\mm$;

\itla{^2}
for any $k,$ 
$\za_{k}=\ang{a_0,b_0}\we\tau_0\we 
\ang{a_1,b_1}\we\tau_1\we\dots\we\ang{a_k,b_k}\we\tau_k$ 
is Cohen generic over $\mm,$ hence so are the subsequences
$\xi_{k}=\ang{a_0}\we\tau_0\we\dots\we\ang{a_k}\we\tau_k$ and 
$\eta_{k}=\ang{b_0}\we\tau_0\we\dots\we\ang{b_k}\we\tau_k$; 

\itla{^3}
for any $k$ and any $z\in\dnd$ such that $\za_{k}\we z$ 
is generic over $\mm$ we have 
$\vt(\xi_{k}\we z)=\vt(\eta_{k}\we z)$.
\een

If this is done then by \ref{^2} choose for any $k$  
a point $z_{k}\in\dnd$ Cohen generic over $\mm[\za_k].$
Then $\za_k\we z_k$ is Cohen generic over $\mm$ by the
product forcing theorem.
It follows by \ref{^3} that $\vt(x_{k})=\vt(y_{k}),$
where $x_{k}=\xi_{k}\we z_{k}$ and 
$y_{k}=\eta_{k}\we z_{k}.$   
Note that $x_{k}\to x$ and $y_{k}\to y$ in $\dnd$ with $k\to\iy,$
and on the other 
hand, all of $x_{k}\yi x\yi y_{k}\yi y$ belong to $D$ because 
of the genericity. 
It follows that $\vt(x)=\vt(y)$ by the choice of $D.$ 
However obviously $\neg\;{x\Ei y},$ so that $\vt$ is not a 
reduction, as required.

To define $a_0\yi b_0\yi \tau_0$ note that  
there exist a perfect set $X\sq\dn$ and a point $z\in\dnd$
such that $\ang{a,b}\we z$ is 
Cohen generic over $\mm$ for any two $a\ne b\in X.$
(Indeed let $\ang{w,z}\in 2^{2\lom}\ti\dnd$ be Cohen generic
over $\mm.$ 
Put $X=\ens{w_a}{a\in\dn},$ where $w_a\in\dn$ is defined by
$w_a(k)=w(a\res k)\zd\kaz k$.)
In particular, $\ang a\we z$ is Cohen generic over $\mm$
for any $a\in X.$
However all points of the form $\ang a\we z$ are 
pairwise \dd\Ei equivalent.
Thus $\vt$ 
sends all of them into one and the same \ddf class, which 
is a countable set by the choice of $\rF.$ 
It follows that there is a pair of $a\ne b$ in $X$ such that 
$\vt(\ang a\we z)\ne\vt(\ang b\we z).$  
This equality is a property of the generic point  
$\ang{a,b}\we z,$ hence, it is forced in the sense that 
there is a number $\ell$ such that  
$\vt(\ang a\we \hat z)=\vt(\ang b\we\hat z)$ whenever
$z\in\dnd,$ $\ang{a,b}\we\hat z$ is Cohen generic over $\mm,$
and $\hat z\res\ell=z\res\ell.$ 
Put $a_0=a\yt b_0=b\yt \tau_0=z\res\ell.$ 

The induction step is carried out by a similar argument.
For instance to define $a_1\yi b_1\yi \tau_1$ we find
points $a'\ne b'\in\dn$ and $z'\in\dnd$ such that
$\ang{a',b'}\we z'$ is Cohen generic over $\mm[a_0,b_0,z]$
and
$\vt(\ang{a_0}\we\tau_0\we\ang{a'}\we z')=
\vt(\ang{a_0}\we\tau_0\we\ang{b'}\we z').$
Yet we have
$\vt(\ang{a_0}\we\tau_0\we\ang{b'}\we z')=
\vt(\ang{b_0}\we\tau_0\we\ang{b'}\we z')$ 
by the choice of $\ell$
(take $\hat z=\tau_0\we\ang{b'}\we z'$). 
Thus
$\vt(\ang{a_0}\we\tau_0\we\ang{a'}\we z')=
\vt(\ang{b_0}\we\tau_0\we\ang{b'}\we z').$
It follows that there is a number $\ell'$
satisfying    
$\vt(\ang{a_0}\we\tau_0\we\ang{a'}\we \hat z)=
\vt(\ang{b_0}\we\tau_0\we\ang{b'}\we \hat z)$
for any $\hat z\in\dnd$ such that  
$\ang{a_0,b_0}\we\tau_0\we\ang{a',b'}\we \hat z$
is Cohen generic over $\mm$ and
$\hat z\res{\ell'}=z'\res{\ell'}.$ 
Put $a_1=a'\yt b_1=b'\yt \tau_1=z'\res{\ell'}.$ 

\ref{eoei'2}
That $\Eo\reb\Ei$ is witnessed by the map 
$f(x)=\ens{\ang{0,n}}{n\in x}.$ 
\epf

\vyk{
\bcot
\lam{eoei"}
$\Eo\rebs\Ei,$ in other words,\/ $\ifi\rebs\Ii$.
\imar{why $\Ei$ and $\Et$ incomparable?}%
\ecot
\bpf
That $\Eo\reb\Ei$ is witnessed by the map 
$f(x)=\ens{\ang{0,n}}{n\in x}.$ 
\epf
}

While $\Ei$ is not countable, the conjunction of 
hypersmoothness and countability characterizes the 
essentially more primitive class of hyperfinite \eqr s. 

\vyk{
The following lemma is a necessary step towards the full 
result below.

\blt
[{{\rm a part of Theorem~5.1 in \cite{djk}}}]
\lam{hf&c=hs}
Suppose that\/ $X\sq\pnd$ is a Borel set and\/ $\Ei\res X$ 
is a countable \er. 
Then\/ $\Ei\res X$ is hyperfinite.
\elt
\bpf

Now, that $\Ei\res X$ is hyperfinite follows from 
Theorem~\ref{thf}.
\epf
}

\punk{The 3rd dichotomy}
\las{hypersm}

The following major result is called the 3rd dichotomy theorem.
 
\bte
[Kechris and Louveau \cite{hypsm}]
\lam{kelu}
Suppose that\/ $\rE$ is a Borel\/ \er\ on some Polish 
space, and\/ $\rE\reb\Ei.$ 
Then\/ \bfei\/ $\rE\reb\Eo$ \bfor\/ $\Ei\reb\rE$.
\ete
\bpf
Starting the proof, we may assume that $\rE$ is a 
$\id11$ \er\ on $\dn,$ and that there is a reduction 
$\rho:\dn\to\dnd$ of $\rE$ to $\Ei,$ of class $\id11.$
In fact it can be assumed that $\rho$ is a bijection.
Indeed define another map $\vpi:\dn\to\dnd$ so that
$\vpi(x)(0)=x$ and $\vpi(x)(n+1)=\rho(x)(n)$ for all
$x\in\dn$ and all $n.$
Then $\vpi$ is a bijection and still a $\id11$ reduction of
$\rE$ to $\Ei$.

Then $R=\ran\rho$ is a $\id11$ subset of $\dnd.$ 
The idea behind the proof is to show that the set $R$ 
is either small enough for $\Ei\res R$ to be 
Borel reducible to $\Eo,$ or otherwise it is big enough 
to contain a closed subset $X$ such that $\Ei\res X$ 
is Borel isomorphic to $\Ei$.

Relations $\cl$ and $\cle$ will denote the inverse 
order relations on $\dN,$ \ie, $m\cle n$ iff $n\le m,$ 
and $m\cl n$ iff $n<m.$  
If $x\in\dnd$ then $x\rec n$ denotes the restriction 
of $x$ (a function defined on $\dN$) 
on the domain ${\cle n},$ \ie, $\iry n.$
If $X\sq\dnd$ then let $X\rec n=\ens{x\rec n}{x\in X}.$ 
Define $x\rc n$ and $X\rc n$ similarly.
In particular, $\dnd\rec n=(\dn)^{\cle n}=(\dn)^{\iry n}.$  

For any sequence $x\in(\dn)^{\cle n},$ let $\glu x$ 
(the {\it depth\/} of $x$) 
be the number (finite or $\iy$) of elements of the set 
$\nab x=\ens{j\cle n}{x(j)\nin \id11(x\rc{j})}.$   
The formula $\glu x\ge d$ 
(of two variables, $d$ running over $\dN\cup\ans\iy$) 
is obviously $\is11$. 

We have two cases:\vtm

{\ubf Case 1:} \ 
all $x\in R=\ran \rho$ satisfy $\glu x<\iy$.\vtm

{\ubf Case 2:} \ 
there exist $x\in R$ with $\glu x=\iy$.\vtm

Case 1 is the easier case. 
The following lemma proves that the Case 1 assumption
implies the \bfei\ case of Theorem~\ref{kelu}.

\ble
\lam{dep3}
Suppose that\/ $X\sq\dnd$ is a\/ $\id11$ set and any\/ 
$x\in X$ satisfies\/ $\glu x<\iy.$ 
Then\/ $\rE_1\res X\reb\rE_0$.
\ele
\bpf
By the choice of $X$ for any $x\in X$ there is a number 
$n$ such that 
${\kaz m\cle n\:\skl x(m)\in\id11(x\rc{m})\skp}.$ 
As the relation between $x$ and $n$ here is clearly $\ip11,$ 
the Kreisel selection theorem (Theorem~\ref{kres})
yields a $\id11$ map 
$\nu:X\to\dN$ such that $x(m)\in\id11(x\rc{n})$ 
holds whenever $x\in X$ and $m\cle\nu(x).$ 
Now define, for each $x\in X,$ $\vt(x)\in\dnd$ as follows: 
$\vt(x)\rec{\nu(x)}=x\rec{\nu(x)},$ but $\vt(x)(j)=\pu$ 
for all $j<\nu(x).$ 
Note that $x\rE_1\vt(x)$ for any $x\in X$.

The other important thing is that 
$\ran\vt\sq Z=\ens{x\in\dnd}{\glu x=0},$ 
where $Z$ is a $\ip11$ set, hence, 
there is a $\id11$ set $Y$ with $\ran\vt\sq Y\sq Z.$ 
In particular $\vt$ reduces $\rE_1\res X$ to 
$\rE_1\res Y.$ 
We observe that $\rE_1\res Y$ is a countable \eqr: 
any \dd{\rE_1}class in $\dnd$ intersects $Y$ by an 
at most countable set 
(as so is the property of $Z,$ a bigger set). 
Thus, $\Ei\res Y$ is hyperfinite by Theorem~\ref{thf}. 
\vyk{
Thus $\rE_1\res Y$ is both countable and hypersmooth. 
It follows, by a theorem of \cite{djk}, that 
\imar{check \cite{djk} !!!}%
$\rE_1\res Y\reb\rE_0,$ as required.
}%
\epf

\punk{Case 2}
\las{3case2}

We are going to prove that then the $\id11$ set
$R=\ran\rho$ contains a $\id11$ subset $X\sq R$ with
$\Ei\reb\Ei\res X.$
This implies the \bfor\ case of Theorem~\ref{kelu}.
Indeed as $\rho$ is a Borel bijection, there exists the
inverse map $\rho\obr,$ and it obviously witnesses
${\Ei\res R}\reb\rE.$ 
On the other hand, $\Ei\reb{\Ei\res X}\reb{\Ei\res R}.$


The required subset $X$ of $R$ will be defined with the
help of a splitting construction developed in \cite{nwf} for 
the study of ``ill''founded Sacks iterations. 

We shall define a map $\vpi:\dN\to\dN,$ which assumes 
infinitely many values and assumes each its value 
infinitely many times (but $\ran\vpi$ may be a proper 
subset of $\dN$), and, for each $u\in\bse,$ a non-empty 
$\is11$ subset $X_u\sq R,$ which satisfy a quite long 
list of properties. 
First of all, if $\vpi$ is already defined at least on 
$\ir0n$ and $u\ne v\in\bse$ then let 
$
\nu_\vpi[u,v]=
\tmin_\cle\ens{\vpi(k)}{k<n\land u(k)\ne v(k)}.
$
(Note that the minimum is taken in the sense of $\cle,$ 
hence, 
it is $\tmax$ in the sense of $\le,$ the usual order). 
Separately, put $\vpi[u,u]=-1$ for any $u$.

Now we give the list of requirements.

\ben
\tenu{{\rm(\roman{enumi})}}
\itla{z1}
if $\vpi(n)\nin\ens{\vpi(k)}{k<n}$ then 
$\vpi(n)\cl \vpi(k)$ for any $k<n$;

\itla{z2}
every $X_u$ is a non-empty $\is11$ subset of $R$;

\itla{z3}
if $u\in2^n\yt x\in X_u,$ and $k<n,$ then 
$\vpi(k)\in \nab x$;

\itla{z4}
if $u\yi v\in2^n$ then 
$X_u\rc{\npi[u,v]}=X_v\rc{\npi[u,v]}$;

\itla{z5}
if $u\yi v\in2^n$ then 
$X_u\rec{\npi[u,v]}\cap X_v\rec{\npi[u,v]}=\pu$;

\itla{z6}\msur
$X_{u\we i}\sq X_u$ for all $u\in\bse$ and $i=0,1$;

\itla{z7}\msur
$\tmax_{u\in2^n}\dia X_u\to0$ as $n\to\iy$ 
(a reasonable 
Polish metric on $\dnd$ is assumed to be fixed);

\itla{z8}
for any $n,$ a certain condition, in terms of the
Gandy -- Harrington forcing,
similar to \ref{d1i,} in \nrf{>smooth} or
\ref{d2i1,} in \nrf{d2split},
related to all sets $X_u\yt u\in 2^n,$ so that, as a 
consequence, $\bigcap_nX_{a\res n}\ne\pu$ 
for any $a\in\dn$.
\een

Let us demonstrate how such a system of sets and a function 
$\vpi$ accomplish Case 2. 
According to \ref{z7} and \ref{z8}, for any $a\in\dn$ 
the intersection $\bigcap_nX_{a\res n}$ contains a single 
point, let it be $F(a),$ and $F$ is continuous and 
$1-1$.

Put $J=\ran\vpi=\ens{j_m}{m\in\dN},$
in the \dd<increasing order; $J\sq\dN$ is infinite. 
Let $n\in\dN.$ 
Then $\vpi(n)=j_m$ for some (unique) $m:$ we put 
$\psi(n)=m.$ 
Thus $\psi:\dN\onto\dN$ and the preimage 
$\psi\obr(m)=\vpi\obr(j_m)$ is an infinite subset of 
$\dN$ for any $m.$ 
This allows us to define a parallel system of 
sets $Y_u\sq\dnd\yt u\in\bse,$ as follows. 
Put $Y_\La=\dnd.$ 
Suppose that $Y_u$ has been defined, $u\in2^n.$ 
Put $j=\vpi(n)=j_{\psi(n)}.$ 
Let $K$ be the number of all indices $k<n$ still 
satisfying $\vpi(k)=j,$ perhaps $K=0.$ 
Put $Y_{u\we i}=\ens{x\in Y_u}{x(j)(K)=i}$ for 
$i=0,1$. 

Each of $Y_u$ is clearly a basic clopen set in $\dnd,$ 
and one easily verifies that conditions 
\ref{z1} -- \ref{z7}, except for \ref{z3},  
are satisfied for the sets $Y_u$ (instead of $X_u$) and 
the map $\psi$ (instead of $\vpi$), in particular, for 
any $a\in\dn,$  $\bigcap_nY_{a\res n}=\ans{G(a)}$ 
is a singleton, and the map $G$ is continuous and $1-1.$ 
(We can, of course, define $G$ explicitly: 
$G(a)(m)(l)=a(n),$ where $n\in\dN$ is chosen so that 
$\psi(n)=m$ and there is exactly $l$ numbers $k<n$ 
with $\psi(k)=m$.) 
Note finally that $\ens{G(a)}{a\in\dn}=\dnd$ since 
by definition $Y_{u\we 1}\cup Y_{u\we 0}=Y_u$ for 
all $u$.

We conclude that the map $\vt(x)=F(G\obr(x))$ is a 
continuous bijection 
(hence, in this case, a homeomorphism by compactness) 
$\dnd\onto X.$ 
We further assert that $\vt$ satisfying the following: 
for each $y\yi y'\in \dnd$ and $m$,
\dm
y\rec m= y'\rec m\quad\hbox{iff}\quad
\vt(y)\rec{j_m}=\vt(y')\rec{j_m}
\,.
\eqno(\ast)
\dm
Indeed, let $y=G(a)$ and $x=F(a)=\vt(y),$ and 
similarly $y'=G(a')$ and $x'=F(a')=\vt(y'),$ where 
$a\yi a'\in\dn.$ 
Suppose that $y\rec m= y'\rec m.$ 
According to \ref{z5} for $\psi$ and the sets $Y_u,$ 
we then have $m\cl\nsi[a\res n,a'\res n]$ for any $n.$ 
It follows, by the definition of $\psi,$ that 
$j_m\cl\npi[a\res n,a'\res n]$ for any $n,$ hence, 
$X_{a\res n}\rec{j_m}=X_{a\res n}\rec{j_m}$ 
for any $n$ by \ref{z4}. 
Assuming now that Polish metrics on all spaces 
$(\dn)^{\cle j}$ are chosen so that 
$\dia Z\ge\dia {(Z\rec j)}$ for all $Z\sq\dn$ and $j,$ 
we easily obtain that $x\rec{j_m}=x'\rec{j_m},$ \ie, 
the right-hand side of $(\ast).$ 
The inverse implication in $(\ast)$ is proved similarly. 

Thus we have $(\ast),$ but this means that $\vt$ is a 
continuous reduction of $\Ei$ to $\Ei\res X,$ thus, 
$\Ei\reb {\Ei\res X},$ as required.\vtm

\qeDD{Theorem~\ref{kelu} modulo the construction\/ 
\ref{z1} -- \ref{z8}}

\punk{The construction}
\las{3constr}

Recall that $R\sq\dnd$ is a fixed non-empty $\is11$ set 
such that $\glu x=\iy$ for each $x\in R.$ 
Set $X_\La=R$.

Now suppose that the sets $X_u\sq R$ with $u\in2^n$ 
have been defined and satisfy the applicable part of 
\ref{z1} -- \ref{z8}.\vom 

{\ubf Step 1.}  
Our 1st task is to choose $\vpi(n).$  
Let $\ans{j_1<\dots<j_m}=\ens{\vpi(k)}{k<n}.$ 
For any $1\le p\le m,$ let $N_p$ be the number of all 
$k<n$ with $\vpi(k)= j_p.$\vom 

{\sl Case 1a\/}.  
If some numbers $N_p$ are $<m$ then  
choose $\vpi(n)$ among $j_p$ with the least $N_p,$ and 
among them the least one.\vom

{\sl Case 1b\/}:  
$N_p\ge m$ (then actually $N_p=m$) for all $p\le m.$ 
It follows from our assumptions, in particular \ref{z4}, 
that $X_u\rc{j_m}=X_v\rc{j_m}$ for all $u\yi v\in2^n.$ 
Let $Y=X_u\rc{j_m}$ for any such $u.$ 
Take any $y\in Y.$ 
Then $\nab y$ is infinite, hence, there is some 
$j\in\nab y$ with $j\cl j_m.$ 
Put $\vpi(n)=j$. 

We have something else to do in this case. 
Let $X'_u=\ens{x\in X_u}{j\in\nab y}$ for any $u\in2^m.$ 
Then we easily have 
$X'_u=\ens{x\in X_u}{x\rc{j_m}\in Y'},$ where 
$Y'=\ens{y\in Y}{j\in\nab y}$ is a non-empty $\is11$ set, 
so that the sets $X'_u\sq X_u$ are non-empty $\is11.$ 
Moreover, as $j_m$ is the \dd\cle least in 
$\ens{\vpi(k)}{k<n},$ we can easily show that the system 
of sets $X'_u$ still satisfies \ref{z4}.
This allows us to assume, without any loss of generality, 
that, in Case~1b, $X'_u=X_u$ for all $u,$ or, 
in other words, that any $x\in X_u$ for any $u\in2^n$ 
satisfies $j=\vpi(n)\in\nab x.$ 
(This is true in Case~1a, of course, because then 
$\vpi(n)=\vpi(k)$ for some $k<n$.)\vom

Note that this manner to choose $\vpi(n)$ implies 
\ref{z1} and also implies that $\vpi$ takes 
infinitely many values and takes each its value 
infinitely many times. 
                              
The continuation of the construction requires the following

\ble
\lam{suz}
If\/ $u_0\in2^n$ and\/ $X'\sq X_{u_0}$ is a non-empty\/ 
$\is11$ set then there is a system of\/ $\is11$ sets\/ 
$\pu\ne X'_u\sq X_u$ with\/ $X'_{u_0}=X',$ which still 
satisfies\/ \ref{z4}. 
\ele
\bpf
For any $u\in2^n,$ let 
$X'_u=\ens{x\in X_u}{x\rc{n(u)}\in X'\rc{n(u)}},$  
where $n(u)=\npi[u,u_0].$ 
In particular, this gives $X'_{u_0}=X',$ because 
$\npi[u_0,u_0]=-1.$ 
The sets $X'_u$ are as required, via a routine verification.
\epF{Lemma}

{\bf Step 2.} 
First of all put $j=\vpi(n)$ and $Y_u=X_u\rc j.$  
(All $Y_u$ are equal to $Y$ in Case~1b, but the argument 
pretends to make no difference between 1a and 1b). 
Take any $u_1\in2^n.$ 
By the construction any element $x\in X_{u_1}$ satisfies 
$j\in\nab x,$ so that $x(j)\nin\id11(x\rc j).$ 
As $X_{u_1}$ is a $\is11$ set, it follows that  
$\ens{x'(j)}{x'\in X_{u_1}\land x'\rc j=x\rc j}$ 
is not a singleton, in fact is uncountable. 
It follows that there is a number $l_{u_1}$ having the 
property that the $\is11$ set 
\dm
Y'_{u_1}\:=\:
\ens{y\in Y_{u_1}}{\sus x,x'\in X_{u_1}\:
\skl x'\rc j=x\rc j=y\land 
l_{u_1}\in x(j)\land l_{u_1}\nin x'(j)\skp}
\dm
is non-empty. 
We now put $X'=\ens{x\in X_{u_1}}{x\rc j\in Y'_{u_1}}$ 
and define $\is11$ sets $\pu\ne X'_u\sq X_u$ as in the 
lemma, in particular,
$X'_{u_1}=X'\yt X'_{u_1}\rc j=Y'_{u_1},$
still \ref{z4} is satisfied, and in addition 
\dm
\kaz y\in X'_{u_1}\rc j\;
\sus x,x'\in X'_{u_1}\:
\skl x'\rc j=x\rc j=y\land 
l_{u_1}\in x(j)\land l_{u_1}\nin x'(j)\skp
\eqno(1)
\dm

Now take some other $u_2\in2^n.$ 
Let $\nu=\npi[u_1,u_2].$ 
If $j\cl \nu$ then $X_{u_1}\rc j=X_{u_2}\rc j,$ so 
that we already have, for $l_{u_2}=l_{u_1},$ that 
\dm
\kaz y\in X'_{u_2}\rc j\;
\sus x,x'\in X'_{u_2}\:
\skl x'\rc j=x\rc j=y\land 
l_{u_2}\in x(j)\land l_{u_2}\nin x'(j)\skp\,,
\eqno(2)
\dm
and can pass to some $u_3\in2^n.$ 
Suppose that $\nu\cle j.$ 
Now things are somewhat nastier. 
As above there is a number $l_{u_2}$ such that 
\dm
Y'_{u_2}\:=\:
\ens{y\in Y_{u_2}}{\sus x,x'\in X_{u_2}\:
\skl x'\rc j=x\rc j=y\land 
l_{u_2}\in x(j)\land l_{u_2}\nin x'(j)\skp}
\dm
is a non-empty $\is11$ set, thus, we can define 
$X''=\ens{x\in X_{u_1}}{x\rc j\in Y'_{u_1}}$ and 
maintain the construction of Lemma~\ref{suz}, getting 
non-empty $\is11$ sets $X''_u\sq X'_u$ still satisfying 
\ref{z4} and $X''_{u_2}=X'',$ therefore, 
we still have $(2)$ for the set $X''_{u_2}.$

Yet it is most important in this case that $(1)$ 
is preserved, \ie, 
it still holds for the set $X''_{u_1}$ instead of 
$X'_{u_1}$! \ 
Why is this ? \ 
Indeed, according to the construction in the proof of 
Lemma~\ref{suz}, we have 
$X''_{u_1}=\ens{x\in X'_{u_1}}{x\rc\nu\in X''\rc \nu}.$ 
Thus, although, in principle, $X''_{u_1}$ is smaller than 
$X'_{u_1},$ for any $y\in X''_{u_1}\rc j$ we have 
\dm
\ens{x\in X''_{u_1}}{x\rc j=y}\;=\;
\ens{x\in X'_{u_1}}{x\rc j=y}\,,
\dm
simply because now we assume that $\nu\cle j.$ 
This implies that $(1)$ still holds.

Iterating this construction so that each $u\in2^n$ is 
eventually encountered, we obtain, in the end, a system 
of non-empty $\is11$ sets, let us call them ``new'' 
$X_u,$ but they are subsets of the ``original'' $X_u,$ 
still satisfying \ref{z4}, still satisfying that 
$\vpi(n)\in\nab x$ for each $x\in\bigcap_{u\in2^n}X_u,$ 
and, in addition, for any $u\in2^n$ there is a number 
$l_u$ such that ${j\cl \npi[u,v]}\limp {l_u=l_v}$ and 
\dm
\kaz y\in X_{u}\rc j\;
\sus x,x'\in X_{u}\:
\skl x'\rc j=x\rc j=y\land 
l_{u}\in x(j)\land l_{u}\nin x'(j)\skp\,.
\eqno(\ast)
\dm

{\ubf Step 3.}  
We define the \dd{(n+1)}th level of sets by 
$X_{u\we 0}=\ens{x\in X_u}{l_u\in x(j)}$ and 
$X_{u\we 1}=\ens{x\in X_u}{l_u\nin x(j)}$ for 
all $u\in2^n,$ where still $j=\vpi(n).$ 
It follows from $(\ast)$ that all these $\is11$ sets 
are non-empty. 

\ble
\lam{n+1}
The system of sets\/ $X_s\yt s\in2^{n+1}$ just defined
satisfies\/ \ref{z4}, 
\ref{z5}.
\ele
\bpf 
Let $s=u\we i$ and $t=v\we i'$ belong to $2^{n+1},$ 
so that $u\yi v\in2^n$ and $i,i'\in\ans{0,1}.$ 
Let $\nu=\npi[u,v]$ and $\nu'=\npi[s,t]$.
\vom

{\it Case 3a\/}: $\nu\cle j=\vpi(n).$ 
Then easily $\nu=\nu',$ 
so that \ref{z5} immediately follows from \ref{z5} 
at level $n$ for $X_u$ and $X_v.$
As for \ref{z4}, we have $X_s\rc\nu=X_u\rc\nu$ 
(because by definition $X_s\rc j=X_u\rc j$), and 
similarly $X_t\rc\nu=X_v\rc\nu,$  
therefore, $X_t\rc{\nu'}=X_s\rc{\nu'}$ since 
$X_u\rc{\nu}=X_v\rc{\nu}$ by \ref{z4} at level $n.$ 
\vom

{\it Case 3b\/}: $j\cl\nu$ and $i=i'.$ 
Then still $\nu=\nu',$ thus we have \ref{z5}.
Further, $X_u\rc\nu=X_v\rc\nu$ by \ref{z4} at 
level $n,$ hence, $X_u\rec j=X_v\rec j,$ hence, 
$l_u=l_v$ (see above). 
Now, assuming that, say, $i=i'=1$ and $l_u=l_v=l,$ 
we conclude that
\dm
X_s\rc{\nu'}=\ens{y\in X_u\rc\nu}{l\in y(j)}
=\ens{y\in X_v\rc\nu}{l\in y(j)}=X_t\rc{\nu'}\,.
\dm

{\it Case 3c\/}: $j\cl\nu$ and $i\ne i',$ say, 
$i=0$ and $i'=1.$ 
Now $\nu'=j.$ 
Yet by definition $X_s\rc j=X_u\rc j$ and 
$X_t\rc j=X_v\rc j,$ so it remains to apply \ref{z4} 
for level $n.$ 
As for \ref{z5}, 
note that by definition 
$l\nin x(j)$ for any $x\in X_s=X_{u\we 0}$ while 
$l\in x(j)$ for any $x\in X_t=X_{v\we 1},$ where 
$l=l_u=l_v$.
\epF{Lemma}

{\it Step 4\/}.
In addition to \ref{z4} and \ref{z5}, we already have 
\ref{z1}, \ref{z2}, \ref{z3}, \ref{z6} at level $n+1.$ 
To achieve the remaining properties \ref{z7} and \ref{z8}, 
it suffices to consider, one by one, all elements 
$s\in2^{n+1},$ 
finding, at each such a substep, a non-empty $\is11$ 
subset of $X_s$ which is consistent with the requirements 
of \ref{z7} and \ref{z8} 
(for instance, for \ref{z7}, just take it so the diameter 
is $\le 2^{-n}$), 
and then reducing all other sets $X_t$ by Lemma~\ref{suz} 
at level $n+1$.
\vtm

\epF{Construction and Theorem~\ref{kelu}}

\punk{Above $\Ei$}
\las{abovei1}

Recall that an embedding is a $1-1$ reduction, and 
an invariant embedding is an embedding $\vt$ such that 
its range is an invariant set, see \nrf{bored}.

\bte
[{{\rm Kechris and Louveau \cite{hypsm}}}]
\lam{abi1}
Suppose that\/ $\Ei\reb\rF,$ where\/ $\rF$ is an 
analytic \er\ on a Polish space\/ $\dY.$ 
Then both\/ $\Ei\emn\rF$ 
and\/ $\Ei\embi\rF$.
\ete
\bpf
To prove the first statement, 
let $\cle$ be the inverted order on $\dN,$ \ie, 
$m\cle n$ iff $n\le m.$ 
Let $\gP$ be the collection of all sets $P\sq\dnd$ 
such that there is a continuous $1-1$ map 
$\eta:\dnd\onto P$ satisfying 
\dm
{x\rec{n}=y\rec{n}}\leqv
{\eta(x)\rec{n}=\eta(y)\rec{n}}
\dm
for all $n$ and $x\yi y\in\dnd,$ where 
$x\rec{n}=\sis{x(i)}{i\cle n}$ for any 
$x\in\dnd.$ 
Clearly any such a map is a continuous embedding of 
$\Ei$ into itself.

This set $\gP$ can be used as a forcing notion to extend the 
universe by a sequence of reals $x_i$ so that each 
$x_n$ is Sacks--generic over $\sis{x_i}{i\curle n}.$ 
This is an example of iterated Sacks extensions with an 
ill-founded ``skeleton'' of iteration, which we 
defined in~\cite{nwf}.
(See \cite{kell} on more recent developments on ill-iterated
forcing.) 
Here, the ``skeleton'' is $\dN$ with the 
inverted order $\cle$.

The method of \cite{nwf} contains a study of 
continuous and Borel functions on 
sets in $\gP.$ 
In particular it is shown there that Borel maps 
admit the following {\it cofinal classification\/} 
on sets in $\gP:$ 
if $\dY$ is Polish, $P'\in\gP,$ and 
$\vt:P'\to\dY$ is Borel then there is a set 
$P\in\gP\yt P\sq P',$ on which $\vt$ is 
continuous, and either a constant or, for some $n,$  
$1-1$ on $P\rec{n}$ in the sense that,   
\dm
\text{for all }\;x\yi y\in P:\quad
{x\rec{n}=y\rec{n}}\leqv{\vt(x)=\vt(y)}\,.
\eqno(\ast)
\dm
We apply this to a Borel map 
$\vt:\dnd\to\dY$ which reduces $\Ei$ to 
$\rF.$ 
We begin with $P'=\dnd$ and find a set $P\in\gP$ 
as indicated. 
Since $\vt$ cannot be a constant on $P$ 
(indeed, any $P\in\gP$ contains many pairwise 
\dd\Ei inequivalent elements), 
we have $(\ast)$ for some $n.$ 
In other words, there is a $1-1$ continuous map 
$f:P\rec n\to \dY$ 
(where ${P\rec n}=\ens{x\rec n}{x\in P}$) 
such that $\vt(x)=f(x\rec n)$ for all $x\in P.$ 
Now, suppose that $x\in\dnd.$ 
Define $\zeta(x)=z\in\dnd$ so that 
$z(i)=\dN\ti\ans0$ for $i<n$ and $z(n+i)=x(i)$ for all $i.$ 
Finally set $\vt'(x)=f(\eta(\zeta(x))\resic n)$ 
for all $x\in\dnd:$ this map turns out to be
a continuous embedding of $\Ei$ in $\rF$. 

Now we prove the second claim. 
We can assume that $\dY=\dn$ and that $\vt:\dnd\to\dn$ 
is already a continuous embedding of $\Ei$ into $\rF.$ 
Let $Y=\ran\vt$ and $Z=[Y]_{\rF}.$ 
Normally $Y\yi Z$ are analytic, but in this case they are 
even Borel. 
Indeed $Z$ is the projection of 
$P=\ens{\ang{z,x}}{z\rF\vt(x)},$ a Borel subset of 
$\dn\ti\dnd$ whose all cross-sections are 
\dd\Ei equivalence classes, \ie, \dds compact sets. 
It is known 
\imar{reference}%
that in this case $Z$ is Borel and, moreover, there is 
a Borel map $f:Z\to\dnd$ such that $f(z)\Ei x$ 
whenever $z\rF\vt(x)$.

We can convert $f$ to a $1-1$ map $g:Z\to\dnd$ 
with the same properties: 
$g(z)(n)=f(z)(n)$ for $n\ge 1,$ but 
$g(z)(0)=z.$ 
Then $\vt:\dnd\to Z\sq\dn$ and $g:Z\to\dnd$ 
are Borel $1-1$ maps 
($\vt$ is even continuous, but this does not matter now), 
and, for any $x\in\dnd,$ $\vt$ maps $[x]_{\Ei}$ into 
$[\vt(x)]_{\rF}\sq Z,$ and $g$ maps $[\vt(x)]_{\rF}$ 
back into $[x]_{\Ei}.$ 
It remains to apply the construction from the Cantor -- 
Bendixson theorem, to get a Borel embedding, say, $f$ 
of $\Ei$ into $\rF$ with $\ran f=Z,$ that is an 
invariant embedding.
\epf

The following theorem shows that orbit \eqr s of  
Polish group actions cannot reduce $\Ei$. 
 
\bte
[{{\rm Kechris and Louveau \cite{hypsm}}}]
\lam{e1pga}
Suppose that\/ $\dG$ is a Polish group and\/ $\dX$ 
is a Borel\/ \dd\dG space. 
Then\/ $\Ei$ is\/ {\ubf not} Borel reducible to\/ $\ergx$.
\ete
\bpf
Towards the contrary, let $\vt:\dnd\to\dX$ 
be a Borel reduction of $\Ei$ to $\rE.$ 
We can assume, by Theorem~\ref{abi1}, that $\vt$ is 
in fact an invariant embedding, \ie, $1-1$ and 
$Y=\ran\vt$ is an \dde invariant set. 
Define, for $g\in\dG$ and $x\in\dnd,$ 
$g\ac x=\vt\obr(g\ac\vt(x)).$ 
Then this is a Borel action of $\dG$ on $\dnd$ such 
that the induced relation $\aer\dG\dnd$ 
coincides with $\Ei$. 

Let us fix $x\in\dnd.$ 

Consider any $y=\sis{y_n}n\in[x]_{\Ei}.$ 
Then $[x]_{\Ei}=\bigcup_nC_n(y),$ where each set 
$C_n(y)=\ens{y'\in\dnd}{\kaz m\ge n\:(y(n)=y'(n))}$  
is Borel (even compact). 
It follows that $\dG=\bigcup_nG_n(y),$ where each 
$G_n(y)=\ens{g\in\dG}{g(x)\in C_n(y)}$ is Borel. 
Thus, as $\dG$ is Polish, there is a number 
$n$ such that $G_n(y)$ is not meager in $\dG$ 
(then this will hold for all $n'\ge n,$ of course). 
Let $n(y)$ be the least such an $n$. 

We assert that for any $n$ the set 
$Y_n(x)=\ens{y\res{\iry n}}{y\in[x]_{\Ei}\land n(x)=n}$ 
is at most countable. 
Indeed suppose that $Y_n(x)$ is not countable. 
Note that if $y_1$ and $y_2$ in $[x]_{\Ei}$ have 
different restrictions $y_i\res{\iry n}$ then the sets 
$C_n(y_1)$ and $C_n(y_2)$ are disjoint, therefore, 
the sets $G_n(y_1)$ and $G_n(y_2)$ are disjoint, 
so we would have uncountably many pairwise disjoint 
non-meager sets in $\dG,$ contradiction. 
Thus all sets $Y_n(x)$ are countable. 

It is most important that $Y_n(x)$ depends on $[x]_{\Ei}$ 
rather than $x$ itself.
More exactly, if $x'\in[x]_{\Ei}$ 
then $Y_n(x)=Y_n(x'):$ this is because any set $G_n(y)$ 
in the sense of $x'$ is just a shift, within $\dG,$ of 
the set $G_n(y)$ in the sense of $x.$ 
Therefore, putting  
$Y(x)=\bigcup_n\ens{\bar u}{u\in Y_n(x)},$ 
where, for $u\in(\dn)^{\iry n},$ $\bar u\in\dnd$ is 
defined by $\bar u\res{\iry n}=u$ and $\bar u(k)(j)=0$
for $k<n$ and all $j,$ we obtain a set
$Y=\bigcup_{x\in\dnd}Y(x)$ 
with the property that $Y\cap[x]_{\Ei}$ is non-empty 
and at most countable for any $x\in\dnd$. 

The other important fact is that the relation $y\in Y(x)$ 
is Borel: this is because it is assembled from Borel 
relations via the Vaught quantifier ``there exists 
nonmeager-many'', known to preserve the borelness. 
\imar{reference}%
It follows that
\dm
Y=\ens{y}{\sus x\:(y\in Y_x)}=
\ens{y}{\kaz x\:(x\in[y]_{\Ei}\imp y\in Y(x)}
\dm
is a Borel subset of $\dnd.$ 
By the uniformization theorem for 
\imar{reference}%
Borel sets with countable sections, there is a Borel 
map $f$ defined on $\dnd$ so that $f(x)\in Y(x)$ for 
any $x.$
This implies $\Ei\reb{\Ei\res Y}.$ 
On the other hand, $\Ei\res Y$ is a countable \eqr\  
by the above, which is a contradiction to Lemma~\ref{eoei'}.
\epf

The theorem just proved allows us to accomplish the proof
of Theorem~\ref{sol} by adding its last condition \ref{s6} 
to the equivalence of its first five conditions established
in \grf{idI1}.
Since $\orb$ implies $\reb,$ the following lemma implies 
the result required.

\ble
\lam{soll}
If\/ $\cI\sq\pn$ is a polishable ideal then\/
$\Ei\not\reb \rei$.
\ele
\bpf
Recall that if $\cI$ is polishable then $\rei$ is induced
by a Polish action of the \dd\sd group of $\cI$ on $\pn.$
It remains to apply Theorem~\ref{e1pga}.
\epf

We are able now to also give another proof of a result
already obtained by different method.
(See Corollary~\ref{eiey}.)

\bcor
\lam{eiey*}
$\Ey\not\reb\Ei$.
\ecor
\bpf
If $\Ey\reb\Ei$ then by Theorem~\ref{kelu} either $\Ey\reb\Eo$
or $\Ey\eqb\Ei.$
The ``either'' case contradicts Theorem~\ref{Eynon}.
The ``or'' case contradicts  
Theorem~\ref{e1pga} since $\Ey$ is induced by a Polish
action of $F_2$.
\epf

\api

\parf{Actions of the infinite symmetric group}
\las{groas}

This Section is connected with the next one (on turbulence). 
We concentrate on a main result in this area, due to 
Hjorth, that turbulent \er s are not reducible to 
those induced by actions of $\isg.$ 
In particular, we shall prove the following:

\ben
\tenu{\Roman{enumi}}
\def\labelenumi{\theenumi.}
\itla{isg1}
Lopez-Escobar: any invariant Borel set of countable models 
is the truth domain of a formula of $\lww$.

\itla{isg2}
Any orbit \er\ of a Polish action of a closed subgroup of 
$\isg$ is classifiable by countable structures (up to isomorphism). 

\itla{isg3}
Any \er, classifiable by countable structures, is Borel reducible 
to isomorphism of countable ordered graphs.

\itla{isg4}
Any \poq{Borel} \er, classifiable by countable structures, 
is Borel reducible to one of \er s $\rT_\xi$.

\itla{isg5}
Any \er, classifiable by countable structures and induced by a 
Polish action (of a Polish group), 
is Borel reducible to one of \er s $\rT_\xi$ on a comeager set.

\itla{isg6}
Any ``turbulent'' \er\ $\rE$ is generically 
\dd{\rT_\xi}ergodic for any $\xi<\omi,$ in particular, $\rE$ 
is not Borel reducible to $\rT_\xi$.

\itla{isg7}
Any ``turbulent'' \er\ is not classifiable by countable 
structures: a corollary of \ref{isg6} and \ref{isg5}. 

\itla{isg8}
A generalization of \ref{isg7}: 
any ``turbulent'' \er\ is not Borel reducible to a \er\ that 
can be obtained from the equality $\rav\dN$ using operations 
defined in \nrf{opeer}.
\een

Scott's analysis, involved in proofs of \ref{isg4} and \ref{isg5}, 
appears only in a rather mild and self-contained version.

\punk{Infinite symmetric group $\isg$}
\las{loac}

Let $\isg$ be the group of all permutations
\index{group!sinfty@$\isg$}
(\ie, 1--1 maps $\dN\onto\dN$) of $\dN,$ with the 
superposition as the group operation. 
Clearly $\isg$ is a $\Gd$ subset of $\dnn,$ hence, a 
Polish group. 
A compatible complete metric on $\isg$ can be defined by 
$D(x,y)=d(x,y)+d(x\obr,y\obr),$ where $d$ is the ordinary 
complete metric of $\dnn,$ \ie, $d(x,y)=2^{-m-1},$ where 
$m$ is the least such that $x(m)\ne y(m).$ 
Yet $\isg$ admits no compatible left-invariant 
complete metric \cite[1.5]{beke}. 
\imar{Proof of $\isg$ not \cli~?}%

For instance isomorphism relations of various kinds of 
countable structures are orbit \er s induced by $\isg.$ 
Indeed, suppose that $\cL=\sis{R_i}{i\in I}$ is a countable 
\index{language!countable relational}%
relational language, \ie, $0<\card I\le\alo$ and each $R_i$ 
is an \dd{m_i}ary relational symbol. 
We put~\footnote
{\ $X_\cL$ is often used to denote $\mox\cL$.}   
$\mox\cL=\prod_{i\in I}\cP(\dN^{m_i}),$ 
\index{zzModL@$\mox\cL$}%
\index{action!logic@logic action $\loa\cL$ of $\isg$}%
\index{zzjL@$\loa\cL$}%
the space of (coded) \dd\cL{\it structures\/} on $\dN.$ 
The {\it logic action\/} $\loa\cL$ of $\isg$ on $\mox\cL$ 
is defined as follows: 
if $x=\sis{x_i}{i\in I}\in\mox\cL$ and $g\in\isg$ 
then $y=\loa\cL(g,x)=g\app x=\sis{y_i}{i\in I}\in\mox\cL,$ 
where we have 
\dm
{\ang{k_1,\dots,k_{m_i}}\in x_i} \leqv 
{\ang{g(k_1),\dots,g(k_{m_i})}\in y_i}  
\dm
for all $i\in I$ and $\ang{k_1,\dots,k_{m_i}}\in\dN^{m_i}.$ 
Then $\stk{\mox\cL}{\loa\cL}$ is a Polish \dd\isg space  
and \dd{\loa\cL}orbits in $\mox\cL$ are exactly the 
isomorphism classes of \dd\cL structures, which is a reason 
to denote the associated equivalence relation 
\index{isomorphism!of st@of structures, $\ism\cL$}%
\index{zzconl@$\ism\cL$}%
$\aer{\loa\cL}{\mox\cL}$ as ${\ism\cL}$. 

If $G$ is a subgroup of $\isg$ then $\loa\cL$ restricted to 
$G$ is still an action of $G$ on $\mox\cL,$ whose orbit \er\ 
\index{isomorphism!of st@induced by $G,$ $\izm\cL G$}%
\index{zzconlg@$\izm\cL G$}%
will be denoted by $\izm\cL G,$ \ie, $x\izm\cL G y$ iff 
${\sus g\in G\:(g\app x=y)}$.

\punk{Borel invariant sets}
\las{BIS}

A set $M\sq\mox\cL$ is {\it invariant\/} if 
\index{set!invariant}%
$\ek M{\ism\cL}=M.$ 
There is a convenient characterization of {\it Borel\/} 
invariant sets, in terms of $\lww,$ an infinitary extension 
of $\cL=\sis{R_i}{i\in I}$ 
\index{language!Lww@$\lww$}%
\index{zzLww@$\lww$}%
by countable conjunctions and disjunctions. 
To be more exact, 
\ben
\tenu{\arabic{enumi})}
\item 
any $R_i(v_0,\dots,v_{m_i-1})$ is an atomic 
formula of $\lww$ (all $v_i$ being variables over $\dN$ 
and $m_i$ is the arity of $R_i$), and 
propositional connectives and quantifiers ${\sus}\yi {\kaz}$ 
can be applied as usual; 

\item
if $\vpi_i,\;i\in\dN,$ are formulas of $\lww$ whose free 
variables 
are among a finite list $v_0,\dots,v_n$ then $\bigvee_i\vpi_i$ 
and $\bigwedge_i\vpi_i$ are formulas of $\lww$.
\een
If $x\in\mox\cL,\msur$ 
$\vpi(v_1,\dots,v_n)$ is a formula of $\lww,$ and 
$i_1,\dots,i_n\in\dN,$ then $x\mo\vpi(i_1,\dots,i_n)$ means that 
$\vpi(i_1,\dots,i_n)$ is satisfied on $x,$ in the 
usual sense that involves transfinite induction on the ``depth'' 
of $\vpi,$ see \cite[16.C]{dst}.

\bte
[{{\rm Lopez-Escobar, see \cite[16.8]{dst}}}] 
\lam{lopes}
A set\/$M\sq\mox\cL$ is invariant and Borel iff\/ 
$M=\ens{x\in\mox\cL}{x\mo\vpi}$ for a 
closed formula\/ $\vpi$ of\/ $\lww$.
\ete
\bpf
To prove the nontrivial direction let $M\sq\mox\cL$ 
be invariant and Borel. 
Put $B_s=\ens{g\in\isg}{s\su g}$ for any injective 
$s\in\dN\lom$ (\ie, $s_i\ne s_j$ for $i\ne j$), this is a 
clopen subset of $\isg$ 
(in the Polish topology of $\isg$ inherited from $\dnn$). 
If $A\sq\isg$ then let $s\forc A(\dog)$ mean that the set 
$B_s\cap A$ is co-meager in $B_s,$ \ie, $g\in A$ holds for \pv\ 
$g\in \isg$ with $s\su g.$ 
The proof consists of two parts:

\ben
\tenu{(\roman{enumi})}
\itla{pari}\msur
$M=\ens{x\in\mox\cL}{\La\forc \dog\app x\in M}$ 
(where $g\app x=\loa\cL(g,x),$ see above);

\itla{pard}
For any Borel $M\sq\mox\cL$ and any $n$ there is a formula 
$\vpi_{M}^n(v_0,\dots,v_{n-1})$ of $\lww$ such that we have, for 
every $x\in\mox\cL$ and every injective $s\in\dN^n:$  
$x\mo\vpi_{M}^n(s_0,\dots,s_{n-1})$ iff 
$s\forc\dog{}\obr\app x\in M$.
\een

\ref{pari} is clear: since $M$ is invariant, we have 
$g\app x\in M$ for all $x\in M$ and $g\in\isg,$ on the other 
hand, if $g\app x\in M$ for at least one $g\in\isg$ then 
$x\in M$. 

To prove \ref{pard} we argue by induction on the Borel 
complexity of $M.$ 
Suppose, for the sake of simplicity, that $\cL$ contains 
a single binary predicate, say, $R(\cdot,\cdot);$ then 
$\mox\cL=\cP(\dN^2).$ 
If $M=\ens{x\sq\dN^2}{\ang{k,l}\nin x}$ for some $k\yi l\in\dN$ 
then take
\dm
\TS
\kaz u_0\,\dots\,\kaz u_m\;
\skl
\bigwedge_{i<j\le m}(u_i\ne u_j)\,\land\,
\bigwedge_{i<n}(u_i=v_i)\,\imp\,\neg\:R(u_k,u_l)
\skp\,,
\dm
where $m=\tmax\ans{l,k,n},$ as $\vpi_{M}^n(v_0,\dots,v_{n-1}).$ 
Further, take 
\dm
\TS
\bay{rl}
\bigwedge_{k\ge n} 
\,\kaz u_0\,\dots\,\kaz u_{k-1}\;
\bigvee_{m\ge k}
\,\sus w_0\,\dots\,\sus w_{m-1}\;
\skl
\bigwedge_{i<j<k}(u_i\ne u_j)\,\land\,
\bigwedge_{i<n}(u_i=v_i)&\\[\dxii]

\limp 
\bigwedge_{i<j<m}(w_i\ne w_j)\,\land\,
\bigwedge_{i<k}(w_i=v_i)\,\land\,\vpi^m_{M}(w_0,\dots,w_{m-1})
\skp
\eay
\dm
as $\vpi_{\neg\,M}^n(v_0,\dots,v_{n-1}).$ 
Finally, if $M=\bigcap_jM_j$ then we take 
$\bigwedge_j\vpi^n_{M_j}(v_0,\dots,v_{n-1})$ 
as $\vpi_M^n(v_0,\dots,v_{n-1})$. 
\epF{Theorem~\ref{lopes}}

\punk{\er s classifiable by countable structures}
\las{erclass}

The classifiability by countable structures means that 
we can associate, in a Borel way, a countable \dd\cL structure, 
say, $\vt(x)$ with any point $x\in\dX=\dom\rE$ so that $x\rE y$ 
iff $\vt(x)$ and $\vt(y)$ are isomorphic. 

\bdf
[{{\rm Hjorth~\cite[2.38]h}}]
\lam{classif}
An \er\ $\rE$ is {\it classifiable by countable structures\/} 
\index{equivalence relation, ER!classifiable by countable structures}%
if there is a countable relational language $\cL$ such that 
$\rE\reb{\ism\cL}$.
\edf

\bre
\lam{classifex}
Any $\rE$ classifiable by countable structures is $\fs11,$ 
of course, and many of them are Borel. 
The \eqr s $\rtd\yt{\Et},$ all countable Borel \er s
(see the diagram on page \pageref{p-p})
are classifiable by countable structures, but 
$\Ei\yt\Ed,$ Tsirelson \er s are not.
\ere

\bte
[{{\rm Becker and Kechris \cite{beke}}}]  
\lam{siy2cc}
Any orbit \er\ of a Polish action of a closed subgroup of 
$\isg$ is classifiable by countable structures. 
\ete

Thus all orbit \er s of Polish actions of $\isg$ and its 
closed subgroups are Borel reducible to a very special kind 
of actions of $\isg.$ 

\bpf
First show that any orbit \er\ of a Polish action of $\isg$ 
itself is classifiable by countable structures.  
Hjorth's simplified argument~\cite[6.19]h is as follows. 
Let $\dX$ be a Polish \dd\isg space with basis 
$\sis{U_l}{l\in\dN},$ and let $\cL$ be the language with 
relations $R_{lk}$ where each $R_{lk}$ has arity $k.$ 
If $x\in \dX$ then define $\vt(x)\in\mox\cL$ by 
stipulation that $\vt(x)\mo R_{lk}(s_0,\dots,s_{k-1})$ 
iff $1)\msur$ $s_i\ne s_j$ whenever $i<j<k,$ and 
$2)\msur$ $\kaz g\in B_s\:(g\obr\app x\in U_l),$ where 
\imar{Hjorth requires $\in\nad{U_l}.$ Why~? 
Also, it seems that $\forall^\ast g\in B_s$ extends the proof to 
Borel actions.}%
$B_s=\ens{g\in\isg}{s\su g}$ and 
$s=\ang{s_0,\dots,s_{k-1}}\in\dN^k.$  
Then $\vt$ reduces $\aer\isg\dX$ to $\ism\cL.$ 

To accomplish the proof of the theorem, it remains to 
apply the following result 
(an immediate corollary of Theorem~2.3.5b in \cite{beke}):

\bpro
\lam{subg2g}
If\/ $\dG$ is a closed subgroup of a Polish group\/ $\dH$ 
and\/ $\dX$ is a Polish\/ \dd\dG space then there is a  
Polish\/ \dd\dH space\/ $\dY$ such that\/ 
$\ergx\reb{\aer\dH\dY}$.
\epro
\bpf
Hjorth~\cite[7.18]h outlines a proof as follows. 
Let $Y=\dX\ti\dH\,;$ define $\ang{x,h}\approx\ang{x',h'}$ if 
$x'=g\app x$ and $h'=gh$ for some $g\in\dG,$ and consider the 
quotient space $\dY=Y/{\approx}$ with the topology induced by 
the Polish topology of $Y$ via the surjection 
$\ang{x,h}\mapsto\ek{\ang{x,h}}{\approx},$ on which $\dH$ 
acts by 
$h'\app\ek{\ang{x,h}}{\approx}=\ek{\ang{x,{hh'}\obr}}{\approx}.$ 
Obviously $\ergx\reb{\aer\dH\dY}$ via the map 
$x\mapsto\ek{\ang{x,1}}{\approx},$ hence, it remains to prove 
that $\dY$ is a Polish \dd\dH space, which is not really 
elementary --- we 
refer the reader to \cite[7.18]h or \cite[2.3.5b]{beke}. 
\epF{Proposition}

To bypass \ref{subg2g} in the proof of Theorem~\ref{siy2cc}, 
we can use a characterization of all closed subgroups of $\isg.$ 
Let $\cL$ be a language as above, and $x\in\mox\cL.$  
Define $\aut x=\ens{g\in\isg}{g\app x=x}:$  
the group of all automorphisms of $x.$ 

\bpro
[{{\rm see \cite[1.5]{beke}}}] 
\lam{closubg}
$G\sq\isg$ is a closed subgroup of\/ $\isg$ iff there is 
an\/ \dd\cL structure\/ $x\in\mox\cL$ of a countable 
language\/ $\cL,$ such that\/ $G=\aut x$. 
\epro 
\bpf
For the nontrivial direction, 
let $G$ be a closed subgroup of $\isg.$ 
For any $n\ge 1,$ let $I_n$ be the set of all \dd Gorbits in 
$\dN^n,$ \ie, equivalence classes of the \er\ 
$s\sim t$ iff $\sus g\in G\:(t=g\circ s),$ thus, $I_n$ is an 
at most countable subset of $\cP(\dN^n).$ 
Let $I=\bigcup_nI_n,$ and, for any $i\in I_n,$ let 
$R_i$ be an \dd nary relational symbol, and 
$\cL=\sis{R_i}{i\in I}.$  
Let $x\in\mox\cL$ be defined as follows: 
if $i\in I_n$ then $x\mo R_i(k_0,\dots,k_{n-1})$ iff 
$\ang{k_0,\dots,k_{n-1}}\in i.$ 
Then $G=\aut x,$ actually, if $G$ is not necessarily closed 
subgroup then $\aut x=\nad G$.
\epF{Proposition}

Now come back to Theorem~\ref{siy2cc}. 
The same argument as in the beginning of the proof shows that 
any orbit \er\ of a Polish action of $G,$ a closed 
subgroup of $\isg,$ is $\reb{\izm\cL G}$ for an appropriate 
countable language $\cL.$ 
Yet, by \ref{closubg}, $G=\aut{y_0}$ where $y_0\in\mox{\cL'}$ 
and $\cL'$ is a countable language disjoint from $\cL.$ 
The map $x\longmapsto\ang{x,y_0}$ witnesses that 
${\izm\cL G}\reb{\ism{\cL\cup\cL'}}$.
 
\epF{Theorem~\ref{siy2cc}}

\punk{Reduction to countable graphs}
\las{gra}

It could be expected that the more complicated a language 
$\cL$ is accordingly the more complicated
isomorphism \eqr\ $\ism\cL$ it 
produces. 
However this is not the case. 
Let $\cG$ be the language of (oriented binary) graphs, \ie, 
\index{language!g@$\cG$ of graphs}%
$\cG$ contains a single binary predicate, say $R(\cdot,\cdot)$.

\bte
\lam{grafs}
If\/ $\cL$ is a countable relational language then\/ 
${\ism\cL}\reb{\ism\cG}.$ 
Therefore, an \er\/ $\rE$ is classifiable by countable 
structures iff\/ $\rE\reb{\ism\cG}.$
In other words, a 
single binary relation can code structures of any countable 
language. 
\ete  

Becker and Kechris \cite[6.1.4]{beke} outline a proof based 
on coding in terms of lattices, unlike the following argument, 
yet it may in fact involve the same idea.

\bpf
Let $\hfn$ be the set of all hereditarily finite sets over the 
set $\dN$ considered as the set of atoms, and $\ve$ be the 
associated ``membership'' (any $n\in\dN$ has no 
\dd\ve elements, $\ans{0,1}$ is different from $2,$ \etc\/.). 
Let $\ihf$ be the $\hfn$ version of $\ism\cG,$ \ie, if 
$P\yi Q\sq\hfn^2$ then $P\ihf Q$ means that there is a 
bijection $b$ of $\hfn$ such that 
$Q=b\app P=\ens{\ang{b(s),b(t)}}{\ang{s,t}\in P}.$ 
Obviously $(\ism\cG)\eqb(\ihf),$ thus, we have to prove that 
${\ism\cL}\reb{\ihf}$ for any $\cL$.

An action of $\isg$ on $\hfn$ is defined as follows. 
If $g\in\isg$ then $g\circ n=g(n)$ for any $n\in\dN,$ 
and, by \dd\ve induction,  
$g\circ\ans{a_1,\dots,a_n}=\ans{g\circ a_1,\dots,g\circ a_n}$ 
for all $a_1,\dots,a_n\in\hfn.$
Clearly the map $a\mto g\circ a$ ($a\in\hfn$) 
is an \dd\ve isomorphism of $\hfn,$ for any fixed $g\in\isg.$ 

\ble
\lam{grafs1}
Suppose that\/ $X\yi Y\sq\hfn$ are\/ \dd\ve transitive subsets 
of\/ $\hfn,$ 
the sets\/ $\dN\dif X$ and\/ $\dN\dif Y$ are infinite, and\/ 
${\ve\res X}\ihf{\ve\res Y}.$ 
Then there is\/ $f\in\isg$ such that\/ 
$Y=f\circ X=\ens{f\circ s}{s\in X}$.
\ele 
\bpf
It follows from the assumption ${\ve\res X}\ism\hfn{\ve\res Y}$ 
that there is an \dd\ve iso\-mor\-phism $\pi:X\onto Y.$  
Easily $\pi\res(X\cap\dN)$ is a bijection of $X_0=X\cap\dN$ 
onto $Y_0=Y\cap\dN,$ hence, there is $f\in\isg$ such that 
$f\res X_0=\pi\res X_0,$ and then we have $f\circ s=\pi(s)$ 
for any $s\in X$.
\epF{Lemma}

Coming back to the proof of Theorem~\ref{grafs}, we first show 
that ${\ism{\cG(m)}}\reb{\ihf}$ for any $m\ge 3,$ where 
$\cG(m)$ is the language with a single \dd mary predicate.  
Note that $\ang{i_1,\dots,i_m}\in\hfn$ whenever  
$i_1,\dots,i_m\in\dN.$ 

Put $\vT(x)=\ens{\vt(s)}{s\in x}$ for every element
$x\in\mox{\cG(m)}=\cP(\dN^m),$ where 
$\vt(s)=\tce{\ans{\ang{2i_1,\dots,2i_m}}}$ for each 
$s=\ang{i_1,\dots,i_m}\in \dN^m,$ and finally, for $X\sq\hfn,$ 
$\tce X$ is the least \dd\ve transitive set $T\sq\hfn$ with 
$X\sq T.$  
It easily follows from Lemma~\ref{grafs} that 
$x\ism{\cG(m)}y$ iff 
${\ve\res\vT(x)}\ihf{\ve\res\vT(y)}.$ 
This ends the 
proof of ${\ism{\cG(m)}}\reb{\ihf}$. 

It remains to show that ${\ism{\cL'}}\reb{\ihf},$ where 
$\cL'$ is the language with infinitely many binary 
predicates. 
In this case $\mox{\cL'}=\cP(\dN^2)^\dN,$ so that we can 
assume that every $x\in\mox{\cL'}$ has the form 
$x=\sis{x_n}{n\ge1},$ with 
$x_n\sq(\dN\dif\ans0)^2$ for all $n.$ 
Let $\vT(x)=\ens{s_n(k,l)}{n\ge1\land\ang{k,l}\in x_n}$ for 
any such $x,$ where 
\dm
s_n(k,l)=\tce{
\ans{\{\dots\{\ang{k,l}\}\dots\}\,\yi 0}}
\,\yi \text{ with }\,n+2\,\text{ pairs of brackets }\,
\{\,\yi \}\,.
\dm
Then $\vT$ is a continuous reduction of $\ism{\cL'}$ to $\ihf$.
\epF{Theorem}

\punk{Borel countably classified \er s: reduction to $\rT_\xi$}
\las{2t}

Equivalence relations $\rT_\xi$ of \nrf{opeer} offer a 
perfect calibration tool for those Borel \er s which admit 
classification by countable structures. 
First of all, 

\bpro
\lam{taccs}
Every \eqr\/ $\rT_\xi$ admits classification by 
countable structures.
\epro
\bpf
$\rT_0,$ the equality on $\dN,$ is the orbit \er\ of the 
action of $\isg$ by $g\app x=x$ for all $g\yi x.$ 
The operation \ref{cdun} of \nrf{opeer} 
(countable disjoint union) 
easily preserves the property of being Borel 
reducible to an orbit \er\ of continuous action of $\isg.$ 

Now consider operation \ref{cp} of countable power. 
Suppose that a \er\ $\rE$ on a Polish space $\dX$ 
is Borel reducible to $\rF,$ the orbit relation of 
a continuous action of $\isg$ on some Polish $\dY.$ 
Let $D$ be the set of all points 
$x=\sis{x_k}{k\in\dN}\in\dX^\dN$ such that either 
$x_k\nE x_l$ whenewer $k\ne l,$ or there is $m$ such 
that $x_k\rE x_l$ iff $m$ divides $|k-l|.$
Then $\ei\reb{(\ei\res D)}$ 
(via a Borel map $\vt:\dX^\dN\to D$ such that 
$x\ei \vt(x)$ for all $x$). 
On the other hand, obviously $(\ei\res D)\reb \rF',$ 
where, for $y\yi y'\in\dY^\dN,$ $y\rF'y'$ means that there 
is $f\in\isg$ such that $y_k\rF y'_{f(k)}$ for all $k.$ 
Finally, $\rF'$ is the orbit \er\ of a 
continuous action of $\isg\ti\isg^\dN,$ which can be realized 
as a closed subgroup of $\isg,$ so it remains to apply 
Theorem~\ref{subg2g}. 
\epf  

The relations $\rT_\al$ are known in different versions, 
which reflect 
the same idea of coding sets of \dd\al th cumulative level 
over $\dN,$ as, \eg, in \cite[\pff1]{sinf}, where results 
similar to Proposition~\ref{taccs} are obtained in much 
more precise form.  

\bte
\lam{2tal}
If\/ $\rE$ is a Borel \er\ classifiable by countable 
structures then\/ $\rE\reb\rT_\xi$ for some\/ $\xi<\omi$. 
\ete
\bpf
The proof (a version of the proof in \cite{frid}) is based on 
Scott's analysis. 
Define, by induction on $\al<\omi,$ a family of Borel \er s
$\rrO\al$ on $\dN\lom\ti\cP(\dN^2):$ 
\bit
\item[$\mtho\ast$]\msur
$\rrq\al stAB$ means $\rro\al stAB$;
\eit  
thus, all $\rrt\al st$ ($s\yi t\in\dN\lom$) 
are binary relations on $\cP(\dN^2),$ and among them all 
relations $\rrt\al ss$ are \er s.
We define them by transfinite induction on $\al$.

\bit
\item\msur
$\rrq0stAB$ \ iff \ 
$A(s_i,s_j)\eqv B(t_i,t_j)$ for all $i\yi j<\lh s=\lh t$;

\item\msur
$\rrq{\al+1}stAB$ \ iff \ 
$\kaz k\:\sus l\:(\rrQ\al{s\we k}{t\we l}AB)$ 
and $\kaz l\:\sus k\:(\rrQ\al{s\we k}{t\we l}AB)$;

\item 
if $\la<\omi$ is limit then: \ 
$\rrq\la stAB$ \ iff \ $\rrq\al stAB$ for all $\al<\la$.
\eit
Easily ${\rrO\ba}\sq{\rrO\al}$ whenever $\al<\ba$. 
%

Recall that, for $A\yi B\sq\dN^2,$ $A\ism\cG B$ means that there 
is $f\in\isg$ with $A(k,l)\eqv B(f(k),f(l))$ for all $k\yi l.$  
Then we have ${\ism\cG}\sq\bigcap_{\al<\omi}{\rrt\al\La\La}$ 
by induction on $\al$ 
(in fact $=$ rather than $\sq,$ see below), 
where $\La$ is the empty sequence. 
Call a set $P\sq \cP(\dN^2)\ti\cP(\dN^2)$ {\it unbounded\/} if 
$P\cap{\rrt\al\La\La}\ne\pu$ for all $\al<\omi$. 

\ble
\lam{2tal1}
Any unbounded\/ $\fs11$ set\/ $P$ contains a pair\/ 
$\ang{A,B}\in P$ such that\/ $A\ism\cG B.$
\ele

It follows that $A\ism\cG B$ iff $\rrq\al\La\La AB$ for all 
$\al<\omi$ (take $P=\ans{\ang{A,B}}$). 

\bpf
Since $P$ is $\fs11,$ there is a continuous map $F:\dnn\onto P.$   
For $u\in\dN\lom,$ let $P_u=\ens{F(a)}{u\su a\in\dnn}.$ 
There is a number $n_0$ such that $P_{\ang{n_0}}$ is still 
unbounded. 
Let $k_0=0.$ 
By a simple cofinality argument, there is $l_0$ such that 
$P_{\ang{n_0}}$ is still unbounded 
{\it over\/} $\ang{k_0}\yi \ang{l_0}$ 
in the sense that there is no ordinal $\al<\omi$ such that 
$P_{\ang{i_0}}\cap{\rrt\al{\ang{k_0}}{\ang{l_0}}}=\pu.$ 
Following this idea, we can define infinite sequences of 
numbers $n_m\yi k_m\yi l_m$ such that both $\sis{k_m}{m\in\dN}$ 
and $\sis{l_m}{m\in\dN}$ are permutations of $\dN$ and, for 
any $m,$ the set $P_{\ang{n_0,\dots,n_m}}$ is still 
unbounded over $\ang{k_0,\dots,k_m}\yi \ang{l_0,\dots,l_m}$ in the 
same sense. 
Note that $a=\sis{n_m}{m\in\dN}\in\dN$ and 
$F(a)=\ang{A,B}\in P.$
(Both $A\yi B$ are subsets of $\dN^2.$) 

Prove that the map $f(k_m)=l_m$ witnesses $A\ism\cG B,$ \ie, 
$A(k_j,k_i)$ iff $B(l_j,l_i)$ for all $j\yi i.$ 
Take $m>\tmax\ans{j,i}$ big enough for the 
following: if $\ang{A',B'}\in P_{\ang{i_0,\dots,i_m}}$ 
then $A(k_j,k_i)$ iff $A'(k_j,k_i),$ and 
similarly $B(l_j,l_i)$ iff $B'(l_j,l_i).$ 
By the construction, there is a pair 
$\ang{A',B'}\in P_{\ang{i_0,\dots,i_m}}$ with 
$\rrq0{\ang{k_0,\dots,k_m}}{\ang{l_0,\dots,l_m}}{A'}{B'},$ 
in particular, $A'(k_j,k_i)$ iff $B'(l_j,l_i),$ as required.
\epF{Lemma}

\bcor
[{{\rm See, \eg, Friedman~\cite{frid}}}] 
\lam{2tal3} 
If\/ $\rE$ is a Borel\/ \er\ and\/ $\rE\reb{\ism\cG}$ then\/ 
$\rE\reb{\rrt\al\La\La}$ for some\/ $\al<\omi$.
\ecor
\bpf
Let $\vt$ be a Borel reduction of $\rE$ to ${\ism\cG}.$ 
Then 
$\ens{\ang{\vt(x),\vt(y)}}{x\nE y}$ is a $\fs11$ subset of 
$\cP(\dN^2)\ti\cP(\dN^2)$ which does not intersect $\ism\cG,$ 
hence, it is bounded by Lemma~\ref{2tal1}. 
Take an ordinal $\al<\omi$ which witnesses the boundedness.
\epf

Now, if $\rE$ is a Borel \er\ classifiable by countable 
structures then $\rE\reb{\ism\cG}$ by Theorem~\ref{grafs}, 
hence, it remains to establish the following: 

\bpro
\lam{r2t}
Any \er\/ $\rrO\al$ is Borel reducible to some\/ $\rT_\xi$.
\epro
\bpf 
We have ${\rrO0}\reb{\rT_0}$ since $\rrO0$ has countably 
many equivalence classes, all of which are clopen sets. 
To carry out the step $\al\mapsto \al+1$ note that the map   
$\ang{s,A}\mapsto \sis{\ang{s\we k,A}}{k\in\dN}$ is a Borel 
reduction of $\rrO{\al+1}$ to $(\rrO\al)^\iy.$ 
To carry out the limit step, let 
$\la=\ens{\al_n}{n\in\dN}$ be a limit ordinal, and 
${\rR}=\bigvee_{n\in\dN}{\rrO{\al_n}},$ \ie, $\rR$ is 
a \er\ on $\dN\ti\dN\lom\ti\cP(\dN^2)$ defined so that 
$\ang{m,s,A}\rR\ang{n,t,B}$ iff $m=n$ and $\rrq{\al_m} stAB.$ 
However the map $\ang{s,A}\mapsto\sis{\ang{m,s,A}}{m\in\dN}$ 
is a Borel reduction of $\rrO\la$ to $\rR^\iy.$
\epF{Proposition}

\epF{Theorem~\ref{2tal}}

\api

\parf{Turbulent group actions}
\las{groat}

This Section accomplishes the proof of irreducibility of the
\eqr s $\Ed$ and ${\fco}$ to $\rtd$ (see Subsection~\ref{brb}).
In fact it will be established that all relations in one family
of \eqr s are Borel irreducible to Borel relations in another
family.
The second family contains all Borel orbit \eqr s which admit
classification by countable structures
(see Section~\ref{groas}), and in fact 
many more equivalences, see below.
The first family consists of orbit equivalences induced by
\rit{turbulent} actions.

\punk{Local orbits and turbulence}
\las{tu:def}

Suppose that a group $\dG$ acts on a space $\dX.$ 
If $G\sq\dG$ and $X\sq\dX$ then let 
\dm
\rr XG=\ens{\ang{x,y}\in X^2}{\sus g\in G\:(x=g\ap y)}
\index{zzrgx@$\rr XG$}%
\dm 
and let $\sym XG$ denote the \er-hull of $\rr XG,$ \ie, 
\index{zzsymgx@$\sym XG$}%
the \dd\sq least \eqr\ on $X$ such that 
${x\rr XG y}\imp {x\sym XG y}.$ 
In particular ${\sym\dX\dG}={\ergx},$ but generally 
we have ${\sym XG}\sneq{\ergx\res X}.$  
Finally, define 
$\lo XGx=\ek x{\sym XG}=\ens{y\in X}{x\sym XG y}$ 
for $x\in X$ -- the {\it local orbit\/} of $x.$ 
\index{orbit!local}%
\index{zzouxx@$\lo XGx$}%
In particular, $\ek x\dG=\ek x{\ergx}=\lo\dX\dG x$ 
is the full \dd\dG orbit of a point $x\in\dX$.

\bdf
[{{\rm This particular version  taken from 
Kechris~\cite[\S~8]{apg}}}]  
\lam{df:durb}
Suppose that $\dX$ is a Polish space and $\dG$ is 
a Polish group acting on $\dX$ continuously. 
\ben
\tenu{(t\arabic{enumi})}
\itla{T1}
A point $x\in\dX$ is {\it turbulent\/} if for any 
\index{turbulent!point}%
non-empty open set $X\sq \dX$ containing $x$ and any nbhd  
$G\sq\dG$ (not necessarily a subgroup) of  
$\ong,$ the local orbit $\lo XGx$ is somewhere 
dense (that is, not a nowhere dense set) in $\dX$.

\itla{T2} 
An orbit $\ek x\dG$ is {\it turbulent\/} if $x$ is such 
\index{turbulent!orbit}%
\index{orbit!turbulent}%
(then all $y\in\ek x\dG$ are turbulent since this notion is
invariant \vrt\ homeomorphisms).

\itla{T3}
The action (of $\dG$ on $\dX$) is 
{\it generically}~\footnote
{\ In this research direction, ``generically'', or, in our 
abbreviation, ``gen.'' (property) intends to 
mean that (property) holds on a comeager domain.}, 
or {\it\gen turbulent\/} and $\dX$ is a 
{\it \gen turbulent\/} Polish \dd\dG space, if 
the union of all
dense (topologically), turbulent, and meager orbits 
$[x]_\dG$ is comeager.\qeD
\index{action!\gen turbulent}%
\index{space!\gen turbulent}%
\index{turbulent!\gen}%
\index{gengen@\gen= generically}%
\een
\eDf

Thus turbulence means that orbits, and even local orbits
of the action considered
behave rather chaotically in some exact sense.
According to the following theorem, this property is
incompatible with the classifiability by countable
structures.

\bte
[{{\rm Hjorth~\cite h}}] 
\lam{turb1}
Suppose that\/ $\dG$ is a Polish group, $\dX$ is a 
\gen turbulent Polish\/ \dd\dG space. 
Then\/ $\ergx$ is\/ {\ubf not} Baire measurable
reducible~\snos
{Reducible via a Baire measurable function. 
This is weaker than the Borel reducibility.}
to a Polish action of\/ $\isg,$ hence, {\ubf not} 
classifiable by countable structures.
\ete

The proof given below is based on general
ideas in \cite[\pff3.2]{h}, \cite[\pff12]{apg}, \cite{frid}. 
Yet it is designed so that only 
quite common tools of descriptive set theory are involved. 
It will also be shown that \lap{turbulent} \eqr s are 
not reducible actually to a much bigger family of relations
than orbit equivalences of Polish actions of $\isg$.

\punk{Shift actions of summable ideals are turbulent}
\las{t:prim}

Quite a lot of examples of turbulent actions is known
(see \eg\ \cite{h}).
The following example will be used in the proof of some
Borel irreducibility results in the end of this Section.
Recall that any summable ideal
$\srn=\ens{x\sq\dN}{\sum_{n\in x}r_n<+\iy}$
(where $r_n\ge 0$ for all $n$)
generates the \eqr\ $\ern=\rE_{\srn}$ on $\pn,$
defined so that $x\ern y$ iff ${x\sd y}\in \srn.$

\bte
\lam{sturb}
If\/ $r_n>0,$ $\ans{r_n}\to0,$ and\/ $\sum_nr_n=+\iy$
then the\/ \dd\sd action of\/ $\srn$ on\/ $\pn$
is Polish and\/ \gen turbulent.
\ete

The condition $\ans{r_n}\to0$ here implies that  
$\srn$ contains some infinite sets.
The condition $\sum_nr_n=+\iy$ means that $\srn$
does not contain co-infinite sets.

\bpf
Show that $\stk{\srn}{\sd}$ is a Polish group 
\imar{move to Sec 2 or 3?}
with the distance $\dpr(a,b)=\vpr(a\sd b),$ where
\dm
\text{
$\vpr(x)=\sum_{n\in x}r_n$ for $x\in\pn,$ 
hence $\srn=\ens{x}{\vpr(x)<\piy}$.}  
\dm
To prove that the operation is continuous, let
$x\yi y\in\pn.$
Fix a real $\da>0,$ and let $\ve=\frac\da2.$
If $x'\yi y'$ belong to the \dd{\ve}nbhds 
of $x\yi y$ in $\srn$ with the distance $\dpr,$ then 
${{(x'\sd y')}\sd{(x\sd y)}}\sq {{(x\sd x')}\cup{(y\sd y')}},$
therefore
$\dpr(x'\sd y',x\sd y)\le \dpr(x,x')+\dpr(y,y')=\da$.

Now prove that the \dd\sd action of $\srn$
on $\pn$ is continuous in the sense of the 
\dd\dpr topology of $\srn$
and the ordinary Polish product topology on $\pn.$
Suppose that $g\in\srn\yt x\in\pn,$ and fix a Polish nbhd
$V=\ens{y\in\pn}{y\cap n= (g\app x)\cap n}$ 
of $g\app x$ in $\pn,$ where $n\in\dN.$
Consider the corresponding nbhd
$U=\ens{x'\in\pn}{x'\cap n= x\cap n}$ of $x.$
Let $\ve=\tmin\ens{r_k}{k< n}.$
Then any element $g'\in\srn$ of the 
\dd\ve nbhd of $g$ in the \dd\dpr topology satisfies 
$g\sd g'\sq\ir{n}\iy,$ therefore $g'\sd x'\in V$
for any $x'\in U.$ 

Finally prove the turbulence of the action. 

Let $x\in\pn.$ 
That $\ek x\srn=\srn\sd x$ is
dense and meager is an easy exercise.
Thus it suffices to check that $x$ is turbulent. 
Consider an open nbhd
$X=\ens{y\in\pn}{y\cap [0,k)=u}$ of $x,$ where $k\in\dN$
and $u=x\cap[0,k),$
and a \dd\dpr nbhd
$G=\ens{g\in\srn}{\vpi(g)<\ve}$ of $\pu$
(the neutral element),
where $\ve>0.$ 
Prove that the local orbit $\lo XGx$
is somewhere dense in $X$. 

Let $l\ge k$ be large enough for $r_n<\ve$ to hold for
all $n\ge l.$
Prove that the orbit $\lo XGx$ is
dense in $Y=\ens{y\in\pn}{y\cap [0,l)=v},$ 
where $v=x\cap [0,l).$
Consider an open set 
$Z=\ens{z\in Y}{z\cap [l,j)=w},$
where $j\ge l\yt w\sq[l,j).$ 
Let $z$ be the only point of $Z$ satisfying 
$z\cap{[j,\piy)}=x\cap{[j,\piy)}.$
Thus
${x\sd z}=\ans{l_1,\dots,l_m}\sq[l,j).$ 
Note that every element of the form $g_i=\ans{l_i}$
belongs to $G$ by the choice of $l$ 
since $l_i\ge l.$  
Moreover, 
$x_i=g_i\sd g_{i-1}\sd\dots\sd g_1\sd x=
\ans{l_1,\dots,l_i}\sd x$ 
belongs to $X$ for each $i=1,\dots,m.$
On the other hand $x_m=z.$ 
It follows that $z\in \lo XGx,$ as required.
\epf

A suitable modification of this argument can be used to
prove the turbulence of the \dd\sd action of some other
ideals including the density ideal $\Zo,$ but as far as
some irreducibility results are concerned,
the turbulence of summable ideals will suffice!

\punk{Ergodicity}
\las{ergo}

The non-reducibility in Theorem~\ref{turb1} will be established 
in a special stronger form. 
Let $\rE\yi \rF$ be \er s on Polish spaces 
resp.\ $\dX\yi \dY.$ 
A map $\vt:\dX\to\dY$ is 

\bit
\item[$-$]\msur
\inw\rE\rF\ \ if \ 
\index{map!efinvariant@\ef invariant}%
${x\rE y}\imp {\vt(x)\rF\vt(y)}$ for all 
$x\yi y\in\dX$;\snos
{Recall that `\gen' means `generic' or `generically'.}

\item[$-$]
{\it \gen}\inw\rE\rF\ \ if \ the implication
${x\rE y}\imp {\vt(x)\rF\vt(y)}$ holds for all $x\yi y$ 
\index{map!genefinvariant@\gen\ef invariant}%
in a comeager subset of $\dX$;

\item[$-$]
{\it \gen reduction of\/ $\rE$ to\/ $\rF$}  \ if \
\index{reduction@reduction!generic}
the equivalence
${x\rE y}\eqv {\vt(x)\rF\vt(y)}$ holds for all $x\yi y$ 
in a comeager subset of $\dX$;

\vyk{
\item\msur 
\ddf{\it constant\/} if we have 
\index{map!fconstant@\ddf constant}%
${\vt(x)\rF\vt(y)}$ for all $x\yi y\in\dX$; 
}

\item[$-$] 
{\it \gen}\ddf constant \ if \ ${\vt(x)\rF\vt(y)}$ 
for all $x\yi y$ 
\index{map!genfconstant@\gen\ddf constant}%
in a comeager subset~of~$\dX.$
\eit 
Finally, following Hjorth and Kechris, say that $\rE$ is 
{\it \gen\ddf ergodic\/} if every Borel   
\index{equivalence relation, ER!genferodic@\gen\ddf ergodic}%
\gen\inw\rE\rF\ map is \gen\ddf constant.

The ergodicity preserves $\reb$ in the sense of the next lemma.

\ble
\lam{sogl*}
If\/ $\rE\yi\rF\yi\rF'$ are Borel \eqr s,
$\rE$ is\/ \gen\ddf ergodic, and\/ $\rF'\reb\rF,$
then\/ $\rE$ is\/ \gen\dd{\rF'}ergodic as well.
\ele
\bpf
Let $\vt$ be a Borel reduction of $\rF'$ to $\rF.$
Given a Borel \gen\/\inw\rE{\rF'} map $f,$ the map
$f'(x)=\vt(f(x))$ is obviously 
\gen\/\inw\rE{\rF}, hence it is a 
\gen\ddf constant --- then easily a \gen\dfp constant,
too.
\epf

The following lemma shows that ergodicity implies
irreducibility. 

\ble
\lam{t2nonr'*}
If an \eqr\/ $\rE$ is \gen\/\ddf ergodic and does not
have co-meager equivalence classes then\/ 
$\rE$ does not admit a Borel\/ \gen reduction to\/ $\rF.$
In addition $\rE$ does not admit a Baire measurable
reduction to\/ $\rF$. 
\ele
\bpf
Suppose towards the contrary that a Borel map
$\vt:\dX\to\dY$
(where $\dX\yi\dY$ are the domains of resp.\ $\rE\yi\rF$) 
is a \gen reduction of $\rE$ to $\rF,$ that is, $\vt$ is a
true reduction on a co-meager set $C\sq\dX.$ 
Then $\vt$ is a \gen\/\ddf constant by the ergodicity,
that is, there exists a co-meager set $C'\sq\dX$ such 
that $\vt(x)\rF\vt(x')$ for all $x\yi x'\in C'.$
Thye set $D=C\cap C'$ is co-meager as well, hence
there exist $x\yi x'\in D$ such that $x\nE x'.$ 
Then $\vt(x)\nF\vt(x')$ holds since $\vt$ is a reduction 
on $C.$
On the other hand, we know that $\vt(x)\rF\vt(x'),$
contradiction.

The additional result follows because it is known that 
any Baire measurable map is continuous on a co-meager set.
\epf

The proof of Theorem~\ref{turb1} consists of the next two 
lemmas.~\snos
{There are slightly different ways to the same goal. 
Hjorth~\cite[3.18]h proves outright and with different 
technique, that any \gen turbulent 
\eqr\ is \gen ergodic \vrt\ any Polish action of $\isg.$ 
Kechris~\cite[\pff12]{apg} proves that 1) any 
\gen  \dd{\rT_2}ergodic equivalence is \gen ergodic 
\vrt\ any 
Polish action of $\isg,$ and 2) any turbulent one is 
\gen \dd{\rT_2}ergodic.}

\ble
\lam{tul1}
If\/ $\dG$ is a Polish group, $\dX$ is a \gen turbulent
Polish\/ \dd\dG space, 
and\/ $\ergx$ is Baire measurable reducible to a Polish
action of\/ $\isg$
then\/ $\ergx$ admits a Borel\/ \gen reduction to an \eqr\  
of the form\/ $\rT_\xi$.
\ele

Saying it differently, any \eqr, Baire measurable reducible
to a Polish action of\/ $\isg,$
is Borel reducible to one of $\rT_\xi$ on a co-meager set.
Note that any \eqr, Borel reducible (in proper sense)
to one of $\rT_\xi,$ is Borel itself.
Yet this cannot be applied to $\ergx$ in the lemma,
since only a generic (on a co-meager set) reduction is
claimed.

 
\ble
\lam{tul2} 
Every \eqr\ induced by a \gen turbulent Polish action
of a Polish group 
is\/ \gen\dd{\rT_\xi}ergodic for all\/ $\xi$. 
\ele

\noi
{\bf Proof}
of Theorem~\ref{turb1} from lemmas \ref{tul1} and \ref{tul2}.
If $\ergx$ is Baire measurable reducible to a Polish action
of $\isg$ then ${\ergx}$ also is Borel \gen reducible to one
of $\rT_\xi$ by Lemma~\ref{tul1}.
On the other hand, ${\ergx}$ is \gen\dd{\rT_\xi}ergodic 
by Lemma~\ref{tul2}.
Thus ${\ergx}$ has a co-meager equivalence class by
Lemma~\ref{t2nonr'*}.
But this contradicts the assumption of \gen turbulence.\vom

\qeDD{Theorem~\ref{turb1} from lemmas \ref{tul1} and \ref{tul2}}%
\vtm\vom

The proof of the lemmas follows below in this Section.

\punk{\lap{Generic} reduction to $\rT_\xi$}
\las{lem3}

Here, we prove Lemma~\ref{tul1}. 
Suppose that $\dG$ is a Polish group, $\dX$ a
\gen turbulent Polish \dd\dG space. 
In particular, the set $W_0$ of all points $x\in\dX$
that belong to
dense turbulent orbits $\ek x G$ is comeager in $\dX.$
It follows that there exists a 
dense $\Gd$ set $W\sq W_0$.

Assume further that the orbit \eqr\ $\rE=\ergx$ is Baire
measurable reducible to a Polish action of $\isg.$
As the latter is Borel reducible to the isomorphism $\ism\cG$ 
of binary relations on $\dN$
according to Theorems \ref{siy2cc} and \ref{grafs},
$\rE$ itself admits a Baire measurable reduction
$\rho:\dX\to\cP(\dN^2)$ to $\ism\cG.$ 
The remainder of the argument borrows notation from the 
proof of Theorem~\ref{2tal}. 

There is a
dense $\Gd$ set $D_0\sq\dX$ such that the restricted map 
$\vt={\rho\res D_0}$ is continuous on $D_0.$ 
%
By definition, we have
\dm
{x\rE y}\imp{\vt(x)\ism\cG\vt(y)}
\quad\text{and}\quad
{x\nE y}\imp{\vt(x)\not\ism\cG\vt(y)}
\dm
for all $x\yi y\in D_0.$ 
We are mostly interested in the second implication,
and the aim is 
to find a 
dense $\Gd$ set $D\sq D_0$ such that, for some $\al<\omi$:
\bit
\item[$\mtho(\ast)$] 
${x\nE y}\imp{\nrq\al\La\La{\vt(x)}{\vt(y)}}$ \ 
holds for all $x\yi y\in D$.
\eit
Recall that $A\ism\cG B$ iff 
$\kaz\al<\omi\:\rrq\al\La\La AB.$
(See a remark after Lemma~\ref{2tal1}.)
It follows that
${x\rE y}\imp{\rrq\al\La\La{\vt(x)}{\vt(y)}}$
holds for all $x\yi y\in D_0.$
Thus $(\ast)$ implies that $\rE\res D$ is Borel reducible
to $\rrq\al\La\La{}{}.$
Now to end the proof of Lemma~\ref{tul1} apply
Proposition~\ref{r2t}.

To find an ordinal $\al$ satisfying $(\ast)$
we make use of {\ubf Cohen forcing}. 
Let us fix a countable transitive model $\mm$ of $\zhc,$ 
\index{theory!zfhc@$\zhc$}%
\index{zfhc@$\zhc$}%
\index{zzzfhc@$\zhc$}%
\ie, $\ZFC$ minus the Power Set axiom but plus the axiom: 
``every set belongs to $\HC
=\ens{x}{x\,\text{ is hereditarily countable}}$''. 

We shall assume that
$\dX$ is coded in $\mm$ in the sense that there is a  
set $D_\dX\in\mm$ which is a 
dense (countable outside of $\mm$)
subset of  $\dX,$ and $d_\dX\res D_\dX$ 
(the distance function of $\dX$ restricted to $D_\dX$) 
also belongs to $\mm.$ 
Further, $\dG,$ the action, $D_0,$ the map $\vt,$ and
the $\Gd$ set $W$ defined above ---  are also 
assumed to be coded in $\mm$ in a similar sense.

In these assumption, the notion of
a point of $\dX$ or of $\dG$
\rit{Cohen generic over $\mm$}
makes sense, 
and, as usual, the set $D$ of all Cohen generic,
over $\mm,$ points of $\dX$ is a
dense $\Gd$ subset of $\dX$ and $D\sq D_0.$ 
We are going to prove that $D$ fulfills $(\ast)$.

Suppose that $x\yi y\in D,$ and $\ang{x,y}$ is a Cohen generic  
pair over $\mm.$  
If $x\ergx y$ is false then   
$\vt(x)\not\ism\cG\vt(y).$ 
Moreover, this fact holds in the extended model
$\mm[x,y]$ by the Mostowski absoluteness.
This allows us to find, arguing in $\mm[x,y]$ 
(which is still a model of $\zhc$),
an ordinal $\al\in\Ord^\mm=\Ord^{\mm[x,y]}$ such that  
$\nrq\al\La\La {\vt(x)}{\vt(y)}.$ 
Moreover, since the Cohen forcing satisfies \ccc, 
there is an ordinal $\al\in\mm$ such that   
$\nrq\al\La\La {\vt(x)}{\vt(y)}$ holds for {\ubf every} 
Cohen generic, over $\mm,$ pair $\ang{x,y}\in D^2$ 
such that $x\ergx y$ is false.  
It remains to show that this also holds provided $x\yi y\in D$ 
(are generic separately, but)  
do not necessarily form a pair Cohen generic over $\mm.$
Now we prove

\ble
\lam{g**}
If\/ $\mn$ is a countable transitive model of\/ $\zhc$
with\/ $\mm\sq\mn,$ a point\/ $x\in \dX\cap\mn$ is 
Cohen generic over\/ $\mm,$ and an element\/ $g\in \dG$
is Cohen generic over\/ $\mn,$ then\/
$x'=g\ap x$ is Cohen generic over\/ $\mn$.   
\ele
\bpf
It follows from the genericity that $x$ belongs to the
set $W$ introduced in the beginning of Subsection~\ref{lem3}.
Thus the \dd\dG orbit $\ens{g'\app x}{g'\in G}$
is turbulent, in particular
dense in $\dX$.

Now consider any
dense open set $X\sq\dX$ coded in $\mn.$  
The set 
$H=\ens{g'\in\dG}{g'\app x\in X}$
is also open and coded in $\mn,$
Moreover $H$ is
dense in $\dG.$
(Indeed otherwise there is an open non-empty set $G\sq\dG$
such that the partial orbit $G\cdot x=\ens{g\app x}{g\in G}$
is nowhere dense.
This leads to a contradiction with the turbulence of $x.$) 
We conclude that $g\in H,$ and further $g\app x\in X,$
as required.
\epf

To make use of the lemma, let $\mn$ be a countable transitive
model of $\zhc$ containing $x\yi y,$ and all sets in $\mm.$
Note that $\mn$ may contain more ordinals than $\mm$ does
since the pair $\ang{x,y}$ is not assumed to be generic over
$\mm.$ 

Fix an element $g\in\dG$ Cohen generic over $\mn.$
Then $x'=g\ap x$ is Cohen generic over $\mn$ by the
lemma, hence over $\mm[y].$
Yet $y$ is generic over $\mm,$ 
thus the pair $\ang{x',y}$ is Cohen generic over $\mm.$
This implies $\nrq\al\La\La {\vt(x')}{\vt(y)}$
by the choice of  $\al.$
On the other hand we have $x'\ergx x$
and hence $\rrq\al\La\La {\vt(x)}{\vt(x')}.$
Thus we finally obtain 
$\nrq\al\La\La {\vt(x')}{\vt(y)},$ as required.\vtm

\qeDD{Lemma~\ref{tul1}}

\punk{Ergodicity of turbulent actions \vrt\ $\rT_\xi$}
\las{lem4}

Here, we prove Lemma~\ref{tul2}.  

We begin with two rather simple technical results
of topological nature, involved in the proof of the lemma. 

\addtocounter{theore}1

\ble
\lam{deno}
Suppose that\/ $\dG$ is a Polish group, and\/ $\dX$ is a 
\gen turbulent Polish\/ \dd\dG space.
Let\/ $\pu\ne X\sq\dX$ be an open set, $G\sq\dG$ be a nbhd 
of\/ $\ong,$ and\/ $\lo XGx$ be
dense in\/ $X$
for\/ \dd Xco\-meager many\/ $x\in X.$ 
Let\/ $U\yi U'\sq X$ be non-empty open and\/ $D\sq X$ 
be comeager in\/ $X.$ 
Then there exist points\/ $x\in D\cap U$ and\/ 
$x'\in D\cap U'$ with\/ $x\sym XG x'$.
\ele
\bpf
Under our assumptions there exist points $x_0\in U$ and 
$x'_0\in U'$ with $x_0\sym XG x'_0,$ that is, there exist
elements 
$g_1,\dots,g_n\in G\cup G\obr$ such that 
$x'_0=g_n\gru g_{n-1}\gru \dots\gru g_1\app x_0$ and in addition 
$g_k\gru \dots\gru g_1\app x_0\in X$ for all $k\le n.$ 
Since the action is continuous, there is a nbhd $U_0\sq U$ of 
$x_0$ such that $g_k\gru\dots\gru g_1\app x\in X$ for all $k$ 
and $g_n\gru g_{n-1}\gru\dots\gru g_1\app x\in U'$ for all  
$x\in U_0.$ 
Since $D$ is comeager, easily there is $x\in U_0\cap D$ such that  
$x'=g_n\gru g_{n-1}\gru\dots\gru g_1\app x\in U'\cap D$. 
\epF{Lemma}

\ble
\lam{kl}
Suppose that\/ $\dG$ is a Polish group, and\/ $\dX$ is a 
\gen turbulent Polish\/ \dd\dG space.
Then for any open non-empty\/ $U\sq\dX$ and\/ 
$G\sq\dG$ with\/ $\ong\in G$ there is an open non-empty\/ 
$U'\sq U$ such that the local orbit\/ $\lo{U'}Gx$ is 
dense in\/ $U'$ for\/ \dd{U'}comeager many\/
$x\in U'$.
\ele
\bpf
Let $\incl X$ be the interior of the closure of $X.$ 
If $x\in U$ and $\lo UGx$ is somewhere dense (in $U$) 
then the set $U_x=U\cap\incl{\lo UGx}\sq U$ is open and\/ 
\dd{\sym UG}invariant 
(an observation made, \eg, in \cite[proof of 8.4]{apg}), 
moreover, ${\lo UGx}\sq U_x,$ hence, ${\lo UGx}={\lo{U_x}Gx}.$ 
It follows from the invariance that the sets $U_x$ are 
pairwise disjoint, and it follows from the turbulence that 
the union of them is
dense in $U.$  
Take any non-empty $U_x$ as $U'.$ 
\epF{Lemma}

The proof of Lemma~\ref{tul2} involves a somewhat stronger 
property than \gen ergodicity in \nrf{ergo}. 
Suppose that $\rF$ is an \er\ on a Polish space $\dX$. 
\bit
\item
An action of $\dG$ on $\dX$ and the induced \eqr{} $\ergx$
are {\it locally generically\/} 
({\it\hgp,} for brevity) \ddf ergodic  
\index{equivalence relation, ER!hgenferodic@\hg\ddf ergodic}%
if the \eqr\ $\sym XG$ is generically \ddf ergodic  
whenever $X\sq \dX$ is a non-empty open set,  
$G\sq\dG$ is a non-empty open set containing $\ong,$ 
and the local orbit $\lo{X}Gx$ is
dense in $X$ for comeager (in $X$) many $x\in X$. 
\eit
This obviously implies \gen \ddf ergodicity of $\ergx$  
provided the action is \gen turbulent. 
Therefore, Lemma~\ref{tul2} is a corollary of the following 
theorem:

\bte
\lam{tut2}
Let\/ $\dX$ be a \gen turbulent Polish\/ \dd\dG space. 
Suppose that an \eqr\ $\rF$ belongs to\/ $\cF_0,$ 
the least collection of \eqr s containing\/ 
$\rav\dN$ 
{\rm(the equality on\/ $\dN$)} 
and closed under operations\/ \ref{oe:beg} -- \ref{oe:end} 
in Section~\ref{opeer}.  
\index{zzFo@$\cF_0$}%
\index{equivalence relation, ER!zzFo@$\cF_0,$ a family of \er s}%
Then\/ $\ergx$ is \hg \ddf ergodic, in particular,\/ $\ergx$ is 
not Borel reducible to\/ $\rF$.
\ete
%

Due to the operation of Fubini product, the family $\cF_0$ 
contains a lot of \eqr s very different from $\rT_\xi\,,$ among  
them some Borel equivalences which do not admit classification 
by countable structures, \eg, all \eqr s of the form $\rE_\cI,$ 
where $\cI$ is one of Fr\'echet ideals, indecomposable ideals, 
or Weiss ideals of \nrf{o}. 
(In fact it is not so easy to show that ideals of the two last 
families produce relations in $\cF_0$.) 
In particular, it follows that 
{\it no \gen turbulent \eqr\ is Borel reducible to a Fr\'echet, 
or indecomposable, or Weiss ideal\/}.

Our proof of Theorem~\ref{tut2} goes on by induction on the number 
of applications of the basic operations, in several following 
subsections. 

Right now, we begin with the initial step: prove that, under the 
assumptions of the theorem, $\ergx$ is \hg\dd{\rav\dN}ergodic. 
Suppose that $X\sq\dX$ and $G\sq\dG$ are non-empty open 
sets, $\ong\in G,$ and $\lo{X}Gx$ is
dense in $X$ for 
\dd Xcomeager many $x\in X,$ and prove that $\sym XG$ is 
generically \dd{\rav\dN}ergodic. 

Consider
a Borel 
\gen\inw{\sym XG}{\rav\dN}{} map $\vt:X\to\dN.$ 
Suppose towards the contrary that $\vt$ is not 
\gen\dd{\rav\dN}constant. 
Then there exist two open non-empty sets $U_1\yi U_2\sq X,$ 
two numbers $\ell_1\ne \ell_2,$ and a comeager set $D\sq X$ 
such that $\vt(x)=\ell_1$ for all $x\in D\cap U_1,$ 
$\vt(x)=\ell_2$ for all $x\in D\cap U_2,$ and $\vt\res D$ 
is ``strictly'' \inw{\sym XG}{\rav\dN}. 
Lemma~\ref{deno} yields a pair of points $x_1\in U_1\cap D$ and 
$x_2\in U_2\cap D$ with $x_1\sym XG x_2,$ contradiction. 

\punk{Inductive step of countable power}
\las{tut:iy}

To carry out this step in the proof of Theorem~\ref{tut2},
suppose that 

\bit
\item\msur 
$\dX$ is a \gen turbulent Polish \dd\dG space, 
$\rF$ is a Borel \er\ on a Polish space $\dY,$ 
and the action of $\dG$ on $\dX$ is \hg\ddf ergodic, 
\eit
and prove that the action is \hg\dd\rfy ergodic. 
Fix a nonempty open set $X_0\sq\dX$ and a nbhd $G_0$ of $\ong$ 
in $\dG,$ such that $\lo{X_0}{G_0}x$ is
dense in $X_0$ for all $x$ in a comeager \dd\Gd set
$D_0\sq X_0.$    
Consider  
a Borel function $\vt:X_0\to\dY^\dN,$ continuous and  
\inw{\sym{X_0}{G_0}}{\rfy} on $D_0,$ so that 
\dm
x\sym{X_0}{G_0} x'\;\limp\;
\kaz k\;\sus l\;\skl\vt_k(x)\rF\vt_l(x')\skp
\quad:\qquad\text{for all }\;x\yi x'\in D_0\,,
\dm
where $\vt_k(x)=\vt(x)(k),$ $\vt_k: X_0\to\dY,$ 
and prove that $\vt$ is \gen \dd\rfy constant. 

\vyk{
Below, let $\px$ be the Cohen forcing for $\dX,$ 
which consists of rational balls with centers in  
a fixed
dense countable subset of $\dX,$ and
let $\pg$ be the Cohen forcing for $\dG$ defined similarly 
(the
dense subset is assumed to be a subgroup). 
Smaller sets are stronger conditions.
}%

Let us fix a countable transitive model 
$\mm$ of $\zhc$ (see above),  
which contains all relevant objects or their codes, in 
particular, codes of the topologies of $\dX\yi \dG\yi \dY,$
of the set $D_0,$ and of the Borel map $\vt.$
Then every point $x\in X_0$ Cohen generic over $\mm$ belongs
to $D_0,$ hence $\vt$ is \inw{\sym{X_0}{G_0}}{\rfy} on
Cohen generic points of $X_0,$ and local orbits $\lo{X_0}{G_0}x$
of Cohen generic points $x\in X_0$ are dense in $X_0$.

\vyk{
\bct
\label{g*}
Suppose that\/ $\ang{x,g}\in \dX\ti\dG$ is\/ 
\dd{\px\ti\pg}generic over\/ $\mm.$  
Then\/ $g\ap x$ is\/ \dd\px generic over\/ $\mm$. \ 
{\rm(Because the action is continuous.)}\qeD
\ect
}

Coming back to the step of countable power, fix $k\in\dN.$ 
Consider any open non-empty set $U_0\sq X_0.$

\ble
\lam{41-}
There exist a number\/ $l$ and open non-empty sets\/
$U\sq U_0$ and\/ $H\sq G_0$ such that both\/
$g\app x\in X_0$ and\/ $\vt_k(x)\rF\vt_l(g\ap x)$
hold for any pair\/ $\ang{g,x}$ in\/ $H\ti U$
Cohen generic over\/ $\mm$.
\ele
\bpf
Consider any point $x_0\in U_0$ Cohen generic over $\mm.$
Note that $\ong\app x_0=x_0\in X_0,$ hence there exist  
a nbhd $U_1\sq U_0$ of $x_0$ and a nbhd
$G_1\sq G_0$ of $\ong$ such that $G_1\app U_1\sq U_0,$
\ie\ $g\app x\in X_0$ for all $g\in G_1$ and $x\in U_1$. 

Consider any pair $\ang{g,x}\in G_1\ti U_1$ Cohen generic
over $\mm.$
Then $g\app x\in U_0.$
In addition, $x$ is Cohen generic over $\mm$ while 
$g$ is Cohen generic over $\mm[x]$ by the forcing product
theorem.
It follows that $g\app x$ is Cohen generic over $\mm[x],$
and hence over $\mm,$ by Lemma~\ref{g**}. 

Furthermore, we have $x\sym{X_0}{G_0}{g\app x}.$   
By the invariance of $\vt$ on generic points this implies
$\vt(x)\rfi\vt(g\ap x).$
It follows that there is an index $l$ such that   
$\vt_k(x)\rF\vt_l(g\ap x).$
Thus there exist Cohen conditions, \ie\ non-empty open
sets $U\sq U_1$ and $H\sq G_1$ such that $x\in U,$
$g\in H,$ and any pair $\ang{g',x'}\in H\ti U$ Cohen generic
over $\mm$ satisfies $g'\app x'\in X_0$ and   
$\vt_k(x')\rF\vt_l(g'\ap x')$.
\epF{Lemma}

Fix $l\yi U\yi H$ as provided by the lemma.
Since $H$ is open, there is $h_0\in H\cap\mm$ and a 
symmetric nbhd $G\sq G_1$ of $\ong$ such that $g_0\,G\sq H$.

\ble
[{{\rm The key point of the turbulence}}]
\lam{41*}
If\/ $x\yi x'\in U$ are Cohen generic over\/ $\mm$  
and\/ $x\sym{U}{G}x'$ then we have\/ $\vt_k(x)\rF\vt_k(x')$.
\ele
\bpf 
We argue by induction on $n(x,x')=$ the least number $n$ 
such that there exist $g_1,\dots,g_n\in G$
(recall: $G=G\obr$) satisfying
\bit
\item[\mtho$(\ast)$]
$x'=g_n\, g_{n-1}\dots g_1\app x,$ \ and \
$g_k\dots g_1\app x\in U'$ \ for all \ $k\le n$.
\eit
Suppose that $n(x,x')=1,$ thus, $x=g\ap x'$ for some
$g\in G.$
Let $\mn$ be any countable transitive model of $\zhc$
containing $x\yi x'\yi g,$ and all sets in $\mm.$ 
Consider any element $h\in H$ Cohen generic over $\mn$
and close enough to $h_0$ for 
$h'=hg\obr$ to belong to $H.$ 
(Note that $h_0g\obr\in H$ by the choice of $G$.) 
Then $h$ is generic over $\mm[x],$ too, and hence
$\ang{h,x}\in H\ti U$ is Cohen generic over $\mm$ 
by the product forcing theorem.
It follows, by the choice of $H,$ that
$h\ap x\in X_0$ and $\vt_k(x)\rF\vt_l(h\ap x).$ 

Moreover, $h'=hg\obr$ also is \dd{\pg}generic over $\mn$
(because $g\in\mn$),   
so that $\vt_k(x')\rF\vt_l(h'\ap x')$ by the same argument. 
Finally $g'\ap x'=gh\obr\ap(h\ap x)={g\ap x},$
and hence $\vt_k(x')\rF\vt_k(x),$ as required. 

As for the inductive step, prove that $(\ast)$ holds for some 
$n\ge 2$ assuming that it holds for $n-1.$ 
Consider an element $g'_1\in G$ close enough to $g_1$ for
$g'_2= g_2\,g_1\,{g'_1}\obr$ to belong to $G$ 
and for $x^\ast=g'_1\app x$ to belong to $U,$ and Cohen generic
over a fixed transitive countable model $\mn$ of $\zhc$
containing $x\yi x'\yi g_1\yi g_2.$ 
Then as above $g'_2$ is Cohen generic over $\mn$ while 
$x^\ast$ is Cohen generic over $\mm,$ and obviously  
$n(x^\ast,x')\le n-1$ because 
$g'_2\app x^\ast=g_2\app g_1\app x.$
It remains to use the induction hypothesis. 
\epF{Lemma}

To summarize, we have shown that for any $k$ and any 
open $\pu\ne U_0\sq X_0$ there 
exist: an open non-empty set $U\sq U_0,$ and an open 
$G\sq G_0$ with $\ong\in G,$ such that the map 
$\vt_k$ is \inw{\sym{U}{G}}{\rF} on $U.$ 
We can also assume that the orbit $\lo{U}Gx$ is
dense in $U$ for \dd{U}comeager many
$x\in U$ by Lemma~\ref{kl}. 
Then, by the \hg\ddf ergodicity of the action considered,
$\vt_k$ is \gen\ddf constant on $U,$ that is, there 
exist a
comeager $\Gd$ set $D\sq U$ and a point $y\in \dY$
such  that $\vt_k(x)\rF y$ for all $x\in D.$ 

We conclude that there exist: an \dd{X_0}comeager set 
$D\sq X_0,$ 
and a countable set $Y=\ens{y_j}{j\in\dN}\sq \dY$ such that, 
for any $k$ and for any $x\in D$ there is $j$ with  
$\vt_k(x)\rF y_j.$ 
Put $\eta(x)=\bigcup_{k\in\dN}\ens{j}{\vt_k(x)\rF y_j}.$ 
Then, for any pair $x\yi x'\in D,$  we have 
$\vt(x)\rfy \vt(x')$ iff $\eta(x)=\eta(x'),$ so that, by 
the invariance of $\vt,$  
\dm
{x\sym{U_0}{G_0}x'}\,\limp \,{\eta(x)=\eta(x')}
\quad:\qquad\text{for all }\;x\yi x'\in D\,.
\eqno(\dag)
\dm
It remains to show that $\eta$ is a constant on a comeager 
subset of $D$.
\imar{check}

Suppose, on the contrary, that there exist two non-empty 
open sets $U_1\yi U_2\sq U_0,$ a number $j\in\dN,$ and a 
comeager set $D'\sq D$ such that 
$j\in\eta(x_1)$ and $j\nin\eta(x_2)$ for all 
$x_1\in D'\cap U_1$ and $x_2\in D'\cap U_2.$ 
Now Lemma~\ref{deno} yields a contradiction to $(\dag),$
as in the end of Subsection~\ref{lem4}.\vtm

\vyk{
Take any $x\in U_1,$ \dd\px generic over $\mm.$ 
Since comeager many orbits are dense, the orbit 
$\lo{U_0}{G_0}x$ intersects $U_2,$ in other words, there exist 
$g_1,\dots,g_n\in G_0$ such that all intermediate points 
$x_k=g_kg_{k-1}\dots g_1\ap x$ belong to $U_0$ and $x_n\in U_2.$ 
As $\mm\cap\dG$ is dense in $\dG,$ we can assume that all 
$g_i$ belong to $\mm,$ and subsequently all $x_i$ are 
\dd\px generic over $\mm,$ so that, in particular 
$y=x_n\in U_2$ is generic. 
It follows, by the choice of $U_1\yi U_2,$ that 
$j\in\eta(x)$ but $j\nin\eta(y).$ 
Yet we have $x\sym{U_0}{G_0} y$ by the construction, 
and $x\yi y\in D$ by the genericity. 
But this contradicts $(\ast),$ as required.\vtm
}

\qeDD{Inductive step of countable power in Theorem~\ref{tut2}}

\punk{Inductive step of the Fubini product}
\las{tut:Fu}

To carry out this step in the proof of Theorem~\ref{tut2},
suppose that 

\bit
\item\msur 
$\dX$ is a \gen turbulent Polish \dd\dG space, 
for any $k,$ $\rF_k$ be a Borel \er\ on a Polish space 
$\dY_k,$ the action of $\dG$ on $\dX$ is \hg\dd{\rF_k}ergodic 
for any $k,$ 
and $\rF=\fps{\ifi}{\rF_k}{k}$ is, accordingly, 
a Borel \er\ on $\dY=\prod_k\dY_k$, 
\eit
and prove that the action is \hg\dd\rF ergodic. 

Fix a nonempty open set $X_0\sq\dX,$ a nbhd $G_0$ of $\ong$ 
in $\dG,$ and a comeager $\Gd$ set $D_0\sq X_0$
such that all local orbits 
$\lo{X_0}{G_0}x$ with $x\in D_0$ are
dense in $X_0.$   
Consider a Borel function $\vt:U_0\to\dY,$ 
\inw{\sym{U_0}{G_0}}{\rF} on $D_0,$ \ie,%
\dm
x\sym{U_0}{G_0} y\;\limp\;
\sus k_0\;\kaz k\ge k_0\;\skl\vt_k(x)\rF_k\vt_k(y)\skp
\quad:\qquad\hbox{for all }\;x\yi y\in D_0\,,
\dm
where $\vt_k(x)=\vt(x)(k),$ and prove that $\vt$ is 
\gen\ddf constant.

Choose a countable transitive model $\mm$ of $\zhc$ as in
\ref{tut:iy}.

Consider an open non-empty set $U_0\sq X_0.$
Similarly to Lemma~\ref{41-},
there exist non-empty open sets $U\sq U_0$ 
and $H\sq G_0,$ and a number $k_0,$ such that
both $g\ap x\in X_0$ and 
$\vt_k(x)\rF_k\vt_k(g\ap x)$ 
hold for all indices $k\ge k_0$ and for all pairs
$\ang{g,x}\in H\ti U$ Cohen generic over $\mm.$ 

As $H$ is open, there exist an element $h_0\in H\cap\mm$ 
and a symmetric nbhd $G\sq G_0$ of $\ong$ such that 
$h_0\,G\sq H$.

\ble
\lam{42}
If\/ $k\ge k_0,$ points\/ 
$x\yi y\in U$ are\/ Cohen generic over\/ $\mm,$  
and\/ $x\sym{U}{G}y,$ then\/ $\vt_k(x)\rF_k\vt_k(y)$. \ 
{\rm(Similarly to Lemma~\ref{41*}.)}\qeD
\ele
\vyk{
\bpf
Similar to the proof of Lemma~\ref{41*}.
\epf
}

Thus, for any open non-empty set $U_0\sq X_0$ there exist: 
a number $k_0,$ an open non-empty $U\sq U_0,$ and a nbhd  
$G\sq G_0$ of $\ong,$ such that $\vt_k(x)$ 
is \gen\inw{\sym{U}{G}}{\rF_k} on $U$ 
for every $k\ge k_0.$ 
We can assume that \dd{U}comeager many orbits $\lo{U}Gx$ 
are
dense in $U,$ by Lemma~\ref{kl}. 
Now, by \dd{\rF_k}ergodicity, any $\vt_k$ with 
$k\ge k_0$ is \gen\dd{\rF_k}constant on such a set $U,$ 
hence, $\vt$ itself is \gen\dd{\rF}constant on $U$ because  
$\rF=\fps{\ifi}{\rF_k}{k}.$ 
It remains to show that these constants are \ddf equivalent  
to each other. 

Suppose, on the contrary, that there exist two non-empty 
open sets $U_1\yi U_2\sq U_0$ and a pair of $y_1\nF y_2$ 
in $\dY$ such that $\vt(x_1)\rF y_1$ and $\vt(x_2)\rF y_2$
for comeager many $x_1\in U_1$ and $x_2\in U_2.$ 
Contradiction follows as in the end of
Subsection~\ref{tut:iy}.\vtm

\vyk{
Take any $x\in U,$ \dd\px generic over $\mm.$ 
Exactly as in the end of \nrf P, there exists a 
\dd\px generic $x'\in U'$ with $x\sym{U'}{G} x'.$ 
Then we have $\vt(x)\rF y$ and $\vt(x')\rF y'$ by the 
genericity, and $\vt(x)\rF\vt(x')$ by the invariance of 
$\vt,$ which contradicts the assumption $y\nF y'$. 
}

\qeDD{Inductive step of Fubini product in Theorem~\ref{tut2}}

\punk{Other inductive steps}
\las{tut:oth}

Here, we accomplish the proof of Theorem~\ref{tut2}, by 
carrying out induction steps, related to operations 
\ref{cun}, \ref{cdun}, \ref{prod} of Subsection~\ref{opeer}.\vtm

{\sl Countable union\/}. 
Suppose that $\rF_1\yi \rF_2\yi \rF_3\yi \dots$ are Borel \eqr s   
on a Polish space\/ $\dY,$ and $\rF=\bigcup_k{\rF_k}$ is 
still an \eqr, and the Polish and \gen turbulent action of 
$\dG$ on $\dX$ is \hg\dd{\rF_k}ergodic for any $k..$ 
Prove that the action remains \hg\ddf ergodic. 

Fix a nonempty open set $X_0\sq\dX,$ a \come $\Gd$ set
$D_0\sq X_0,$ and a nbhd $G_0$ of $\ong$ in $\dG$ such that
all local orbits $\lo{X_0}{G_0}x$ with $x\in D_0$ are
dense in $U_0.$   
Consider a Borel function $\vt:U_0\to\dY,$  
\inw{\sym{U_0}{G_0}}{\rF} on $D_0.$  
It follows from the invariance that for any open
non-empty $U\sq U_0$ there exist: a number $k$ and open 
non-empty sets $U\sq U_0$ and $H\sq G_0$  
such that both $g\ap x\in X_0$ and
$\vt(x)\rF_k\vt(g\ap x)$ hold 
for any pair $\ang{x,g}\in U\ti H$
Cohen generic over a fixed countable transitive model $\mm$
of $\zhc$ chosen as above.
Further, there exist $h_0\in H\cap\mm$ 
and a symmetric nbhd $G\sq G_0$ of $\ong$ such that
$h_0\,G\sq H.$

Similarly to Lemmas~\ref{41*} and \ref{42},   
$\vt(x)\rF_k\vt(x')$ holds for any pair of  
elements $x\yi x'\in U$ Cohen generic over $\mm$ 
and satisfying $x\sym{U}{G}x'.$
It follows, by the ergodicity, that 
$\vt$ is \dd{\rF_k}constant, hence, \ddf constant,  
on a comeager subset of $U.$ 
It remains to show that these \ddf constants are 
\ddf equivalent to each other, which is demonstrated exactly 
as in the end of \nrf{tut:iy}.\vtm

{\sl Disjoint union\/}. 
Let\/ $\rF_k$ be Borel \er s on Polish spaces\/ 
$\dY_k,\msur$ $k=0,1,2,\dots\;.$ 
By definition, $\bigvee_k{\rF_k}=\bigcup_k{\rF'_k},$
where each $\rF'_k$ is a Borel \eqr\ defined on the space
$\dY=\bigcup_k{\ans k\ti\dY_k}$ as follows: 
$\ang{l,y}\rF'_k \ang{l',y'}$ iff either $l=l'$ and 
$y=y'$ or $l=l'=k$ and $y\rF_k y'$.\vtm

{\sl Countable product\/}. 
Let\/ $\rF_k$ be \eqr s on a Polish spaces $\dY_k.$ 
Then $\rF=\prod_k{\rF_k}$ is an \eqr\ on the space 
$\dY=\prod_k\dY_k.$ 
For any map $\vt:\dX\to\dY,$ to be \gen\inw{\rE}{\rF}  
(where $\rE$ is any equivalence on $\dX$) 
it is necessary and sufficient that every 
co-ordinate map $\vt_k(x)=\vt(x)(k)$ 
is \gen\inw{\rE}{\rF_k}. 
This allows to easily accomplish this induction step.\vtm

\qeDD
{Theorem~\ref{tut2}, Lemma~\ref{tul2}, Theorem~\ref{turb1}}

\punk{Applications to the shift action of ideals}
\las{turbI}

We are going to apply the results of this Section in order
to prove that \eqr s generated by many Borel ideals
(in particular almost all polishable ideals) 
are not Borel reducible to Borel actions of the permutation  
group $\isg,$ and hence not classifiable by countable
structures.
The difficult problem of verification of the turbulence
can fortunately be circumpassed by reference to   
theorems \ref{tut2} and \ref{sturb} (the turbulence of
summable ideals).  

Say that a Borel ideal\/ $\cZ\sq\pn$ is {\it special\/} 
if there is a sequence of reals\/ $r_n>0$ with\/ 
$\ans{r_n}\to0,$ such that\/ $\sui{r_n}\sq\cZ.$ 
{\it Nontrivial\/} in the next theorem means: containing 
\index{ideal!special}%
\index{ideal!nontrivial}%
no cofinite sets.
In the context of summable ideals the nontriviality means
simply that $\sum_n r_n=\piy$.

\bte
\lam{speid}
Suppose that\/ $\cZ$ is a nontrivial Borel special ideal, 
and\/ $\rF$ belongs to the family\/ $\cF_0$ of Theorem~\ref{tut2}. 
Then\/ $\rez$ is generically \ddf ergodic, hence, is 
not Borel reducible to\/ $\rF$.  
\ete
\bpf
The \lap{hence} statement follows because by the nontriviality 
all \dd\rez equivalence classes are meager subsets of $\pn$.
 
As $\cZ$ is special, let $\sis{r_k}{}\to 0$ be a sequence of 
positive reals such that $\sum_n r_n=\piy$ and $\srn\sq\cZ.$
Note that $x\ern y$ implies $x\rez y,$ and hence any  
\gen\inw\rez\rF{} map is \gen\inw\ern\rF{} as well
(on the same \come set).
Thus it suffices to prove that $\ern=\rE_{\srn}$ is 
\gen\ddf ergodic.

Recall that the shift action of $\srn$ on $\pn$ is Polish
and \gen turbulent by Theorem~\ref{sturb}.
Thus $\rE_{\srn}$ is \gen\ddf ergodic by Theorem~\ref{tut2},
as required. 
\epf

\vyk{
}

The next corollary returns us to the discussion in
the end of Subsection~\ref{brb}.

\bcor
\lam{-t2}
The equivalence relations\/ $\fco$ and\/ $\Ed$ are
not Borel reducible to any ideal\/ $\rF$ in 
the family\/ $\cF_0,$ in particular,  
are not Borel reducible to\/ $\rtd$.
\ecor
\bpf
According to lemmas \ref{co=d} and \ref{l1=s}, it suffices to
prove that the ideals $\zo$ (density $0$) and $\sui{1/n}$
are special.
(Their nontriviality is obvious.)
The ideal $\sui{1/n}$ is special by definition.
As for $\zo,$ it suffices to prove that $\sui{1/n}\sq\zo.$
Consider a set $x\sq\dN\yt x\nin\zo.$
There is a real $\ve>0$ such that
$\frac{\#(x\cap\ir0n)}n>2\ve$ for infinitely many numbers $n.$ 
One easily defines an increasing sequence
$n_0<n_1<n_2<\dots$ such that $n_{i+1}\ge 2n_i$ and 
$\frac{\#(x\cap\ir{n_i}{n_{i+1}})}{n_{i+1}-n_i}>\ve$ for all $i.$
Then
$\sum_{n\in x}\frac1n\ge\ve\sum_i\frac{n_{i+1}-n_i}{n_{i+1}}=
\piy,$
hence $x\nin\sui{1/n}$.
\epf

The next theorem shows that, with three exceptions, there
exist no polishable ideals Borel reducible to \eqr s in $\cF_0.$
(Note that $\cF_0$ contains various \eqr s of the form $\rei,$
generated by non-polishable ideals $\cI,$ for instance,
by Fr\'echet ideals.) 
Kechris~\cite{rig} proved a similar theorem, with the
assumption of reducibility to a relation in $\cF_0$ the
reducibility to a Borel action of $\isg$ is considered.
Recall that $\cI\cong\cJ$ means isomorphism via bijection 
between the ground sets of the ideals.
\vyk{
The following theorem 
An ideal $\cJ\sq\pn$ is a
\rit{trivial variation of\/ $\ifi$\/}
(the ideal of all finite subsets of $\dN$)
iff there is an infinite
and co-infinite set $d\sq\dN$ such that
$\cJ=\ens{x}{x\cap d\,\text{ is finite}}$.
}

\bte
\lam{rig2t}
If\/ $\cI$ is a nontrivial Borel polishable ideal on\/ $\dN,$
$\rF$ an \eqr\ in\/ $\cF_0,$ and\/ $\rei\reb\rF,$ then\/
$\cI$ is isomorphic to one of the following three ideals$:$
$\It\yt\ifi\yt\fsm\ifi\pn$.
\vyk{
either\/
$\cI\cong\It,$ or\/ $\cI\cong\ifi,$
or\/ $\cI$ is a trivial variation of\/ $\ifi.$
}
\ete

Note that in each of the three cases $\rei\reb\Et$ holds. 

\bpf
It follows from Theorem~\ref{sol} that $\cI=\Exh_\vpi$
for a \lsc\ submeasure $\vpi$ on $\dN.$ 
We can assume that $\vpi(x)\le 1$ for all $x\in\pn.$
(Otherwise put $\vpi'(x)=\tmin\ans{1,\vpi(x)}.$)
Consider the sets $U_n=\ens{k}{\vpi(\ans k)\le\frac1n}$
and $U_\iy=\ens k{\vpi(\ans k)=0}.$
Clearly $U_{n+1}\sq U_n$ and $\vpi(U_\iy)=0,$ therefore
$\vpi(x)=\vpi(x\dif U_\iy)$ for all $x\in\pn$.

We claim that $\tlim_{n\to\iy}\vpi(U_n)=0.$

Suppose towards the contrary that there exists $\ve>0$
such that $\vpi(U_n)>\ve$ for all $n.$
By definition for every $m$ there is $n\ge m$ satisfying
$U_n\sq\iry m\cup U_\iy$ --- 
then $\vpi(U_n\dif m)>\ve$ as well.
Moreover there exists $n'\ge n$ satisfying
$\vpi(U_n\cap{\ir m{n'}})>\ve.$ 
This leads to a sequence $n_1<n_2<n_3<\dots$ of numbers
and a sequence of finite sets  
$w_j\sq U_{n_j}\dif U_{n_{j+1}}$ such that $\vpi(w_j)>\ve.$ 
The sets $w_j$ are pairwise disjoint, hence every \lap{tail}
$W\cap\iry n$ of their union $W=\bigcup_jw_j$
includes at least one of $w_j$ as a subset.
It follows that $W\nin\cI=\Exh_\vpi.$
The ideal $\cJ=\cI\cap\pws W$ on $W$ is then nontrivial.
We also have $\ans{\vpi(\ans k)}_{k\in W}\to 0$
and $\sum_k\vpi(\ans k)=\piy$ since 
for any $n$ all but finite sets $w_l$ satisfy $w_l\sq W.$ 
Finally the equivalence 
${x\sd y\in \cI}\eqv{x\sd y\in \cJ}$
holds for all $x\yi y\sq W.$
It follows that $\rej\reb\rei$ by means of the identity map.

Since $\vpi$ is a \lsc\ submeasure, we have  
$\vpi(y)\le \sum_{k\in y}\vpi(\ans k)$
for all $y\sq\dN.$  
It follows that every set $x\sq W$ satisfying
$\sum_{k\in x}\vpi(\ans k)<+\iy$ belongs to $\cI,$
hence to $\cJ$ as well.
Thus $\cJ$ is isomorphic to a special ideal via a bijection
of $W$ onto $\dN.$ 
We conclude that $\rej,$ and hence $\rei,$ are Borel
irreducible to relations in the family $\cF_0$ by
Theorem~\ref{speid}, contradiction.

Thus $\vpi(U_n)\to 0.$ 
It follows that for any set $x\in\pn$ to belong to $\cI$
it is necessary and sufficient that   
$x\cap (U_n\dif U_{n+1})$ is finite for every $m.$
This observation allows us to accomplish the proof:
if the difference $U_n\dif U_{n+1}$ is infinite for
infinitely many indices $n$ then $\cI\cong\It$;
if there exist only finitely many 
infinite differences
$U_n\dif U_{n+1}$ and their union is co-finite in $\dN$
then $\cI\cong\ifi$;
and finally $\cI\cong\fsm\ifi\pn$ iff 
there exist only finitely many (but $>0$) infinite differences
$U_n\dif U_{n+1}$ but their union is co-infinite in $\dN$.
\epf

\bcor
\lam{t2i}
There is no Borel ideal\/ $\cI$ such that\/
$\rei\eqb\rtd$. 
\ecor
\bpf
Suppose towards the contrary that $\cI$ is such an ideal.
Then $\cI$ is polishable.
(Indeed otherwise $\Ei\reb\cI$ by Theorem~\ref{sol},
and hence $\Ei\reb\rtd.$
But this contradicts Theorem~\ref{e1pga} since $\rtd$
is easily Borel reducible to a Polish action.)
Thus $\rei\reb\Et$ by Theorem~\ref{rig2t}.
On the other hand, recall that recall that the ideal 
$\It=\ofi$ is a P-ideal (Example~\ref{exh:e}), hence it is
polishable by Theorem~\ref{sol}.
Thus $\rtd\not\reb\Et$ by Theorem~\ref{h+orb},
which is applicable in this case because the \dd\sd group
of $\It$ (basically, of any ideal) is abelian.
Therefore  $\rtd\not\reb\rei,$ as required.
\epf

The next application of Theorem~\ref{speid} is related
to the structure of ideals Borel reducible to $\Et.$
The result is similar to Theorem~\ref{rig1}.
We begin with the following irreducibility lemma:

\ble
\lam{0vs3}
$\Eo\rebs \Et.$
\Eqr s\/ $\Et$ and\/ $\Ei$ are\/ \dd\reb incomparable.
\Eqr s\/ $\Ed$ and\/ $\Ei$ are\/ \dd\reb incomparable
as well.
\ele
\bpf
It is quite obvious that $\Eo\reb\Et$ and $\Eo\reb\Ei.$  
Thus $\Eo\rebs \Et$ strictly since we have
$\Et\not\reb\Ei$ by Corollary~\ref{ne<e1}.
To prove $\Ei\not\reb\Et$ recall that the ideal 
$\It$ is polishable (see above).
Now $\Ei\not\reb\Et$ follows from Theorem~\ref{sol}.

The proof of the second claim is similar.
\epf

The following result of Kechris~\cite{rig}
should be compared with Theorem~\ref{rig1}.

\bcor
\lam{rig2'}
If\/ $\cI$ is a nontrivial Borel ideal on\/ $\dN$
and\/ $\rei\reb\Et$ then\/
$\cI$ is isomorphic to one of the following three ideals$:$
$\It\yt\ifi\yt\fsm\ifi\pn$.
\ecor
\bpf
We have $\Ei\not\reb\rei$ by Lemma~\ref{0vs3}.
Therefore $\cI$ is a polishable ideal by Theorem~\ref{sol}.
It remains to apply Theorem~\ref{rig2t}.
\epf

\newpage

\api

\parf{Ideal $\It$ and the \eqr\ $\Et$}
\las{It:id}

This Chapter is devoted to the ideal $\It$ 
and the corresponding \eqr\ $\Et.$
Recall that $\It$
(also denoted by or $\ofi$)
consists of all sets $x\sq\pnn$ such that all 
cross-sections $\seq xn=\ens{k}{\ang{n,k}\in x}$ are finite.
Accordingly  the relation $\Et=\rE_{\It}$ is defined on 
$\pnn$ by $x\Et y$ iff $x\sd y\in\It.$ 
But we'll rather consider $\Et$ as an equivalence on $\dnd$
defined 
\index{equivalence relation, ER!E3@$\Et$}%
\index{zzE3@$\Et$}%
so that $x\Et y$ iff $x(n)\Eo y(n)$ for all $n:$ 
here $x\yo y$ belong to $\dnd.$
More detaily, $x\Et y$ holds iff
\dm
\kaz n\;\sus k_0\;\kaz k\ge k_0\;(x(n,k)=y(n,k)).
\dm

The main goal of this Section will be the proof of
the following theorem of
Hjorth and Kechris \cite{hk:nd,hk:rd} known as the 6th
dichotomy theorem. 

\bte
\lam{6dih}
If\/ $\rE\reb\Et$ is a Borel \eqr\ then\/ \bfei\/ 
$\rE\reb\Eo$ \bfor\/ $\rE\eqb \Et$.
\ete

Thus similarly to $\Ei,$ the \er\ $\Et$ is an immediate
successor of $\Eo$ in a rather strong sense.
Let us mention an immediate corollary.

\bcor
\lam{etey}
$\Ey\not\reb\Et$.
\ecor
\bpf
If $\Ey\reb\Et$ then by Theorem~\ref{6dih} either $\Ey\reb\Eo$
or $\Ey\eqb\Et.$
The \bfei\ case contradicts Theorem~\ref{Eynon}.
To derive a contradiction from the \bfor\ case recall that 
$\Ey\reb{\bel\iy}$ by Theorem~\ref{eiy} but on the other hand
$\Et\not\reb{\bel\iy}$ by Lemma~\ref{e3}.
\epf

The proof of Theorem~\ref{6dih} employs the Gandy -- Harrington
topology in a manner rather similar to the proof of
Theorem~\ref{kelu}
(3rd dichotomy).
The scheme of the proof given here is designed on the base of
the proofs of Theorems 7.2 and 7.3 in \cite{hk:rd}.
The first of them contains a different dichotomy
while the second theorem contains a result
thal allows to derive Theorem~\ref{6dih} from the first theorem.
To present Theorem 7.2 in \cite{hk:rd}
recall that $\Ey$ is a \dd\reb largest countable Borel \eqr,
realized in the form of a certain equivalence on the Polish
space $2^{F_2},$ where $F_2$ is the free group with two generators.
(Here it is essential only that $F_2$ is a countable set.)
Let $\Eya$ denote the \eqr\ on $(2^{F_2})^\dN,$ defined so
that $x\Eya y$ iff $x(n)\Ey y(n)$ for all $n.$
Thus $\Eya$ is related to $\Ey$ just as $\Et$ to $\Eo.$
Theorem 7.2 in \cite{hk:rd} asserts that any Borel \eqr\
$\rE$ such that $\rE\reb\Eya$ satisfies either $\rE\reb\Ey$
or $\Et\reb \rE$.
                   
\vyk{
\punk{Ideals below $\It$}
\las{<it}

\ble
\lam{III}
$\ifi\rebs \It.$
$\It$ and\/ $\Ii$ are\/ \dd\reb incomparable.
\ele
\bpf
To see that $\ifi\rebs \It$ take 
$\vt(x)=\ens{\ang{n,0}}{n\in x}.$ 
That $\It\not\reb\Ii$ can be shown as follows: otherwise 
by Theorem~\ref{rig1} $\It$ would be isomorphic either 
to one of $\ifi,\msur$ $\Ii,$ or to a trivial variation 
of $\ifi,$ which can be easily shown to be not the case. 
To see that $\Ii\not\reb\It$ recall that $\It=\ofi$ is 
of the form $\Exh_\psi$ for a \lsc\ submeasure $\psi$ 
(Example~\ref{exh:e}) and apply Theorem~\ref{sol}. 
\epf

The following theorem is analogous to Theorem~\ref{rig1}, 
yet the method of its proof is absolutely different.

\bte[{{\rm Kechris~\cite{rig}}}]
\lam{rig2}
If\/ $\cI\reb\It$ is a Borel (nontrivial) ideal on\/ $\dN$ 
then either\/ $\cI\cong\It$ or\/ $\cI$ is a trivial 
variation of\/ $\ifi$. 
\ete
\bpf
First of all we make use of Theorem~\ref{sol}: 
$\Ii\not\reb \cI$ according to Lemma~\ref{III}, therefore, 
$\cI=\Exh_\vpi$ for a \lsc\ submeasure $\vpi$ on $\dN.$ 
We can \noo\ suppose that $\vpi(x)\le 1$ for any $x\in\pn.$ 
Now put $U_n=\ens{k}{\vpi(\ans k)\le\frac1n}$.

We assert that $\tlim_{n\to\iy}\vpi(U_n)=0.$ 
Indeed, otherwise $\vpi(U_n)>\ve$ for some $\ve>0$ and all $n.$  
As $\vpi$ is \lsc\ we can choose a sequence 
of numbers $n_1<n_2<n_3<\dots$ and for any $l$ a finite set 
$w_l\sq U_{n_l}\dif U_{n_{l+1}}$ with $\vpi(w_l)>\ve.$ 
Then $W=\bigcup_lw_l\nin\cI$ and obviously 
$\ans{\vpi(\ans k)}_{k\in W}\to 0.$ 
Note that the Borel ideal $\cZ=\cI\res W$ satisfies 
$\cZ\reb\cI$ (via the identity map), because $W\nin\cI.$ 
On the other hand, $\cZ$ is isomorphic to a special ideal 
(see \nrf{turbI})
via the order preserving bijection of $W$ onto $\dN.$ 
It follows from Theorem~\ref{speid} that $\rez$ is not Borel 
reducible to any \eqr\ in $\cF_0,$ hence, neither is $\rei.$ 
But $\rE_{\It}=\Et$ obviously belongs to $\cF_0,$ which is 
a contradiction because $\cI\reb\It$.

Thus $\vpi(U_n)\to 0.$ 
Then clearly a set $x\in\pn$ belongs to $\cI$ iff 
$x\cap (U_n\dif U_{n+1})$ is finite for any $m,$ which 
easily implies that $\cI$ is as required. 
\imar{check the proof}
\epf
}

\punk{Continual assembling of equivalence relations}
\las{asser}

The next theorem  
will be used in the proof of Theorem~\ref{6dih}.
The result is somewhat similar to Theorem~\ref{cudf} in that
it evaluates the type of an \eqr\ $\rE$ on the base of the
types of certain fragments of $\rE.$
But in this case the number of fragments can be continual.

\bte
\lam{t:ass} 
Suppose that\/ $\dX\yi \dY$ are 
Polish spaces, $P\sq\dX\ti\dY$ is a Borel set, 
$\rE$ is a Borel \eqr\ on\/ $P,$ and\/ 
$\dG$ is a countable group acting on\/ $\dX$ in a Borel 
way, and\/ $\ang{x,y}\rE\ang{x',y'}$ implies\/ 
$x\ergx x'.$ 

Finally, assume that\/ $\rE\res{P(x)}$ is smooth for 
each\/ $x\in \dX,$ where\/ $P(x)=\ens{\ang{x',y}\in P}{x'=x}.$
Then\/ $\rE$ is Borel-reducible to a Borel action of\/ 
$\dG$.
\ete
\bpf
We can assume that $\dX=\dY=\dn$ and both $P$ and $\rE$ 
are $\id11.$ 
We can also assume that the action of $\dG$ 
(a countable group) is $\id11.$ 
Then clearly ${x\ergx x'}\imp{\id11(x)=\id11(x')}.$ 
Define 
$\tP(x)=\bigcup_{a\in\dG}P(a\app x)$ for $x\in\dX$.

\bcy
\lam{ass1}
Suppose that pairs\/ $\ang{x,y}$ and\/ $\ang{x',y'}$ belong 
to\/ $P$ and\/ $x\ergx x'.$ 
Then\/ $\ang{x,y}\rE\ang{x',y'}$ iff the equivalence\/ 
$\ang{x,y}\in U\eqv\ang{x',y'}\in U$ holds for any\/  
\dd{\rE\res{\tP(x)}}-invariant $\id11(x)$ 
set\/ $U\sq\tP(x)$.
\ecy
\bqf
Note that $\rE\res{\tP(x)}$ is still smooth by 
Corollary~\ref{cudc} because $\dG$ is countable. 
In addition $\rE\res{\tP(x)}$ is $\id11(x).$ 
This observation yields the result.
Indeed otherwise the \eqr, defined om $\tP(x)$ by intersections 
with \dd{\rE\res{\tP(x)}}invariant $\id11(x)$ sets, 
is coarser than $\rE\res{\tP(x)}.$
It follows
(see the proof of the 2nd dichotomy theorem, Theorem~\ref{2dih}) 
that $\Eo\reb{\rE\res{\tP(x)}},$ 
a contradiction with the smoothness.
\eqF{Claim}

In the continuation of the proof of Theorem~\ref{t:ass}
we make use of 
a standard enumeration of $\id11$ sets.
It follows from Theorem~\ref{penur} that
there exist $\ip11$ sets $C\sq\dX\ti\dN$ and  
$W\sq\dX\ti\dN\ti\dX\ti\dY$ and a $\is11$ set  
$W'\sq\dX\ti\dN\ti\dX\ti\dY$ such that the sets 
\dm
W_{xe}= \ens{\ang{x',y'}}{\ang{x,e,x',y'}\in W}
\;\;\text{ and }\;\;
W'_{xe}= \ens{\ang{x',y'}}{\ang{x,e,x',y'}\in W'}
\dm
coincide whenever $\ang{x,e}\in C,$ and  
for any $x\in\dX$ a set $R\sq\dX\ti\dY$ is\/ $\id11(x)$ iff
there is $e\in C_x=\ens{e}{\ang{x,e}\in C}$
such that $\ang{x,e}\in C$ and $X=W_{xe}=W'_{xe}.$

Let $\inva(x,e)$ be the formula 
\dm
x\in\dX\lland e\in C_x\lland W_{xe}\sq\tP(x)\lland
\text{$W_{xe}$ is \dd{\rE\res{\tP(x)}}invariant}\,. 
\dm

\bcoy
\lam{ass3}
Let\/ $\ang{x,y}\yi \ang{x',y'}$ be as in 
Claim~\ref{ass1}. 
Then\/ $\ang{x,y}\rE\ang{x',y'}$ iff the equivalence\/ 
$\ang{x,y}\in W_{xe}\eqv\ang{x',y'}\in W_{xe}$ 
holds for any\/ $e$ with\/ $\inva(x,e)$.\qfD
\ecoy

Let us change ``iff'' here to $\nmp.$
Such a reduced claim can be formally represented in the
form
$(P\ti P)\cap{\ergx}\sq U\cap{\rE},$
where $U=\bigcup_{e\in\dN}U_e$ and\pagebreak[0] 
\dm
\textstyle
U_e=\ens{\ang{\ang{x,y},\ang{x',y'}}}
{\ang{x,e}\in J\land\neg\:
(\ang{x,y}\in W_{xe}\eqv\ang{x',y'}\in W_{xe})}.
\dm
As $J\sq C,$ we can re-write the negation of $\eqv$ in
the last formula as follows:
\dm
\big(
\ang{x,y}\in W_{xe} \land \ang{x',y'}\nin W'_{xe}
\big)
\land
\big(
\ang{x,y}\nin W'_{xe} \land \ang{x',y'}\in W_{xe}
\big).
\dm
Thus the inclusion $(P\ti P)\cap{\ergx}\sq U\cap{\rE}$ as a
property of a $\ip11$ set $J$ is $\ip11$ in the codes.
It follows by Theorem~\ref{refl} (\Refl) that there is a
$\id11$ set $J'\sq J$ such that $(P\ti P)\cap{\ergx}\sq U'\cap{\rE}$
holds, where $U'$ is defined in terms of $J'$ similarly to the
definition of $U$ in terms of $J.$

\vyk{
Implication $\nmp$ of the ``iff'' in this Corollary can 
be considered as a property of the $\ip11$ set 
$J=\ens{\ang{x,e}}{\inva(x,e)},$ \ie, the property that 

\bit 
\item 
for all pairs $\ang{x,y}$ and $\ang{x',y'}$ in $P$ 
with $x\ergx x':$  
if the equivalence  
$\ang{x,y}\in W_{xe}\eqv\ang{x',y'}\in W_{xe}$ holds
for every $e\in J_x=\ens{e}{\ang{x,e}\in J}$ 
then $\ang{x,y}\rE\ang{x',y'}$.
\eit
This is easily a $\ip11$ property, hence, 
by the $\ip11$ Reflection, there is a $\id11$ set 
$B\sq J$ satisfying the same property, that is, we have

For any $x\in\dX$ let $E(x)$ be the set of all $e\in\dN$ 
which code a $\id11(x)$ subset of $P,$ and, for 
$e\in E(x),$ let $\dc ex$ be the $\id11(x)$ subset of $P$ 
coded by $e,$ in the sense of Theorem~\ref{penur}.
(It is known that $\ens{\ang{x,e}}{e\in E(x)}$ is 
$\ip11.$) \ 
}

\vyk{
Unfortunately the formula $\inva$ is a $\ip11$ formula, 
which makes it difficult to use Corollary~\ref{ass3} 
directly if one wants to define a Borel reduction. 
Yet there is a typical trick which allows to extract a 
sufficient $\id11$ part. 
Note that  
\dm
Q=\ens{\ang{x,y,x',y'}}{\ang{x,y}\in P\lland
\ang{x',y'}\in P\lland {x\ergx x'} \lland 
{\ang{x,y}\nE\ang{x',y'}}} 
\dm
is a $\id11$ subset of $P\ti P$ while 
\dm
S=\ens{\ang{x,y,x',y',e}\in Q\ti\dN}
{\inva(x,e)\lland \neg\:
(\ang{x,y}\in \dc ex\eqv\ang{x',y'}\in \dc ex)}
\dm
is a $\ip11$ set whose projection on $Q$ coincides with 
$Q.$ 
Then there is a $\id11$ map $\ve:Q\to\dN$ 
such that $\ang{x,y,x',y',\ve(x,y,x',y')}\in S$ holds 
whenever $\ang{x,y,x',y'}\in Q.$ 
Now
$A=\ens{\ang{x,\ve(x,y,x',y')}}{\ang{x,y,x',y'}\in Q}$ 
is a $\is11$ subset of the $\ip11$ set 
$C=\ens{\ang{x,e}}{\inva(x,e)}\sq\dX\ti\dN,$ 
hence, there is a $\id11$ set $B$ with $A\sq B\sq C$. 
}

\bcoy
\lam{ass4}
Let\/ $\ang{x,y}\yi \ang{x',y'}$ be as in 
Claim~\ref{ass1}. 
Then\/ $\ang{x,y}\rE\ang{x',y'}$ iff the equivalence\/ 
$\ang{x,y}\in W_{xe}\eqv\ang{x',y'}\in W_{xe}$ 
holds for any\/ $e$ with\/ $\ang{x,e}\in J'$.\qfD
\ecoy

To continue the proof of the theorem, define, for 
any $\ang{x,y}\in P$, 
\dm
D_{xy}=\ens{\ang{a,e}}{a\in \dG\lland 
\ang{a\app x,e}\in J'
\lland \ang{x,y}\in W_{a\app x,e}}\,.
\dm
Clearly $\ang{x,y}\mapsto D_{x,y}$ is a $\id11$ map 
$P\to\pws{\dG\ti\dN}$.
 
If $D\sq\dG\ti\dN$ and $b\in\dG$ then put 
$b\circ D=\ens{\ang{ab\obr,e}}{\ang{a,e}\in D}$.

\bcy
\lam{ass2}
Suppose that\/ $\ang{x,y}$ and\/ $\ang{x',y'}$ belong 
to\/ $P,$ $b\in\dG,$ and\/ $x'=b\app x.$ 
Then\/ $\ang{x,y}\rE\ang{x',y'}$ iff\/ 
$b\circ D_{xy}=D_{x'y'}$.
\ecy
\bqf
Assume that $b\circ D_{xy}=D_{x'y'}.$ 
According to Corollary~\ref{ass4}, to establish  
$\ang{x,y}\rE\ang{x',y'}$ it suffices to prove that 
${\ang{x,y}\in W_{x e}}\eqv{\ang{x',y'}\in W_{x e}}$ 
holds whenever $\ang{x,e}\in J'.$ 
We have 
\dm
{\ang{x,y}\in W_{x e}}\seqv{\ang{\La,e}\in D_{xy}} 
\seqv{\ang{b\obr,e}\in D_{x'y'}}\seqv
{\ang{x',y'}\in W_{b\obr\app x',\,e}=W_{x e}}\,,
\dm
as required. 
Conversely, let $\ang{x,y}\rE\ang{x',y'}.$ 
If $\ang{a,e}\in D_{xy}$ then $\ang{a\app x,e}\in J'$ 
and $\ang{x,y}\in W_{a\app x,\,e},$ hence, 
$\ang{x',y'}\in W_{a\app x,\,e},$ too, because 
the set $W_{a\app x,\,e}$ is invariant and 
$\ang{x,y}\rE\ang{x',y'}.$  
Yet $a\app x=ab\obr\app x',$ therefore, 
by definition, $\ang{ab\obr,e}\in D_{x'y'}.$ 
The same argument can be carried out in the opposite 
direction, so that $\ang{a,e}\in D_{xy}$ 
iff $\ang{ab\obr,e}\in D_{x'y'},$ that means 
$b\circ D_{xy}=D_{x'y'}$.
\eqF{Claim}

To end the proof of the theorem, consider 
$\dS=\dX\ti\pws{\dG\ti\dN},$ a Polish space.  
Define a Borel action   
$b\app\ang{x,D}=\ang{b\app x, b\circ D}$ 
of $\dG$ on $\dS.$  
We assert that $\vt(x,y)=\ang{x,D_{xy}}$ 
is a Borel reduction of $\rE\res P$ to 
the action $\aer\dG\dS.$ 
Indeed, let $\ang{x,y}$ and $\ang{x',y'}$ belong 
to $P.$ 
Suppose that $\ang{x,y}\rE\ang{x',y'}.$ 
Then $x\aer\dG\dX x',$ so that $x'=b\app x$ for some 
$b\in \dG.$ 
Moreover, $b\circ D_{xy}=D_{x'y'}$ by Claim~\ref{ass2}, 
hence, $\vt(x',y')=b\app\vt(x,y).$  
Let, conversely, $\vt(x',y')=b\app\vt(x,y),$ so that 
$x'=b\app x$ and $D_{x'y'}=b\circ D_{xy}.$ 
Then $\ang{x,y}\rE\ang{x',y'}$ by Claim~\ref{ass2}, 
as required. 
\epf

\punk{The two cases}
\las{Et:er}

Here we begin the proof of Theorem~\ref{6dih}.

We may assume that $\rE$ is a $\id11$ \eqr\
on the Cantor space $\dn,$ and there is 
a $\id11$ reduction $\vt:\dn\to\dnd$ of $\rE$ to $\Et.$
In this case, it can be \noo\ assumed that in fact $\vt$
is a $\id11$ bijection.
Indeed,
define $\vpi:\dX\to\dnd$ so that for any $x\in\dn:$
$\vpi(x)(n)(0)=x(n)$ for all $n$ and
$\vpi(x)(n)(k+1)=\vt(x)(n)(k)$ for all $n$ and $k.$
The map $\vpi$ is a bijection and still a Borel reduction
of $\rE$ to $\Et$.

Define $R=\ran \vt,$ a $\id11$ subset of $\dnd.$
(That $R$ is $\id11$ follows from the assumption that $\vt$
is a Borel bijection.)
 
For $x\yi y\in\dnd$ and $n\in\dN,$ define 
$x\eqn n y$ iff $x\Et y$ and $x\rnq n=y\rnq n$
\index{zzequivn@$\eqn n$}%
(the latter requirement means $x(k)=y(k)$ for all $k\le n$).
Put
\dm
\ah{n}{kp}=\ens{A\sq\dnd}{A\,\hbox{ is }\,\is11\lland
\index{zzAnkp@$\ah n{kp}$}%
\kaz x,y\in A\;
\skl{x\eqn n y}\imp
{x(k)\ap y(k)\sq\ir0p}
\skp}
\dm
for all $n,k,p\in\dN,$\snos
{Hjorth and Kechris~\cite{hk:rd} define $\ah{n}{kp}$ 
with $\kaz x,y\in R\cap A$ instead of $\kaz x,y\in A.$ 
Let us use $\cA'_{nkp}$ to denote their 
version, thus, $\ah{n}{kp}\sq\cA'_{nkp}.$ 
However if Case 1 holds in the sense of $\cA'_{nkp}$ 
then it also holds in the sense of $\ah{n}{kp}$ because 
$A\in \cA'_{nkp}$ iff $A\cap R\in \ah{n}{kp}$.}
where $a\ap b\in\dn$ is defined for any pair of
\index{zzadotb@$a\ap b$}%
$a\yi b\in\dn$ so that $(a\ap b)(k)=0$ whenever $a(k)=b(k)$
and $(a\ap b)(k)=1$ otherwise --- for all $k\in\dN.$
Thus for a $\is11$ set $A\sq\dnd$ to belong to
$\bigcup_p\ah{n}{kp}$
it is necessary and sufficient that for any $x\in A$ the set
$
\ens{y(k)}{y\in A\land y\eqn n x}
$
is finite.
To strengthen here finiteness to being a singleton with
the help of something like Theorem~\ref{cenu} (\Cenu) is
hardly possible.

\ble
\lam{utv1}
If\/ $A\in\ah{n}{kp}$ then there is a\/ $\id11$ set\/ 
$B\in\ah{n}{kp}$ with\/ $A\sq B$.
\ele
\bpf
The definition of $A\in\ah{n}{kp}$ as a property of a $\is11$
set $A$ is obviously $\ip11$ in the codes.
Therefore Theorem~\ref{refl} (\Refl) implies the result
required.
There is a more pedestrian proof based on the \Sepa\ rather
than \Refl\ theorem.
First consider the $\ip11$ set
\dm
P=\ens{y\in\dnd}{\kaz x\in A\;
\skl{x\eqn n y}\imp
{x(k)\ap y(k)\sq\ir0p}
\skp}\,.
\dm
Then $A\sq P$ since $A\in \ah{n}{kp}.$
Take a $\id11$ set $D$ such that $A\sq D\sq P.$
Now consider the $\ip11$ set
\dm
P'=\ens{x\in D}{\kaz y\in D\;
\skl{x\eqn n y}\imp
{x(k)\ap y(k)\sq\ir0p}
\skp}\,.
\dm
Then $A\sq P'$ since $A\sq D\sq P.$
Any $\id11$ set $B$ such that $A\sq B\sq P$
is as required.
\epf

\bcor
\lam{ankp}
The sets\/ $\ar{n}{kp}=\bigcup\ah{n}{kp}$\snos
{That is, $\ar{n}{kp}=\bigcup\ens{A}{A\in\ah{n}{kp}},$
the union of all sets in $\ah{n}{kp}$.}
belong to\/ $\ip11$ uniformly on\/ $n\yi k\yi p.$
Therefore the set\/
$\wA=\bigcup_n\bigcap_{k>n}\bigcup_p\:\ar{n}{kp}$
also belongs to\/ $\ip11$.
\ecor
\bpf
The result follows from Lemma~\ref{utv1} by standard
computations based on the coding of $\id11$ sets
(see \grf{i})
and Theorem~\ref{penu}.  
\epf
 
This leads us to the following partition onto cases.
\vtm

{\ubf Case 1:} $R\sq \wA$.\vtm

{\ubf Case 2:}
otherwise. \vtm

\punk{Case 1}
\las{6case1}

We are going to prove that in this case $\rE\reb\Eo.$
The proof of the next theorem shows that the Case 1 condition
makes all \dd\Et classes inside the domain $R=\ran\vt$
looking in a sense similar to \dd\Eo classes.
This will allow to employ Theorem~\ref{t:ass} to obtain the
result required.

Thus we are going to prove:

\bte
\lam{.key} 
In all asssumptions above,\/ ${\Et\res R}\reb\Eo$.
\ete
\bpf
By \Kres\ (Theorem~\ref{kres}) there exists a $\id11$ map 
$\nu:R\to\dN$ such that for any $x\in R$ we have
\dm
\kaz k>\nu(x)\;\sus p\;\sus B\in \ah{\nu(x)}{k,p}\;
(x\in B\in \id11)\,.
\dm
Let $R_n=\ens{x\in R}{\nu(x)\le n},$ these are increasing 
$\id11$ subsets of $R,$ and $R=\bigcup_nR_n.$ 
According to Corollary~\ref{cudc}, it suffices to prove 
that ${\Et\res R_n}\reb\Eo$ for any $n.$ 
Thus let us fix $n.$ 
Then by definition 
\pagebreak[0]%
\dm
\kaz x\in R_n\;\kaz k>n\;\sus p\;
\sus B\in \ah{n}{kp}\;
(x\in B\in \id11)\,.
\eqno(\ast)
\dm 

Recall that $\fC$ is the least class of sets containing 
all open sets and closed under the A-operation and the 
complement. 
A map $f$ is called {\it \fC-measurable\/} iff all 
\dd fpreimages of open sets belong to $\fC$.

\ble
\lam{.k1}
For any\/ $n$
there is a \dd\fC measurable map\/ $f:R_n\to\dnd$ such 
that\/ $f(x)=f(y)\eqn n x$ whenever\/ $x\yi y\in R_n$ 
satisfy\/ $x\eqn n y$.
\ele
\bpf
Let $C\sq\dN$ be the $\ip11$ set of all codes of 
$\id11$ subsets of $\dnd,$ and let 
$W_e\sq\dnd$ be the $\id11$ set coded by $e\in C.$ 
We have, by $(\ast)$,
\dm
\kaz x\in R_n\;
\kaz k>n\;\sus p\;\sus e\in C\;
(x\in W_e\in \ah{n}{kp})\,.
\dm

Now a straightforward application of \Kres\
(Theorem~\ref{kres})
yields a pair of $\id11$ maps 
$\pi\yi \ve:R_n\ti\dN\to\dN$ such that  
$\ve(x,k)\in C$ and 
$x\in W_{\ve(x,k)}\in\ah{n}{k,\pi(x,k)}$ 
hold whenever $x\in R_n$ and $k>n.$ 
Let $\pit(x,k)$ and $\tve(x,k)$ 
be the least, in the sense of any fixed recursive 
\dd\om long wellordering of $\dN\ti\dN,$ of all 
possible pairs $\pi(x',k)$ and $\ve(x',k)$ with 
$x'\in R_n\cap[x]_{\eqn n}.$ 
Then $\pit$ and $\tve$ are \dd{\eqn n}invariant in the 
1st argument. 
In addition, we have  
$W_{\tve(x,k)}\in\ah{n}{k,\pit(x,k)}$ and the set 
$Z_{xk}=R_n\cap[x]_{\eqn n}\cap W_{\tve(x,k)}$ 
is nonempty, 
whenever $x\in R_n$ and $k>n$. 

Let $x\in R_n.$ 
For any $k>n,$ the set 
$Y_{xk}=\ens{y(k)}{y\in Z_{xk}}\sq\dn$ is finite 
(and nonempty) by the definition of $\ah{n}{kp}.$ 
Let $f_k(x)$ be the least member of $Y_{xk}$ 
in the sense of the lexicographical order of $\dn$. 
Define $f(x)\in\dnd$ so that $f(x)(k)=x(k)$ for $k\le n$ 
and $f(x)(k)=f_k(x)$ for $k>n$. 

That $f(x)=f(y)$ whenever $x\eqn n y$ follows from the 
invariance of $\ve$ and $\pi.$ 
To see that $f(x)\eqn n x$ note that by definition 
$f_k(x)\Eo x_k$ for $k>n:$ 
indeed, $f_k(x)=y_k$ for some $y\in[x]_{\eqn n},$ but 
$x\eqn n y$ implies $x_k\Eo y_k$ for all $k.$ 
Finally, the \dd\fC measurability needs a routine check.
\epF{Lemma}

For any $u\in(\dn)^n$ define  
$R_n(u)=\ens{x\in R_n}{x\rnq n=u}$.

\ble
If\/ $u\in(\dn)^n$ then\/ $\Et\res R_{n}(u)$ is 
smooth.
\lam{utv4}
\ele
\bpf 
As $\Et$ and $\eqn n$ coincide on $R_n(u),$ the relation 
$\Et\res R_n(u)$ is smooth by means of a \dd\fC measurable,
hence, a Baire-measurable map. 
Suppose, towards the contrary, that it is not really 
smooth, \ie, smooth by means of a Borel map. 
Then, by the 2-nd dichotomy theorem, we have 
$\Eo\reb{\Et\res R_n(u)},$ hence, $\Eo$ turns out to 
be smooth by means of a Baire-measurable map, which is easily 
impossible.
\epF{Lemma}

To complete the proof of the theorem, let 
$\dG$ denote the group ${\pwf\dN}^n,$ that is, the product
of $n$ copies of $\stk{\pwf\dN}{\ap}.$ 
Let $\dG$ act on $\dX=(\dn)^n$ componentwise 
and by $\ap$ on each of the $n$ co-ordinates.
(Recall that $(a\ap b)(k)=0$ iff $a(k)=b(k)$ whenever
$a\yi b\in\dn$ and $k\in\dN.$)
Then, for any $u\yi v\in\dX,$ $u\ergx v$ is equivalent to 
$u(k)\Eo v(k)$ for all $k<n.$  
Let us apply Theorem~\ref{t:ass} with $\dG$ and $\dX$ 
as indicated, and $P=R_n$ and $\rE={\Et\res R_n},$  
Lemma~\ref{utv4} witnesses the principal requirement. 
Thus $\Et\res R_n$ is Borel reducible to an \eqr\ 
induced by a Borel action of $\dG.$ 
Yet $\dG$ is the increasing union of a countable 
sequence of its finite subgroups, hence any \er\ 
induced by a Borel action of $\dG$ is hyperfinite,
therefore Borel reducible to $\Eo$.\vtm

\epF{Theorem~\ref{.key} and Case 1 in Theorem~\ref{6dih}}

\punk{Case 2}
\las{6case2}

Then the $\is11$ set $H=R\dif \wA$ is non-empty.
A rather typical example is
\dm
R=\ens{x\in\dnd}{\kaz n,k,l\:(x(\bie nk)=x(\bie nl))},
\dm
where $n,k\mto\bie nk$ is a recursive pairing function on $\dN,$
see below.
Thus members of $R$ are those infinite
sequences of elements of $\dn$
in which every term is duplicated in infinitely many copies.
It can be verified that the intersection
$\wA\cap R$ consists of all sequences $x\in R$ that
contain a finite number of terms $x(0),\dots,x(n)$ such that any
other term is $\id11$ in $x(0),\dots,x(n).$
Obviously the difference $R\dif \wA$ is non-empty.

We are going to prove

\bte
\lam{,key} 
In all asssumptions above, including the Case 2 assumption,
there exists a Borel subset\/ $X$ of\/ $H$ such 
that $\Et\reb{\Et\res X}$.
\ete

This result leads to the \bfor\ case of Theorem~\ref{6dih}.
Indeed the Borel map $\vt$ that reduces $\rE$ to $\Et,$
actually to ${\Et}\res R,$ is a bijection
(see the beginning of \nrf{Et:er}),
therefore there is an inverse map
$\vpi=\vt\obr:R\to\dX=\dom\rE,$ also Borel, of course.
The map $\vpi$ then witnesses ${{\Et}\res R}\reb\rE.$
On the other hand, $\Et\reb{\Et\res X}\reb{\Et\res R}$.

\bpf
By definition, 
$H=\bigcap_n\bigcup_{k>n} \hh nk,$ where
$\hh nk=H\dif \bigcup_p\ar{n}{kp}.$    
Note that%
\jmar{hnk}%
\bus
\label{hnk}
\bay{rcl}
\hh nk 
&=& \ens{x\in H}
{\kaz p\:\kaz A\in\ah{n}{kp}\:(x\nin A)} \\[\dxii]

&=&
\ens{x\in H}
{\kaz p\:\kaz A\in\id11\:(x\in A\imp A\nin\ah{n}{kp})}
\eay
\eus
by Lemma~\ref{utv1}, and hence $\hh nk$ is $\is11$ by
rather elementary computation.
Note also that for any $\is11$ set $A$ and any
$n,k,p$ the following holds:
\jmar{ninnkp}%
\bus
\label{ninnkp}
A\nin\ah{n}{kp}\leqv
\sus y,z\in A\:\sus j\ge p\;
({y\eqn n z}\land y(k)(j)\ne z(k)(j))\,.
\eus

To prove the theorem, we are going to define a rather
complicated splitting system of non-empty $\is11$ subsets
of $\dnd.$ 
Let us take some space for technical notation 
involved in the construction of the splitting system.

Put $\bie rq=2^r(2q+1)-1$ for all $r\yi q\in\dN.$
\index{zzprecrq@$\bie rq$}%
Thus $\ang{r,q}\mto\bie rq$ is a recursive bijection
$\dN^2\onto\dN,$ increasing in each argument.
Put
\dm
L(n)=\tmax\ens{r}{\sus q\:(\bie rq\le n)}=\ens{r}{2^r-1\le n}
\dm
\index{zzLn@$L(n)$}%
for any $n$
--- for instance $L(0)=0$ and
$L(1)=L(2)=1.$
For any  $r\le L(n)$ define
$(n)_r=\ens{q}{\bie rq\le n}$ --- this is a natural number
$\ge1$ (assuming $r\le L(n)$).
\index{zznr@$(n)_r$}%
For instance $(0)_0=1$ (since $\bie 00=0$),
$(1)_0=2,$ and $(1)_1=1.$
Obviously $n=\sum_{r=0}^{L(n)-1}(n)_r$.

Suppose that $n\in\dN$ and $s\in 2^n$
(a dyadic sequence of length $n$).
For any $r<L(n)$ define $(s)_r\in 2^{(n)_r}$ so that
\index{zznr@$(s)_r$}%
$(s)_r(q)=s(\bie rq)$ for all $q<(n)_r.$
Thus the original sequence $s\in\bse$ of length $\lh s=n$
is split into a \dd{L(n)}sequence of dyadic sequences of
lengths $\lh{(s)_r}=(n)_r.$
Formally this secondary sequence $\sis{(s)_r}{r<L(n)}$
belongs to the product set $\prod_{r=0}^{L(n)-1}2^{(n)_r}$.

We consider $\dn$ as a group with the componentwise
operation,
that is, if $a,b\in\dn$ then $a\ap b\in\dn$ and
$(a\ap b)(k)=a(k)+_2 b(k)\yt\kaz k,$ where $+_2$ is the addition
modulo $2.$
The neutral element is the constant-$0$ sequence
$\fo=\dN\ti\ans0$
(that is, $\fo(k)=0\yd \kaz k$),
clearly $\fo\ap a=a$ for all $a\in\dn$.

Accordingly consider $\dnd$ as the product of \dd\dN many
copies of $\dn,$ a group with the componentwise operation
still denoted by $\ap,$ so that
$(f\ap g)(n)(k)=f(n)(k)+_2g(n)(k)$ for all $n,k.$
The neutral element is the constant-$\fo$ sequence
$\fo^\dN\in\dnd.$
Define $\supp g=\ens{n\in\dN}{g(n)\ne\fo},$ the domain of
\index{zzsuppg@$\supp g$}%
non-triviality of $g\in\dnd$.

The group $\dnd$ contains the subgroups
\dm
\xH = \ens{g\in\dnd}
{\kaz n\:\sus k_0\:\kaz k\ge k_0\:(f(n)(k)=0)},
\dm
essentially the ideal $\It,$ acting on $\dnd$ by the group
operation $\ap$, and
\dm
\bay{rclcl}
\xh n &=& \ens{g\in\xH}{\supp g\sq\iv n\iy} &=&
\ens{g\in\xH}{\kaz k\le n\:(g(k)=\fo)}\,,\\[0.8\dxii]

\zh n &=& \ens{g\in\xH}{\supp g\sq\ir n\iy} &=&
\ens{g\in\xH}{\kaz k<n\:(g(k)=\fo)}\,,\\[0.8\dxii]

\yh n &=& \ens{g\in\xH}{\supp g\sq\ix 0n} &=&
\ens{g\in\xH}{\kaz k>n\:(g(k)=\fo)}\,.
\eay
\dm
for any $n.$ 
Obviously $x\Et y$ iff $y\in \xH\ap x,$ and
$x\eqn n y$ iff $y\in \xh n \ap x,$

\vyk{
If $s\yi t\in 2^n$ for some $n$
then define $s\ap t\in 2^n$ so that
$(s\ap t)(i)=s(i)+_2 t(i)$ for all $i<n,$ 
where $+_d$ is the addition modulo $2.$
If $a\in\dn$ and $s\in\bse$ then define $s\ap a\in\bse$ so
that $(s\ap a)(k)=a(k)+_2 s(k)$ whenever
$k<\lh s,$ where $(s\ap a)(k)=a(k)$ whenever $k\ge\lh{s}.$
}

\vyk{
If $a\in\dn$ and $u\in\pwf\dN$ then define
$u\ap a\in\dN$ so that
\dm
(u\ap a)(k)=
\left\{
\bay{rcl}
a(k)    &\text{whenever}& k\nin u\,,\\[0.8\dxii]

1- a(k) &\text{whenever}& k\in u\,.
\eay
\right.
\dm

Assume that $g\in\xH=\pwf{\nn}.$
Put $(g)_n=\ens{k}{\ang{n,k}\in g}$ for any $n,$
\index{zzgn@$(g)_n$}%
and $\supp g=\ens{n}{(g)_n\ne\pu}.$ 
\index{zzsuppg@$\supp g$}%
If $x\in\dnd$ then define $g\ap x\in\dnd$ so that 
$(g\ap x)(n)=(g)_n\ap x(n)$ for all $n.$
Note that $(g\ap x)(n)=x(n)$ for all $n\nin\supp g.$
}%

Finally if $X\sq\dnd$ then put $g\ap X=\ens{g\ap x}{x\in X}$.

\vyk{
Now suppose that $x\yi y\in\dnd$ and
$n\in\dN\yt g\in\xH\yt \supp g\sq \ir0n.$
We write $x\eqo ng y$ iff
\index{zzequivsg@$\eqo ng$}%
$x\Et y$ and $x(i)= (g)_i\app y(i)$ for all $i<n.$\snos
{This generalizes the definition of $x\eqn n y$ in \nrf{Et:er}:
indeed, $x\eqn n y$ is equivalent to $x\eqo n\pu y$.}
Obviously $y= g\ap x$ implies $x\eqo ng y$ for any $n$ such that
$\supp g\sq n,$ but not conversely.

For $X,Y\sq\dnd$ define $X\eqo n g Y$ iff
for any $x\in X$
there exists $y\in Y$ satisfying $x\eqo n g y$ and \vva\
for any $y\in Y$
there exists $x\in X$ satisfying $x\eqo n g y$.
}

The splitting system used here will contain non-empty
$\is11$ sets $X_s\sq \dnd\yt s\in\dln,$
the increasing sequence of numbers $k_0<k_1<k_2<\dots\in\dN,$
a collection of natural
numbers $p_{mj}\yt m,j\in\dN,$ and elements
$\gs st\in\xH,$ where $s,t\in\bse\yt \lh s =\lh t,$ 
satisfying the following requirements \ref{6i} -- \ref{6last}:

\ben
\tenu{(\roman{enumi})}
\itlm{6i}\msur
$X_\La\sq H,$ $X_{s\we i}\sq X_s,$ $\dia{X_s}\le2^{-\lh s}$.

\itlm{6igen}%
A certain condition similar to \ref{z8} in \nrf{3case2}
holds, connecting each $X_{s\we i}$ with $X_s$ so that, as a 
consequence, $\bigcap_nX_{a\res n}\ne\pu$ 
for any $a\in\dn$.

\itlm{6ii}%
If $s\in 2^{n+1}$ then 
$X_{s}\sq \bigcap_{r\le L(n)}\hh{r}{k_r}$.

\itlm{6iii}%
If $s,t\in2^{n+1}$ then $\supp{\gs st}\sq\ix0{k_{L(n)}},$
that is, $\gs st\in\yh{k_{L(n)}}$.

\itlm{6iii.}\msur
$k_0<k_1<k_2<\dots,$
and $p_{m0}<p_{m1}<p_{m2}<\dots$ for any $m$.

\itlm{6iii+}\msur
$\gs su=\gs tu\ap \gs st$ for all $s, t, u\in2^{n+1}.$
It easily follows that $\gs ss=\fo^\dN\yd\kaz s$.


\itlm{6iv}%
For any $s\yi t\in2^{n+1},$ we have
$\gs st\ap X_s\eqn {k_{L(n)}} X_t$.
\setcounter{enui}{\value{enumi}}
\een
We define $X\eqn m Y$ (for any sets $X,Y\sq\dnd$) iff
\index{zzequivn@$\eqn n$}%
$\ek X{\eqn m}=\ek Y{\eqn m},$ that is,
for any $x\in X$ there exists $y\in Y$ satisfying
$x\eqn m y$ and \vva\
for any $y\in Y$
there exists $x\in X$ satisfying $x\eqn m y.$
This is equivalent to $\xh m\ap X=\xh m\ap Y.$

\ben
\tenu{(\roman{enumi})}
\setcounter{enumi}{\value{enui}} 
\itlm{6v}%
For any $s,t\in 2^{n+1},$ if $\ell\le L(n),$
$n'\le n,$ and $s',t'\in 2^{n'}$   
satisfy $s'\sq s\yt t'\sq t,$ and
the equality $(s)_r(q)=(t)_r(q)$ holds whenever 
$r\le\ell$ and $q\in\dN$ satisfy $n'\le \bie rq\le n,$
then $\gs st(i)=\gs{s'}{t'}(i)$ for any $i\le\ell.$

\itlm{6vi}%
\label{6last}%
For any $s,t\in 2^{n+1},$ if $s(n)=0\ne1=t(n),$ and
$n=\bie mj,$ then $x(k_m)(p_{mj})=0$ for all $x\in X_s$ but 
$y(k_m)(p_{mj})=1$ for all $y\in X_t$.
\een

\punk{The embedding}
\las{6embe}

Suppose that a system of sets $X_s,$ elements $\gs st,$ and
numbers $k_m$ and $p_{mj}$ satisfying \ref{6i} -- \ref{6last}
has been defined.
Let us show that this leads to the proof of Theorem~\ref{,key}.

As usual it follows from \ref{6i} and \ref{6igen} that
for any $a\in\dn$ the intersection $\bigcap_n X_{a\res n}$ is
a singleton.
Let us denote by $\vt(a)=\sis{\vt_n(a)}{n\in\dN}$ its only
element.
Thus $a\mto\vt(a)$ is a map $\dn\to\dnd$ while each $\vt_n$
is a map $\dn\to\dn.$
In addition both $\vt$ and all $\vt_n$ are continuous
(in the Polish product topology).

On the other hand for any $a\in\dn$ there is a unique point
$\wa=\sis{(a)_n}{n\in\dN}\in\dnd$
such that $(a)_n(k)=a(\bie nk)$ for all $n,k.$
The map $a\mto \wa$ is a homeomorphism of $\dn$ onto $\dnd,$
while each $a\mto{(a)_n}$ is a continuous map $\dn\to\dn.$
Thus the following lemma suffices to prove Theorem~\ref{,key}:

\ble
\lam{keyl}
For any\/ $a,b\in\dn,$ we have$:$
${\wa\Et\wb}$ iff\/ ${\vt(a)\Et\vt(b)}$.
\ele
\bpf
Assume that ${\wa\Et\wb},$ take an arbitrary $\ell\in\dN$
and prove that $\vt_\ell(a)\Eo\vt_\ell(b).$
In our assumptions there exists a number $n'$ such that
$\ell\le L(n')$ and for any $r\le\ell$ and $q,$ if
$\bie rq\ge n'$ then $a(\bie rq)=b(\bie rq).$
Put $s'=a\res {n'}$ and $t'=b\res n'.$
Then $g'=\gs{s'}{t'}\in\xH.$ 
Our goal is to prove that $\vt_\ell(b)=(g')_\ell\ap\vt_\ell(a),$ 
that obviously implies $\vt_\ell(a)\Eo\vt_\ell(b)$.

It suffices to show that
$\gs{s'}{t'}\app X_{s}\eqn \ell X_{t}$
holds for any $n>n',$ where $s=a\res n$ and $t=b\res n.$
We observe that $\gs st\ap X_s\eqn {\ell} X_t$
by \ref{6iv} because $\ell\le L(n')\le k_{L(n')}\le k_{L(n)}.$
On the other hand, $\gs st(i)=\gs{s'}{t'}(i)$ for any
$i\le\ell$ by \ref{6v} and the choice of $n'.$  
It follows that
$\gs{s'}{t'}\app X_{s}=\gs st\ap X_s\eqn {\ell} X_t,$
as required.

To prove the converse suppose that ${\wa\Et\wb}$ fails, and hence
there is at least one index $m$ such that $(a)_m\Eo (b)_m$ fails
as well, meaning that $a(\bie mj)\ne b(\bie mj)$ holds for
infinitely many numbers $j\in\dN.$
Then by \ref{6vi} we obtain
$\vt_{k_m}(a)(p_{mj})=0\ne 1=\vt_{k_m}(b)(p_{mj})$
for all $j,$ and hence $\vt_{k_m}(a)\Eo\vt_{k_m}(b)$
fails since the numbers $p_{mj}\yt j\in\dN,$ form a strictly
increasing sequence by \ref{6iii.}.\vtm

\epF{Lemma~\ref{keyl}}

\punk{The construction of a splitting system: warmup}
\las{6sis}

Now to prove Theorem~\ref{,key} it remains to carry out
the construction of a system of sets $X_s$ and $\gs st$ and
numbers $k_m$ and $p_{mj}$ satisfying
conditions \ref{6i} -- \ref{6last} of \nrf{6case2}.
The construction goes on by induction on $n,$ so that at each
step $n$ we define the sets $X_s\yt s\in 2^n$ and elements
$\gs st\yt s,t\in2^n.$
Here we present only the transition from $0$ to $1$ as a
warmup.

Put $X_\La= H$ and by default $\gs\La\La=\fo^\dN$ for the only
sequence $\La$ of length $0.$

At the next stage, we have to define $\is11$ sets
$X_{\ang0}\yi X_{\ang1}\sq X_\La,$ an element
$\gs{\ang0}{\ang1}=\gs{\ang1}{\ang0}\in\xH,$ and numbers
$k_0$ and $p_{00}$ such that a relevant fragment of
\ref{6i} -- \ref{6last} is satisfied.
Note that $L(0)=0$.\vom

{\ubf Stage 1.}
We shrink $X_\La$ to make sure that conditions
\ref{6i} and \ref{6igen} are satisfied; the resulting
$\is11$ set is still denoted by $X_\La$.\vom

{\ubf Stage 2.}
Consider any $x\in X_\La.$
Then $x\in\bigcap_{k>0}\hh{0}{k}$
(see the beginning of the proof of Theorem~\ref{,key}).
Fix a number $k=k_0>0$ such that $x\in \hh{0}{k_0}.$
The set $X'_\La=X_\La\cap \hh{0}{k_0}$ is still of class
$\is11,$ and for any $p$ it does not belong to the
family $\ah{0}{k_0,p}$ by \eqref{hnk} in \nrf{6case2}.
Thus by \eqref{ninnkp} there exist points
$y_0\yi z_0\in X'_\La$ satisfying $y_0\eqn {L(0)} z_0$
and numbers $k_0>0=L(0)$ and $p_{00}$ such that
$y_0(k_0)(p_{00})=0 \ne 1=z_0(k_0)(p_{00}).$
The $\is11$ sets
\dm
\bay{rcll}
Y &=&\ens{y\in X'_\La}
{y\eqn {L(0)} y_0\land y(k_0)(p_{00})=0}\,,&\text{and}\\[0.8\dxii]

Z &=&\ens{z\in X'_\La}
{z\eqn {L(0)} z_0\land z(k_0)(p_{00})=1}\,
\eay
\dm
still contain resp.\ $y_0,z_0,$ therefore so do the $\is11$ sets
\dm
Y' =\ens{y'\in Y}
{\sus z\in Z\;(y'\eqn {L(0)} z)}
\quad\text{and}\quad
%
Z' =\ens{z'\in Z}
{\sus y\in Y\;(y\eqn {L(0)} z')}\,.
\dm
Finally define
$\gs{\ang0}{\ang1}=\gs{\ang1}{\ang0}\in\xH$
so that $\gs{\ang0}{\ang1}(k_0)(p_{00})=1$ and
$\gs{\ang0}{\ang1}(m)(j)=0$ for any other pair of $m,j.$
Then easily $\gs{\ang0}{\ang1}\ap y_0\eqn{k_0} z_0,$ hence 
$\gs{\ang0}{\ang1}\ap Y'\eqn{k_0} Z'.$
Thus we get a pair of sets $X_{\ang0}=Y'$ and $X_{\ang1}=Z'$
compatible with \ref{6iv}.
This ends the construction for $n=1$.

\punk{The construction of a splitting system: the step}
\las{6sig}

Now suppose that $n=\bie mj\ge1,$ and the construction has been
accomplished up to the level $n,$ that is, there exist
sets $X_s\sq H$ and elements $\gs st\in\xH,$
where $s,t\in 2^{n'}\yt n'\le n,$
and numbers $k_0,\dots,k_{L(n-1)}$ and $p_{m'j'},$ where
$\bie {m'}{j'}<n,$ such that conditions \ref{6i} -- \ref{6last}
are satisfied in this domain.
The goal is to define $X_s$ and $\gs st,$ where
$s,t\in 2^{n+1},$
and numbers $k_n$ and $p_{mj},$ such that conditions
\ref{6i} -- \ref{6last} are satisfied in the extended domain.

The numbers $n,m,j$ are fixed in the course of the arguments
in this Section.

\ble
\lam{shr}
Suppose that collections of\/ $\is11$ sets\/
$P_s\sq\dnd\yt s\in2^{n}$
and elements\/ $\gs st\in\xH\yt s,t\in2^n,$ satisfy\/
both \ref{6iii+} and \ref{6iv} for a fixed\/ $k,$ that is, 
$\gs su=\gs tu\ap \gs st$ and\/ 
$\gs st\ap P_s\eqn{k} P_t$
for all\/ $s, t, u\in2^{n+1}.$ 

If\/ $\sg\in2^n$ and\/ $P\sq P_\sg$ is a non-empty\/
$\is11$ set then the sets\/
\dm
{P'_s=\ens{x\in P_s}
{\sus y\in P\,(\gs \sg s\ap y\eqn{k}x)}}
,\quad s\in2^n,
\dm
are non-empty\/ $\is11$ sets still satisfying\/
\ref{6iv}, \ie\ $\gs st\ap P'_s\eqn{k} P'_t$
for all\/ $s, t\in2^{n+1}.$
\ele
\bpf
Fix $s\yi t\in2^{n}.$
To show that
$\gs st\ap X'_s\eqn {k} X'_t$
consider any $x\in X'_s,$ so that
$\gs st\ap x\in \gs st\ap X'_s.$
By definition there exists $y\in X$ satisfying
$\gs \sg s\ap y\eqn{k}x.$
It follows from \ref{6iii} that $\gs st\in\yh{k},$
therefore
$\gs st\ap\gs \sg s\ap y\eqn{k}\gs st\ap x,$
that is, $\gs \sg t\ap y\eqn{k}\gs st\ap x$ by
\ref{6iii+}.
However by definition $\gs \sg t\ap y\in X'_t,$ as required.

The converse is similar.
\epF{Lemma}

It follows from \ref{6iii} (in the domain $2^n$)
that there is a number $\mu\in\dN$
such that $\gs st(r)(q)=1$ holds only in the case when
both $r\le k_{L(n-1)}$ and $q\le \mu.$
We proceed with several stages of successive reduction
and splitting of the $\is11$ sets $X_s\yt s\in2^n.$
These further stages depend on whether the number
$n=\bie mj$ considered opens a \lap{new} axis $k_m$
of splitting.\vom

{\bf Case A\/}: $j>0.$\vom

Then $n'=\bie{m}{j-1}<n,$ thus $m$ is \lap{old}.
Moreover, $L(n)=L(n-1).$
We have to define $p_{mj}$ but needn't to define any
new $k_r.$\vom  

{\ubf Stage 1.}
Fix an arbitrary sequence $\sg\in2^n;$
this can be \eg\ the sequence $0^n$ of $n$ zeros.
Consider any $x\in X_\sg.$
Then $x\in \hh{m}{k_m}$ by \ref{6ii}, and hence there exist
points $y_0,z_0\in X_\sg$ and a number $p_{mj}>\mu$ such that
$y_0\eqn{m-1} z_0$ and $y_0(k_m)(p_{mj})=0$
but $z_0(k_m)(p_{mj})=1.$
Easily $p_{mj} >p_{m,j-1}:$ indeed $p_{m,j-1}\le \mu$ by the
choice of $\mu.$\vom  

{\ubf Stage 2.}
Define $g\in\dnd$ so that $g(r)(q)=1$ iff
both $m\le r\le k_{L(n)}$ and $y_0(r)(q)\ne z_0(r)(q).$
Then $g\in\xH$ since $y\Et z.$
Moreover we have 
$\supp g\sq\ix {m}{k_{L(n)}},$ in other words,
$g\in\zh{m}\cap\yh{k_{L(n)}}.$
In addition $g(k_m)(p_{mj})=1.$

We observe that by definition $g\app y_0\eqn{k_{L(n)}}z_0.$
Thus the $\is11$ sets
\dm
\bay{rcl}
Y&=&\ens{y\in X_\sg}
{y(k_m)(p_{mj})=0\land \sus z\in Z\:
(z(k_m)(p_{mj})=1\land g\app y\eqn{k_{L(n)}}z)}\,,\\[0.8\dxii]

Z&=&\ens{z\in X_\sg}
{z(k_m)(p_{mj})=1\land \sus y\in Y\:
(y(k_m)(p_{mj})=0\land g\app y\eqn{k_{L(n)}}z)}
\eay
\dm
are still non-empty (contain resp.\ $y_0\yi z_0$) and
satisfy $g\ap Y\eqn{k_{L(n)}}Z;$
in addition $y(k_m)(p_{mj})=0$ and $z(k_m)(p_{mj})=1$ for all
$y\in Y$ and $z\in Z.$

As a matter of fact we can \noo\ assume
that $Y\cup Z=X_\sg$:
indeed otherwise put $P=Y\cup Z$ and apply Lemma~\ref{shr}.\vom  

{\ubf Stage 3.}
Put $X_{\sg\we 0}=Y$ and $X_{\sg\we 1}=Z,$
thus
\bus
\label{Kk}
g\ap X_{\sg\we 0}\eqn{k_{L(n)}}X_{\sg\we 1},
\eus
and then
\dm
X_{s\we \xi}=\ens{x\in X_s}
{\sus y\in X_{\sg\we \xi}\:(\gs{\sg}s(y)\eqn{k_{L(n)}}x)}
\dm
for all $s\in2^n$ and $\xi=0,1.$
It follows, by \ref{6iv} at the level $n,$ that
\bus
\label{KK}
X_{s\we\xi}\eqn{k_{L(n)}}\gs\sg s \ap X_{\sg\we\xi}
\quad\text{for all $s\in2^n$ and $\xi=0,1$}\,.
\eus
Put $\gs{s\we \xi\,,\,}{t\we \xi}=\gs st$ but
$\gs{s\we \xi\,,\,}{t\we {(1-\xi)}}=\gs st\ap g$ for all
$s,t\in2^n$ and $\xi=0,1,$\snos
{\label{nonab}%
In the definition of $g_{st}$ we make use of the fact that
$\stk{\dnd}{\ap}$ is an abelian group.
In the non-abelian case we would have to define
$\gs{s\we i\,,\,}{t\we{(1-i)}}=\gs\sg t\ap g\ap\gs s\sg$
and accordingly render some other related definitions in
somewhat more complicated way.}
or saying it differently
\bus
\label{kkk}
\gs{s\we \xi\,,\,}{t\we {\eta}}=\gs st\ap g^{\xi-\eta}
\quad\text{for all $s,t\in2^n$ and $\xi,\eta=0,1$}
\eus
where $g^1=g^{-1}=g$ while $g^0=\fo^\dN$
is the neutral element in $\stk\dnd\ap$.\vom  

{\ubf Stage 4.}
Lemma~\ref{shr} allows us to reduce the sets
$X_s\yt s\in2^{n+1},$
in several rounds to make sure that conditions
\ref{6i} and \ref{6igen} are satisfied at level $n+1;$
the resulting $\is11$ sets are still denoted by $X_s$. 

This ends the transition from $n$ to $n+1.$
It remains to show that conditions \ref{6i} -- \ref{6last}
are satisfied in the extended \dd{(\le n+1)}domain.\vom

{\ubf Verification.}
As \ref{6i} and \ref{6igen} are explicitly fulfilled,
\ref{6ii} in Case 1 is vacuous,
and \ref{6iii}, \ref{6iii.} clearly hold by definition,
we begin with \ref{6iii+}.
We have to prove that
\dm
\gs{s\we\xi\,,\,}{u\we\za}=\gs{t\we\eta\,,\,}{u\we\za}
\ap\gs{s\we\xi\,,\,}{t\we\eta}
\dm
for all $s, t, u\in2^{n}$ and $\xi,\eta,\za=0,1.$
By definition this equality is equivalent to
$
\gs{s}{u}\ap g^{\xi-\za}=\gs{t}{u}\ap g^{\eta-\za}
\ap\gs{s}{t}\ap g^{\xi-\eta}.
$
However obviously
$g^{\xi-\za}=g^{\eta-\za}\ap g^{\xi-\eta},$ 
and on the other hand in our assumptions
$\gs{s}{u}=\gs{t}{u}\ap\gs{s}{t}$ by \ref{6iii+} at level $n.$\vom

Let us check \ref{6iv}, that is,
$\gs{s\we\xi\,,\,}{t\we\eta}\ap X_{s\we\xi}\eqn {k_{L(n)}}
X_{t\we\eta}$
for all $s, t\in2^{n}$ and $\xi,\eta=0,1.$
It follows from \eqref{KK} that 
the left-hand side is \dd{\eqn{k_{L(n)}}}equivalent to
$\gs st\ap g^{\xi-\eta}\ap\gs\sg s\ap X_{\sg\we\xi}$ 
while the right-hand side is \dd{\eqn{k_{L(n)}}}equivalent to
$\gs\sg t\ap X_{\sg\we\eta}.$
On the other hand it follows from \eqref{Kk} that  
$g^{\xi-\eta}\ap X_{\sg\we\xi}\eqn{k_{L(n)}}X_{\sg\we\eta}.$
This allows to easily get the result required.\vom

Let us check \ref{6v}.
Suppose that $s\yi t\yi \ell\yi n'\yi s'\yi t'$ are as indicated
in \ref{6v}.
Then $s=\ws\we\xi$ and $t=\wt\we\eta,$ where $\ws\yi\wt\in2^n$
while $\xi=s(n)$ and $\eta=t(n)$ are numbers in $\ans{0,1}.$ 
Then $\gs \ws\wt(i)=\gs{s'}{t'}(i)$ for any $i\le\ell$ by
\ref{6v} in the domain $2^n.$
Thus if $\xi=\eta$ then the result holds immediately because
then $\gs st=\gs \ws\wt$ by \eqref{kkk}.
Assume that \eg\ $\xi=0$ and $\eta=1.$
Then $\ell<m$ in the assumptions of \ref{6v}, and hence
the set $\supp g$ does not contain numbers $i\le\ell,$
in other words, $g(i)=\fo$ for any $i\le\ell.$
It follows that $\gs st(i)=\gs \ws\wt(i)$ for any $i\le\ell,$
as required.\vom

We finally check \ref{6vi}.
Suppose that $s\we\xi$ and $t\we \eta$ belong to $2^{n+1}$ and
$\xi\ne\eta,$ say $\xi=0\ne1=\eta.$
We have to prove that $x(k_m)(p_{mj})=\xi$ for all
$\xi=0,1$ and $x\in X_{s\we\xi}.$
First of all note that by definition  
$x(k_m)(p_{mj})=\xi$ for all $x\in X_{\sg\we\xi}.$
On the other hand $\gs\sg s(k_m)(p_{mj})=0$ since
$p_{mj}>\mu$ by the construction.
Thus $(\gs\sg s\ap x)(k_m)(p_{mj})=\xi$  for all
$x\in X_{\sg\we\xi}.$
It remains to use \eqref{KK}.
\vtm

{\bf Case B\/}: $j=0.$\vom

Then there is no number $n'=\bie{m'}{j'}<n$ such that $m'=m$
--- in other words, $m$ is \lap{new}.
Obviously $m=L(n-1)+1=L(n)$ in this case.\vom

{\ubf Stage 1.}
The first goal is to appropriately choose a number $k_m.$
Let us fix an arbitrary $\sg\in 2^n.$
Consider any $x\in X_\sg.$
As $X_\sg\sq X_\La\sq H=\bigcap_n\bigcup_{k>n}\hh{n}{k},$
it follows from \eqref{hnk} in
\nrf{6case2} that $x\in \hh{k_{L(n-1)}+1}{k_m}$ for some
$k_m>k_{L(n-1)}+1.$
In particular $k_m>k_{m-1},$ $k_m>L(n),$ and $x\in \hh{L(n)}{k_m}.$

It can be \noo\ assumed that $X_\sg\sq \hh{L(n)}{k_m}.$
(Indeed otherwise we can replace the set $X_\sg$ by
$X'_\sg=X_\sg\cap \hh{L(n)}{k_m},$ still a non-empty $\is11$ set,
and apply Lemma~\ref{shr} to shrink all sets $X_s\yt s\in2^n,$ 
accordingly.)

\ble
\lam{ita}
In this assumption, $X_{s}\sq \hh{L(n)}{k_m}$ for all\/ $s\in2^n.$
\ele
\bpf
Consider an arbitrary point $x_0\in X_{s}$ and prove that
$x_0\in \hh{L(n)}{k_m}.$
Fix any number $p$ and a $\id11$ set $A\sq\dnd$
containing $x_0;$ we have to show that
$A\nin\ah{L(n)}{k_m,p}$.

Recall that $X_s\eqn {k_{L(n-1)}} \gs \sg s\ap X_\sg$
by \ref{6iv} in the domain $2^n,$ therefore there is
a point $y_0\in X_\sg$ satisfying
$x_0\eqn {k_{L(n-1)}} \gs \sg s\ap y_0.$
Then there exists an element $g\in\xH$ with
$\supp g\sq\ix0{k_m}$ such that $x_0\eqn{k_m}g\ap y_0.$
And it is clear that $g$ extends $\gs\sg s$ in the
sense that $g(r)=\gs\sg s(r)$ for all $r\le L(n-1).$

The pre-image
$B=\ens{y\in \dnd}{\sus x\in A\:(g\ap y\eqn {k_m}x)}$
is a $\is11$ set containing $y_0.$
But in our assumptions $y_0\in X_\sg\sq \hh{L(n)}{k_m},$
and hence there exist points $y\yi y'\in B$ such that
$y\eqn {{L(n)}} y'$ but $y'(k_m)\ap y(k_m)\not\sq\ix0p.$
In other words, there is a number $j>p$ with
$y'(k_m)(j)\ne  y(k_m)(j).$
By definition there exist poits $x\yi x'\in A$ such that
$g\ap y\eqn {k_m}x$ and $g\ap y'\eqn {k_m}x'.$
In particular $x(r)=g(r)\ap y(r)$ and
$x'(r)=g(r)\ap y'(r)$ for all $r\le k_n.$
We conclude that $x\eqn {{L(n)}} x'$ but 
$x'(k_m)(j)\ne  x(k_m)(j).$
It follows that $A\nin\ah{L(n)}{k_m,p},$
as required.
\epF{Lemma}

{\ubf Stage 2.}
It follows from \eqref{ninnkp} in
\nrf{6case2} that there exist points $y_0,z_0\in X_\sg$
and a number $p_{m0}\in\dN$ such that
$y_0\eqn{L(n)} z_0$ and
$y_0(k_m)(p_{m0})=0\ne 1=z_0(k_m)(p_{m0}).$
Following the construction in Case A, 
define $g\in\zh{m}\cap\yh{k_{L(n)}}$ so that
$g\app y_0\eqn{k_{L(n)}}z_0,$  in particular,
$g(k_m)(p_{m0})=1.$
Then the $\is11$ sets
\dm
\bay{rcl}
Y&=&\ens{y\in X_\sg}
{y(k_m)(p_{m0})=0\land \sus z\in Z\:
(z(k_m)(p_{m0})=1\land g\app y\eqn{k_{L(n)}}z)}\,,\\[0.8\dxii]

Z&=&\ens{z\in X_\sg}
{z(k_m)(p_{m0})=1\land \sus y\in Y\:
(y(k_m)(p_{m0})=0\land g\app y\eqn{k_{L(n)}}z)}
\eay
\dm
are still non-empty sets containing resp.\ $y_0\yi z_0$ and
satisfying $g\ap Y\eqn{k_{L(n)}}Z;$
in addition $y(k_m)(p_{m0})=0$ and $z(k_m)(p_{m0})=1$ for all
$y\in Y$ and $z\in Z.$
And still we can \noo\ assume that $Y\cup Z=X_\sg.$\vom

{\ubf Stage 3.}
We define the sets $X_{s\we \xi}\sq X_s$ and elements
$\gs{s\we \xi\,,\,}{t\we {\eta}}$
($s,t\in2^n$ and $\xi=0,1$)
exactly as on Stage 3 of Case A.
Conditions \eqref{Kk}, \eqref{KK}, \eqref{kkk}
still hold and by the same reasons.\vom  

{\ubf Stage 4.}
Shrink the sets $X_s\yt s\in2^{n+1},$ with the help 
of Lemma~\ref{shr}, in several rounds, so that
the resulting $\is11$ sets, still denoted by $X_s,$ 
satisfy \ref{6i} and \ref{6igen} in the domain $2^{n+1}.$
This completes the transition from $n$ to $n+1.$\vom

{\ubf Verification.}
A new feature here in comparison to Case A is the non-vacuous
character of condition \ref{6ii}.
It suffices to show that $X_{s\we\xi}\in \hh{L(n)}{k_m}$ for all
$s\in2^n$ and $\xi=0,1,$ or, that is sufficient,
$X_{s}\in \hh{L(n)}{k_m}$ for all $s\in2^n$ --- but this follows
from Lemma~\ref{ita}.
The verification of \ref{6iii} -- \ref{6last} is quite similar
to the verification in Case 1, we leave it to the reader.\vtm

\epF{Theorem~\ref{,key} and Case 2 in Theorem~\ref{6dih}}

\qeDD{Theorem~\ref{6dih}}

\newpage

\api

\parf{Summable ideals and \eqr s}
\las{summI}

Given a sequence of nonnegative reals $r_n$ with
$\sum_{n=0}^\iy r_n=\piy,$ the summable
ideal $\srn$ consists of all sets $x\sq\dN$ such that
$\mrn(x)=\sum_{n\in x}r_n<+\iy.$
The corresponding \eqr\ $\ern$ is defined on $\pn$ so that
$x\ern y$ iff $x\sd y\in\srn.$
Equivalently $\ern$ is defined on $\dn$ the same way, with
$a\sd b=\ens{n}{a(n)\ne b(n)}$ for $a\yi b\in\dn$.

Farah~\cite[\pff 1.12]{aq} gives the following classification 
of summable ideals   
based on the distribution of reals $r_n$:
\ben
\tenu{(S\arabic{enumi})}
\itla{S1}
{\it Atomic\/} ideals: there is $\ve>0$ such that 
the set $A_{\ve}=\ens{n}{r_n\ge\ve}$ is infinite and 
satisfies $\mrn(\dop{A_{\ve}})<\piy.$ 
In this case $\sui{r_n}=\ans{x:x\cap A_{\ve}\in\ifi};$
so this is what 
Kechris~\cite{rig} called  
{\it trivial variations of\/ $\ifi$},
see Footnote~\ref{trvar} in \gla{idI1}.

\itla{S2}
{\it Dense\/} (summable) ideals: $r_n\to 0$.

\itla{S3}
There is a decreasing sequence of 
positive reals $\ve_n\to0$ sich that all sets 
$D_n=A_{\ve_{n+1}}\dif A_{\ve_n}$ are infinite.

\itla{S4}
Ideals of the form $\ifi\oplus\hbox{dense}:$ 
there is a real $\ve>0$ such that 
\imar{define $\oplus$ somewhere}%
the set $A_{\ve}$ is infinite, 
$\mrn(\dop{A_{\ve}})=\piy,$ and 
$\tlim_{n\to\iy\,,\;n\in\dop{A_{\ve}}}r_n=0$.
\een

In the sense of $\reb,$ all ideals of types \ref{S2}, 
\ref{S3}, \ref{S4} are equivalent to each other, 
and all ideals of type \ref{S1} are equivalent to each other, 
so that we have just $2$ summable ideals modulo $\eeb,$ 
namely $\ifi$ and, say, $\sui{1/n}.$ 
The structure under $\orb$ (which we don't consider here)
is much more complicated.

This \gla\ is mainly devoted to the following theorem
of Hjorth~\cite{h-ban},
often called the 4th dichotomy theorem.

\bte
\lam h
Let\/ $\rE$ be a Borel \er\ on a Polish space\/ $\dX,$ 
and\/ $\rE\reb \esn.$ 
Then either\/ $\rE\eqb\esn$ or\/ $\rE$ is essentially 
countable.
\ete

\punk{Grainy sets}
\las S

We begin the proof of Theorem~\ref{h} with a few definitions.


For $a\yi b\in \dn$ put $a\sd b=\ans{n:a(n)\ne b(n)}$ 
(identified with the function $c(n)=1$ iff $a(n)\ne b(n)$) 
and $\sm(a,b)=\sum_{n\in a\sd b,\:n\ge1}\frac1{n}$ --- 
this can be a nonnegative real or $\piy.$  
Generally, we define 
$\sm_k^m(a,b)=\sum_{n\in a\sd b\,,\;k\le n\le m}\frac1{n}$
for $1\le k\le m,$
and accordingly
$\sm_k^\iy(a,b)=\sum_{n\in a\sd b\,,\;k\le n<\iy}\frac1{n}.$ 

Define $\sm(a)=\sum_{a(n)=1,\:n\ge1}\frac1{n}$ 
and similarly $\sm^m_k(a)$ and $\sm^\iy_k(a)$. 


Recall that the {\it summable ideal\/} is defined as 
\dm
\sun=\ans{a\in\dn:\sm(a)<\piy}\,.
\dm
(The notation $\cI_2$ and $\cI_0$ is also used.) 
$\esn$ will denote the associated Borel \er\ on $\dn,$ 
\ie, $a\esn b$ iff $\sm(a,b)<\piy$. 

Suppose that $\vt:\dX\to\dn$ is a Borel reduction of $\rE$ to 
$\esn.$ 
We can assume that $\vt$ is in fact continuous.
Indeed 
it is known that there is a stronger Polish topology 
on $\dX$ which makes $\vt$ continuous but does not add 
new Borel subsets of $\dX.$ 
Moreover, as any Polish space $\dX$ is a $1-1$ continuous
image of a closed subset of $\bn,$ we can assume that 
$\dX=\bn$. 

Finally, we can assume that $\vt$ is $\id11,$ not 
merely Borel.

If $a\in A\sq\dn$ and $q\in\dqp$ then let $\gal qAa$ 
be the set of all $b\in A$ such that there is a finite 
chain $a=a_0\yi a_1\yi \dots \yi a_n=b$ of reals $a_i\in A$ 
such that $\sm(a_i,a_{i+1})<q$ for all $i,$ the 
\dd q{\it galaxy of\/ $a$ in\/ $A$.\/} 

\bdf
\lam{d:good}
A set $A\sq\dn$ is {\it\gra q\/}, where 
$q\in\dqp,$ iff $\sm(a,b)<1$ for all $a\in A$
and $b\in\gal qAa.$ 
A set $A$ is {\it\grap\/} if it is \gra q for some $q\in\dqp.$ 
(In other words it is required that the galaxies are 
rather small.) 
\edf


\bcl
\lam{cl1}
Any\/ \gra q $\is11$ set\/ $A\sq\dn$ is covered  
by a\/ \gra q $\id11$ set. 
\ecl 
\bpf\snos
{\label{rp}%
The result can be achieved as a routine application 
of a reflection principle, yet we would like to show 
how it works with a low level technique.} 
The set 
$D_0=\ans{b\in\dn:A\cup\ans b\,\hbox{ is \gra q}}$ 
is $\ip11$ and $A\sq D_0,$ hence, there is a $\id11$ set 
$B_1$ with $A\sq B_1\sq D_0.$ 
Note that $A\cup\ans a$ is \gra q for any $a\in B_1.$ 
It follows that the $\ip11$ set 
\dm
D_1=\ans{b\in B_1:A\cup\ans{a,b}\,
\hbox{ is \gra q for any }\,a\in B_1}
\dm
still contains $A,$ hence, there is a $\id11$ set 
$B_2$ with $A\sq B_2\sq D_1\sq B_1.$ 
Note that $A\cup\ans{a_1,a_2}$ is \gra q for any 
$a_1\yi a_2\in B_2.$ 
In general, as soon as we have got a $\id11$ set $B_n$ 
with $A\sq B_n$ and such that $A\cup\ans{a_1,\dots ,a_n}$ is 
\gra q for any $a_1,\dots ,a_n\in B_n,$ then 
the $\ip11$ set 
\dm
D_{n}=\ans{b\in B_n:A\cup\ans{a_1,\dots ,a_n,b}\,
\hbox{ is \gra q for any }\,a_1,\dots ,a_n\in B_n}
\dm
contains $A,$ hence, there is a $\id11$ set 
$B_{n+1}$ with $A\sq B_{n+1}\sq D_{n}\sq B_n$.

As usual in similar cases, the choice of the sets $B_n$ 
can be made effective enough for the set 
$B=\bigcap_nB_n$ to be still $\id11,$ not merely Borel. 
On the other hand, $A\sq B$ and $B$ is \gra q.
\epF{Claim}

Coming back to the proof of the theorem, 
let $C$ be the union of all \grap\ $\id11$ sets. 
An ordinary computation shows that $C$ is $\ip11.$ 
We have two cases.\vtm

{\ubf Case 1:} \ \msur$\ran\vt\sq C.$\vtm 

{\ubf Case 2:} \ otherwise.

\punk{Case 1}
\las{4case1}

We are going to prove that, in this case, $\rE$ is 
essentially countable.
First note that, 
by \Sepa, there is a $\id11$ set $\aH\sq\dn$ 
with $\ran \vt\sq \aH\sq C$. 

Fix a standard enumeration $\sis{W_e}{e\in E}$ of all 
$\id11$ subsets of $\dn,$ where, as usual, $E\sq\dN$ is 
a $\ip11$ set. 
By \Kres, there exist $\id11$ functions $a\longmapsto e(a)$ 
and $a\longmapsto q(a),$ defined on $\aH,$ 
such that for any $a\in\aH$ the $\id11$ set 
$W(a)=W_{e(a)}$ contains $a$ and is \gra{q(a)}. 
The final point of our argument will be an application 
of Lemma~\ref{KT}, where $\rho$ will be a derivate of 
the function $G(a)=\gal{q(a)}{W(a)}a.$ 
We prove 

\bcl
\lam{cl3}
If\/ $a\in\aH$  then\/ 
$\ga_a=\ans{G(b):b\in [a]_{\esn}\cap\aH}$ 
is at most countable.
\ecl
\bpf
Otherwise there is a pair of $e\in E$ and $q\in\dqp$ 
and an uncountable set $B\sq [a]_{\esn}\cap\aH$ 
such that $q(b)=q$ and $e(b)=e$ for any $b\in B$ and 
$G(b')\ne G(b)$ for any two different $b,b'\in B.$ 
Note that any $G(b),\msur$ $b\in B,$ is a \dd qgalaxy 
in one and the same set $W(a)=W(b)=W_e,$ therefore, 
if ${b\ne b'}\in B$ then $b'\nin G(b)$ and 
$\sm(b,b')\ge q.$ 
On the other hand, as $B\sq [a]_{\esn},$ we have 
$\sm(a,b)<\piy$ for all $b\in B,$ hence, there is  
$m$ and a still uncountable set $B'\sq B$ such that 
$\smy m (a,b)<q/2$ for all $b\in B'.$ 
Now take a pair of ${b\ne b'}\in B'$ with 
${b\res \ir0m}={b'\res\ir0m}:$ then  
$\sm(b,b')<q,$ contradiction.
\epF{Claim}

It follows that $x\mapsto G(\vt(x))$ maps any 
\dde class into a countable set of galaxies $G(a).$ 
To code the galaxies by single points, let 
$S(a)=\bigcup_m\ens{b\res m}{b\in G(a)}.$ 
Thus $S(a)\sq\bse$ codes the Polish topological
closure of the galaxy $G(a)$.

\bcl
\lam{cl2}
If\/ $a\yi b\in\aH$ and\/ $\neg\;a\esn b$ then\/ $b$ 
does not belong to the (topological) closure of\/ 
$G(a),$ in particular,\/ $b\res m\nin S(a)$ for some $m$.
\ecl
\bpf 
Take $m$ big enough for $\smp0m(a,b)\ge 2.$ 
Then $s={b\res m}$ does not belong to $S(a)$ because 
any $a'\in G(a)$ satisfies $\sm(a,a')<1$.
\epF{Claim}

Elementary computation shows that the sets 
\dm
{\bf G}=\ans{\ang{a,b}:a\in\aH\land b\in G(a)}
\quad\hbox{and}\quad
{\bf S}=\ans{\ang{a,s}:a\in\aH\land s\in S(a)}\,.
\dm
belong to $\is11,$ but this is not enough to claim 
that $a\mapsto S(a)$ is a Borel map. 
Yet we can change it appropriately to get a Borel map 
with similar properties. 
First of all define the following $\is11$ \er\ 
on $\aH$:
\dm
a\rF b\qquad\hbox{iff}\qquad
e(a)=e(b)\land q(a)=q(b)\land G(a)=G(b)\,.
\dm
(To see that $\rF$ is $\is11$ note that here $G(a)=G(b)$ 
is equivalent to $b\in G(a),$ and that $\bf G$ is $\is11.$)
It follows from Claim~\ref{cl2} and \Kres\  
that there is a $\id11$ function $\mu:\aH\ti\aH\to\dN$ 
such that for any pair of $a\yi b\in\aH$ with $a\nsn b$ 
we have $b\res \mu(a,b)\nin S(a).$ 
Then the set
\dm
R(a)=\ans{b\res \mu(a',b):
a',b\in\aH \land {a\rF a'}\land {a'\nsn b})}\,\sq\,\bse
\dm
does not intersect $S(a),$ for any $a\in\aH,$ hence, 
the $\is11$ set 
\dm
{\bf R}=\ans{\ang{a,s}:a\in\aH\land s\in R(a)}
\dm 
does not intersect ${\bf S}.$ 
Note that by definition $\bf R$ is \ddf invariant \vrt\ 
the 1st argument, \ie, if $a\yi a'\in\aH$ satisfy 
$a\rF a'$ then $R(a)=R(a').$ 
It follows from Lemma~\ref{d20} that there is a 
$\id11$ set ${\bf Q}\sq\aH\ti\bse$ with 
${\bf S}\sq{\bf Q}$ but ${\bf R}\cap{\bf Q}=\pu,$ 
\ddf invariant in the same sense. 
Then the map $a\mapsto Q(a)=\ans{s:{\bf Q}(a,s)}$ is $\id11$. 

\bcl
\lam{cl4}
Suppose that\/ $a\yi b\in\aH.$ 
Then$:$ ${a\rF b}$ implies\/ ${Q(a)=Q(b)}$ and\/ ${a\nsn b}$ 
implies\/ ${Q(a)\ne Q(b)}$.
\ecl
\bpf 
The first statement holds just because $Q$ is 
\ddf invariant. 
Now suppose that $a\nsn b.$ 
Then by definition $s=b\res \mu(a,b)\in R(a),$ hence, 
$s\nin Q(a).$ 
On the other hand, $s\in S(b)\sq Q(b)$.
\epF{Claim}

Define $\tau(x)=Q(\vt(x))$ for $x\in\bn,$ so that $\tau$ 
is a $\id11$ map $\bn\to\cP(\bse)$. 

\bcl
\lam{cl5}
If\/ $x\in\bn$ then\/ 
$T_a=\ans{\tau(y):y\in\eke x}$ is at most countable. 
\ecl
\bpf 
Suppose that $y\yi z\in\eke x.$ 
Then $a=\vt(x),\msur$ $b=\vt(y),$ and $c=\vt(z)$ belong to 
$\aH,$ and $b\yi c\in[a]_{\esn}.$ 
It follows from Claim~\ref{cl4} that if $G(b)=G(c),\msur$ 
$e(b)=e(c),$ and $q(b)=q(c),$ then $Q(b)=Q(c).$ 
It remains to note that $G$ takes only countably many 
values on $\aH\cap \ek a{\esn}$ by Claim~\ref{cl3}. 
\epF{Claim}

Finally note that, if ${x\nE y}\in\bn$ then 
$\vt(x)\yi \vt(y)$ belong to $\aH$ and satisfy 
$\vt(x)\nsn \vt(y),$ hence, $\tau(x)\ne \tau(y)$ by 
Claim~\ref{cl4}. 
Thus, the Borel map $\tau$ witnesses that the given 
\er\ $\rE$ is 
essentially countable by Lemma~\ref{KT}.

\punk{Case 2}
\las{4case2}
 
Thus we suppose that the $\is11$ set 
$\aB=\ran\vt\dif C$ is non-empty. 
Note that, by Claim~\ref{cl1}, there is no non-empty 
$\is11$ \grap\ set $A\sq\aB$.\vom 

Let $\cB_s=\ans{a\in\dn:s\su a}$ for $s\in\bse$ and 
$\cN_u=\ans{x\in\dnn:u\su x}$ for $u\in\nse$ 
(basic open nbhds in $\dn$ and $\dnn$). 

If $A\yi B\sq\dn$ and $m\yi k\in\dN,$ then 
$A\xE km B$ will mean that for any $a\in A$ there is $b\in B$ 
with $\smy{k}(a,b)<2^{-m},$ and conversely, for any $b\in B$ 
there is $a\in A$ with $\smy{k}(a,b)<2^{-m}.$ 
This is not a \er, of course, yet the conjunction of 
$A\xE km B$ and $B\xE km C$ implies $A\xE k{m-1} C$.

$0^m$ will denote the sequence of $m$ zeros. 

To prove that $\esn\reb\rE$ in Case 2, we define 
an increasing sequence of natural numbers 
$0=\yk0<\yk1<\yk2<\dots ,$ and also objects 
$A_s,\:g_s,\:v_s$ for any $s\in\bse,$ 
which satisfy the following list of requirements 
\ref{j1} -- \ref{jx}.  

\ben
\tenu{(\roman{enumi})}
\itla{j1}\msur
if $s\in 2^m$ then $g_s\in 2^{\yk m},$ and 
${s\su t}\imp {g_s\su g_t}$; 

\itla{j2}\msur
$\pu\ne A_s\sq\aB\cap\cB_{g_s},$ $A_s$ is $\is11,$ 
and $s\su t\imp A_t\sq A_s$.

\itla{j3}
if $s\in 2^n$ then $A_{0^n}\xE{\yk n}{n+2} A_s$;

\itla{j5}
if $s\in 2^n,\msur$ $m<n,\msur$ $s(m)=0,$ then 
$\smp{\yk{m}}{\yk{m\plo}}(g_s,g_{0^m})<2^{-m-1}$;

\itla{j6}
if $s\in 2^n,\msur$ $m<n,\msur$ $s(m)=1,$ then 
$|\smp{\yk{m}}{\yk{m\plo}}(g_s,g_{0^m})-\frac1{m+1}|<2^{-m-1}$;

\itla{j7}
if $s\yo t\in 2^n\yt m<n\yt s(m)=t(m),$ then 
$|\smp{\yk{m}}{\yk{m\plo}}(g_s,g_{t})|<2^{-m}$;

\vyk{
\itla8 
sets $A_s$ are enough generic;
}

\itla{j8}
for any $n,$ a certain condition, in terms of the
Gandy -- Harrington forcing,
similar to \ref{d1i,} in \nrf{>smooth} or
\ref{d2i1,} in \nrf{d2split},
related to all sets $A_s\yt s\in 2^n,$ so that, as a 
consequence, $\bigcap_nA_{a\res n}\ne\pu$ 
for any $a\in\dn$;

\itla{j9}
if $s\in2^n$ then $v_s\in\dN^n,$ and $s\su t\imp v_s\su v_t$;

\itla{jx}\msur
$A_s\sq 
\ans{a\in\aB:\vt\obr(a)\cap\cN_{v_s}\ne\pu}$.
\een


We can now accomplish Case~2 as follows. 
For any $a\in\dn$ define 
$F(a)=\bigcup_n g_{a\res n}\in\dn$ 
(the only element satisfying $g_{a\res n}\su F(a)$ 
for all $n$) 
\vyk{
(this is where \ref8 is applied to provide the 
non-emptiness of the intersection)
}%
and $\rho(a)=\bigcup_n v_{a\res n}\in\bn.$ 
It follows, by \ref{jx} and the continuity of $\vt,$ 
that $F(a)=\vt(\rho(a))$ for any $a\in\dn.$ 
Thus the next claim proves that $\rho$ is a 
Borel (in fact, here continuous) 
reduction $\esn$ to $\rE$ and ends Case~2. 

\ble
\lam{cl8'}
The map\/ $F$ reduces\/ $\esn$ to\/ $\esn,$ that is, 
the equivalence\/ 
${a\esn b}\eqv{F(a)\esn F(b)}$ holds for all\/ 
$a\yi b\in\dn$.
\ele
\bpf
By definition 
$\sm(F(a),F(b))=
\tlim_{n\to\iy}\smp0{\yk n}(g_{a\res n},g_{b\res n}).$ 
\vyk{
As $\ga_s=g_{0^n}\sd g_s$ for $s\in 2^n,$ we can 
re-write this as follows:
\dm
\sm(F(a),F(b))=
{\textstyle\tlim_{n\to\iy}}\:\sm(\ga_{a\res n},\ga_{b\res n})\,.
\dm
}
However it follows from \ref{j5}, \ref{j6}, \ref{j7} that 
\dm
|\smp0{\yk n}(g_{a\res n},g_{b\res n})-
\smp0{n}(a\res n,b\res n)|
\,\le\,
{\textstyle\sum_{m<n}}2^{-m}\,<\,2\,.
\dm
We conclude that 
$|\sm(F(a),F(b)) - \sm(a,b)|\le2,$ as required.
\epF{Lemma}

\punk{Construction}
\las{4cons}

The construction of a system of sets satisfying \ref{j1}
-- \ref{jx} goes on by induction. 
To begin with we set $\yk0=0,\msur$ $g_\La=\La$ and 
$A_\La=\aB.$ 
Suppose that, for some $n,$ we have the objects as 
required for all $n'\le n,$ and extend the construction 
on the level $n+1$.

As $A_{0^n}$ is not \grap\ (see above), there is a 
pair of elements $a^0\yi a^1\in A_{0^n}$ such that 
$|\sm(a^0,a^1)-\frac1{n+1}|<2^{-n-2}.$ 
Note that $a^0\res {\yk n}=a^1\res {\yk n}$ by \ref{j1} and 
\ref{j2}, therefore there is $\yk{n\plo}>\yk n$ such that 
$|\smp{\yk n}{\yk{n\plo}}(a^0,a^1)-\frac1{n+1}|<2^{-n-2}.$ 
According to \ref{j3}, for any $s\in2^n$ there exist 
$b^0_s\yi b^1_s\in A_s$ such that  
and $\smy{\yk n}(a^i,b^i_s)<2^{-n-2}$ for $i=0,1;$ 
we can, of course, assume that $b^i_{0^n}=a^i.$ 
Moreover, the number $\yk{n+1}$ can be chosen big enough 
for the following to hold: 
\dm
\smy{\yk{n\plo}}(b^i_s,a^0)< 2^{-n-3}
\quad\text{--- \ for all}\quad s\in 2^n
\quad\text{and}\quad i=0,1.
\eqno(1)
\dm

We let $g_{s\we i}=b^i_s\res {\yk{n\plo}}$ for all  
$s\we i\in 2^{n+1}.$ 
This definition preserves \ref{j1}. 
To check \ref{j5} for $s'=s\we 0\in2^{n+1}$ and $m=n,$ 
note that
\dm
\smp{\yk{n}}{\yk{n\plo}}(g_{s'},g_{0^{n+1}})=
\smp{\yk{n}}{\yk{n\plo}}(b^0_s,a^0)<2^{-n-2}.
\dm 
To check \ref{j6} for $s'=s\we 1\in2^{n+1}$ and $m=n,$ 
note that
\dm
|\smp{\yk{n}}{\yk{n\plo}}(g_{s'},g_{0^{n+1}})-
{\textstyle\frac1{n+1}}| \,\le \,
\smp{\yk{n}}{\yk{n\plo}}(b^1_s,a^1)
+|\smp{\yk n}{\yk{n\plo}}(a^0,a^1)-{\textstyle\frac1{n+1}}|
\,<\,2^{-n-1}.
\dm 

To fulfill \ref{j9}, choose, for any $s\we i\in 2^{n+1},$ 
a sequence $v_{s\we i}\in\dN^{n+1}$ so that $v_s\su v_{s\we i}$ 
and there is $\cN_{v_{s\we i}}\cap \vt\obr(b^i_{s})\ne\pu$. 

Let us finally define the sets $A_{s'}\sq A_s,$ 
for all $s'=s\we i\in 2^{n+1}$ 
(so that $s\in2^n$ and $i=0,1$).  
To fulfill \ref{j2} and \ref{jx}, we begin with 
\dm
A'_{s\we i}=\ans{a\in A_s\cap\cB_{g_{s\we i}}:
\vt\obr(a)\cap\cN_{v_{s\we i}}\ne\pu}\,.
\dm 
This is a $\is11$ subset of $A_s,$ containing $b^i_{s}.$ 
To fulfill \ref{j3}, we define $A_{0^{n+1}}$ to be the set of 
all $a\in A'_{0^{n+1}}$ such that 
\dm
\kaz s'=s\we i\in 2^{n+1}\;\sus b\in A'_{s'}\;
\skl\smy{\yk{n\plo}}(a,b)<2^{-n-3}\skp\;;
\dm
this is still a $\is11$ set containing $b^0_{0^n}=a^0$ by 
(1). 
It remains to define, for any $s\we i\ne 0^{n+1},$ 
$A_{s\we i}$ to be the set of 
all $b\in A'_{s\we i}$ such that 
\dm
\sus b\in A_{0^{n+1}}\;
\skl\smy{\yk{n\plo}}(a,b)<2^{-n-3}\skp\;.
\dm
This ends the definition for the level $n+1$.

\imar{\ref{j8}?}

\vyk{

According to \ref3, for any $s\in2^n$ there exist 
$h^0_s\yi h^1_s\in H_{\yk{n}}(2^{-n-2})$ such that 
\dm
b^0_{s}=\ovp s\sd h^0_s\sd a_0
\quad\hbox{ and }\quad
b^1_{s}=\ovp s\sd h^1_s\sd a_m=
\ovp s\sd h^1_s\sd h\sd a_0
\dm
belong to $A_s.$ 
We can assume that $h^0_{0^n}=h^1_{0^n}=0^\dN$ 
(the constant $0$), so that 
$b^0_{0^n}=a_0,\msur$ $b^1_{0^n}=a_m$. 
Now, as $h^0_s$ and $h^1_s$ belong to 
$H_{\yk{n}}(2^{-n-2}),$ 
there is a number $\yk{n\plo}>\yk{n}$ big enough for 
the following to hold~\footnote
{\ Here Hjorth also adds inequalities 
$\sm_{\ge \yk{n+1}}(h^0_s) < 2^{-n-3}$ and 
$\sm_{\ge \yk{n+1}}(h^1_s) < 2^{-n-3},$ which does 
not seem to be of any use.-- ?}:
\pagebreak[0]%
\ben
\tenu{(\roman{enumi})}
\addtocounter{enumi}8
\itla y
$|\sm(h\res\ir{\yk n}{\yk{n\plo}})-\frac1{n+1}| 
< 2^{-n-2}$.
\een

\vyk{
\dm 
\left. 
\bay{rcll}
\sm_{\ge \yk{n+1}}(h^0_s) &<& 2^{-n-3} &;\\[1ex]

\sm_{\ge \yk{n+1}}(h^1_s) &<& 2^{-n-3} &;\\[1ex]

|\sm(h\res\ir{\yk n}{\yk{n\plo}})-\frac1{n+1}| 
&<& 2^{-n-2} &.
\eay
\right\}
\eqno(\ast)
\dm
}

We let $g_{s\we 0}=b^0_s\res \ir0{\yk{n\plo}}$ and 
$g_{s\we 1}=b^1_s\res \ir0{\yk{n\plo}}.$ 
This definition preserves \ref{7i}. 
To check \ref5 for $s'=s\we 0\in2^{n+1}$ and $m=n,$ 
note that
\dm
\sm(\ga_{s'}\res\ir{\yk n}{\yk{n+1}}) =
\sm(h^0_s\res\ir{\yk n}{\yk{n+1}})<2^{-n-2}
\dm 
since $h^0_s\in H_{\yk{n}}(2^{-n-2}).$ 
To check \ref6 for $s'=s\we 1\in2^{n+1}$ and $m=n,$ 
note that
\dm
|\sm(\ga_{s'}\res\ir{\yk n}{\yk{n+1}})-
{\textstyle\frac1{n+1}}| \,\le \,
\sm(h^1_s\res\ir{\yk n}{\yk{n+1}}) 
+|\sm(h\res\ir{\yk n}{\yk{n+1}})-
{\textstyle\frac1{n+1}}|
\,<\,2^{-n-1}
\dm 
because $h^0_s\in H_{\yk{n}}(2^{-n-2})$ 
and by \ref y.  
The verification of \ref7 for $s'\yi t'\in 2^{n+1}$ and 
$m=n$ is similar. 
}

\vtm

\qeDD{Construction and Theorem~\ref{h}}

\api

\parf{\protect\co equalities}
\las{co}

Recall that the \eqr\ $\fco$ is defined on $\rn$ as follows:
$x\fco y$ iff $x(n)-y(n)\to0$ with $n\to\iy.$
This definition admits a straightforward generalization.

\bdf
[Farah~\cite{f-co}]
\lam{genC}
Suppose that $K$ is a non-empty index set, and 
$\stk{X_k}{d_k}$ is a metric space for any index $k\in K.$
An \eqr\snos
{The letter $\rD$ in this context is due to Farah~\cite{f-co}.
One has to suppress any association with the diagonal, \ie, 
the true equality.}
${\td}=\qqd{X_k}{d_k}{k\in K}$ on 
\index{zzcxkdkk@$\qqd{X_k}{d_k}{k\in K}$}
\index{equivalence relation, ER!zzcxkdkk@$\qqd{X_k}{d_k}{k\in K}$}%
the cartesian product $X=\prod_kX_k$ is defined so that 
$x\td y$ iff $\tlim\, d_n(x(n),y(n))=0,$ where the limit is 
associated with the filter of all finite subsets of $K.$\snos
{Thus $\tlim\, d_n(x(n),y(n))=0$ iff for any $\ve>0$
there exist only finitely many indices $k\in K$ such
that $d_n(x(n),y(n))>\ve$.}

If $K=\dN$ (the most typical case below) then we'll 
write $\xd{X_k}{d_k}$ instead of
\index{equivalence relation, ER!zzcxkdk@$\xd{X_k}{d_k}$}%
\index{zzcxkdkk@$\xd{X_k}{d_k}$}
$\qd{X_k}{d_k}{k\in\dN}$ for the sake of brevity.

We'll be mostly interested in the case when
\bit
\item[\text{\mtho$(\ast)$}]
$X_k$ are Borel sets in Polish spaces $\dX_k,$
and the distance
functions $d_k$ are Borel maps $X_k\ti X_k\to \dR^+,$ not
necessarily equal to the restrictions of Polish metrics
of $\dX_k$.
\eit
Then $\xd{X_k}{d_k}$ is obviously a Borel \eqr\ on
$X=\prod_kX_k$.

The \eqr\ $\xd{X_k}{d_k}$ is {\it nontrivial\/} 
if $\tlis_{k\to\iy}\dia(X_k)>0.$ 
(Otherwise $\xd{X_k}{d_k}$ obviously makes everything 
equivalent.)

A \rit{\co equality} is any \eqr\ of the form
\index{c0equality@\co equality}%
\index{equivalence relation, ER!c0equality@\co equality}%
$\qd{X_k}{d_k}k,$ where all sets $X_k$ are finite.
\edf

Every \co equality is easily a Borel \eqr, more exactly, 
of type $\fp03.$ 
The \eqr\ $\fco$ itself is essentially a \co equality
(see below) ---
this explains the meaning of the term \lap{\co equality}.

The \dd\reb properties of these \er s are largely unknown,
except for the case of \dd\fsg compact metric spaces
$\stk{X_k}{d_k},$
easily reducible to the case of $X_k$ finite (= \co equalities).
This case is presented in this \gla.
We prove that Borel reducibility of a \co equality to another
one implies a stronger additive reducibility of an infinitely
generated \co subequality (Theorem~\ref{ab4}),
show that $\fco$ is a \dd\reb maximal \co equality
(Theorem~\ref{comax}),
prove Theorem~\ref{coT} that shows the turbulence of
\co equalities except those \dd\eqb equivalent to $\Eo$ and $\Et,$
and finally show that the \dd\reb structure of \co equalities
includes a substructure similar to $\stk\pn{\sqa}$
(Theorem~\ref{coAB}).

\punk{Some examples and simple results}
\las{co1}

The following examples show that many typical \eqr s can
be defined in the form of \co equalities.

\bex
\ben
\renu
\itla{coex1}
\lam{co1ex}
Let $X_k=\ans{0,1}$ with $d_k(0,1)=1$ for all $k.$ 
Then clearly the relation 
$\xd{X_k}{d_k}$ on $\dn=\prod_k\ans{0,1}$ 
is just $\Eo$.

\itla{coex2}
Let $X_{kl}=\ans{0,1}$ with $d_{kl}(0,1)=k\obr$ for all 
$k\yi l\in\dN.$ 
Then the relation $\qd{X_{kl}}{d_{kl}}{k,l}$  
on $2^{\dN\ti\dN}=\prod_{k,l}\ans{0,1}$ 
is exactly $\Et$.

\vyk{
\snot{The transformation of $\qd{X_{kl}}{d_{kl}}{k,l}$
to the canonical form of Definition~\ref{genC} can be
carried out by means of an arbitrary bijection of $\dN\ti\dN$
onto $\dN$.}
}

\itla{coex3}
Generally, if $0=n_0<n_1<n_2<\dots$ and $\vpi_i$ is a 
submeasure on $\il i{i+1},$ then let $X_i=\pws{\il i{i+1}}$ 
and $d_i(u,v)=\vpi_i(u\sd v)$ for $u\yi v\sq \il i{i+1}.$ 
Then $\xd{X_{i}}{d_{i}}$ is isomorphic to $\rei,$ where 
\dm
\cI=\Exh(\vpi)=
\ens{x\sq\dN}{\tlim_{n\to\iy}\vpi(x\cap\ir n\iy)=0}
\dm
and $\vpi(x)=\tsup_i\vpi_i(x\cap\il i{i+1})$.

\itla{coex4}
Let, for all $k,$ $X_k=\dR$ with $d_k$ being 
the usual distance on $\dR.$
Then the relation $\qd{X_k}{d_k}k$ on $\rn$ 
is just $\fco$.\qed
\een
\eeX

\ble
[ Farah~\cite{f-co} with a reference to Hjorth]
\lam{coact}
\sloppy
Every \co equality\/ ${\td}={\xd{X_k}{d_k}}$ 
is induced by a continuous action of a Polish group.
\ele

The domain $\dX=\prod_kX_k$ of $\td$ 
is considered with the product topology.

\bpf[sketch] 
For any $k$ let $S_k$ be the (finite) group of all 
permutations of $X_k,$ with the distance  
$\rho_k(s,t)=\tmax_{x\in X_k}\,d_k(s(x),t(x)).$ 
Then
\dm
\dG=\ens{g\in{\TS\prod_kS_k}}
{\tlim_{k\to\iy}\rho_k(g_k,e_k)=0}\,,
\quad\text{where $e_k\in S_k$ is the identity}\,,
\dm
is easily a subgroup of $\prod_kS_k.$
Moreover, the distance  
$d(g,h)=\tsup_k\rho_k(g_k,h_k)$ converts $\dG$ into a 
Polish group, the natural action of which on $\dX,$ 
that is, $(g\app x)_k=g_k(x_k)\zd\kaz k,$  
is continuous and induces $\td$.
\epf

Let us finally show that the case of \dd\fsg compact spaces
$X_k$ does not give anything beyond the case of \co equalities.

\ble
\lam{sc2co}
Suppose that in the assumptions of\/ \ref{genC}(\mast)
$\stk{X_k}{d_k}$ are \dd\fsg compact spaces.
Then\/ $\xd{X_k}{d_k}$ is\/
\dd\eqb equivalent to a\/ \co equality.
\ele
\bpf
Suppose that all spaces $X_k$ are compact.
Then for any $k$ there exists a finite \dd{\frac1k}net
$X'_k\sq X_k.$  
Given $x\in X=\prod_kX_k,$ we define $\vt(x)\in X'=\prod_kX'_k$
so that $\vt(x)(k)$ is the \dd{d_k}closest to $x(k)$
element of $X'_k$
(or the least, in the sense of a fixed ordering of $X'_k,$ of
such closest elements, whenever there exist two or more of them)
for each $k.$ 
Then $\vt$ is obviously a Borel reduction of $\xd{X_k}{d_k}$
to the \co equality $\xd{X'_k}{d_k}$.

The general \dd\fsg compact case can be reduced to the compact
case by the same trick as in the beginning of the proof of
Lemma~\ref{redx}.
\epf

\punk{\co equalities and additive reducibility}
\las{aco}

The structure of \co equalities tend to be  
connected more with the additive reducibility $\rea$ 
than with the general Borel reducibility.\snos 
{See \nrf{b-c} on $\rea$ and the associated relations
$\reas$ and $\eqa$.}
In particular, we have 

\ble
\lam{<co}
For any\/ \co equality\/ ${\td}=\qd{X_k}{d_k}k,$ 
if\/ $\td'$ is a Borel \eqr\ on\/ a set of the form\/ 
$\prod_kX'_k$ with finite nonempty factors\/ $X'_k,$  
and\/ ${\td'}\rea{\td}$ then\/ $\td'$ itself is a\/ 
\co equality.
\ele
\bpf
Let a sequence $0=n_0<n_1<n_2<\dots$ and a collection 
of maps $H_i:X'_i\to\prod_{n_i\le k<n_{i+1}}X_k$ 
witness ${\td}'\rea{\td}.$ 
For $x'\yi y'\in X'_i$ put 
\dm
d'_i(x',y')=
\tmax_{n_i\le k<n_{i+1}}d_k(H_i(x')_k,H_i(y')_k)\,.
\dm
Then easily ${\td'}=\qd{X'_k}{d'_k}k$. 
\epf

It is perhaps not true that ${\td}\reb{\td'}$ implies
${\td}\rea{\td'}$ for any pair of \co equal\-ities.
Yet a somewhat weaker statement holds by
the next theorem of Farah~\cite{f-co}.

\bte
\lam{ab4}
If\/ $\td=\qd{X_k}{d_k}k$ and\/ $\td'=\qd{X'_k}{d'_k}k$
are\/ \co equalities and\/ ${\td}\reb{\td'}$
then there is an infinite set\/ $A\sq\dN$ such that the\/
\co equality\/ $\td_A=\qqd{X_k}{d_k}{k\in A}$ satisfies\/ 
${\td_{A}}\rea{\td'}$.
\ete
\bpf
Define $X_C=\prod_{k\in C}X_k$ and $X'_C=\prod_{k\in C}X'_k$
for any set $C\sq\dN,$ and
$d'_C(x,y)=\tsup_{k\in C}d'_k(x(k),y(k))$ for all $x,y\in X'.$ 
Suppose that
\dm
\TS
\vt:X=\prod_{k\in\dN}X_k\to X'=\prod_{k\in\dN}X'_k
\dm
is a Borel reduction of $\td$ to $\td'.$
Then there exists an infinite set $A'\sq\dN$ such that
$\qqd{X_k}{d_k}{k\in A'}\ren {\td'}$
(via a continuous reduction) --- this can be proved analogously
to the second claim of Lemma~\ref{l:bc}.
Thus it can be accumed from the beginning that $\vt$ is a
continuous reduction of $\td$ to $\td'$.

To extract an additive reduction, we employ a version of
the stabilizers construction used in the proof of
Theorem~\ref{t:>s}\ref{t>s1}.
In fact our task here is somewhat simpler because the
given countinuity of $\vt$ allows us to avoid the Cohen
genericity arguments.

Put $[s]=\ens{x\in X}{x\res u=s}$ 
for any $u\sq\dN$ and $s\in X_u.$
Consider the closed set $W=\bigcap_{i\in\dN}[s_i]$
of all points $x\in X$ such that
$x\res{\iv{n_i}{n_{i+1}}}=s_i$ for all $i.$ 
Arguing approximately as in the proof of
Theorem~\ref{t:>s}\ref{t>s1}, we can define an increasing
sequence $0=k_0=n_0<k_1<n_1<k_2<n_2<\dots$
and elements $s_i\in X_{\iv{n_{i}}{n_{i+1}}}$ such that 
for all $u,v\in X_{\ix0{n_i}}$
and all $x,y\in X_{\iry{n_{i+1}}}$ satisfying 
$x\res{\iv{n_j}{n_{j+1}}}=y\res{\iv{n_j}{n_{j+1}}}=s_j$
for all indices $j>i$ and
$u\res{\iv{n_j}{n_{j+1}}}=v\res{\iv{n_j}{n_{j+1}}}=s_j$
for all indices $j<i,$\snos
{Note that under this assumption the points
$u\cup s_i\cup x\yt u\cup s_i\cup y\yt v\cup s_i\cup x$
mentioned in \ref{ss3=a}, \ref{ss3=b}, belong to $W.$}
the following holds:
\ben
\tenu{(\alph{enumi})}
\itla{ss3=a}
$\vt(u\cup s_i\cup x)\res \ir0{k_{i+1}}=
\vt(u\cup s_i\cup y)\res \ir0{k_{i+1}}$, \ 
and

\itla{ss3=b}
$d_{\iry{k_{i+1}}}(\vt(u\cup s_i\cup x),
\vt(v\cup s_i\cup x))<\frac1i$.
\een

Put $A=\ens{n_i}{i\in\dN}$
and fix any $z\in X_A.$
For any $i,$ if $\xi\in X_{n_i}$ then define $z^{i\xi}\in W$ so
that $z^{i\xi}(n_i)=\xi,$ $z^{i\xi}(n_j)=z(n_j)$ for all
$j\ne i,$ and $z^{i\xi}\res {\iv{n_j}{n_{j+1}}}=s_j$ for all $j.$
If $x\in X_{A}$ then define $H(x)\in X'$ as follows: 
\dm
H(x)\res{\ir{k_i}{k_{i+1}}}=\vt(z^{i,x(n_i)})\res{\ir{k_i}{k_{i+1}}}
\quad\text{for every}\quad i\in\dN.
\eqno(1)
\dm
Clearly $H$ is a continuous map from $X_{A}$ to $X'$
(in the sense of the Polish product topologies).
Moreover for any $i$ the value $H(x)\res{\ir{k_i}{k_{i+1}}}$
obviously depends only on $x(n_i).$
Thus to accomplish the proof of the theorem we need
only to prove that $H$ is a reduction of $\td_A$ to $\td'.$

For any $x\in X_A$ define $f(x)\in W$ so that $f(x)\res A=x$
and $f(x)\res {\iv{n_j}{n_{j+1}}}=s_j$ for all $j.$
Then $f$ is a reduction of $\td_A$ to $\td,$ therefore it
suffices to prove that $\vt(f(x))\td' H(x)$ for every $x\in X_A.$
For an arbitrary $i\ge1,$ let us show that  
\dm
d'_{\ir{k_i}{k_{i+1}}}(\vt(f(x)),H(x))\le 1/i\,.
\eqno(2)
\dm
The key fact is that by the construction
the elements $a=f(x)$ and $b=z^{i,x(n_i)}$ of $W$ satisfy
$a\res {\iv{n_j}{n_{j+1}}}=b\res {\iv{n_j}{n_{j+1}}}=s_j$
for all $j$ and in addition $a(n_i)=b(n_i)=x(n_i).$
Define an auxiliary element $c\in W$ by
\dm
c\res{\ix0{n_i}}=a\res{\ix0{n_i}}\quad\text{and}\quad
c\res{\iry{n_{i+1}}}=b\res{\iry{n_{i+1}}}.
\dm
Then   
$
d'_{\ir{k_i}{k_{i+1}}}(\vt(b),\vt(c))\le \frac1i  
$
by \ref{ss3=b},
and $\vt(a)\res{\ir{k_i}{k_{i+1}}}=\vt(c)\res{\ir{k_i}{k_{i+1}}}$ 
by \ref{ss3=a}.
(Note that \ref{ss3=b} is applied in fact for the value
$i-1$ instead of $i.$)
It follows that 
$
d'_{\ir{k_i}{k_{i+1}}}(\vt(a),\vt(b))\le \frac1i\,. 
$
However
$H(x)\res{\ir{k_i}{k_{i+1}}}=\vt(b)\res{\ir{k_i}{k_{i+1}}}$
by (1).
This proves (2) as required.
\epf

\punk{A maximal \co equality}
\las{co-max}

We define $\rdm=\xd{X_k}{d_k},$ where
\index{zzDmax@$\rdm$}%
\index{equivalence relation, ER!zzDmax@$\rdm$}%
$X_k=\ans{0,\frac1k,\frac2k,\dots,1}$ and  
$d_k$ is the distance on $X_k$ inherited from the real line $\dR.$
The next theorem says that
$\rdm$ is \dd\reb maximal among all \co equalities.
The proof will show that in fact ${\td}\rea\rdm$ in \ref{comax2},
in the sense of the additive reducibility.\pagebreak[2]

\bte
[{{\rm Farah~\cite{f-co} with a reference to Oliver}}]
\label{comax}
\ben
\tenu{{\rm(\roman{enumi})}}
\itla{comax1} 
$\rdm\eqb{\fco}\,;$  
\imar{comax}%

\itla{comax2} 
if\/ $\td$ is a\/ \co equality then\/ 
${\td}\reb\rdm$. 
\een
\ete

It follows from \ref{comax1} and Lemma~\ref{co=d} that
$\rdm\eqb{\rzo}$. 

\bpf
\ref{comax1} 
It is clear that $\rdm$ is the same as ${\fco}\res\dX,$ 
where $\dX\sq\rtn$ is defined as in the proof of 
Lemma~\ref{co=d}, where it is also shown that  
${\fco}\eqb{{\fco}\res\dX}$. 

\ref{comax2} 
To prove ${\td}\reb\rdm,$ it suffices by \ref{comax1}
to show that ${\td}\reb{\fco}.$
The proof is based on the following: 

\bcl
\lam{comax'}
Any finite\/ \dd nelement metric space\/ $\stk Xd$ is 
isometric to an\/ \dd nelement subset of\/ 
$\stk{\dR^n}{\rho_n},$ where\/ $\rho_n$ is the distance 
on\/ $\dR^n$ 
defined by\/ $\rho_n(x,y)=\tmax_{i<n}|x(i)-y(i)|$.
\ecl
\bpc
Let $X=\ans{x_1,\dots,x_n}.$ 
It suffices to prove that for any $k\ne l$ there is a 
set of reals $\ans{r_1,\dots,r_n}$ 
such that $|r_k-r_l|=d(x_k,x_l)$ and 
\bit
\item[$(\ast)$]
$|r_i-r_j| \le d_{ij}=d(x_i,x_j)$ for all $i\zi j$. 
\eit
We can assume that $k=1$ and $l=n.$ 

{\sl Step 1\/}. 
There is a least number $h_1\ge0$ such that 
$(\ast)$ holds for the reals 
$\ans{r_i}=\ans{\underbrace{0,0,\dots,0}_{n-1\;\,\text{times}},h}$ 
for any $0\le h\le h_1.$ 
Then, for some index $k\yt 1\le k<n,$ we have   
$h_1-0=d_{kn}$ exactly. 
Suppose that $k\ne1;$ 
then it can be assumed that $k=n-1.$\vom 

{\sl Step 2\/}. 
Similarly, there is a least number $h_2\ge0$ such that 
$(\ast)$ holds for the reals 
$\ans{r_i}=
\ans{\underbrace{0,0,\dots,0}_{n-2\;\,\text{times}},h,h_1+h}$ 
for any $0\le h\le h_2.$ 
(For example, $h_2=0$ in the case when on step 1 we have one
more index $k'\ne k$ such that $h_1=d_{k'\,n}.$)
Then, for some $k,\nu\yt 1\le k<n-1\le\nu\le n,$ 
we have $h_2-0=d_{k\nu}$ exactly. 
Suppose that $k\ne1;$ 
then it can be assumed that $k=n-2.$\vom 

{\sl Step 3\/}. 
Similarly, there is a least number $h_3\ge0$ such that 
$(\ast)$ holds for the reals 
$\ans{r_i}=\ans{\underbrace{0,0,\dots,0}_{n-3\;\,\text{times}},
h,h_2+h,h_1+h_2+h}$ 
for any $0\le h\le h_3.$ 
Then again, for some $k,\nu\yt 1\le k<n-2\le\nu\le n,$ 
we have $h_3-0=d_{k\nu}$ exactly. 
Suppose that $k\ne1;$ 
then it can be assumed that $k=n-3.$\vom 

{\sl Et cetera\/}.\vom

This process ends, after a number $m$ ($m<n$) steps, 
in such  a way that the index $k$ obtained at the final 
step is equal to $1.$
Then $(\ast)$ holds for the numbers 
$\ans{\underbrace{0,0,\dots,0}_{n-m\;\,\text{times}},
r_{n-m+1},r_{n-m+1},\dots,r_n},$ 
where $r_{n-m+j}=h_m+h_{m-1}+\dots+h_{m-j+1}$ for each 
$j=1,\dots,m.$ 
Moreover it follows from the construction that there is 
a decreasing sequence $n=k_0>k_1>k_2>\dots>k_\mu=1$ 
($\mu\le m$) 
such that $r_{k_{i}}-r_{k_{i+1}}=d_{k_{i+1},k_{i}}$ 
exactly for any $i.$ 
Then $d_{1n}\le \sum_i r_{k_{i}}-r_{k_{i+1}}$ by the 
triangle inequality. 
But the right-hand side is a part of the sum 
$r_n=h_1+\dots+h_m,$ and hence $r_n\ge d_{1n}.$
On the other hand we have $r_n\le d_{1n}$ by \mast.
We conclude that $r_n=d_{1n},$
\vyk{
It follows that, cutting the construction at an 
appropriate step $m'\le m$  
(and taking an appropriate value of $h\le h_{m'}$), 
we obtain a sequence of numbers 
$r_1=0\le r_2\le\dots\le r_{n-1}\le r_n$ 
still satisfying $(\ast)$ and satisfying 
$r_n=r_n-r_0=d_{1n}.$ 
This ends the proof.
}%
as required.
\epc

We come back to the proof of \ref{comax2}, that is,
${\td}\reb{\fco}$ for
an arbitrary \co equality ${\td}=\xd{X_k}{d_k}$ on
$\dX=\prod_{k\in\dN}X_k,$ where each 
$\stk{X_k}{d_k}$ is a finite metric space.
Let $n_k$ be the number of elements in $X_k.$ 
Let, by the claim, $\eta_k:X_k\to \dR^{n_k}$ be an 
isometric embedding of $\stk{X_k}{d_k}$ into 
$\stk{\dR^{n_k}}{\rho_{n_k}}.$ 
It easily follows that the map 
$\vt(x)= \eta_0(x_0)\we\eta_1(x_1)\we\eta_2(x_2)\we\dots$  
(from $\dX$ to $\rtn$) 
reduces ${\td}$ to $\fco$.
\epF{Theorem~\ref{comax}}

\punk{Classification}
\las{coclass}

Recall that for a metric space $\stk Ad,$ a rational 
$q>0,$ and $a\in A,$ the galaxy $\gal qAa$ is the set of all 
$b\in A$ which can be connected with $a$ by a finite 
chain $a=a_0,a_1,\dots,a_n=b$ with $d(a_i,a_{i+1})<q$ 
for all $i.$ 
Define, for $r>0,$ 
\dm
\da(r,A)=\tinf{
\ens{q\in\dqp}{\sus a\in A\:(\dia{(\gal qAa)}\ge r)}}
\dm 
(with the understanding that here $\tinf\pu=\piy$), 
and
\dm
\Da(A)=\ens{d(a,b)}{a\ne b\in A}\,,\quad
\text{so that}\quad 
\dia A=\tsup(\Da(A)\cup\ans0)\,.
\dm

Now suppose that ${\td}=\xd{X_k}{d_k}$ is a \co equality 
on $\dX=\prod_{k\in\dN}X_k.$ 
The next theorem of Farah~\cite{f-co} 
shows that basic properties of ${\td}$ in the
\dd\reb structure of Borel \er s are determined
by the following two conditions:

\ben
\tenu{{\rm(co\arabic{enumi})}}
\itla{*1}
$\tlii_{k\to\iy}\da(r,X_k)=0$ for some $r>0$.\enuci

\itla{*2}
$\kaz\ve>0\;\sus \ve'\in(0,\ve)\;\susi k\;
\skl\Da(X_k)\cap\ir{\ve'}\ve\ne\pu\skp$.
\een
Clearly \ref{*1} implies both the nontriviality 
of $\xd{X_k}{d_k}$ and \ref{*2}.

\bte
\lam{coT}
Let\/ ${\td}=\qd{X_k}{d_k}k$ be a nontrivial 
\co equality. 
Then\/$:$
\imar
{Comment upon turbulent in \ref{c03}.}%
\ben
\tenu{{\rm(\roman{enumi})}}
\itla{c01}
if\/ \ref{*2} fails
{\rm(then \ref{*1} also fails)}
then\/ ${\td}\eqb\Eo\;;$

\itla{c02}
if\/ \ref{*1} fails but\/ \ref{*2} holds then\/ 
${\td}\eqb\Et\;;$

\itla{c03}
if\/ \ref{*1} holds
{\rm(then \ref{*2} also holds)}
then there exists a turbulent \co equality\/ ${\td}'$
satisfying\/ $\Eo\rebs{\td}'$ and\/ ${\td}'\reb{\td}$.%
\een
\ete

Thus any nontrivial 
\co equality ${\td}$ \dd\reb contains a turbulent 
\co equality $\td'$ with $\Et\rebs{\td'}$   
unless $\td$ is \dd\eqb equivalent to either $\Eo$ or $\Et.$
In addition, \ref{*1} is necessary for the turbulence 
of $\td$ itself and sufficient for a turbulent 
\co equality ${\td'}\reb{\td} $ to exist.
The proof will show that in fact $\reb$ can be improved
to $\rea$ in the theorem.

\bpf
\ref{c01}
To show that $\Eo\reb{\td}$ note that, by the 
nontriviality of $\td,$ there exist: a number 
$p>0,$ an increasing sequence $0=n_0<n_1<n_2<\dots\;,$ 
and, for any $i,$ a pair of elements 
$x_{n_i}\yi  y_{n_i}\in X_{n_i}$ with 
$d_{n_i}(x_{n_i},y_{n_i})\ge p.$ 
For $n$ not of the form $n_i$ fix an arbitrary 
$z_n\in X_n.$ 
Now, if $a\in\dn,$ then define $\vt(a)\in\prod_kX_k$ 
so that $\vt(a)(n)=z_n$ for $n$ not of the form $n_i,$ 
while $\vt(a)(n_i)=x_{n_i}$ or $=y_{n_i}$ if resp.\ 
$a_i=0$ or $=1.$ 
This map $\vt$ witnesses $\Eo\reb{\td}$.

Now prove that ${\td} \reb\Eo.$ 
As \ref{*2} fails, there is $\ve>0$ such that for each 
$\ve'$ with $0<\ve'<\ve$ we have only finitely many $k$ 
with the propery that $\ve'\le d_k(\xi,\eta)<\ve$ 
for some $\xi\yi\eta\in X_k.$ 
Let $G_k$ be the (finite) set of all 
\dd{\frac\ve2}galaxies in $X_k,$ and let 
$\vt:\dX=\prod_kX_k\to G=\prod_kG_k$ be defined as 
follows:
for every $k,$ $\vt(x)(k)$ is that galaxy in $G_k$
to which $x(k)$ belongs. 
Let $\rE$ be the \dd Gversion of $\Eo,$ that is, if 
$g\yi h\in G$ then 
$g\rE h$ iff $g(k)=h(k)$ for all but finite $k.$ 
As easily $\rE\reb\Eo,$ it suffices to demonstrate that 
${\td}\reb\rE$ via $\vt.$ 

Suppose that $x\yi y\in \dX$ and $\vt(x)\rE\vt(y)$ and 
prove $x\td  y$ (the nontrivial direction). 
Suppose towards the contrary that $x\ntd y,$
so that there is a number 
$p>0$ with $d_k(x(k),y(k))>p$ for infinitely many $k.$ 
We can assume that $p<\frac\ve2.$ 
On the other hand, as $\vt(x)\rE\vt(y),$ there is 
$k_0$ such that $x(k)$ and $y(k)$ belong to one and the 
same \dd{\frac\ve2}galaxy in $X_k$ for all $k>k_0.$ 
Then, for any $k>k_0$ with $d_k(x(k),y(k))>p$ 
(and hence for infinitely many indices $k$) 
there exists an element $z_k\in X_k$ in the same 
galaxy such that $p<d_k(x(k),z_k)<\ve,$ but this is a 
contradiction to the choice of $\ve$ 
(indeed, take $\ve'=p$).\vom

\ref{c02} 
First prove that if \ref{*2} holds then 
$\Et\reb{\td}.$
It follows from \ref{*2} that there exist: 
an infinite sequence $\ve_1>\ve_2>\ve_3>\dots>0,$ 
for any $i$ an infinite set $J_i\sq\dN,$ 
and for any $j\in J_i$ a pair of elements 
$x_{ij}\yi y_{ij}\in X_j$ with 
$d_j(x_{ij},y_{ij})\in\ir{\ve_{i+1}}{\ve_i}.$ 
We may assume that the sets $J_i$ are pairwise 
disjoint. 
Then the \co equality 
${\td'}=\qqd{\ans{x_{ij},y_{ij}}}{d_j}
{i\in\dN,\:j\in J_i}$ 
satisfies both ${\td'}\reb{\td} $ and ${\td'}\isi\Et$ 
(isomorphism via a bijection between the underlying sets). 

Now, assuming that, in addition, \ref{*1} fails, 
we show that ${\td} \reb\Et.$ 
For all $k\yi n\in\dN$ let $G_{kn}$ be the (finite) 
set of all \dd{\frac1n}galaxies in $X_k.$ 
For any $x\in \dX=\prod_iX_i$ define 
$\vt(x)\in G=\prod_{k,n}G_{kn}$ so that for any $k,n$ 
$\vt(x)(k,n)$ is that \dd{\frac1n}galaxy in 
$G_{kn}$ to which $x(k)$ belongs (for all $k\yi n$). 
The \eqr\ $\rE$ on $G,$ defined so that 
\dm
g\rE h\quad\text{ iff }\quad
\kaz n\:\kazi k\:(g(k,n)=h(k,n))\qquad
(g\yo h\in G)
\dm
is obviously $\reb\Et,$ so it suffices 
to show that ${\td} \reb\rE$ via $\vt.$ 
Suppose that $x\yi y\in \dX$ and $\vt(x)\rE\vt(y)$ 
and prove $x\td  y$ (the nontrivial direction). 
Otherwise there is some $r>0$ with $d_k(x(k),y(k))>r$ 
for infinitely many indices $k.$ 
As \ref{*1} fails for this $r,$ there is $n$ big 
enough for $\da(r,X_k)>\frac1n$ to hold for almost 
all $k.$ 
Then, by the choice of $r,$ we have 
$\vt(x)(k,n)\ne\vt(y)(k,n)$ for infinitely many $k,$ 
hence, $\vt(x)\nE\vt(y),$ contradiction.\vom

\ref{c03} 
Fix $r>0$ with $\tlii_{k\to\iy}\da(r,X_k)=0.$ 
For any increasing sequence $n_0<n_1<n_2<\dots$ 
we have $\qd{X_{n_k}}{d_{n_k}}k\reb{\td}.$ 
Therefore it can be assumed that $\tlim_{k}\da(r,X_k)=0,$ 
and further that $\da(r,X_k)<\frac1k$ for all $k.$
(Otherwise choose an appropriate subsequence.)
Then every set $X_k$ contains a \dd{\frac1k}galaxy 
$Y_k\sq X_k$ such that $\dia Y_k\ge r.$ 
As easily $\xd{Y_k}{d_k}\reb{\td},$ 
the following lemma suffices to prove \ref{c03}.

\ble
\lam{coturb}
Suppose that\/ $r>0$ and each\/ $X_k$ is a\/ 
\dd{\frac1k}galaxy and\/ $\dia(X_k)\ge r.$  
Then the\/ \co equality\/ ${\td} =\qd{X_k}{d_k}k$
is turbulent and satisfies\/ $\Et\reb{\td}$.
\ele
\bpf
We know from the proof of \ref{c03} above 
that $\Et\reb{\td} .$
Now prove that the natural action of the Polish group 
$\dG$ defined as in the proof of Lemma~\ref{coact} is 
turbulent under the assumptions of the lemma. 

That every \dd\td  class is dense in $\dX=\prod_kX_k$ 
(with the product topology on $\dX$) is an easy exercise. 
To see that every \dd\td  class $[x]_{\td }$ also is 
meager in $\dX,$ note that by the assumptions of the lemma 
any $X_k$ contains a pair of elements $x'_k\yi x''_k$ with 
$d_k(x'_k,x''_k)\ge r.$ 
Let $y_k$ be one of $x'_k\yi x''_k$ which is 
\dd{d_k}fahrer than $\frac r2$ from $x_k.$ 
The set $Z=\ens{z\in \dX}{\susi k\;(z(k)=y_k)}$ 
is comeager in $\dX$ and disjoint from $[x]_{\td }.$ 

It remains to prove that local orbits are somewhere 
dense.
Let $G$ be an open nbhd of the neutral element in $\dG$ and 
$\pu\ne X\sq\dX$ be open in $\dX.$ 
We can assume that, for some $n,$ $G$ is the 
\dd{\frac1n}ball around the neutral element in $\dG$ while 
$X=\ens{x\in\dX}{\kaz k<n\;(x(k)=\xi_k)},$ where 
elements $\xi_k\in X_k\yt k<n,$ are fixed. 
It is enough to prove that all local orbits, \ie\
equivalence classes of $\ler GX,$ are dense subsets of $X.$ 
Consider an open set 
$Y=\ens{y\in\dX}{\kaz k<m\;(y(k)=\xi_k)}\sq X,$ where 
$m>n$ and elements $\xi_k\in X_k\yt n\le k<m,$ are fixed 
in addition to the above. 

Let $x\in X.$ 
Then $x(k)=\xi_k$ for all $k<n.$ 
Let $n\le k<m.$ 
The elements $\xi_k$ and $x(k)$ belong to $X_k,$ 
which is a \dd{\frac1k}galaxy, therefore, there is 
a chain, of a length $\ell(k),$ 
of elements of $X_k,$ which connects $x(k)$ 
to $\xi_k$ so that every step within the chain has 
\dd{d_k}length $<\frac1k.$ 
Then there is a permutation $g_k$ of $X_k$ 
such that $g_k^{\ell(k)}(x(k))=\xi_k,$ 
$g_k(\xi_k)=x(k),$ and 
$d_k(\xi,g_k(\xi))<\frac1k$ for all $\xi\in X_k.$

In addition let $g_k$ be the identity on $X_k$
whenever $k<n$ or $k\ge m.$ 
This defines an element $g\in\dG$ which obviously 
belongs to $G.$
Moreover, the set $X$ is \dd ginvariant 
and $g^\ell(x)\in U,$ where 
$\ell=\prod_{k=n}^{m-1}\ell(k).$ 
It follows that $x\ler GX g(x),$ as required.
\epF{Lemma}

\epF{Theorem~\ref{coT}}

\punk{LV-equalities}
\las{lv}

By Farah, an {\it\lv equality\/} is a \co equality 
${\td} =\qd{X_k}{d_k}k$ satisfying 
\dm
\kaz m\;\kaz\ve>0\;\kazi k\;
\kaz x_0,\dots,x_m\in X_k\;
\skl
d_k(x_0,x_m)\le\ve+\tmax_{j<m}d_k(x_j,x_{j+1}) 
\skp.
\eqno(\ast)
\dm
In other words, the metrics involved are postulated 
to be \rit{asymptotically close} to ultrametrics. 
This sort of \co equalities was first considered by 
Louveau and Velickovic \cite{lv}.

\bup
\lam{lvex}
Put $X_k=\ans{1,2,\dots,2^k}$ and
$d_k(m,n)=\frac{\log_2(|m-n|+1)}{k}$ for all $k$ and
$1\le m,n\le2^k.$
Prove that $\xd{X_k}{d_k}$ is an \lv equality and 
satisfies \ref{*1} of \nrf{coclass}. 
\eup

\vyk{
The following simple fact is analogous to 
Lemma~\ref{<co}.

\ble
\lam{<lv}
For any\/ \lv equality\/ $\td ,$ 
if\/ $\td'$ is a Borel \er\ on\/ a set\/ 
$\prod_kX'_k$ (with finite nonempty\/ $X'_k$) 
and\/ ${\td'}\rea{\td} $ then\/ $\td'$ is an\/ 
\lv equality.\qeD
\ele
}

The next theorem of Louveau and Velickovic \cite{lv} 
is a major application of \co equalities.
One of its corollaries is that there exist big 
families of mutually irreducible Borel \eqr s, see
below. 

\bte
\lam{coAB}
Let\/ ${\td} =\qd{X_k}{d_k}k$ be an\/ 
\lv equality satisfying\/ \ref{*1} of \nrf{coclass}. 
Then we can associate, with each infinite set\/ $A\sq\dN,$ 
an\/ \lv equality\/ ${\td_A}\rea{\td} $ such that for 
all\/ $A\yi B\sq\dN$ the following are equivalent$:$
\ben
\tenu{{\rm(\roman{enumi})}}
\itla{AB1} 
$A\sqa B$ \ {\rm(that is, $A\dif B$ is finite)}$;$

\itla{AB2} 
${\td_A}\rea{\td_B}$ \ {\rm(the additive reducibility)}$;$


\itla{AB3}
${\td_A}\reb{\td_B}$.
\een
\ete 
\bpf 
Since $\td $ is turbulent, the necessary turbulence 
condition \ref{*1} of \nrf{coclass} holds.
Moreover, 
as in the proof of Theorem~\ref{coT} (part~\ref{c03}), 
we can assume that it takes the following special form 
for some $r>0$: 
%
\ben
\tenu{(\arabic{enumi})}
\itla{lv+}
Each $X_k$ is a 
\dd{\tmin\ans{\frac r{2},\frac1{k+1}}}galaxy 
and $\dia(X_k)\ge4r$.
\een
The intended transformations
(reduction to a certain 
infinite subsequence of spaces $\stk{X_k}{d_k},$ 
and then of each $X_k$ to a suitable galaxy $Y_k\sq X_k$) 
preserve the \lv condition (\mast), of course.
Moreover, we can assume that (\mast) 
holds in the following special form: 
%
\ben
\tenu{(\arabic{enumi})}
\addtocounter{enumi}1
\itla{lv5}\msur
$d_k(x_0,x_{\mu_k})\le
\frac1{k+1}+\tmax_{i<\mu_k}d_k(x_i,x_{i+1})$ 
whenever $x_0,\dots,x_{\mu_k}\in X_k,$ 
where $\mu_k=\prod_{j=0}^{k-1}\#(X_j)$ and $\#X$ is the number
of elements in a finite set $X.$
\een
(For if not then take a suitable subsequence once again.)

We can derive the following important consequence:

\ben
\tenu{(\arabic{enumi})}
\addtocounter{enumi}2
\itla{lv4}\msur
For any 
$k$ there is a set $Y_k\sq X_k$ having exactly $\#(Y_k)=\mu_k$
elements and
such that $d_k(x,y)\ge r$ for all $x\ne y$ in $Y_k$.
\een

To prove this note that by \ref{lv+} there is a 
set $\ans{x_0,\dots,x_m}\sq X_k$ such that 
$d_k(x_0,x_m)\ge 4r$ but 
$d_k(x_i,x_{i+1})<r$ for all $i.$ 
We may assume that $m$ is the least possible length 
of such a sequence $\ans{x_i}.$ 
Define a subsequence $\ans{y_0,y_1,\dots,y_n}$ 
of $\ans{x_i},$ the number $n\le m$ will be specified 
in the course of the construction. 
\ben
\tenu{\alph{enumi})}
\item
Put $y_0=x_0.$

\item
If $y_j=x_{i(j)}$ has been defined, and there is an index
$l>i(j)\yt l\le m,$ such that $d_k(y_j,x_l)\ge r,$ 
then let $y_{j+1}=x_l$ for the least such $l.$

Note that in this case $d_k(y_j,y_{j+1})<2r,$
for otherwise $d_k(y_j,x_{l-1})>r$ because
$d_k(x_{l-1},x_l)<r$.

\item
Otherwise put $n=j$ and stop the construction.
\een

By definition $d_k(y_j,y_{j+1})\ge r$ for all $j<n,$ 
moreover, $d_k(y_{j'},y_{j+1})\ge r$ for any 
$j'<j$ by the minimality of $m.$ 
Thus $Y_k=\ens{y_j}{j\le n}$ satisfies 
$d_k(x,y)\ge r$ for all $x\ne y$ in $Y_k.$ 
It remains to prove that $n\ge\mu_k.$
Suppose otherwise.
Add $y_{n+1}=x_m$ as an extra term.
Then $d_k(x_0,x_m)=d_k(y_0,y_{n+1})\le 3r$ by \ref{lv5}
because $d_k(y_j,y_{j+1})<2r$ (see above).
However we know that $d_k(x_0,x_m)\ge 4r,$ contradiction.
This proves \ref{lv4}.

In continuation of the proof of the theorem,
define $\td_A=\qqd{X_k}{d_k}{k\in A}$ for any $A\sq\dN.$
Thus $\td_A$ is essentially a \co equality on
$\prod_{k\in A}X_k.$
The direction $\ref{AB1}\imp\ref{AB2}\imp\ref{AB3}$
is routine.
Thus it remains to prove $\ref{AB3}\imp \ref{AB1}.$
In view of Theorem~\ref{ab4}, it is enough to prove the
following  lemma.

\ble
\lam{ab21}
If\/ $A,B\sq\dN$ are infinite and disjoint then\/
${\td_A}\rea{\td_B}$ fails.
\ele
\bpf
Suppose, towards the contrary, that ${\td_A}\rea{\td_B}$ 
holds, and let this be witnessed by a reduction 
$\Psi$ defined (as in \nrf{b-c}) from an increasing 
sequence $\tmin B=n_0<n_1<n_2<\dots$ of numbers 
$n_k\in B$ and a collection of maps 
$H_k:X_k\to \prod_{m\in\il k{k+1}\cap B}X_m\yt k\in A.$ 
We put 
\dm
f_k(\da)=
\tmax_{\xi\yi \eta\in X_k,\:d_k(\xi,\eta)<\da}
\;\;\;
\tmax_{m\in\il k{k+1}\cap B}
\;\;\;
d_m(H_k(\xi)(m),H_k(\eta)(m))\,,
\dm 
for  $k\in\dN$ and $\da>0$  
(with the understanding that $\tmax\pu=0$ if 
applicable).
Then $f(\da)=\tsup_{k\in A}f_k(\da)$ 
is a nondecreasing map $\dR^+\to\ir0\iy$.

We claim that
$\tlim_{\da\to0}f(\da)=0.$ 
Indeed otherwise there is $\ve>0$ such that $f(\da)\ge\ve$ 
for all $\da.$ 
Then the numbers
\dm
s_k=
\TS\tmin_{\xi\yi \eta\in X_k\,,\;\xi\ne\eta}\,
d_k(\xi,\eta)
\quad\hbox{(all of them are $>0$)}
\dm
must satisfy $\tinf_{k\in A}s_k=0.$ 
This allows us to define a sequence 
$k_0<k_1<k_2<\dots$ of numbers $k_i\in A,$ and, 
for any $k_i,$ a pair of elements $\xi_i,\eta_i\in X_{k_i}$ 
with $d_{k_i}(\xi_i,\eta_i)\to0,$ and also a number 
$m_i\in\il{k_i}{k_i+1}\cap B$ such that 
$d_{m_i}(H_{k_i}(\xi_i)(m_i),H_{k_i}(\eta_i)(m_i))
\ge\ve.$ 
Let $x\yi y\in\prod_{k\in A}X_k$ satisfy 
$x(k_i)=\xi_i$ and $y(k_i)=\eta_i$ for all $i$ and 
$x(k)=y(k)$ for all $k\in A$ not of the form $k_i.$  
Then easily $x\td_A y$ holds but 
$\Psi(x)\td_B\Psi(y)$ fails, which is a contradiction.  
Thus in fact $\tlim_{\da\to0}f(\da)=0$.

Let $k\in A,$ and let $Y_k\sq X_k$ be as in \ref{lv4}. 
Then there exist elements $x_k\ne y_k$ 
in $Y_k$ such that $H_k(x_k)\res k=H_k(y_k)\res k.$ 
By \ref{lv+} there is a chain 
$x_k=\xi_0,\xi_1,\dots,\xi_n=y_k$ of elements 
$\xi_i\in X_k$ with $d_k(\xi_i,\xi_{i+1})\le\frac1{k+1}$ 
for all $i<n.$ 
Now  
$H_k(\xi_i)\in\prod_{m\in\il k{k+1}\cap B}X_m$ 
for each $i\le n.$

Suppose that $m\in\il k{k+1}\cap B,$
and hence $m\ge n_k\ge k.$
The elements $y^m_i=H_k(\xi_i)(m)\yt i\le n,$
satisfy $d_m(y^m_i,y^m_{i+1})\le f_k(\frac1{k+1}).$
Note that $m\ne k$ because $k\in A$ while $m\in B.$ 
Thus we have $m>k$ strictly.
It follows that $n\le\mu_m,$ therefore, by \ref{lv5}, 
\ben
\tenu{(\arabic{enumi})}
\addtocounter{enumi}3
\itla{lv=}\msur
$d_m(H_k(x_k)(m),H_k(y_k)(m)) 
\le f_k(\frac1{k+1})+\frac1{m+1}\le
f(\frac1{k+1})+\frac1{k+1}$ \ 
\een 
for all \ $m\in\il k{k+1}\cap B$.

Both $x=\sis{x_k}{k\in A}$ and 
$y=\sis{y_k}{k\in A}$ are elements of 
$\prod_{k\in A}X_k,$ and $x\td_A y$ fails because 
$d_k(x_k,y_k)\ge r$ for all $k.$ 
On the other hand, we have $\Psi(x)\td_B\Psi(y)$ by 
\ref{lv=}, because $\tlim_{\da\to0}f(\da)=0.$ 
This is a contradiction to the assumption that $\Psi$ 
reduces ${\td_A}$ to ${\td_B}$.
\epF{Lemma~\ref{ab21}}

\epF{Theorem~\ref{coAB}}

\punk{Non-\dd\fsg compact case}
\las{nesk}

For any metric space $\dX=\stk X d,$ let
$\ld \dX$ denote the \eqr\ $\xd{X_k}{d_k}$ on $X^\dN,$ where
\index{zzDX@$\ld{\dX}$}
\index{equivalence relation, ER!zzDX@$\ld{\dX}$}%
$\stk{X_k}{d_k}=\stk\dX d$ for all $k.$ 
Thus $\fco$ is equal to $\ld\dR.$
One may ask what is the place of \eqr s of the form $\ld \dX,$
where $\dX$ is a Polish space, in the global \dd\reb structure of
Borel \eqr s?

The case of \dd\fsg compact Polish spaces here can be reduced 
to the case of finite spaces, \ie\ to \co equalities,
by Lemma~\ref{sc2co}.
Thus in this case we obtain a family of Borel \er s 
situated \dd\reb between the relations $\Et$ and $\fco$
by Theorems \ref{coT} and \ref{comax}, and this family 
has a rather rich \dd\reb structure by Theorem~\ref{coAB}.

The case of non-\dd\fsg compact spaces is much less studied.

\bex
\lam{C}
Let $\dX=\bn$ be the Baire space,
with the standard distance $d(a,b)=\frac1{m(a,b)+1},$
where $m(a,b)$ (for $a\ne b\in\bn$) is the largest integer $m$
such that $a\res m=b\res m.$~\snos
{Note that the relation $\ld \dX$ depends on the metric
rather than topological structure of a space $\dX,$ and hence
it is, generally speaking, essential to specify a concrete
distance compatible with the given topology.}
If $x\in\dnn$ and $n\zi k\in\dN$ then $x(n)\res k$ is
a finite sequence of $k$ integers.
It follows from the fact that $\bn$ is 0-dimensional that
$x \mathbin{\ld \bn} y$ is equivalent to
\dm
\kaz n\; \sus k_0\; \kaz k\ge k_0\;(x(n)\res k=y(n)\res k)\,.
\dm
for any $x\yi y\in\bn.$   
{\ubf Exercise:}
use this to show that ${\ld \bn}\eqb\Et$.
\eex

\bqe
\lam{c?}
Now let $\dX$ be the Polish space $C[0,1]$ of all continuous
maps $f:{[0,1]}\to\dR,$ with the distance 
$d(f,g)=\tmax_{0\le x\le1}|f(x)-g(x)|.$
(This space is not \dd\fsg compact, of course.
What is the position of $\ld{C[0,1]}$ in the global
\dd\reb structure of Borel \eqr s and what are its
\dd\reb connections with such better known \eqr s as
$\rE_i\yd i=0,1,2,3,$ and ${\bel p}\yd{\fco}$?
\eqe

This question
(see, \eg, Su Gao~\cite{su9as})
remains open.
The question is also connected with \co equalities,
in particular, with $\fco$ itself from another side.
Let us consider the following continual version 
$\fCo$ of the \eqr\ ${\fco}.$
If $f\yi g$ are continuous maps from $[0,\piy)$ to $\dR$
then we
define $x\fCo y$ iff 
$\tlim_{x\to\piy}|f(x)-g(x)|=0.$
\index{equivalence relation, ER!C0@$\fCo$}%
\index{zzC0@$\fCo$}%

It is clear that any continuous map $f:[0,\piy)\to\dR$
can be identified with the sequence of its restrictions 
to intervals of the form $\il n{n+1}\zt n\in\dN,$ that is,
with a certain point of the Polish product space
$C[0,1]^\dN.$
With such an identification, the domain of $\fCo$ is
naturally identified with a certain Borel set in $C[0,1]^\dN,$
while $\fCo$ itself is identified with a Borel \eqr{},
equal to $\ld{C[0,1]}$ on that set.
(The domain of $\ld{C[0,1]}$ is the whole space 
$C[0,1]^\dN.$)
Question~\ref{c?} also can be addressed to $\fCo$.

Su Gao proved in \cite{su9as} that  $\fCo$
(there defined as $E_K$)
satisfies ${\fCo}\reb {\uoa},$ where
\index{equivalence relation, ER!u0a@$\uoa$}%
\index{zzu0a@$\uoa$}%
$\uoa$ is an \eqr\ on $\dR^{\dN\ti\dN}$ defined as follows:
\dm
x\uoa y\qquad\text{iff}\qquad
\kaz\ve>0\;\sus m_0\;\kaz m\ge m_0\;\kaz n\;
(|x(m,n)-y(m,n|<\ve)\,.
\dm
In addition, a more complicated Borel \er\ $\uo$ on
$\dR^{\dN\ti\dN}\ti\dN^\dN$ is defined in
\cite{su9as} such that ${\fCo}\eqb {\uoa}.$
Investigations of ${\uo}\yd{\uoa}\yd{\fCo}\yd{\ld{C[0,1]}}$
remain work in progress.

\api

\parf{Pinned \eqr s}
\las{rtd}

In this \gla\ we consider a class of \eqr s $\rE$ characterized
by the property that if $\rE$ has an equivalence class in a
generic extension $\dvp$ of the ground set universe $\dV,$
definable in $\dvp$ in certain way in terms of sets in $\dV$
as parameters then this equivalence class contains an element
in $\dV.$
We call them \rit{pinned} \er s.

The main goal will be to prove that certain families of Borel
\er s are pinned, while on the other hand the \eqr\ $\rtd$ 
of equality of countable sets of the reals is not pinned, and
hence not Borel reducible to any pinned \eqr.
The class of pinned \er s includes, for 
instance, continuous actions of \cli\ groups and some 
ideals, not necessarily Polishable, and is closed under the 
Fubini product modulo $\ifi.$   

Recall that $\rtd$ is defined on $\nnn$ as follows: 
$x\rtd y$ iff $\ran x=\ran y$.

\bdf
\lam{pinnm}
$\dV$ will denote the ground set universe.
In this \gla\/ we'll consider forcing extensions of
\index{zzV@$\dV$}%
$\dV.$\snos
{Basically, a more rigorous treatment would be either to
consider boolean-valued extensions of the universe, or to
to assume that in fact $\dV$ is a countable model in a wider
universe.}

Suppose that $X$ is $\fs11$ or $\fp11$ in the universe 
$\dV,$ and an extension $\dvp$ of $\dV$ is considered. 
In this case, let $\di X$ denote what results by the
\index{zzXsharp@$\di X$}%
definition of $X$ applied in $\dvp.$ 
There is no ambiguity here by the Shoenfield absoluteness
theorem,
and easily $X=\di X\cap\dV$.
\edf

\punk{The definition of pinned \eqr s}
\las{defpin}

For instance, if, in the universe $\dV,$ $\rE$ is a  
$\fs11$ \eqr\ on a fixed 
polish space $\dX,$ then, still by the Shoenfield 
absoluteness, $\drE$ is a $\fs11$ \er\ on $\di\dX.$  
If now $x\in \dX$ (hence, $x\in\dV$) then the \dde class 
$\ek x\rE\sq\dX$ of $x$(defined in $\dV$) is included in 
a unique \dd\drE class $\ek x{\drE}\sq\di\dX$ (in $\dvp$).  
Classes of the form $\ek x{\drE}\zT x\in\dX,$ belong to a 
wider category of \dd\drE classes which admit a description 
from the point of view of the ground universe $\dV$. 

\vyk{
Fix a Polish space $\dX$ and let $\sis{B_n}{n\in\dN}$ 
be a base of its topology.  
By a {\it Borel code\/} for $\dX$ we shall understand 
\index{Borel code}%
a pair $p=\ang{T,f}$ of a wellfounded tree 
$\pu\ne T=T_p\sq\Ord\lom$ (then $\La\in T$) 
and a map $f:\Max T\to\dN,$ where $\Max T$ 
is the set of all \dd\sq maximal elements of $T.$  
\index{zzMaxT@$\Max T$}%
We define a set $\bk p(t)\sq\dnn$ for any $t\in T$ 
by induction on the rank of $t$ in $T,$ so that 
\bit
\item\msur
$\bk p(t)=\bk{f(t)}$ for all $t\in\Max T,$ and 

\item\msur
$\bk p(t)=\dop\bigcup_{t\we\xi\in T}\bk p(t\we\xi)$ for 
\index{zzcf@$\bk\tau$}%
$t\in T\dif\Max T$
($\dop$ denotes the complement);

\item
finally, put $\bk p=\bk p(\La)$. 
\eit

For a Borel code $p=\ang{T,F},$ let $\tsup p=\tsup T$ 
be the least ordinal $\ga$ with $T\sq\ga\lom.$ 
A code $p$ is {\it countable\/} 
\index{Borel code!countable}%
if $\tsup p<\omi,$ 
in this case the coded set $\bk p$ is a Borel subset 
of $\dX$.

\bdf
[based on an argument in Hjorth~\cite{h:orb}]
\lam{t2like}
A $\fs11$ \eqr\ $\rE$ is {\it\PP\/} if,
\index{equivalence relation, ER!pinned@\PP}%
\index{pointp@\PP}%
for any (perhaps, uncountable) Borel code $p,$ 
{\ubf if} $\bk p$ is  
pairwise \dd{\drE}equivalent in any generic extension   
of $\dV$ 
and non-empty in some generic extension of $\dV,$  
{\ubf then} there is a point $x\in\dom\rE,$ ``pinning'' 
$p$ in the sense that $\bk p\sq \ek{x}{\drE}$ in any  
extension of $\dV$. 
\edf
}

\bdf
[based on an argument in Hjorth~\cite{h:orb}]
\lam{virt}
Assume that $\rE$ is a $\fs11$ \eqr\ on a Polish space
$\dX$ and $\dP$ is a notion of forcing in $\dV.$ 
A {\it virtual\/ \dde class\/} is any \ddp term $\xib$
\index{virtual class}
such that $\dP$ forces $\xib\in\di\dX$ and
\index{zzxileft@$\tal$}%
\index{zzxiri@$\tar$}%
$\dPP$ forces $\tal\drE \tar$.~\footnote
{\label{xilr}\ $\tal$ and $\tar$ are \dd\dPP terms  
meaning $\xib$ associated with the resp.\ left and right 
factors $\dP$ in the product forcing. 
Formally, 
$\tal[U\ti V]=\xib[U]$ and $\tar[U\ti V]=\xib[V]$ for any 
\dd\dPP generic set $U\ti V,$ where 
$\xib[U]$ is the 
interpretation of a term $\xib$ via a generic set $U$.}  

A virtual class is {\it\PP\/} if there is, in $\dV,$ a 
point $x\in\dX$ which pins it, in the sense that 
$\dP$ forces $x\drE \xib.$  
Finally, $\rE$ is {\it\PP\/} if, 
for any forcing notion $\dP\in\dV,$ 
all virtual \dde classes are \PP.
\edf

If $\xib$ is a virtual \dde class then, in any extension 
$\dvp$ of $\dV,$ if $U$ and 
$V$ are generic subsets of $\dP$ then $x=\xib[U]$ 
and $y=\xib[V]$ belong to $\di\dX$ and satisfy $x\drE y,$ 
hence $\xib$ induces a \dd\drE class in the 
extension. 
If $\xib$ is \PP\ then this class contains an element in the 
ground universe $\dV$ --- in other words, \PP\ virtual 
classes induce \dd\drE equivalence classes of the form 
$\ek x\drE\zT x\in\dV,$ in the extensions of the universe 
$\dV.$


The following theorem
(originally \cite{calt1,kr:izv})
is the main result in this \gla.
Part \ref{horb2} here is from \cite{h:orb}.
Part \ref{horb3} also belongs to Hjorth and is published
with his permission.

\bte
\lam{h+orb}
The class of all pinned\/ $\fs11$ \eqr s$:$
\ben
\renu
\itla{horb1}
is closed 
under Fubini products modulo\/ $\ifi\:;$
\een
and contains the following equivalences$:$
\ben
\addtocounter{enumi}1
\renu
\itla{horb2}
all orbit \er s of Polish actions of (Polish)
\cli\ groups on a Polish space$;$\snos
{Recall that a Polish group $\dG$ is 
{\it complete left-invariant\/}, \cli\ for brevity, 
if $\dG$ admits a compatible  
\index{group!cli@\cli}%
\index{group!complete left-invariant}%
\index{cli@\cli}%
left-invariant complete metric.}

\itla{horb3}
all Borel \er s, all of whose equivalence classes are\/ 
$\Gds\:;$

\itla{horb4}
all \er s of the form\/
$\Exh_{\sis{\vpi_i}{}}=\ens{X\sq\dN}{\vpy(X)=0},$ where\/ 
$\vpi_i$ are lower semicontinuous (\lsc) submeasures on\/
$\dN$.
\een
On the other hand, $\rtd$ is not pinned and hence\/ $\rtd$
in Borel irreducible to any pinned \eqr.
\ete

Quite recently, Thompson~\cite{tho} proved that for a Polish
group $\dG$ to be \cli\ it is not only necessary
(which is by \ref{horb2})
but also sufficient that all orbit \eqr s of Polish actions
of $\dG$ are pinned.

\punk{$\rtd$ is not \PP}
\las{t2nonp}

Here we prove the last claim of Theorem~\ref{h+orb}.

\bcl
\lam{t2np}
$\rtd$ is \poq{not} \PP.
\ecl
\bpf
To prove that $\rtd$ is not \PP, consider, in $\dV,$ 
the forcing notion $\dP=\coll(\dN,\dn)$ to produce a 
generic map $f:\dN\onto\dn.$ 
($\dP$ consists of all functions $p:u\to\dn$ where  
$u\sq\dN$ is finite.)
The \ddp term $\xib$ for the set 
$\ran f=\ans{f(n):n\in\dN}$ is obviously a virtual 
\dd\rtd class, but it is not \PP\ because $\dn$ is 
uncountable in the ground universe $\dV$.  
\epf

\ble
\lam{p<p}
\vyk{
{\rm(Hjorth~\cite{h:orb})} \ 
\PPP\/ $\fs11$ \er s\/ $\rE$ do not satisfy\/ 
$\rtd\reb\rE$.
}%
If\/ $\rE,\,\rF$ are\/ $\fs11$ \er s, $\rE\reb\rF,$  and\/ 
$\rF$ is \PP, then so is\/ $\rE$. 
\ele
\bpf
Suppose that, in $\dV,$ $\vt:\dX\to\dY$ is a Borel
reduction of $\rE$ to $\rF,$ where $\dX=\dom\rE$ and 
$\dY=\dom\rF.$    
We can assume that $\dX$ and $\dY$ are just two copies 
of $\dn.$ 
Let $\dP$ be a forcing notion and a \ddp term $\xib$ be a 
virtual \dde class.  
By the Shoenfield absoluteness, $\di\vt$ is a 
reduction of $\di\rE$ to $\di\rF$ in any extension of $\dV,$ 
hence, $\sgb,$ a \ddp term for $\di\vt(\xib),$ is also a 
virtual \ddf class. 
Since $\rF$ is \PP, there is $y\in\dY$ such that 
$\dP$ forces $y\drF\sgb.$ 
Note that it is true in the \ddp extension that 
$y\drF \di\vt(x)$ for some $x\in\di\dX,$ hence, by the 
Shoenfield theorem, 
in the ground universe there is $x\in\dX$ with $y\rF\vt(x).$ 
Clearly $\dP$ forces $x \drE \xib$.
\epf

\punk{Fubini product of \PP\ \er s is \PP}
\las{fbbsbs}

Here we prove part \ref{horb1} of Theorem~\ref{h+orb}.
Recall that the Fubini product $\rE=\fps{\ifi}{\rE_k}{k\in\dN}$ 
of \er s $\rE_k$ on $\dnn$ modulo $\ifi$ is an \eqr\ on $\nnn$ 
defined as follows: $x\rE y$ if $x(k)\rE_k y(k)$ for all 
but~finite~$k$.

Suppose that $\fs11$ \eqr s $\rE_k$ on Polish spaces 
$\dX_k$ are \PP. 
Prove that the Fubini product 
$\rE=\fps{\ifi}{\rE_k}{k\in\dN}$ is a \PP\ \er\
(on the Polish space $\dX=\prod_k\dX_k$). 
Consider a forcing notion $\dP$ and a \ddp term $\xib.$ 
Assume that $\xib$ is a virtual \dde class.  
There is a number $k_0$ and conditions $p,\,q\in \dP$ 
such that $\ang{p,q}$ \dd\dPP forces 
$\tal(k)\mathbin{\di{\rE_k}}\tar(k)$ for all $k\ge k_0.$ 
As all $\rE_k$ are \er s, we conclude that the condition  
$\ang{p,p}$ also forces 
$\tal(k)\mathbin{\di{\rE_k}}\tar(k)$ for all $k\ge k_0.$
Therefore, since $\rE_k$ are \PP, there is in $\dV$ 
a sequence of points $x_k\in\dX_k$ such that $p$ \ddp forces 
$x_k \mathbin{\di{\rE_k}}\xib(k)$ for any $k\ge k_0.$ 
Let $x\in\dX$ satisfy $x(k)=x_k$ for all $k\ge k_0.$ 
(The values $x(k)\in\dX_k$ for $k<k_0$ can be arbitrary.)  
Then $p$ obviously \ddp forces $x\drE \xib.$ 

It remains to show that just every $q\in\dP$ also forces 
$x\drE \xib.$ 
Suppose otherwise, that is, some $q\in\dP$ forces that 
$x\drE \xib$ \poq{fails}. 
Consider the pair $\ang{p,q}$ as a condition in $\dPP:$ 
it forces $x\drE\tal$ and $\neg\;{x\drE\tar},$ as 
well as $\tal\drE\tar$ by the choice of $\rE$ and 
$\xib,$ which is a contradiction.

\vyk{
Define $x\rF_k y$ iff $x(k)\rE_k y(k):$ $\rF_k$ are 
$\fs11$ \er s on $\nnn$ and $x\rE y$ iff $x\rF_k y$ for 
almost all~$k$.

\bcl
\label{bsfub'}
Each\/ $\rF_k$ is \PP.
\ecl
\bpf
Consider a Borel code $p$ for a subset of $\nnn,$ 
satisfying \poq{if} of Definition~\ref{t2like} \vrt\ $\rF_k.$  
By the same argument as in the proof of Lemma~\ref{t2unlike}, 
there is a Borel code $q$ for a subset of $\dnn,$ such that 
$\bk q\ne\pu$ in some extension of $\dV$ and 
$\bk q\sq\ens{x(k)}{x\in\bk p}$ in any extension of $\dV,$ 
hence, $q$ satisfies \poq{if} of Definition~\ref{t2like} 
\vrt\ $\rE_k.$ 
As $\rE_k$ is \PP, there is $a\in\dnn$ such that 
$\bk q\sq\ek{a}{\di\rE_k}$ in any extension, but then easily 
$\bk p\sq\ek{x}{\di\rF_k}$ in any extension, where 
$x\in\nnn\cap\dV$ has only to satisfy $x(k)=a$ for the given 
$k$.
\epF{Claim}

In continuation of the proof of \ref{horb1} of
Theorem~\ref{h+orb}, 
consider a Borel code $p$ for a subset of $\nnn,$ 
satisfying \poq{if} of Definition~\ref{t2like} \vrt\ $\rE.$ 
Our plan is to find another Borel code $\bap$ with 
$\bk\bap\sq\bk p$ everywhere, which satisfies 
\poq{if} of Definition~\ref{t2like} for almost all $\rE_k.$ 
This involves a forcing by Borel codes.

Let, in $\dV,$ $\la=\tsup p$ and $\ka=\la^+,$ 
thus, $\tsup p<\ka.$  
Let $\dP$ be the set of all Borel codes $q\in\dV$ for  
subsets of $\nnn$ such that $\tsup q<\ka$ 
and $\bk q\ne\pu$ in a generic extension of the 
universe $\dV.$ 
$\dP$ is considered as a forcing, with   
$q\lef p$ ($q$ is stronger) iff $\bk q\sq \bk p$ 
in all generic extensions of $\dV.$  
It is known that $\dP$ forces a point of $\nnn,$ so that  
$\bigcap_{q\in G}\bk q=\ans{x_G}$ for any \dd\dP generic, 
over $\dV,$ set $G\sq\dP.$ 
Let $\dox$ be the name of the generic element of $\nnn$.

By the choice of $p,$ $\ang{p,p}$ \dd{\dP\ti\dP}forces 
$\doxl\drE\doxr,$ hence, there are codes $q,\,r\in\dP$ 
and a number $k_0$ such that $\ang{q,r}$ 
\dd{\dP\ti\dP}forces 
$\doxl\mathbin{\di{\rF_k}}\doxr$ for any $k\ge k_0.$ 
By a standard argument, we have 
$x\mathbin{\di{\rF_k}} y$ for all $k\ge k_0$ 
in any extension of $\dV$ for any two \dd\dP generic, over 
$\dV,$ elements $x,\,y\in\bk q.$ 
We can straightforwardly define in $\dV$ a Borel code 
$\bap$ (perhaps, not a member of $\dP$!) 
such that, in 
any extension of $\dV,$ $\bk\bap$ is the set of all 
\dd\dP generic, over $\dV,$ elements of $\bk q.$ 
Then $\bap$ satisfies \poq{if} of Definition~\ref{t2like} 
\vrt\ any $\rF_k$ with $k\ge k_0.$ 
Hence, by the claim, there is, in $\dV,$ a sequence of points 
$x_k\in\nnn$ such that $\bk\bap \sq\ek{x_k}{\di\rF_k}$ 
in any generic extension of $\dV,$ for any $k\ge k_0.$ 
Define $x\in\nnn\cap\dV$ so that $x(k)=x_k(k)$ for any 
$k\ge k_0,$ then, by the definition of $\rF_k,$ we have  
$\bk\bap\sq\ek{x}{\di\rF_k}$ for all $k\ge k_0$ 
in any extension of $\dV.$ 
Yet 
$\bigcap_{k\ge k_0}\ek{x}{\di\rF_k}\sq\ek{x}{\drE}$.
%
}

\punk{Complete left-invariant actions induce 
\PP\ \er s}
\las{bs:cli}

Here we prove part \ref{horb2} of Theorem~\ref{h+orb}.
Suppose that $\dG$ is a Polish \cli\ group
continuously  acting on a Polish space $\dX.$
By definition $\dG$ admits a compatible  
left-invariant complete metric. 
Then easily $\dG$ also admits a compatible 
{\ubf right}-invariant 
complete metric, which will be practically used.

Let $\dP$ be a forcing notion and $\xib$ be a virtual 
\dde class. 
Let $\le$ denote the partial order of $\dP\,;$ we assume, 
as usual,  
that $p\le q$ means that $p$ is a stronger condition.
Let us fix a compatible complete right-invariant 
metric $\rho$ on $\dG.$ 
For any $\ve>0,$ put $G_\ve=\ans{g\in\dG:\rho(g,1_\dG)<\ve}.$ 
Say  that $q\in\dP$ {\it is of size\/ $\le\ve$\/} 
if $\ang{q,q}$ \dd{\dPP}forces the existence of 
$g\in \di{G_\ve}$ such that $\tal=g\app \tar$. 

\ble
\label{l4}
If\/ $q\in\dP$ and\/ $\ve>0,$ then there is a condition\/ 
$r\in\dP\zT r\le q,$ of size\/ $\le \ve$.
\ele
\bpf
Otherwise for any $r\in\dP\zT r\le q,$ there 
is a pair of conditions $r',\,r''\in\dP$ stronger than 
$r$ and such that $\ang{r',r''}$ \dd{\dPP}forces that 
there is no $g\in \di{G_\ve}$ with $\tal=g\app\tar.$ 
Applying an ordinary splitting construction in such 
a generic extension $\dvp$ of $\dV$ where 
$\cP(\dP)\cap\dV$ is countable, we find an uncountable 
set $\cU$ of generic sets $U\sq\dP$ with $q\in U$ such 
that any pair $\ang{U,V}$ with $U\ne V$ in $\cU$ is 
\dd\dPP generic (over $\dV$), hence, there is no 
$g\in \di{G_\ve}$ with 
$\xib[U]=g\app \xib[V].$~\footnote
{\label{tauu}\ $\xib[U]$ is the interpretation of the 
\ddp term $\xib$ obtained by taking $U$ as the generic 
set.} 
Fix $U_0\in \cU.$
We can associate in $\dvp$ with each $U\in \cU,$ an 
element $g_U\in\di G$ such that 
$\xib[U]=g_U\app \xib[U_0];$ then $g_U\nin\di{G_\ve}$ 
by the above. 
Moreover, we have $g_Vg_U\obr\app \xib[U]=\xib[V]$ 
for all $U,\,V\in \cU,$ 
hence $g_Vg_U\obr\nin\di{G_\ve}$ whenever $U\ne V,$ which 
implies $\rho(g_U,g_V)\ge\ve$ by the right invariance. 
But this contradicts the separability of $G$.
\epF{Lemma}

Coming back to the proof of \ref{horb3} of
Theorem~\ref{h+orb}, suppose towards the contrary that 
a condition $p\in\dP$ forces that there is no $x\in\dX$ 
(in the ground universe $\dV$) 
satisfying $x\drE\xib.$ 
According to Lemma~\ref{l4}, there is, in $\dV,$ a sequence 
of conditions $p_n\in\dP$ of size $\le 2^{-n},$ 
and closed sets $X_n\sq\dX$ with \dd\dX diameter 
$\le 2^{-n},$ such that 
$p_0\le p\zT\,p_{n+1}\le p_n\zT\,X_{n+1}\sq X_n,$ 
and $p_n$ forces $\xib\in\di{X_n}$ for any $n.$ 
Let $x$ be the common point of the sets $X_n$ in $\dV.$ 
We claim that $p_0$ forces $x\drE \xib$. 

Indeed, otherwise there is $q\in\dP\zT q\le p_0,$ which 
forces $\neg\;{x\drE \xib}.$
Consider an extension $\dvp$ of $\dV$ rich enough to 
contain, for any $n,$ a generic set $U_n\sq\dP$ with 
$p_n\in U_n$ such that each pair $\ang{U_n,U_{n+1}}$ 
is \dd\dPP generic 
(over $\dV$), and, in addition, $q\in U_0.$    
Let $x_n=\xib[U_n]$ (an element of $\di\dX$),  
then $\sis{x_n}{}\to x.$ 
Moreover, for any $n,$ both $U_n$ and $U_{n+1}$ 
contain $p_n,$ hence, as $p_n$ has size 
$\le 2^{-n-1},$ there is $g_{n+1}\in \di{\dG_\ve}$ with 
$x_{n+1}=g_{n+1}\app x_n.$ 
Thus, $x_n=h_n\app x_0,$ where $h_n=g_{n}\dots g_1.$ 
However  
$\rho(h_n,h_{n-1})=\rho(g_n,1_\dG)\le 2^{-n+1}$ by 
the right-invariance of the metric, thus, 
$\sis{h_n}{n\in\dN}$ is a Cauchy sequence in $\di\dG.$ 
Let $h=\tlim_{n\to\iy}h_n\in\di\dG$ be its limit. 
As the action considered is continuous, we have 
$x=\tlim_nx_n=h\app x_0.$ 
It follows that $x\drE x_0$ holds in $\dvp,$ hence 
also in $\dV[U_0].$ 
However $x_0=\xib[U_0]$ while $q\in U_0$ forces 
$\neg\;{x\drE \xib},$ which is a contradiction. 

Thus $p_0$ \ddp forces $x\drE \xib.$ 
Then any $r\in\dP$ also forces $x\drE \xib:$ 
indeed, if some $r\in\dP$ forces $\neg\;{x\drE \xib}$ 
then the pair $\ang{p_0,r}$ \dd\dPP forces  that 
$x\drE \tal$ and $\neg\;{x\drE \tar},$ which 
contradicts the fact that $\dPP$ forces 
$\tal\drE\tar$.

\punk{All \er s with $\Gds$ classes are \PP}
\label{s02}

Here we prove part \ref{horb3} of Theorem~\ref{h+orb}.
Suppose that $\rE$ is a Borel \eqr\ on $\bn$ and all
\dde equivalence classes are $\Gds.$
Prove that then $\rE$ is \PP.

It follows from a theorem of Louveau~\cite{louv} that 
there is a Borel map $\ga,$ defined on $\dnn,$ so that 
$\ga(x)$ is a \dd\Gds code of $\ek x\rE$ for any 
$x\in\dnn,$ that is, for instance, 
$\ga(x)\sq\dN^2\ti\dN\lom$ and  
\dm
\ek x\rE=
\bigcup_i\bigcap_j\bigcup_{\ang{i,j,s}\in\ga(x)}B_s,
\quad\text{where}\quad 
B_s=\ans{a\in\dnn:s\su a}
\;\text{ for all }\;s\in\dN\lom.  
\dm
We consider a forcing notion $\dP=\stk\dP\le$ and a  
virtual \dde class $\xib.$ 
Then $\dPP$ forces $\tal\drE\tar;$ hence there is a 
number $i_0$ and a condition $\ang{p_0,q_0}\in\dPP$ 
which forces $\tal\in\di\vt(\tar),$ where 
$\vt(x)=\bigcap_j\bigcup_{\ang{i_0,j,s}\in\ga(x)}B_s$ for 
all $x\in\dnn.$ 

The key idea of the proof is to substitute $\dP$ by the 
Cohen forcing. 
Let $\dS$ denote the set of all $s\in\dN\lom$ such that 
$p_0$ does \poq{not} \ddp force that $s\not\su \xib.$ 
We consider $\dS$ as a forcing, and $s\sq t$ 
(that is, $t$ is an extension of $s$) means that $t$ is a 
stronger condition; $\La,$ the empty sequence, 
is the weakest condition in $\dS.$  
If $s\in \dS$ then obviously there is at least one $n$ 
such that $s\we n\in \dS;$ hence $\dS$ forces an element 
of $\dnn,$ whose \dd\dS name will be $\fa$. 

\ble
\label{s0}
The pair\/  $\ang{\La,q_0}$ \dd{\dSP}forces\/ 
$\fa\in\di\vt(\xib).$ 
\ele
\bpf
Otherwise some condition $\ang{s_0,q}\in\dSP$ with  
$q\le q_0$ forces $\fa\nin\di\vt(\xib).$ 
By the definition of $\vt$ we can assume that 
\dm
\ang{s_0,q}
\quad \text{\dd{\dSP}forces}\quad
\neg\;\sus s\:\skl
\ang{i_0,j_0,s}\in\ga(\xib)\land s\su\fa
\skp 
\eqno(\ast)
\dm
for some $j_0.$ 
Since $s_0\in\dS,$ there is a condition 
$p'\in\dP\zq p'\le p_0,$ which \ddp forces $s_0\su\xib.$ 
By the choice of $\ang{p_0,q_0}$ we can assume that  
\dm
\ang{p',q'}
\quad \text{\dd{\dPP}forces}\quad
\ang{i_0,j_0,s}\in\ga(\tar)\land s\su\tal\,.
\dm
for suitable $s\in\dS$ and $q'\in\dP\zq q'\le q.$ 
This means that 
$1)\msur$ $p'$ \ddp forces $s\su\xib$ and 
$2)\msur$ $q'$ \ddp forces $\ang{i_0,j_0,s}\in\ga(\xib).$ 
In particular, by the above, $p'$ forces both $s_0\su\xib$ 
and $s\su\xib,$ therefore, 
either $s\sq s_0$ -- then let $s'=s_0,$ 
or $s_0\su s$ -- then let $s'=s.$ 
In both cases, $\ang{s',q'}$ \dd\dSP forces 
$\ang{i_0,j_0,s}\in\ga(\xib)$ and $s\su\fa,$ 
contradiction to $(\ast)$.
\epF{Lemma} 

Note that $\dS$ is a subforcing of the Cohen forcing 
$\dC=\dN\lom,$ therefore, by Lemma~\ref{s0}, there is a 
\dd\dC term $\sgb$ such that $\ang{\La,q_0}$ \dd\dCP 
forces 
$\sgb\in\di\vt(\xib),$ hence, forces $\sgb\drE\xib.$ 
It follows, by consideration of the forcing $\dCP\ti\dP,$ 
that generally $\dCP$ forces $\sgb\drE\xib.$ 
Therefore, by ordinary arguments, first,  
$\dCC$ forces $\sgbl\drE\sgbr,$ and second, to prove the 
theorem it suffices now to find $x\in\dnn$ in $\dV$ such 
that $\dC$ forces $x\drE\sgb.$ 
This is our next goal.

Let $\fa$ be a \dd\dC name of the Cohen generic element 
of $\dnn.$ 
The term $\sgb$ can be of complicated nature, but we can 
substitute it by a term of the form $\di f(\fa),$ where 
$f:\dnn\to\dnn$ is a Borel map in the ground universe $\dV.$ 
It follows from the above that 
$\di f(\fa)\drE \di f(\fb)$ for any \dd\dCC generic, 
over $\dV,$ pair $\ang{\fa,\fb}\in\dnn\ti\dnn.$ 
We conclude that $\di f(\fa)\drE \di f(\fb)$ 
also holds even for any pair 
of separately Cohen generic $\fa,\,\fb\in\dnn.$ 
Thus, in a generic extension of $\dV,$ where there are  
comeager-many Cohen generic reals, there is a comeager 
$\Gd$ set $X\sq\dnn$ such that $\di f(a)\drE \di f(b)$ 
for all $a,\,b\in X.$ 
By the Shoenfield absoluteness theorem, the statement of 
existence 
of such a set $X$ is true also in $\dV,$ hence, in $\dV,$ 
there is $x\in\dnn$ such that we have $x\rE f(a)$  
for comeager-many $a\in\dnn.$ 
This is again a Shoenfield absolute property of $x,$ 
hence, $\dC$ forces $x\drE \di f(\fa),$ as required.

\punk{A family of \PP\ ideals}
\las{bsti}

Here we prove part \ref{horb4} of Theorem~\ref{h+orb}.

Let us say that a Borel ideal $\cI$ is {\it\PP\/} 
if the induced \er\ $\rei$ is such. 
It follows from Theorem~\ref{h+orb}\ref{horb2} that any 
P-ideal is \PP\ because Borel P-ideals are polishable
by Theorem~\ref{sol} while all Polish abelian groups are \cli.
\imar{reference}
Yet there are non-{\it P\/} \PP\ ideals. 

Suppose that $\sis{\vpi_i}{i\in\dN}$ is a sequence of 
lower semicontinuous (\lsc) submeasures on $\dN.$ 
Define the exhaustive ideal of the sequence,
\dm
\Exh_{\sis{\vpi_i}{}}\,=\,
\ens{X\sq\dN}{\vpy(X)=0}\,,
\quad\hbox{where}\quad
\vpy(X)=\tlis_{i\to\iy}\vpi_i(X)\,.
\dm
It follows from Theorem~\ref{sol} that for any Borel
P-ideal $\cI$ there is a \lsc\ submeasure $\vpi$ such that  
$\cI=\Exh_{\sis{\vpi_i}{}}=\Exh_\vpi,$ where 
$\vpi_i(x)=\vpi(x\cap\iry i).$
On the other hand, the 
non-polishable ideal $\Ii=\fio$ also is of the form 
$\Exh_{\sis{\vpi_i}{}},$ where for $x\sq\dN^2$ we define 
$\vpi_i(x)=0\,\hbox{ or }\,1$ if resp.\ 
$x\sq\,\hbox{ or }\,\not\sq\ans{0,\dots,n-1}\ti\dN$. 

Thus suppose that $\vpi_i$ is a \lsc\ submeasure on $\dN$
for each $i\in\dN.$  
The goal is to prove that the ideal
$\cI=\Exh_{\sis{\vpi_i}{}}$ is \PP.

We can assume that the submeasures $\vpi_i$ decrease, 
that is $\vpi_{i+1}(x)\le\vpi_i(x)$ for any $x,$ for if not 
then consider the \lsc\ submeasures 
$\vpi'_i(x)=\tsup_{j\ge i}\vpi_j(x).$ 

\vyk{
Let $\wtau$ be a Borel code, for a subset of $\pn,$ 
satisfying \poq{if} of Definition~\ref{t2like} \vrt\ 
the induced \er\ $\rei$ on $\pn,$ 
thus, $\wtau\in\dP,$ where $\dP$ is a forcing 
defined as in \nrf{fbbsbs} 
($\dP$ forces a subset of $\pn$). 

Using the same arguments as above, we see that for any 
$p\in\dP,\msur$ $p\lef\wtau,$ and $n\in\dN,$ 
there are $i\ge n$ and codes $q,\,r\in\dP$ with 
$q,\,r\lef p,$ 
such that $\ang{q,r}$ \dd{\dP\ti\dP}forces that 
$\vpi_i(\doxl\sd\doxr)\le 2^{-n-1},$ hence, any two 
\dd\dP generic, over $\dV,$ elements $x,\,y\in\bk q$ 
satisfy $\vpi_i(x\sd y)\le 2^{-n}.$ 
It follows that, in $\dV,$ there is a sequence of numbers 
$i_0<i_1<i_2<\dots,$ a sequence 
$\wtau\gef p_0\gef p_1\gef p_2\gef\dots$ of codes 
in $\dP,$ and, for any $n,$ a set $u_n\sq\ir0n,$ such that 
\ben
\tenu{(\arabic{enumi})}
\itla{bs1}
each $p_n$ \dd\dP forces $\dox\cap\ir0n=u_n$;

\itla{bs2}
any \dd\dP generic, over $\dV,$ $x,\,y\in\bk{p_n}$ 
satisfy $\vpi_{i_n}(x\sd y)\le 2^{-n}.$ 
\een
Let, in $\dV,$ $a=\bigcup_nu_n,$ then  
$a\cap\ir0n=u_n$ for all $n.$ 
Prove that $a$ pins $\bk\wtau,$ \ie, 
$\bk\wtau\sq\ek a{\di\rE_\cI}$ in any extension of $\dV$.

We can assume that, in the extension, for any $n$ there is 
a \dd\dP generic, over $\dV,$ element $x_n\in\bk{p_n}.$  
Then we have, by \ref{bs2}, $\vpi_{i_n}(x_n\sd x_m)\le 2^{-n}$ 
whenever $n\le m.$ 
It follows that $\vpi_{i_n}(x_n\sd a)\le 2^{-n},$ because 
$a=\tlim_mx_m$ by \ref{bs1}. 
However we assume that the submeasures $\vpi_j$ decrease, 
hence, $\vpy(x_n\sd a)\le 2^{-n}.$ 
On the other hand, $\vpy(x_n\sd x_0)=0$ because 
all elements of $\bk{p_0}$ are pairwise 
\dd{\di\rE_\cI}equivalent. 
We conclude that $\vpy(x_0\sd a)\le 2^{-n}$ 
for any $n,$ 
in other words, $\vpy(x_0\sd a)=0,\msur$  
$x_0\mathbin{\di\rE_\cI} a,$ and 
$\bk\wtau\sq \ek a{\di\rE_\cI},$ as required.

*****
}%

Suppose towards the contrary that the equivalence $\rE=\rei$
is not \PP.
Then there is a forcing notion $\dP,$ a virtual 
\dde class $\xib,$ and a condition $p\in \dP$ which 
\ddp forces $\neg\;x\drE\xib$ for any $x\in\pn$ in 
$\dV.$ 
By definition, for any $p'\in\dP$ and $n\in\dN$ there 
are $i\ge n$ and conditions $q,\,r\in\dP$ with 
$q,\,r\le p',$ such that $\ang{q,r}$ \dd{\dPP}forces  
the inequality $\vpi_i(\tal\sd\tar)\le 2^{-n-1},$ 
hence, $\ang{q,q}$ \dd{\dPP}forces  
$\vpi_i(\tal\sd\tar)\le 2^{-n}.$ 
It follows that, in $\dV,$ there is a sequence of numbers 
$i_0<i_1<i_2<\dots,$ and a sequence 
$p_0\ge p_1\ge p_2\ge\dots$ of conditions in $\dP,$ and, 
for any $n,$ a set $u_n\sq\ir0n,$ such that $p_0\le p$ and  
\ben
\tenu{(\arabic{enumi})}
\itla{bs1}
each $p_n$ \dd\dP forces $\xib\cap\ir0n=u_n$;

\itla{bs2}
each $\ang{p_n,p_n}$ \dd\dPP forces  
$\vpi_{i_n}(\tal\sd \tar)\le 2^{-n}.$ 
\een
Arguing in the universe $\dV,$ put $a=\bigcup_nu_n;$ then  
$a\cap\ir0n=u_n$ for all $n.$ 
We claim that $p_0$ forces $a\drE\xib.$ 
This contradicts the assumption above, ending the
proof of \ref{horb4} of Theorem~\ref{h+orb}. 

To prove the claim, note that otherwise there is a condition
$q_0\le p_0$ which forces $\neg\;a\drE\xib.$ 
Consider a generic extension $\dvp$ of the universe, where 
there exists a sequence of \ddp generic sets $U_n\sq\dP$ such 
that for any $n,$ the pair $\ang{U_n,U_{n+1}}$ is 
\dd\dPP generic, $p_n\in U_n,$ and in addition 
$q_0\in U_0.$ 
Then, in $\dvp,$ the sets $x_n=\xib[U_n]\in\pn$ satisfy
$\vpi_{i_n}(x_n\sd x_m)\le 2^{-n}$ by \ref{bs2}, 
whenever $n\le m.$ 
It follows that $\vpi_{i_n}(x_n\sd a)\le 2^{-n},$ because 
$a=\tlim_mx_m$ by \ref{bs1}. 
However we assume that the submeasures $\vpi_j$ decrease, 
therefore $\vpy(x_n\sd a)\le 2^{-n}.$ 
On the other hand, $\vpy(x_n\sd x_0)=0$ because $\xib$ 
is a virtual \dde class.
We conclude that $\vpy(x_0\sd a)\le 2^{-n}$ 
for any $n.$
In other words, $\vpy(x_0\sd a)=0,$ that is,  
$x_0\drE a,$ which is a contradiction with the choice of 
$U_0$ because $x_0=\xib[U_0]$ and $q_0\in U_0$.\vtm

\qeDD{Theorem~\ref{h+orb}}\vtm\vom

One might ask whether all Borel ideals are \PP.
This question answers in the negative.
Indeed it will be proved in the next \gla\ that for 
every Borel \eqr\ $\rE$ there exists a Borel ideal
$\cI$ such that $\rE\reb\rei.$
In particular this is true for the \er\ $\rtd,$ non-pinned
by Theorem~\ref{h+orb}.
It follows, still by Theorem~\ref{h+orb}, that any Borel 
ideal $\cI$ satisfying $\rtd\reb\rei$ is non-pinned as well.

\bqe
[Kechris]
\lam{ppq2}
Is it 
true that $\rtd$ is the \dd\reb least non-\PP\ Borel \eqr\,?
\eqe

\api

\parf
[Borel equivalences vs. Borel ideals]
{Reduction of Borel equivalences
to Borel ideals}
\las{rosen}

The main goal of this \gla\ is to show that any Borel \eqr\
is Borel reducible to a relation of the form $\rei$ for some
Borel ideal $\cI,$ and moreover, there is a \dd\reb cofinal
\dd\omi sequence of Borel ideals in the sense of the
next theorem:

\bte
\label{M}
There is a\/ \dd\sq decreasing sequence of Borel ideals\/
$\cI_\xi\,\;(\xi<\omi)$ on\/ $\dN,$ \dd\reb cofinal in the
sense that every Borel \eqr\ is Borel reducible
to one of the relations $\rE_{\cI_\xi}$.
\ete

The proof (due to Rosental~\cite{roz}) of this important result
involves a universal analytic equivalence generated by an
analytic ideal, followed by 
a well-known construction of upper Borel approximations
of $\fs11$ sets.
Note that this theorem, together with Corollary~\ref{t2i},
accomplishes the proof of Theorem~\ref{troz}.

In the end we briefly outline the results of subsequent
study~\cite{kl5}: the ideals $\cI$ and
the corresponding relations $\rE_{\cI_\xi}$ as above can be   
explicitly and meaningfully defined on the base of a certain
game.

\punk{Trees}
\label{n+}

We begin with a review of basic notation related to trees
of finite sequences. 
Recall that for any set $X,$ $X^n$ denotes the set of all
sequences, of length $n,$ of elements of $X,$ and 
\index{sequence}
$X\lom=\bigcup_{n\in\dN}X^n$ -- the set of all finite sequences 
of elements of $X.$
Regarding product sets, note that any 
$s\in {(X_1\ti\cdots\ti X_n)}\lom$ is formally a finite sequence 
of \dd ntuples $\ang{x_1,\dots,x_n},$ where 
$x_i\in X_i\yd\kaz i.$ 
We identify such a sequence $s$ with the \dd ntuple 
$\ang{s_1,\dots,s_n},$ 
where all $s_i\in {X_i}\lom$ have the same length as $s$ 
itself, and $s(i)=\ang{s_1(i),\dots,s_n(i)}$ for all $i$.

$\lh s$ is the {\it length\/} of a sequence $s.$
\index{sequence!length}
$\La,$ {\it the empty sequence\/}, is the only one of length $0.$
If $s$ is a finite sequence and $x$ any set then by $s\we x,$ resp.,
$x\we s$ we denote the result of adjoining $x$ as the new
rightmost, resp, leftmost term to $s.$
If $s,t$ are sequences then $s\sq t$ means that $t$ is an
{\it extension\/}
of $s,$ that is, $s=t\res m$ for some $m\le \lh t$. 

A {\it tree\/} on a set $X$ is any subset $T\sq X\lom$ closed
\index{tree}
under restrictions --- that is, if $t\in T,$ $s\in X\lom,$ and 
$s\sq t,$ then $s\in T.$
Note that $\La,$ the empty sequence, belongs to any 
tree $\pu\ne R\sq X\lom.$
An {\it infinite branch \/} in a tree $T\sq X\lom$ is any
infinite sequence $b\in X^\om$ such that $b\res m\in T\zd\kaz m.$ 
A tree $T$ is {\it well-founded\/} iff it has no infinite branches.
\index{tree!well-founded}%
Otherwise $T$ is {\it ill-founded\/}.
\index{tree!ill-founded}%

The following transformations of trees on $\dN$ 
preserve in this or another way
the properties of well- and ill-foundednes.\vom

{\ubf Finite union.}
If $S,T$ are trees then so is $W=S\cup T,$ and clearly 
$S\cup T$ is ill-founded iff so is at least one of $S,T.$
\vom

{\ubf Contraction.} 
Let $S\sq\tlw$ be a tree. 
Fix once and for all a bijection $b:\dN^2\onto\dN.$ 
For any sequence $s=\ang{k_0,k_1,\dots,k_n}\in\tlw$ with 
$\lh s=n+1\ge 2$ define a sequence 
$\cntr s=\ang{b(k_0,k_1),k_2,\dots,k_n}$ of length $n.$
The {\it contracted tree\/}
\dm
\cntr S=\ans\La\cup\ens{\hs}{s\in S\land \lh s\ge 2}
\dm
is ill-founded iff so is $S$ itself.\vom

{\ubf Countable sum.} 
Countable unions do not preserve well-foundedness. 
Yet there is another useful operation. 
For any sequence of trees $T_n\sq\nlw,$ we let 
$\sum_n^\ast T_n$ denote the tree
$T=\ans\La\cup\ens{n\we t}{t\in T_n}.$ 
Clearly $T$ is ill-founded iff so is {\ubf at least one}
of the trees $T_n$.\vom

{\ubf Countable product.} 
Let 
$\prod_n^\ast T_n$ denote the set $T$ of all finite sequences 
of the form 
$t=\ang{t_0,\dots,t_n},$ where $t_k\in T_k$ and $\lh{t_k}=n$ 
for all $k\le n.$ 
We put $\ang{t_0,\dots,t_n}\lef \ang{s_0,\dots,s_m}$ iff 
$n\le m$ and $t_k\sq s_k$ (in $\nlw$) for all $k\le n.$ 
In addition, let $\La$ belong to $T,$ with $\La\lef t$ for any 
$t\in T.$ 
Obviously $\stk T\lef$ is an at most countable tree, order 
isomorphic to a tree in $\nlw.$ 
Moreover $T=\prod_n^\ast T_n$ 
is ill-founded iff so is {\ubf every} tree $T_n$.\vom

{\ubf Componentwise addition.}
This is a less trivial operation.
First of all, if $s\yi t\in \tlw$ then $s\lcw t$
(the {\it componentwise\/} ordering)
means that $\lh s=\lh t$ and $s(i)\le t(i)$ for all 
$i<\lh s.$
Similarly,
then $s\pcw t$ denotes the componentwise addition of finite
sequences $s,t$ of equal length.
%
We now define 
\dm
S\pcw T=\ens{s\pcw t}{s\in S\land t\in T\land\lh s=\lh t}
\dm
for any trees $S\yi T\sq\nlw.$
The following lemma
shows that the componentwise addition of trees behaves 
somewhat like the ``equal-length''  cartesian product 
$S\ti T=\ens{\ang{s,t}}{s\in S\land t\in T\land \lh s=\lh t}$.

\ble
\label{r}
Let\/ $S\yi T\sq\nlw$ be any trees.
The tree\/ $W=S\pcw T$ is ill-founded iff so are both\/ $S$
and\/ $T.$
\ele
\bpf
In the nontrivial direction,
suppose that $\ga\in\nw$ is an infinite branch in $W,$ \ie,
$\ga\res n\in W$ for all $n.$ 
Then, for each $n,$ there exist $s_n\in S$ and $t_n\in T$ of 
length $n$ such that $s_n\pcw t_n=\ga\res n.$ 
The sequences $s_n\yi t_n$ then belong to 
$\ens{t\in\nlw}{t\lcw \ga\res\lh t},$ a finite-branching tree. 
Therefore, by K\"onig's lemma, there exist infinite sequences 
$\al\yi\ba\in\nw$ such that 
$\kaz m\:\sus n\ge m\:
({\al\res m=s_n\res m} \land {\ba\res m=t_n\res m}).$
Then $\al\yi\ba$ are infinite brances in resp.\ $S\yi T,$   
as required.
\epf

\punk{Louveau -- Rosendal transform}
\label{LR}

Suppose that $A$ is a $\fs11$ subset of $\tw\ti\tw.$ 
It is known from elementary topology of Polish spaces that any 
$\fs11$ subset of a Polish space $S$ is is equal to the projection 
of a closed subset of $S\ti\nw$ on $S.$
Thus there exists a closed set $P\sq \tw\ti\tw\ti\nw$ 
satisfying 
$A=\dom P=\ens{\ang{x,y}}{\sus z\:P(x,y,z)}.$ 
Further, there is a tree $R\sq{(2\ti2\ti\dN)}\lom$ 
(a tree on $2\ti2\ti\dN$) such that 
$P=[R]=\ens{\ang{x,y,\ga}}
{\kaz n\:R(x \res n, y \res n,\ga\res n)},$ 
and hence 
\bus
\label{e1}
\ang{x,y}\in A \; \eqv \; 
\ceq Rxy=\ens{s\in\nlw}{R(x\res\lh s,y\res\lh s,s)}
\,\text{ is ill-founded.}
%
%
\eus
(Obviously $\ceq Rxy$ is a tree in $\nlw.$) 
If $A$ is an arbitrary $\fs11$ set then, perhaps, not much can 
be established regarding the structure of a tree $R$ which 
generates $A$ in the sense of \eqref{e1}. 
However, assuming that $A={\rE}$ is an equivalence relation on 
$\tw,$ we can expect a nicer behaviour of $R.$ 
This is indeed the case.

The following key definition goes back to \cite{lr,roz}.
 
\bdf
\label{U1+}
A tree $T$ on a set of the form $X\ti\dN$ is {\it normal\/}
\index{tree!normal}%
if for any $u\in X\lom$ and $s\yi t\in \nlw$ such that 
$\lh u=\lh s=\lh t$ and $s\lcw t,$  we have 
$\ang{u,s}\in T\imp \ang{u,t}\in T.$ 
\edf

Thus normality means that the tree is \dd\lcw closed
upwards \vrt\ the second component.
$X=2\ti2$ in the next theorem, 
and the case $X=2=\ans{0,1}$ will also be considered.
But in all cases ${(X\ti\dN)}\lom$ itself is a normal tree.

\vyk{
\bdf
\label{rt}
Suppose that $\xi<\omi.$ 
A tree $R\sq{(2\ti2\ti\dN)}\lom$ is \dd\xi{\it restricted\/} 
if and only if for any $x\yi y\in\tw$ we have 
$\rah{\ceq Rxy}=\iy$ or $\rah{\ceq Rxy}<\xi$. 
\edf
}

\bte
\label{q2r}
Suppose that\/ $Q\sq{(2\ti2\ti\dN)}\lom$ is a tree and 
the set\/ 
\bus
\label{QR}
{\rE} \,=\, 
\ens{\ang{x,y}\in\tw\ti\tw} 
{\ceq Qxy\,\text{\rm\ is ill-founded}}
\eus
is an \eqr{} on\/ $\tw.$
Then there is a  
tree\/ $R\sq{(2\ti2\ti\dN)}\lom$ satisfying the following  
requirements\/\, \ref{1+} -- \ref{last}$:$
\ben
\def\theenumi{{\rm(\roman{enumi})}}
\def\labelenumi{{{\rm(\roman{enumi})}}}
\vyk{\itla0
still 
${x \rE y } \,\eqv\, 
\text{the set }\,R_{x  y }\;\,\text{is ill-founded}$;
}

\itla{1+}
symmetry$:$ $R(u,v,s)\eqv R(v,u,s),$ hence\/ 
$\ceq Rxy=\ceq Ryx$ for all $x \yo y\,;$

\itla{2+}
if\/ $u\in\tw\yt s\in\nw\yt\lh s=\lh u$ then\/ $R(u,u,s)\,;$

\itla{3+}
normality$:$ 
if\/ $R(u,v,s)\yt t\in\nw,$ and\/ $s\lcw t,$ then\/ $R(u,v,t)\,;$

\itla{4+}
transitivity$:$ 
if\/ $R(u,v,s)$ and\/ $R(v,w,t)$ then\/ $R(u,w,s\pcw t)\,;$

\itla{5+}
\label{last}
for any\/ $x\yi y\in\tw,$ 
${\ceq Rxy}$ is ill-founded iff so is\/ ${\ceq Qxy}$ --- 
and hence\/ \eqref{QR} holds for the tree\/ $R$ 
instead of\/ $Q\,;$ 
\een
\ete

This theorem is equal to Theorem 4 in \cite{lr}. 
 
\bpf
{\it Part 1\/}.
We observe that the tree   
\dm
\wQ=Q\cup
\ens{\ang{u,u,s}}
{u\in\tw\land s\in\nw\land\lh s=\lh u}
\cup
\ens{\ang{u,v,s}}{Q(v,u,s)}.
\dm 
satisfies $\ceq \wQ xy=\ceq Qxy\cup\ceq Qyx\cup D_{xy},$ 
where $D_{xy}=\nlw$ provided $x=y$ and $D_{xy}=\pu$ otherwise.
It easily follows that \eqref{QR} still holds for $\wQ.$   
In addition, $\wQ$ obviously satisfies both \ref{1+} and \ref{2+}.
Thus {\it we can assume, from the beginning, that\/ $Q$ 
satisfies both\/ \ref{1+} and \ref{2+}\/}.\vom

{\it Part 2\/}.
In this assumption, to fulfill \ref{3+}, we define 
\dm
\wQ=\ens{\ang{u,v,t}\in {(2\ti2\ti\dN)}\lom}
{\sus\ang{u,v,s}\in Q\:(s\lcw t)}.
\dm 
This is still a tree on $2\ti2\ti\dN,$ containing $Q$ 
and satisfying \ref{1+}, \ref{2+}, \ref{3+}. 
In addition, we have $\ceq \wQ xy=\ceq Qxy\pcw\tlw$ 
for any $x\yi y\in\tw,$ therefore the trees 
$\ceq Qxy$ and $\ceq\wQ xy$ are ill-founded simultaneously 
by Lemma~\ref r.
It follows that \eqref{QR} still holds for $\wQ.$ 
Thus, 
{\it we can assume that\/ $Q$ itself
satisfies\/ \ref{1+}, \ref{2+}, \ref{3+}\/}.\vom 

{\it Part 3\/}.
It is somewhat more difficult to fulfill \ref{4+}. 
A straightforward plan would be
to define a new tree $R$ containing all triples 
of the form $\ang{u_0,u_{n+1},s_0\pcw \cdots\pcw s_k},$ where 
$\ang{u_i,u_{i+1},s_i}\in Q$ for all $i=0,1,\dots,k.$ 
However, to work properly, such a construction has to be 
equipped with a kind of counter for the number $k$ of steps
in the finite chain. 
This idea can be realized as follows. 

Working in the assumption that $Q$ 
satisfies \ref{1+}, \ref{2+}, \ref{3+} (see Part 2), we define a 
tree $R\sq{(2\ti2\ti\dN)}\lom$ as follows. 
Suppose that $n\in\dN,$ $u\yi v\in2^n,$ $s\in\dN^n,$ 
$k\in\dN,$ and $i\yi j\in 2=\ans{0,1}.$ 
We put $\ang{u\we i,v\we j,k\we s}\in R$ iff
\bus
\label{e3+}
\sus u_0,u_1,\dots,u_k\in 2^n\:
(u_0=u\land u_k=v\land \kaz\ell<k\:
Q(u_\ell,u_{\ell+1},s)).
\eus
In addition, we put $\ang{\La,\La,\La}\in R,$ of course.
($\La$ is the empty sequence.) 
Note that $R$ is a tree on $2\ti2\ti\dN$ because so is $Q.$ 

We claim that, in our assumptions, the tree\/ $R$ 
satisfies all of \ref{1+} -- \ref{last}. 

\ref{1+}
If $u_0,\dots,u_k$ witness $R(u\we i,v\we j,k\we s)$ 
then the reversed sequence $u_k,\dots,u_0$ witnesses 
$R(v\we j,u\we i,k\we s)$ in the sense of \eqref{e3+}, 
because the tree $Q$ satisfies \ref{1+}.

\ref{3+}
Suppose that $\ang{u\we i,v\we j,k\we s}\in R,$ and let 
$u_0,\dots,u_k$ witness \eqref{e3+}. 
Let $n=\lh u=\lh v=\lh s=\lh u_\ell\yd\kaz \ell.$ 
Suppose that $k\le k'$ and $s\lcw s'$ (still $\lh {s'}=n$). 
Put $u_\ell=v$ whenever $k<\ell\le k'.$ 
Note that $Q(u_\ell,u_{\ell+1},s)$ also holds for 
$k<\ell< k'$ by \ref{2+} for $Q.$ 
(Indeed, in this case $u_\ell=u_{\ell+1}.$)  
Thus, 
$Q(u_\ell,u_{\ell+1},s')$ holds for all $\ell<k'$ 
by \ref{3+} for $Q.$ 
By definition, this witnesses 
$\ang{u\we i,v\we j,k'\we s'}\in R,$ as required.

\ref{2+}
If $k=0$ and $u=v$ then \eqref{e3} obviously holds 
(with the empty list of intermediate sequences $u_1,\dots,u_{k-1}$), 
and hence $R(u\we i,u\we j,0\we s)$ holds for all 
$u\in\tw\yt s\in\nw$ of 
equal length, in particular, $R(u,u,0^n)$ for all $n$ and 
$u\in\nw$ with $\lh u=n.$  
It remains to apply property \ref{3+} just proved.

\ref{4+}
Suppose that the triples 
$\ang{u\we i,v\we j,k\we s}$ and 
$\ang{v\we j,w\we \rho,\ka\we\sg}$ belong to $R,$ 
and $n$ is the length of all sequences 
$u\yi v\yi s\yi w\yi t.$ 
Let $R(u\we i,v\we j,k\we s)$ be witnessed, in the sense 
of \eqref{e3+}, by $u_0,\dots,u_k$ and, accordingly, 
$R(v\we j,w\we \rho,\ka\we\sg)$ be witnessed by 
$v_0,\dots,v_\ka.$
(All $u_\ell$ and $v_\ell$ belong to $2^n.$)
Since $Q$ satisfies \ref{3+}, the same sequences also witness 
$R(u\we i,v\we j,k\we t)$ and  
$R(v\we j,w\we \rho,\ka\we t),$
where $t=s\pcw\sg$ (componentwise). 
It easily follows that the concatenated complex  
$u_0,\dots,u_{k-1},u_k=v_0,v_1,\dots,v_\ka$
witnesses $R(u\we i,w\we \rho,(k+\ka)\we t),$
as required. 

\ref{5+}
We observe that, by definition, 
$Q(u,v,s)\imp R(u\we i,v\we j,1\we s)$ for any $i\yo j=0,1.$
It follows that, for any $x\yo y\in\tw,$ 
$s\in \ceq Qxy\imp 1\we s\in \ceq Rxy,$
and hence $\ceq Rxy$ is ill-founded provided so is 
$\ceq Qxy.$
\vyk{
\dm
\rah{\ceq Rxy}\ge \rag{\ceq Rxy}{\ang1}+1\ge 
\rag{\ceq Qxy}{\La}+1=\rah{\ceq Qxy}+1.
\dm
Thus $\rah{\ceq Rxy}\ge \rah{\ceq Qxy}+1.$
}%
The inverse implication in \ref{5+} needs more work. 
This argument belongs to Louveau and Rosendal \cite{lr}.
Assume that $\ceq Rxy$ is ill-founded, that is, there exists 
an infinite sequence $\da\in\nw$ such that 
$\kaz n\:R(x\res n,y\res n,\da\res n).$ 
Let $k=\da(0)$ and $\ga(m)=\da(m+1)$ for all $m,$ 
so that $\da=k\we\ga.$ 
By definition, for any $n$ there exist sequences 
$u^n_0,\dots,u^n_k\in 2^n$ such that 
$u^n_0= x \res n\yt u^n_k= y \res n,$ and 
$Q(u^n_\ell,u^n_{\ell+1},\ga\res n)$ 
for all $\ell<k.$ 
Each \dd{k+1}tuple $\ang{u^n_0,\dots,u^n_k}\in (2^n)^{k+1}$ 
can be considered as an \dd ntuple in $(2^{k+1})^n.$ 
By K\"onig's lemma, there exist infinite 
sequences $x_0,\dots,x_{k}\in\tw$ such that for any $m$ 
there is a number $n\ge m$ with 
$x _\ell\res m=u^n_\ell\res m$ for all $\ell\le k.$  
It follows that $x_0=x \yt x_k= y,$ and, as $Q$ is a 
tree, $Q(x_\ell\res m, x_{\ell+1}\res m,\ga\res m)$ holds 
for all $\ell<k$ and all $m.$ 
We conclude that $x_\ell\rE x_{\ell+1}$ for all $\ell<k$ 
by \eqref{QR} for $Q,$ 
therefore, $x \rE y$ because $\rE$ is an equivalence. 
Finally, 
$\ceq Qxy$ is ill-founded still by \eqref{QR} for $Q$.
\epf

\punk{Embedding and equivalence of normal trees}
\label{cae}

Let $\nt$ denote the set of all non-empty 
normal trees $T\sq (2\ti\dN)\lom.$
\index{zzNT@$\nt$}%
Suppose that $S,T\in\nt.$ 
The set of all finite sequences 
$f\in\nlw$ such that   
${\ang{u,s}\in S}\imp{\ang{u,s\pcw (f\res n)}\in T}$ 
for all $n\le \lh f$ and $u\in 2^n\yt s\in\dN^n,$ 
will be denoted by $\lip ST.$  
Obviously $\lip ST$ is a tree in $\nlw$ containing $\La$.
 
We proceed with the following key definition of \cite{lr}.

\bdf
\label{U2}
Define $S\lnt T$ iff the tree $\lip ST$ is ill-founded, 
that is,
\dm
\sus\ga\in\nw\:\kaz n\:\kaz u\in2^n\:\kaz s\in\dN^n\:
({\ang{u,s}\in S}\imp{ \ang{u,s\pcw {\ga\res n}}\in T }).
\dm
Define $S\ent T$ iff $S\lnt T$ and $T\lnt S$.~\snos
{$\lnt$ and $\ent$ are denoted in \cite{roz} by, resp., 
$\le^\ast_{\text{max}}$ and $E^\ast_{\text{max}}$.}\qed%
\eDf

Thus $S\lnt T$ indicates the existence of a
certain  shift-type embedding of $S$ into $T.$ 
The relation $\lnt$ is a partial order on
the set $\nt,$
To check that $\lnt$ is transitive,
suppose that $R\lnt S$ and $S\lnt T,$ where $R\yi S\yi T$ are 
normal trees in $(2\ti\dN)\lom.$ 
Then the trees $U=\lip RS$ and $V=\lip ST$ (trees in $\nlw$) 
are ill-founded, and hence so is $W=U\pcw V$ 
by Lemma~\ref r. 
On the other hand, easy verification shows that $W\sq\lip RT.$ 
Thus $\lip RT$ is ill-founded, as required.
It follows that $\ent$ is an \eqr\ on $\nt.$

Moreover, applying the componentwise addition to the sequences 
$\ga$ that witness $\lnt,$ one proves that $S\ent T$ is 
equivalent to the existence of $\ga\in\nw$ such that for all 
$n$ and all $u\in2^n\yt s\in\dN^n$ the following holds 
simultaneously: 
\dm
{\ang{u,s}\in S}\imp{ \ang{u,s\pcw {\ga\res n}}\in T}
\:\text{ and }\:
{\ang{u,s}\in T}\imp{ \ang{u,s\pcw {\ga\res n}}\in S}.
\dm

\bcor
\label{e0}
If\/ $S,T\in\nt$ then\/ $S\ent T$ iff the tree\/
$\lip ST\cap\lip TS$ is ill-founded.\qed
\ecor

Note that any tree $T\in\nt$ is, by definition, 
a subset of the countable set ${(2\ti\dN)}\lom,$ 
Thus, $\nt$ is a subset of the Polish space 
$\pws{{(2\ti\dN)}\lom},$ 
identified, as usual, with the product space 
$2^{{(2\ti\dN)}\lom}.$  
(Elementary computations show that in fact $\nt$ is a 
closed set.)
Therefore, the relations $\lnt$ and $\ent$ are, 
formally, subsets of 
$\pws{{(2\ti\dN)}\lom}\ti\pws{{(2\ti\dN)}\lom}$.

\ble
\label{U2l}
$\lnt$ and\/ $\ent$ are\/ $\fs11$ relations.
\ele
\bpf
Straightforward estimations. 
The principal quantifier expresses the existence of 
$\ga\in\nw$ with certain properties.
\epf

It occurs that $\ent$ belongs to a special
type of $\fs11$ \eqr s.

\bdf
\label{uer}
An $\fs11$ \eqr\  
$\rU$ is {\it universal\/}, or {\it complete\/},
\index{equivalence relation, ER!complete}%
\index{equivalence relation, ER!universal}%
if and only if $\rF\reb\rU$ holds 
for any other $\fs11$ \eqr\ $\rF.$
\edf

There is a simple construction that yields a
universal $\fs11$ \eqr. 

\bex
\label u
We begin with a $\fs11$ set $U\sq{(\dnn)}{}^3,$ 
universal in the sense that for any $\fs11$ set $P\sq{(\dnn)}{}^2$ 
there is an index $x\in\dnn$ such that $P$ is equal to the 
cross-section $U_x=\ens{\ang{y,z}}{\ang{x,y,z}\in U}.$ 
Define a set $P\sq {(\dnn)}{}^3$ so that every cross-section 
$P_x$ is equal to the {\it equivalence hull\/} of $U_x,$ 
that is, to the least equivalence relation containing $U_x.$ 
Formally, $\ang{y,z}\in P_x$ iff there is a finite chain 
$y=y_0,y_1,y_2,\dots y_n,y_{n+1}=z$ such that, for any 
$k\le n,$ either $\ang{y_k,y_{k+1}}$ belongs to $U_x,$ or 
$\ang{y_{k+1},y_{k}}$ belongs to $U_x,$ or just $y_k=y_{k+1}.$ 

Clearly $P$ is still a $\fs11$ subset of ${(\dnn)}{}^3,$
with each $P_x$ being a $\fs11$ \eqr. 
Moreover, if $U_x$ is an \eqr\ then $P_x=U_x.$ 
Thus the family of all cross-sections $P_x\yd x\in\dnn,$ is 
equal to the family of all $\fs11$ \eqr s on $\dnn.$ 
We claim that the \eqr\ $\rU$ on ${(\dnn)}{}^2,$ defined 
so that $\ang{x,y}\rU\ang{x',y'}$ iff 
$x=x'$ and $\ang{y,y'}\in P_x,$ is universal.
For take any $\fs11$ equivalence $\rF$ on $\dnn.$ 
Then $\rF=P_x$ for some $x$ by the above, therefore, the map 
$\vt(y)=\ang{x,y}$ is a continuous reduction of $\rF$ to $\rU,$ 
as required.  
\eex

\bte[{{\rm Theorem~5 in \cite{lr}}}]
\label{max}
$\ent$ is a universal\/ $\fs11$ equivalence on\/ $\nt.$
\ete
\bpf
Consider any $\fs11$ \eqr\ $\rE$ on $\tw.$ 
Then $\rE$ is a $\fs11$ subset of $\tw\ti\tw,$ and hence 
there is a tree $Q\sq{(2\ti2\ti\dN)}\lom$ 
(a tree on $2\ti2\ti\dN$) such that, for all $x\yo y\in\tw,$ 
\bus
\label{,}
{x\rE y}\leqv \,
\text{the cross-section tree $\ceq Qxy$ is ill-founded}.
\eus
It follows from Theorem~\ref{q2r} that it can be assumed that  
$Q$ satisfies requirements 
\ref{1+} -- \ref{last} of Theorem~\ref{q2r}. 
%
We claim that the map
\bus
\label{vt}
x\lmto \vt(x )=\ens{\ang{u,s}\in {(2\ti\dN)}\lom}
{Q(u,x\res\lh u,s)}
\qquad
(x\in\tw)
\eus 
is a Borel reduction of $\rE$ to $\ent.$ 
That $\vt$ is a Borel, even continuous map, is rather easy. 
That $\vt(x)\in\nt$ immediately follows from \ref{3+}.
The reduction property follows from the next lemma.

\ble
\label j
If a tree\/ $Q\sq{(2\ti2\ti\dN)}\lom$ satisfies requirements\/ 
\ref{1+} -- \ref{4+} of Theorem~\ref{q2r}, and\/ $x\yo y\in\tw,$ 
then\/ $\lip{\vt(x)}{\vt(y)}=\ceq Qxy$. 
\ele
\bpf 
Suppose that $f\in\lip{\vt(x)}{\vt(y)},$ $m=\lh f.$ 
Then, by definition, we have   
$Q(u,x\res m,s)\imp R(u,y\res m,s\pcw f)$ 
for all $u\in2^m$ and $s\in\dN^m.$ 
Take here $u=x\res m$ and $s=0^m$ 
(the sequence of $m$ 0s); 
then 
$Q(x\res m,x\res m,0^n)\imp Q(x\res m,y\res m,f).$ 
Yet the left-hand side holds by \ref{2+}. 
Therefore, the right-hand side holds, thus $f\in \ceq Qxy$. 

To prove the converse let $f\in \ceq Qxy,$ 
that is, $Q(x\res m,y\res m,f),$ 
where $m=\lh f$ --- and hence $Q(x\res n,y \res n,f\res n)$ 
for any $n\le m$ as $Q$ is a tree. 
Assume that $n\le m$ and $u\in2^n\yt s\in \dN^n.$ 
We have to prove 
$Q(u,x\res n,s)\imp Q(u,y\res n,s\pcw (f\res n)).$
So suppose $Q(u, x \res n,s).$ 
In addition, $Q(x\res n,y\res n,f\res n)$ holds by the 
above. 
Then $Q(u, y \res n,s\pcw (f\res n))$ holds by \ref{4+}, 
as required.
\epF{Lemma}

To accomplish the proof of Theorem~\ref{max}, suppose that 
$x \yi y \in\tw.$ 
Then $x\rE y$ iff the tree $\ceq Rxy$ is ill-founded, 
iff (by the lemma) $\lip{\vt(y)}{\vt(x)}$ is ill-founded, 
iff $\vt(x)\ent\vt(y)$ 
(by Definition~\ref{U2}).

\epF{Theorem~\ref{max}}

\punk{Reduction to Borel ideals: first approach}
\label{nti1}

We present two different proofs of Theorem~\ref{M}.
The first one, due to Rosendal~\cite{roz}, involves 
the ideal $\ndi$ on ${(2\ti\dN)}\lom$ finitely 
generated by all sets of the form $S\sd T,$ where 
$S,T\sq{(2\ti\dN)}\lom$ are normal trees and $S\ent T.$ 
Thus $\ndi$ consists of all subsets of ${(2\ti\dN)}\lom,$ 
covered by unions of finitely many symmetric differences $S\sd T$ 
of the type just indicated.

\bte
\label{roi}
The ideal\/ $\ndi$ is\/ $\fs11$ as a subset of the Polish space\/ 
$\pws{{(2\ti\dN)}\lom}.$  
Furthermore, the \eqr\/ $\ent$ is equal to\/ 
${\rE_{\cI}}\res\nt$ --- 
this means that for any\/ $S,T\in\nt,$ the following holds$:$
$S\ent T$ if ant only if\/ $S\sd T\in\ndi$. 
\ete
\bpf
That $\ndi$ is $\fs11$ is quite clear: 
the principal quantifier expresses the existence of a 
finite collection of elements of $\nt,$ whose properties are 
expressible still by a $\fs11$ relation because $\ent$ is $\fs11$.

Suppose that $S\sd T\in\ndi,$ and prove $S\ent T$ 
(the nontrivial direction).
By definition  
$S\sd T\sq \bigcup_{i=1}^k(S_i\sd T_i),$  where 
$S_i,T_i\in\nt$ and $S_i\ent T_i.$ 
Then the trees
$R_i=\lip{S_i}{T_i}\cap\lip{T_i}{S_i}$
are ill-founded by Corollary~\ref{e0}.  
We have to prove that $\lip ST$ and $\lip TS$
are ill-founded trees, too. 
To check the ill-foundedness of $\lip ST,$ 
note that the tree $R=R_1\pcw \cdots\pcw R_k$
is ill-founded by Lemma~\ref r. 
Thus it remains to prove that $R\sq \lip ST.$ 

Consider any 
$r=r_1\pcw \cdots\pcw r_k\in R,$ where all sequences 
$r_i\in R_i\yt i=1,\dots,k,$ 
have one and the same length, say $m.$ 
Suppose towards the contrary that 
$r\nin\lip ST,$ \ie\ there exists a pair $\ang{u,s}\in S$ 
such that $\ang{u,s\pcw (r\res n)}\nin T,$
where $n=\lh u=\lh s\le m.$
Then $(\ast)\msur$ $\ang{u,s\pcw r'}\nin T$ whenever
$r'\in 2^n\yt r'\lcw{r\res n}.$ 
In particular, $\ang{u,s}\nin T$ by the normality,
and hence $\ang{u,s}\in S\sd T,$ 
thus $\ang{u,s}\in S_{i_1}\sd T_{i_1}$ for some $1\le i_1\le k.$ 
This implies $\ang{u,s_1}\in S_{i_1}\cap T_{i_1},$ where 
$s_1=s\pcw (r_{i_1}\res n).$
(Indeed we have $\ang{u,s}\in S_{i_1}\cup T_{i_1}$ by the
choice of $i_1.$ 
If say $\ang{u,s}\in S_{i_1}$ then $\ang{u,s_1}\in T_{i_1}$
because $r_{i_1}\in R_{i_1}\sq\lip{S_{i_1}}{T_{i_1}}.$
In addition $\ang{u,s_1}\in S_{i_1}$ by the normality of $S_{i_1}$.)

Once again, $\ang{u,s_1}\in S\bez T$ by $(\ast)$ above.
It follows that 
$\ang{u,s_1}\in S_{i_2}\sd T_{i_2}$ for some $1\le i_2\le k$
by the same argument. 
This implies $\ang{u,s_2}\in S_{i_2}\cap T_{i_2},$ where 
$s_2=s_1\pcw (r_{i_2}\res n),$ because $r_{i_2}$ belongs to 
$R_{i_2}.$
Note that $i_2\ne i_1$ as $\ang{u,s_1}\in S_{i_1}\cap T_{i_1},$
and still $\ang{u,s_2}\in S_{i_1}\cap T_{i_1}$
since $S_i$ and $T_i$ are normal trees.

After $k$ steps of this construction, all indices $1\le i\le k$ 
will be considered, and the final sequence 
$s_k=s\pcw (r\res n)$ will satisfy $\ang{u,s_k}\in S_{i}\cap T_{i}$ 
for all $i=1,\dots,k.$ 
It follows that $\ang{u,s_k}\nin S\sd T.$ 
However $\ang{u,s_k}\in S$ since $\ang{u,s}\in S$ and $S$ is 
a normal tree. 
Thus $\ang{u,s_k}$ belongs to $T,$ contrary to the above.
\epf

Theorems \ref{roi} and \ref{max} imply

\bcor
\label{roic}
$\rE_{\ndi}$ is a universal\/ $\fs11$ equivalence relation.\qed 
\ecor

\vyk{

\punk{Borel approximations}
\label b

Thus we have defined a $\fs11$ ideal $\cI=\ndi$ on
a countable set $X= (2\ti\dN){}\lom,$
which generates a universal  $\fs11$ \eqr\ $\rE_{\ndi}.$
}

Let us show now that these properties of $\ndi$ suffice to prove
Theorem~\ref{M}.

We begin with a very general fact of basic descriptive set theory:
as any $\fs11$ set, $\ndi$ can be presented in the form
$\ndi=\bigcap_{\xi<\omi}\ndi^\xi,$ where $\ndi^\xi$ are Borel
subsets of $\cP((2\ti\dN){}\lom),$
${\xi<\eta}\imp \ndi^\eta\sq\ndi^\xi,$
and for any $\fp11$
set $X$ in the same Polish space containing $\ndi$ there
is an ordinal $\xi<\omi$ such that $\ndi^\xi\sq X.$~\snos
{This {\it index restriction\/} property was first established
by Lusin and Sierpi\'nski \cite{ls}, essentially in the dual
form saying that the canonical representation of any $\fp11$
set $C$ in the form $C=\bigcap_{\xi<\omi}C_\xi$ of a union of
\dd\sq increasing Borel approximations has the property
that for any $\fs11$ set $X\sq C$ there is an index $\xi<\omi$
with $X\sq C_\xi.$ 
The shortest proof consists of observation that otherwise the
relation $x\preccurlyeq y$ iff $x$ appears in sets $C_\xi$ not later
than $y$ on $X$ is a $\fs11$ prevellordering of uncountable
length, contrary to the Kunen -- Martin prewellordering theorem  
(see, \eg, \cite[2G.2]{mDST}).}
The sets $\ndi^\xi$ are called
{\it (upper) Borel approximations\/} of $\ndi$.

The following lemma is the key fact.

\ble
\label{kl.f}
For any\/ $\xi<\omi$ there exists an ordinal\/ 
$\nu\yt\xi<\nu<\omi,$ such that the Borel approximation\/
$\ndi^\nu$ is still an ideal. 
\ele
\bpf
{\it Step 1\/}: we claim that for any $\xi<\omi$ there is 
an ordinal $\eta=\eta(\xi)\yt\xi<\eta<\omi,$ such that 
${y\sq x\in \ndi^\eta}\imp{y\in \ndi^\xi}.$ 
Indeed the set 
$P=\ens{x\in\ndi^\xi}{\kaz y\sq x\:(y\in \ndi^\xi)}$ 
is a $\fp11$ superset of $\ndi$   
(since $\ndi$ is an ideal). 
It follows that there is an ordinal $\eta>\xi$ with 
$\ndi^\eta\sq P.$

{\it Step 2\/}: we claim that for any $\xi<\omi$ there is 
an ordinal $\za=\za(\xi)\yt\xi<\za<\omi,$ such that 
${x,y\in \ndi^\za}\imp{x\cup y\in \ndi^\xi}.$ 
The argument contains two substeps. 
First, the set 
$
X=\ens{x\in\ndi^\xi}{\kaz y\in\ndi\:
(x\cup y\in\ndi^\xi)} 
$       
is a $\fp11$ superset of $\ndi$ since $\ndi$ is an ideal. 
Thus there is an ordinal $\al>\xi$ with $\ndi^\al\sq X.$ 
Then we have $x\cup y\in \ndi^\xi$ whenever $x\in \ndi^\al$ 
and $y\in \ndi.$ 
It follows that the $\fp11$ set
$
Y=\ens{y\in\ndi^\al}{\kaz x\in\ndi^\al\: 
(x\cup y\in\ndi^\xi)} 
$ 
is a superset of $\ndi,$    
and hence there is an an ordinal 
$\eta>\al$ such that $\ndi^\eta\sq Y.$ 
Obviously $\eta$ is as required.

{\it Final argument\/}. 
Put $\xi_0=\xi$ and $\xi_{n+1}=\eta(\za(\xi_n))$ for all $n.$ 
The ordinal $\nu=\tsup_n\xi_n$ is as required. 
\epf

It follows that the set
$\Xi=\ens{\xi<\omi}{\ndi^\xi\,\text{ is an ideal}}$
is unbounded in $\omi.$
We also note that $\rE_{\ndi^\xi}$ is a Borel \eqr\ on
$\cP((2\ti\dN){}\lom)$ for any $\xi\in\Xi,$ and the sequence
of these equivalences is \dd\sq decreasing and satisfies 
$\rE_{\ndi}=\bigcap_{\xi\in\Xi}\rE_{\ndi^\xi}.$
The proof of Theorem~\ref{M}, our main result
here, is accomplished with the following lemma.

\ble
\label{w+}
If\/ $\rE$ is a Borel \eqr\ on a Polish space\/ $X$ then there 
is an ordinal\/ $\xi\in\Xi$ such that\/ $\rE\reb\rE_{\ndi^\xi}$.
\ele
\bpf
It follows from Corollary~\ref{roic} that $\rE\reb\rE_{\ndi},$  
that is, there exists a Borel map $\vt:X\to\cP((2\ti\dN){}\lom)$
such that ${x\rE y}\eqv{\vt(x)\sd\vt(y)\in\ndi}.$ 
Thus the full \dd\vt image $\vt\imb P$ of the set
$P={(X\ti X)}\bez{\rE}$ is a $\fs11$ set disjoint from $\ndi.$ 
Then by Lemma~\ref{kl.f} there is an ordinal $\xi\in \Xi$
such that $\vt\imb P$ does not intersect $\ndi^\xi,$ too. 
Thus $\vt$ reduces $\rE$ not only to $\rE_{\ndi}$ but also to 
the approximating Borel \eqr{} $\rE_{\ndi^\xi}$.
\epf

\qeDD{Theorem~\ref{M}, first proof}

\punk{Reduction to Borel ideals: second approach}
\label{nti2}

Is there any method to prove Theorem~\ref{M} by a sequence
of more ``effective'' and mathematically meaningful upper
Borel approximations of a \dd\reb maximal analytic ideal?
Paper~\cite{kl5} suggested a suitable definition. 

First of all recall that any tree $T\sq X\lom$ admits  
{\it the rank function\/},   
\index{rank function}%
\index{zzrnkR@$\raf R$}%
a unique map $\raf R:R\to\Ord\cup\ans\iy,$ where
$\iy$ denote a formal element bigger than any ordinal,
satisfying the following requirements:
\ben
\tenu{(\alph{enumi})}
\itla{iy0}\msur 
$\rag Rr=-1$ whenever $r\nin R$;

\itla{iy2}\msur
$\rag Rr=\sup_{r\we n\in R}\rag R{r\we n}$ for any 
$r\in R.$~\snos 
{\label{sup}%
We define $\sup \Om,$ for $\Om\sq\Ord,$ to be the least ordinal 
strictly bigger than all ordinals in $\Om.$ 
We also define $\sup\Om=\iy$ provided $\Om$ contains $\iy.$} 
In particular, $\rag Rr=0$ if and only if $r\in R$ is 
a \dd\sq maximal element of $R\,;$ 

\itla{iy1}\msur 
$\rag Rr=\iy$ if and only if $R$ has an infinite branch 
containing $r,$ 
\ie, there exists $\ga\in X^\om$ such that $\ga\res n\in R$ 
for all $n,$ and $\ga\res{\lh r}=r$.
\een
In addition, put $\rah \pu=-1$ for the empty tree $\pu,$
and $\rah R=\rag R\La$ for any non-empty tree $R.$ 
($\La,$ the empty sequence, belongs to any 
tree $\pu\ne R\sq X\lom.$)  
Obviously any tree $R$ is
well-founded iff $\rah R<\iy.$

\bdf
\label V
Suppose that $S,T\in\nt$ and $\xi<\omi.$ 

Define $S\lmt\xi T$ iff the tree $\lip ST$  
satisfies $\rah{\lip ST}\ge\xi.$~\snos
{The inequality $\rah{\lip ST}\ge\xi$ means that either 
$\lip ST$ (a tree in $\nlw$) is ill-founded
(then $\rah{\lip ST}=\iy$)
or it is well-founded and its rank is an ordinal 
$\ge\xi.$}

Define $S\emt\xi T$ iff both $S\lmt\xi T$ and 
$T\lmt\xi S$.
\edf

It is demonstrated in \cite{kl5} by simple and rather
straightforward arguments that all relations $\emt\xi$
are Borel \eqr s on $\nt,$ of certain explicitly defined
Borel ranks.
A notable part of this result is the proof of transitivity
of $\lmt\xi$ and $\emt\xi,$ based on the following
generalization of Lemma~\ref{r}.

\ble[\rm Lemma 4 in \cite{kl5}]
\label{rx}
We have\/ $\rah {S\pcw T}=\tmin\ans{\rah S,\rah T}$
for any trees\/ $S,T\sq\dN\lom,$
well- or ill-founded independently of each other.\qed
\ele

In addition,
${\ent}=\bigcap_{\xi<\omi}{\emt\xi},$ and this intersection
has the same restriction property as above: if $P$ is a
$\fp11$ subset of $\nt\ti\nt$ containing $\ent$ then there
is an ordinal $\xi<\omi$ such that $\emt\xi\sq P.$ 

It follows, essentially by the same arguments as above,
that the sequence of Borel relations $\emt\xi$ is \dd\reb
cofinal among all Borel \eqr s.


The following construction of Borel ideals that generate the
\eqr s $\emt\xi$ is a modification of a construction in \cite{kl5}.

Consider a set $X\sq {(2\ti\dN)}\lom.$
Suppose that $f\in\nlw\yt u\in \tlw\yt n=\lh u\le\lh f.$ 
Let $\gam fuX$ be the game in which \pli\ plays
$s_1,s_2,\dots\in \dN^n,$
\pld\ plays $t_1,t_2\dots\in \dN^n$ so that
$t_1\pcw\cdots\pcw t_m\lcw f\res n$ for all $m,$
and \pli\ wins if and only if 
$\ang{u,\ts_k}\in X$ for all $k,$ where
$\ts_k=s_1\pcw t_1\pcw\cdots\pcw s_{k-1}\pcw t_{k-1}\pcw s_{k}.$

Define $\vid X$ to be the tree of all $f\in\nlw$ such that
for any $n\le\lh f$ and $u\in2^n$ 
\pld\ has a winning strategy in $\gam fuX.$ 
Thus, informally, $f\in \vid X$ can be seen as a statement of 
the possibility to leave $X$ for good in finitely many steps, 
the \dd\pcw total length of which is at most $f.$ 
%
%
Let $\nti$ be the collection of all sets 
$X\sq {(2\ti\dN)}\lom$ 
such that $\vid X$ is ill-founded.
For $\xi<\omi,$ let $\mti\xi$ be the collection of all 
sets $X\sq {(2\ti\dN)}\lom$ 
with $\rah{\vid X}\ge\xi$.

\ble
\label{jjid}
$\nti$ and all sets\/ $\mti\xi$ are ideals on\/ ${(2\ti\dN)}\lom$.
\ele
\bpf
Suppose that sets $X\zi Y\sq (2\ti\dN)\lom$ belong to $\nti,$ 
and hence the trees $F=\vid X$ and $G=\vid Y$ are ill-founded.
Then the tree $F\pcw G$ is ill-founded by Lemma~\ref r
(to be replaced by Lemma~\ref{rx} for the ideals $\mti\xi$),
and hence it
suffices to prove that $F\pcw G\sq\vid{X\cup Y}.$ 

Take any $f\in F$ and $g\in G$ with $\lh f=\lh g.$
To prove that $h=f\pcw g$ belongs to $\vid{X\cup Y}$ fix
any $u\in2^n\yt n\le\lh f,$    
and a pair of winning strategies $\xi\yi\eta$
for \pld\ in games resp. $\gam fuX$ and $\gam guY.$
To describe a winning strategy for \pld\ in
$\gam hu{X\cup Y},$ let
$s_1,t_1,s_2,t_2,\dots$ be a full sequence of moves.
Put $K=\ens k{\ts_k\in X}$ and $K'=\ens k{\ts_k\in Y\bez X}.$
Let $K=\ans{k_1,k_2,\dots}$ and $K'=\ans{k'_1,k'_2,\dots},$ in 
the increasing order. 

For any $k,$ if $k=k_j\in K$ then \pld\ plays 
$t_k=\xi(\sg_1,\tau_1,\dots,\sg_{j-1},\tau_{j-1},\sg_j),$
where
$\tau_i=t_{k_i}$ and, for any $1\le i\le j$, 
\dm
\sg_i=s_{k_{i-1}+1}\pcw t_{k_{i-1}+1}\pcw s_{k_{i-1}+2}\pcw t_{k_{i-1}+2}
\pcw\cdots\pcw s_{k_i-1}\pcw t_{k_i-1}\pcw s_{k_i}.
\dm
Accordingly if $k=k'_j\in K'$ then
$t_k=\eta(\sg'_1,\tau'_1,\dots,\sg'_{j-1},\tau'_{j-1},\sg'_j),$
where
\dm
\sg'_i=s_{k'_{i-1}+1}\pcw t_{k'_{i-1}+1}\pcw s_{k'_{i-1}+2}\pcw
t_{k'_{i-1}+2}
\pcw\cdots\pcw s_{k'_i-1}\pcw t_{k'_i-1}\pcw s_{k'_i}
\dm
and $\tau'_i=t'_{k_i}$ for any $1\le i\le j.$
If to the contrary \pli\ wins then $K\cup K'=\dN.$
Let, say, $K=\ans{k_1,k_2,\dots}$ be infinite.
Then \pld\ must win the auxiliary play
$\sg_1,\tau_1,\sg_2,\tau_2,\dots$ in $\gam fuX,$ hence one of
the finite sums
$\tsg_j=\sg_1\pcw \tau_1\pcw\cdots
\pcw \sg_{j-1}\pcw\tau_{j-1}\pcw \sg_{j}$
satisfies $\tsg_j\nin X.$
But obviously $\tsg_j=\ts_{k_j},$ which is a contradiction with
$k_j\in K.$
\epf

Thus $\nti$ is a $\fs11$ ideal while 
each $\mti{\xi}$ is a Borel ideal. 

\bte
\label{il}
The \eqr\/ $\ent$ is equal to\/ ${\rE_{\nti}}\res\nt,$ 
while for any\/ $\xi,$ $\emt\xi$ is equal to\/ 
${\rE_{\mti\xi}}\res{\nt}$.
%
\ete
\bpf
Consider any $S,T\in\nt.$ 
Assume that $S\ent T.$ 
Then the trees $F=\lip ST$ and $G=\lip TS$ are ill-founded, and hence 
so is $H=F\pcw G$ by Lemma~\ref r.
(Lemma~\ref{rx} is used in the case of $\mti\xi.$)
Note that $H\sq G\cap F$ since both $S$ and $T$ are
\dd\lcw transitive to the right.
Thus it suffices to prove that $G\cap F\sq\vid{S\sd T}.$ 
Consider any $f\in G\cap F.$ 
By definition, for any $\ang{u,s}\in S\cup T\yt \lh u=\lh s=n\le\lh f,$ 
we have $\ang{u,s\pcw{(f\res n)}}\in S\cap T.$ 
In particular, 
${\ang{u,s}\in S\sd T}\imp {\ang{u,s\pcw{(f\res n)}}\nin S\sd T},$ 
and easily $f\in\vid{S\sd T}.$ 

To prove the converse, suppose that $S\sd T\in\nti,$ thus 
$\vid {S\sd T}$ is ill-founded. 
It suffices to prove that $\vid {S\sd T}\sq \lip ST.$
Suppose, towards the contrary, that $f\in \vid {S\sd T}$ but 
$f\nin \lip ST.$
The latter means that there exists a pair 
$\ang{u,s}\in S\yt \lh u=\lh s=n\le\lh f,$ such that 
$\ang{u,s\pcw{(f\res n)}}\nin T.$
Then also $\ang{u,s}\nin T,$ and hence both $\ang{u,s}$ and 
$\ang{u,s\pcw{(f\res n)}}$ belong to $S\bez T.$
It follows that 
\bit
\item[$(\ast)$]
$\ang{u,s+g}\in S\bez T$ for any  
$g\in\dN^n\yt g\lcw {(f\res n)}.$ 
\eit
Now consider a play in $\gam fu{S\sd T}$ in which \pld\ follows its 
winning strategy (which exists because $f\in \vid {S\sd T}$) 
while \pli\ plays $s_1=s$ and
$s_k=0^n$ (the sequence of $n$ zeros) on every move 
$k\ge2.$ 
Let $t_1,t_2,\dots$ be the sequence of \pld's moves. 
Then $t_1\pcw\cdots \pcw t_k\lcw {(f\res n)}$ for all $k,$ 
and hence, by $(\ast),$ the sum
$\ts_k=s\pcw t_1\pcw\cdots \pcw t_k$ satisfies
$\ang{u,\ts_k}\in {S\sd T},$ 
which contradicts the choice of the strategy.
\epf

\qeDD{Theorem~\ref{M}, second proof}

\vyk{
Rosendal proved in \cite{roz} that any Borel \eqr{} is Borel reducible 
to an \eqr\ of the form $\rE_\cJ,$ where $\cJ$ is a Borel ideal. 
The ideals $\cJ$ is his argument were taken as suitable upper Borel 
approximations of a certain $\fs11$ ideal $\cI.$ 
(The ideal $\cI$ was described in the beginning of this Subsection.) 
The following corollary shows that Borel ideals in Rosendal's 
result can be chosen in more direct and effective way.

\bcor
\label{rcor}
For any Borel \eqr\/ $\rE$ on a Polish space there is 
an ordinal\/ $\xi<\omi$ such that\/ ${\rE}\reb{\rE_{\mti\xi}}$.
\ecor
\bpf
According to Corollary~\ref{w}, there is an ordinal $\xi<\omi$ 
such that that $\rE\reb{\emt\xi}.$ 
On the other hand, ${\emt\xi}\reb{\rE_{\mti\xi}}$ by 
Theorem~\ref{il}\ref{il2}.
\epf

Lemma~\ref{nn}, Corollary~\ref w, and Theorem~\ref{il}  
enable us to immediately accomplish the proof of Theorem~\ref{M}: 
the sequence of 
\eqr s $\rE_\nu=\emt{\om\nu+2},$ 
where $\nu<\omi,$ is as required.
(Well, we have to additionally apply any homeomorphism of 
$\pws{(2\ti\dN)\lom}$ onto the space $\tw,$ of course,
which sends the domain $\nt$ of all \eqr s $\emt{\om\nu+2}$ 
onto a closed set $D\sq\tw.$)
}

\vyk{

\punk{Are the approximations complete in their Borel classes\,?}
\label{rem}

The \eqr s $\emt\xi$ lead to further questions, 
the most interesting of which is whether they can be complete 
for those Borel classes indicated in Lemma~\ref{nn}\ref{nn2}, 
analogously to the completeness of $\ent$ by Theorem~\ref{max}. 
We don't know the answer yet, but there is a promising approach.

\bdf
\label{rt'}
For any tree $R\sq(2\ti2\ti\dN)\lom,$ we define 
\dm
\bay{lcl}
\use R &=& \ens{\ang{x,y}\in\tw\ti\tw}
{\rah{\ceq Rxy}=\iy}\,;
\\[\dxii]

\uap\xi R &=& \ens{\ang{x,y}\in\tw\ti\tw}{\rah{\ceq Rxy}\ge\xi}\,,
\eay
\dm
where $\ceq Rxy=\ens{s\in\nlw}{R(x\res{\lh s},y\res{\lh s},s)},$ 
as in \eqref{e1}. 
A tree $R\sq{(2\ti2\ti\dN)}\lom$ is \dd\xi{\it restricted\/}, 
where $\xi<\omi,$ if and only if $\use R=\uap\xi R$. 
\edf

Thus $\use R$ is a $\fs11$ subset of $\tw\ti\tw,$ 
any $\fs11$ set $A\sq\tw\ti\tw$ can be presented in such a form, 
and the intersection $\use R=\bigcap_{\xi<\omi}\uap\xi R$ 
satisfies the canonical intersection property \cip. 
To establish \cip, note that the sets $A(R)$ and $A^\xi(R)$ are 
preimages of the sets, resp., $\cT^\iy$ and $\cT^\xi,$ with respect 
to the continuous map $\ang{x,y}\lmto \ceq Rxy.$ 
It remains to recall that $\cT^\iy=\bigcap_{\xi}\cT^\xi$ 
satisfies \cip\ by Proposition~\ref{18}. 

The following lemma is, essentially, a classical result. 
We refer to \cite{lus:cl} for quite similar results and 
constructions.
Let a \dd m{\it Lipschitz set\/} be any set $X\sq\tw\ti\tw$ such 
that we have ${\ang{x,y}\in X}\eqv{\ang{x',y'}\in X}$ whenever 
$x\yo x'\yo y\yo y'\in\tw$ satisfy $x\res m=x'\res m$ and 
$y\res m=y'\res m.$
Any such a set is obviously clopen in $\tw\ti\tw.$ 

\ble
\label{z}
If\/ $R\sq(2\ti2\ti\dN)\lom$ is a tree, $\nu<\omi,$ and\/ 
$k\ge1,$  then\/ 
$\uap{k}R$ is a\/ \dd kLipschitz set,   
$\uap{\om(1+\nu)}R$ is\/ $\fp0{1+2\nu},$ and\/ 
$\uap{\om(1+\nu)+k}R$ is\/ $\fs0{1+2\nu+1}$.

Conversely, if\/ $\nu<\omi\yt k\ge1,$ and\/ $A\sq\tw\ti\tw$ 
is a\/ \dd kLipschitz set, resp., 
a set in\/ $\fp0{1+2\nu},$ resp., 
a set in\/ $\fs0{1+2\nu+1},$
then there is a tree $R\sq(2\ti2\ti\dN)\lom$ such that\/ 
$A=\uap{k}R=\use R,$ resp.,
$A=\uap{\om(1+\nu)}R=\use R,$ resp., 
$A=\uap{\om(1+\nu)+k}R=\use R$.~\snos
{Note that $1+\nu=\nu$ and $1+2\nu=2\nu$ provided $\nu\ge\om$.}
\ele
\bpf
The proof of the first claim is somewhat similar to the estimation 
of Borel classes of the relations $\emt\xi$ in the proof of 
Lemma~\ref{nn}\ref{nn2}. 
Put 
\dm
\cer Rxys=\ens{t\in\nlw}{s\we t\in\ceq Rxy}
\quad\text{and}\quad
\uaq\xi Rs = \ens{\ang{x,y}\in\tw\ti\tw}{\rah{\cer Rxys}\ge\xi}
\dm 
for all 
$x\yo y\in\tw\yt s\in\nlw\yt\xi<\omi.$ 
We claim that, in the conditions of the lemma,  
$\uaq{k}Rs$ is a\/ \dd{(\lh s+ k)}Lipschitz set,   
$\uaq{\om(1+\nu)}Rs$ is\/ $\fp0{1+2\nu},$ and\/ 
$\uaq{\om(1+\nu)+k}Rs$ is\/ $\fs0{1+2\nu+1}.$    
To prove the first claim, it suffices to note that $\uaq{k}Rs$ 
is equal to the set of all pairs $\ang{x,y}$ such that 
$\ceq Rxy$ contains at least one sequence of the form 
$s\we t,$ where $t\in\dN^k.$ 
The other two claims are verified by induction with the help 
of obvious equalities 
\dm
\TS
\uaq{\xi+1}Rs=\bigcup_{j\in\dN}\uaq\xi R{s\we j}  
\quad\text{and}\quad
\uaq\la RS=\bigcap_{\xi<\la}\bigcup_{j\in\dN}\uaq\xi R{s\we j}
\quad\text{for limit $\la$}.
\dm
Finally, we observe that $\uap \xi R=\uaq\xi R\La$.

The proof of the converse in the lemma goes on by induction on 
$\nu$ and $k.$ 

{\it Case of Lipschitz sets.\/} 
Consider a \dd{k}Lipschitz set $A\sq\tw\ti\tw.$ 
There is a set $P\sq 2^k\ti2^k$ such that 
\dm
A=\ens{\ang{x,y}\in\tw\ti\tw}
{\sus \ang{u,v}\in P\:(u\su x\land v\su y)}.
\dm
Let $R$ be the tree of all triples 
$\ang{u,v,s}\in\tlw\ti\tlw\ti\nlw$ such that 
$\lh u=\lh v=\lh s$ and there is a pair $\ang{u',v'}\in P$ 
such that $u$ is compatible with $u'$ 
(that is, $u\sq u'$ or $u'\sq u$) and  
$v$ is compatible with $v'.$ 
Then obviously $\ceq Rxy=\nlw$ provided $\ang{x,y}\in A,$ 
and $\ceq Rxy$ does not contain any sequence of length $k$ 
otherwise. 
It follows that 
$A=\uap{k}R=\uap{\om}R=\use R$.  

{\it Case\/ $\fs0{1+2\nu+1}.$\/} 
Any $\fs0{1+2\nu+1}$ set $A\sq\tw\ti\tw$ is equal to a union 
$A=\bigcup_n A_n,$ where each $A_n$ is a set in 
$\fp0{1+2\nu}.$ 
Then, by the induction hypothesis, for any $n$ there is a tree 
$R^n\sq(2\ti2\ti\dN)\lom$ such that 
$A_n=\use{R^n}=\uap{\om(1+\nu)}{R^n}.$ 
Define a tree $R\sq(2\ti2\ti\dN)\lom$ such that 
$R(u\we i,v\we j,n\we s)$ iff $R^n(u,v,s)$ for all 
$i\yo j=0,1$ and $n\in\dN$ 
(and separately, $\ang{\La,\La,\La}\in R$).
Then 
$\ceq Rxy=\ans\La\cup \ens{n\we s}{s\in \ceq {R^n}xy}$ 
for any $x\yo y\in\tw,$  and hence the cross-section tree 
$\ceq Rxy$ is equal to $\sum^\ast_n{\ceq {R^n}xy}$ 
(see Subsection~\ref T on the countable sum operation). 
This implies $\rah{\ceq Rxy}=\tsup \rah{\ceq {R^n}xy}.$~\snos
{Recall that $\tsup X$ denotes the lest ordinal strictly larger 
than any ordinal in $X,$ or $\iy$ whenever $X$ contains $\iy$.} 
In particular, $\rah{\ceq Rxy}=\iy$ iff 
$\rah{\ceq {R^n}xy}=\iy$ for some $n,$ and hence we have 
$A=\use R$ by the choice of the trees $R^n.$ 
Furthermore, if $\ang{x,y}\nin A$ then $\ang{x,y}\nin A_n$ 
for all $n,$ therefore $\rah{\ceq {R^n}xy}<\om(1+\nu)\zd\kaz n,$ 
and this implies $\rah{\ceq Rxy}<\om(1+\nu)+1.$
We conclude that 
$A=\use R=\uap{\om(1+\nu)+1}R=\uap{\om(1+\nu)+k}R$ 
for any $k\ge1,$ as required.  

{\it Case\/ $\fp0{1+2\nu}.$\/} 
Here $\nu<\omi$ can be either a limit ordinal or a successor.
Any $\fp0{1+2\nu}$ set $A\sq\tw\ti\tw$ is equal to an intersection  
$A=\bigcap_n A_n,$ where each $A_n$ is a set in 
$\fs0{1+2\eta_n+1}\yt\eta_n<\nu,$ 
or just a \dd{k_n}Lipschitz set for some $k_n$ in the case when 
$\nu=0$ 
(then $\fp01=$ closed sets).
Then, by the induction hypothesis, for any $n$ there exists 
a tree $R^n\sq(2\ti2\ti\dN)\lom$ such that 
$A_n=\use{R^n}=\uap{\om(1+\eta_n)+1}{R^n},$
or $A_n=\use{R^n}=\uap{k_n}{R^n}$ whenever $\nu=0.$ 
%
A tree $R\sq(2\ti2\ti\dN)\lom$ can be easily defined so that, 
for all $x\yi y\in\tw,$ the cross-section tree $\ceq Rxy$ 
is isomorphic to 
$\prod^\ast_{n\in\dN}\ceq{R^n}xy.$ 
(See Subsection~\ref T on the countable product operation. 
The tree $R$ can be defined as essentially the set of 
all triples $\ang{u,v,\sg},$ where, for some 
$m,$ $u$ and $v$ belong to $2^m$ while $\sg=\ang{s_0,\dots,s_m},$ 
each $s_k$ belongs to $\dN^m$ and $\ang{u,v,s_k}\in R^k,$
with the order 
$\ang{u,v,\sg}\lef\ang{u',v',\sg}$ whenever $u\sq  u'\yt v\sq v',$ 
and $s\lef s'$ in the sense defined in Subsection~\ref T.) 

By Lemma~\ref{protl}, we have $\rah{\ceq Rxy}=\iy$ if and 
only if $\rah{\ceq{R^n}xy}=\iy$ for all $n.$
It follows, by the choice of the trees $R^n,$ that $A=\use R.$ 
Moreover, if $\ang{x,y}\nin A$ then $\ang{x,y}\nin A_n,$ 
and hence $\rah{\ceq{R^n}xy}\le\om(1+\eta_n)<\om(1+\nu)$ 
(or $<k_n<\om$ whenever $\nu=0$)
for at least one $n.$ 
On the other hand, we have
$\rah{\ceq Rxy}\le\tmin_{n}\rah{\ceq{R^n}xy}+n$ 
by Lemma~\ref{protl}. 
Thus $\rah{\ceq Rxy}\le\om(1+\nu).$  
It follows that  $A=\use R=\uap{\om(1+\nu)}R,$ as required.
\epf

In view of the lemma, we could try to prove the completeness of, 
\eg, $\emt{\om(1+\nu)+2}$ for the Borel class $\fs0{1+2\nu+1}$ as 
follows.
Consider any $\fs0{1+2\nu+1}$ \eqr{} $\rE$ on $\tw.$ 
Thus $\rE,$ as a set of pairs, is $1+\fs0{2\nu+1}$ in $\tw\ti\tw.$ 
It follows from Lemma~\ref z that there is a tree 
$Q\sq(2\ti2\ti\dN)\lom$ such that 
${\rE}=\uap{\om(1+\nu)+1}Q=\use Q,$ 
or, saying it differently, for any $x\yo y\in\tw,$
\dm
{x\rE y} \leqv {\rah{\ceq Qxy}\ge\om(1+\nu)+1}
\leqv {\rah{\ceq Qxy}=\iy}.
\dm
In particular, $Q$ is an \dd{(\om(1+\nu)+1)}restricted tree, 
in the sense of Definition~\ref{rt'}. 

\bcj
\label J
In this case, there exists a tree $R\sq(2\ti2\ti\dN)\lom$ 
satisfying \ref{1+} -- \ref{last} of Theorem~\ref{max}, together 
with the following additional reqirement: $R$ 
is an \dd{(\om(1+\nu)+2)}restricted tree, and hence, 
for all $x\yo y\in\tw,$
\dm
{x\rE y} \leqv {\rah{\ceq Rxy}\ge\om\nu+2}
\leqv {\rah{\ceq Rxy}=\iy}\,.\eqno\qed 
\dm
\ecj

Taking this on trust, we obtain 
\dm
{x\rE y} \leqv {\vt(x)\emt{\om(1+\nu)+2}\vt(y)}
\leqv {\vt(x)\ent\vt(y)}, 
\dm
where the map $\vt$ is defined by \eqref{vt} in Subsection~\ref{cae}. 
(Indeed the equality  
$\lip{\vt(x)}{\vt(y)}=\ceq Rxy$ holds for all $x\yo y$ 
by Lemma~\ref j.)
Thus, as in the proof of Theorem~\ref{max}, 
the map $\vt$ is a Borel reduction 
of $\rE$ to $\emt{\om(1+\nu)+2}.$ 
As $\rE$ is an arbitrary \eqr{} in $\fs0{1+2\nu+1},$ 
the relation $\emt{\om(1+\nu)+2}$ turns out to be complete  
in the class of all $\fs0{1+2\nu+1}$ \eqr s.
Thus, since $\emt{\om(1+\nu)+2}$ itself belongs to  
$\fs0{1+2\nu+1}$ by Lemma~\ref{nn}\ref{nn2}, we conclude that  
the class $\fs0{1+2\nu+1}$ contains a complete \eqr{}  
provided $\nu\ge\om$ (then $1+\nu=\nu$).

However this argument involves Conjecture~\ref J, which,   
unfortunately, we cannot yet confirm or reject.

}

\punk{Some questions}
\label Q


It can be reasonably conjectured that 
$\emt\eta\rebs\emt{\om\nu}\rebs\emt{\om\nu+n}$ 
whenever $\eta<\om\nu$ and $n\ge1.$ 
The background idea here is that there is no \dd\reb largest Borel 
\eqr{} (noted in \cite{sinf}), therefore, the sequence  of \eqr s 
$\emt\xi$ has uncountably many indices of  
\dd\rebs increase (in strict sense).
On the other hand, it seems plausible that 
$\emt{\om\nu+n}\eqb \emt{\om\nu+n+1}$ provided $n\ge1$.

Few more interesting questions. 

Which Borel classes contain complete \eqr s\,?

A related problem can be discussed here. 
It was once considered a viable conjecture 
(see, \eg, \cite{kr:izv}) 
that the \eqr{} $\rtd$
called {\it the equality of countable sets of reals\/}~\snos
{$\rtd$ is defined on $\dR^\om,$ the set of countable sequences of 
reals, so that $x\rtd y$ iff the sets $\ens{x(n)}{n\in\om}$ and 
$\ens{y(n)}{n\in\om}$ (countable sets of reals) are equal to each 
other.}
is not Borel reducible to any \eqr\/ $\rei$ induced by a Borel 
ideal $\cI\sq\pws{\dN}.$
It follows from Theorem~\ref{M} that this is not the case,
in fact there is an ordinal $\xi<\omi$ such that
$\rtd\reb\emt\xi.$ 
What is the least ordinal $\xi$ satisfying this statement\,?

Finally, it should be stressed that all evaluations of Borel class 
of equivalence relations in this paper were related to the actual 
Borel class in Cantor's discontinuum-like spaces. 
A somewhat deeper approach of ``potential'' Borel classes of 
equivalence relations in \cite{sinf} may require corresponding 
adjustment of arguments.


\backmatter

\addcontentsline{toc}{chapter}{\bibname}
\bibliographystyle{plain}
{\api
}

\mtho
\addcontentsline{toc}{chapter}{\indexname}
{\small\printindex}

\end{document}